\documentclass[11pt]{article}
\usepackage[T1]{fontenc}
\usepackage{amsmath, amsthm, amssymb}
\usepackage{mathrsfs}
\usepackage{enumitem}
\usepackage{bm}
\usepackage{tikz}
\usepackage{xcolor}
\usepackage{xypic}
\usepackage{verbatim}
\usepackage{cite}
\usepackage{hyperref}
\usepackage{nameref}
\usepackage{blindtext}
\usepackage{comment}
\usepackage{titling}
\usepackage{titlesec}
\usetikzlibrary{arrows, decorations.pathmorphing}
\allowdisplaybreaks

\titleformat{\section}[block]{\color{black}\normalsize\bfseries\filcenter}{\thesection}{1em}{}

\titleformat{\subsection}[runin]{\color{black}\normalsize\bfseries}{\thesubsection}{1em}{}

\titleformat{\subsubsection}[runin]{\color{black}\normalsize\bfseries}{\thesubsubsection}{1em}{}

\titleformat{\paragraph}[runin]{\color{black}\normalsize\bfseries}{\theparagraph}{1em}{}

\usetikzlibrary{quotes,angles}

\usepackage[
top    = 2.00cm,
bottom = 1.50cm]{geometry}

\hypersetup{
    colorlinks=true,
    linkcolor=black,
    filecolor=magenta,
    urlcolor=cyan,
}

\newtheoremstyle{mystyle}%                % Name
  {}%                                     % Space above
  {}%                                     % Space below
  {\itshape}%                                     % Body font
  {}%                                     % Indent amount
  {\bfseries}%                            % Theorem head font
  {.}%                                    % Punctuation after theorem head
  { }%                                    % Space after theorem head, ' ', or \newline
  {\thmname{#1}\thmnumber{ #2}\thmnote{ (#3)}}%                                     % Theorem head spec (can be left empty, meaning `normal')

\theoremstyle{mystyle}

\newtheorem{theorem}{Theorem}[section]

\newtheorem{conjecture}[theorem]{Conjecture}

\newtheorem{definition}[theorem]{Definition}

\newtheorem{lemma}[theorem]{Lemma}
\newtheorem{proposition}[theorem]{Proposition}
\newtheorem{remark}[theorem]{Remark}
\newtheorem*{assumption*}{Assumption}

\newcommand{\bkt}[1]{\biggl( #1 \biggr) }
\newcommand{\sbkt}[1]{\bigl( #1 \bigr) }
\newcommand{\bigbkt}[1]{\biggl[ #1 \biggr] }
\newcommand{\sbigbkt}[1]{\bigl[ #1 \bigr] }

\newcommand{\bignorm}[1]{\biggl|\biggl| #1 \biggr| \biggr|}
\newcommand{\abs}[1]{\biggl| #1 \biggr| }

\newcommand{\set}[1]{\biggl \{ #1 \biggr\} }

\newcommand{\biggbkt}[1]{\biggl\{ #1 \biggr\} }
\numberwithin{equation}{section}

\title{\textbf{Asymptotics and the sub-limit at $L^{2}$-criticality of higher moments for the SHE in
dimension $d\geq 3$}}
\author{
  Te-Chun Wang 
  \thanks{Department of Mathematics and Statistics, University of Victoria, British Columbia, Canada.}
  \thanks{Email: \href{mailto:techun@uvic.ca}{techun@uvic.ca}}
}

\begin{document}
\maketitle

\vspace{-1em}
\begin{abstract}
In this article, we consider the $d$-dimensional mollified stochastic heat equation (SHE) when the mollification parameter is turned off. Here, we concentrate on the high-dimensional case $d \geq 3$. Recently, the limiting higher moments of the two-dimensional mollified SHE have been established. However, this problem in high dimensions remains unexplored to date. The main theorems of this article aim to answer this question and prove some related properties: (1) Our first main result, based on the spectral theorem for the unbounded operator, proves the divergence of the higher moments of the high-dimensional mollified SHE even when the system is strictly inside the $L^{2}$-regime. This phenomenon is completely opposite to its two-dimensional counterpart; (2) To further differentiate the nature of the high-dimensional case from the case in two dimensions, our second main result proves the unboundedness of the sub-limiting higher moments of the three-dimensional mollified SHE at the $L^{2}$-criticality. Here, the sub-limiting higher moment is a natural limit of the higher moment of the three-dimensional mollified SHE at the $L^{2}$-criticality; (3) As an application, we provide partial results for the conjecture about the high-order critical regimes of the continuous directed polymer. The other byproduct of the above results gives a proper estimate for the critical exponent of the polymer in the $L^{2}$-regime.

%where the sub-limiting higher moment is a non-Gaussian formal limit of the higher moment of the three-dimensional mollified SHE. Here, the evaluation of the formal limit is as the non-Gaussian limiting moment of the two-dimensional mollified SHE studied in \cite{chen}, but the contact interactions between each pair of Brownian motions are governed by the three-dimensional delta-Bose gas \cite[Chapter I.1]{albeverio2012solvable}; (3) As an application, we provide a partial result for the conjecture about the high-order critical regimes of the continuous directed polymer. Moreover, if the polymer is in the $L^{2}$-regime, we give a proper estimate for its critical exponent.

$\newline$
{\textsl{Keywords}: stochastic heat equation; Kardar–Parisi–Zhang equation; directed polymers;  high-order critical regimes; critical exponents; delta-Bose gas; Efimov effect.}
$\newline$
{\textsl{Mathematics Subject Classification}: 60K37, 60K35, 35R60,  35Q82,  60H15, 82D60, 82D60.}

\end{abstract}

\fontdimen2\font=0.29em

\tableofcontents
\section{Introduction} 

\subsection{Background.}
In this article, we consider the $d$-dimensional stochastic heat equation (SHE) with multiplicative noise, which is formally written as follows:
\begin{equation} \label{SHE}
    \partial_{t} \mathscr{U}(x,t) = \frac{1}{2} \Delta \mathscr{U}(x,t) + \beta \mathscr{U}(x,t) 
    \xi(x,t), 
    \quad  (x,t) \in \mathbb{R}^{d} \times \mathbb{R}_{+}.
\end{equation}
Here, $\xi(x,t)$ is a space-time white noise and $\beta$ is a positive coupling constant playing the role of the strength of the noise. The SHE is by now a crucial model since it lies at the heart of two areas of intense research in recent years. On the one hand, equation (\ref{SHE}) enjoys a long history of study due to its connection with the Kardar-Parisi-Zhang (KPZ) equation, which was introduced by Kardar, Parisi and Zhang in order to describe the dynamic of random growth interfaces in $d + 1$ dimensions \cite{KPZ}. On the other hand, $\mathscr{U}$ shares close analogies to the partition function of the continuous directed polymer \cite{UU,Alberts_2013,Feynman-Kac,NK_KPZ}, which is a fundamental model in statistical physics. See \cite{comets2017directed} for a review of the directed polymer model. 

The overall goal of this article is to investigate the asymptotics of the higher moments of a regularization scheme of equation (\ref{SHE}) in high dimensions with various coupling constants. The first objective is to prove the divergences of the higher moments \emph{inside the $L^{2}$-regime}. The main theorem is in contrast with the convergences in two dimensions. This phenomenon partially explains why the singularity of the high-dimensional SHE is stronger than the two-dimensional case. To further illustrate this phenomenon, we consider a natural three-dimensional limit of the higher moment at the $L^{2}$-criticality. Then, the second goal is to reveal its unboundedness. This fact is also opposite to its two-dimensional counterpart. Thanks to these results, the third goal is to illustrate their applications to the continuous directed polymer. Finally, the connections between our main results and various results in quantum physics will be discussed.

In the remainder of the introduction, we first introduce an approximation scheme of equation (\ref{SHE}) in the next subsection, and then we will give a brief overview the main results and their background in Section \ref{section_overview_main_result}. Throughout this article, we assume that $d \geq 3$.

\subsection{Mollified stochastic heat
equation.} \label{section_the_model}
From the mathematical perspective, equation (\ref{SHE}) is meaningless since the noise is too irregular for the multiplicative term on the right-hand side of equation (\ref{SHE}) to be well-defined. To rigorously study equations (\ref{SHE}), we then introduce in this subsection a renormalization of this equation. 

This method starts with mollifying the white noise in equation (\ref{SHE}). For this purpose, let $\sbkt{\xi(x,t)}_{(x,t)\in \mathbb{R}^{d}\times \mathbb{R}}$ be a space-time white noise. To be more specific, $\{\xi(h):h\in\mathscr{S}(\mathbb{R}^{d} \times \mathbb{R}_{+})\}$ is a zero-mean Gaussian random field with covariance 
\begin{align} \label{definition_white}
    \mathbb{E}[\xi(h_{1}) \xi(h_{2})] = \int_{\mathbb{R}^{d}\times \mathbb{R}_{+}} h_{1}(x,t) h_{2}(x,t) dtdx, \quad\text{where}\quad \xi(h) := \int_{\mathbb{R}^{d} \times \mathbb{R}_{+}} h(x,t) \xi(x,t) dxdt.
\end{align}
Here, $\mathscr{S}(\mathbb{R}^{d} \times \mathbb{R}_{+})$ denotes the Schwartz space on $\mathbb{R}^{d}\times \mathbb{R}_{+}$. With this notation, we consider the spatially mollified Gaussian noise as given below:
\begin{equation*}
    \xi_{\varepsilon} (x,t)
    := \phi_{\varepsilon} * \xi (x,t)
    = \int_{\mathbb{R}^{d}} \phi_{\varepsilon}(x-y) \xi(y,t) dy.
\end{equation*}
Here, $\phi_{\varepsilon}(x) := \varepsilon^{-d} \phi(x/\varepsilon)$ is an approximation of the identity. In the sequel, we work with a fixed mollifier $\phi: \mathbb{R}^{d} \mapsto \mathbb{R}_{+}$ which satisfied the following properties: \emph{$\phi \in C_{c}^{\infty}(\mathbb{R}^{d};\mathbb{R}_{+})$ and $\phi$ is symmetric-decreasing}. To be more accurate, $\phi$ is spherically symmetric, and for any unit vector $e\in \mathbb{R}^{d}$, the function $r\in \mathbb{R}_{+} \mapsto \phi(r \cdot e)$ is decreasing. As a result, $(\xi_{\varepsilon}(x,t))_{(x,t)\in \mathbb{R}^{d}\times \mathbb{R}}$ is a centered Gaussian random field with the following covariance:
\begin{equation} \label{definition_R_epsilon}
    \mathbb{E}[\xi_{\varepsilon}(x,t)\xi_{\varepsilon}(y,s)] = \delta_{0}(t-s) R_{\varepsilon}(x-y), \quad \text{where} \;  R_{\varepsilon}(x) := \phi_{\varepsilon} * \phi_{\varepsilon}(x).
\end{equation}

With the above notations at hand, we now clarify the core of the renormalization. The main idea is to find a proper constant $\beta_{\varepsilon}$ such that $\mathscr{U}_{\varepsilon}(x,t)$ has a non-trivial limit. Here, $\mathscr{U}_{\varepsilon}(x,t)$ solves the \emph{mollified stochastic heat equation} as follows:
\begin{equation} \label{Mollified_SHE}
    \partial_{t} \mathscr{U}_{\varepsilon}(x,t) = \frac{1}{2} \Delta \mathscr{U}_{\varepsilon}(x,t) + \beta_{\varepsilon} \mathscr{U}_{\varepsilon}(x,t) 
    \xi_{\varepsilon}
    (x,t), \quad \mathscr{U}_{\varepsilon}(0,x) = U_{0}(x).
\end{equation}
where $U_{0} \in C_{b}(\mathbb{R}^{d};\mathbb{R}_{+})$. The authors in \cite{weakandstrong} observed that, in view of the generalized Feynman-Kac formula, the solution of equation (\ref{Mollified_SHE}) admits the following representation:
\begin{equation} \label{representation_FK}
    \mathscr{U}_{\varepsilon}(x,t) 
    = \mathbf{E}^{B}_{x}\bigbkt{\exp\bkt{\beta_{\varepsilon} \int_{0}^{t} ds\int_{\mathbb{R}^{d}} dy\phi_{\varepsilon}(B(t-s) - y) \xi(y,s) - \frac{\beta_{\varepsilon}^{2} R_{\varepsilon}(0) t }{2}} U_{0}(B(t))}.
\end{equation}
Here, $\mathbf{P}^{B}_{x}$ is the law of $d$-dimensional Brownian motion $B$ starting at $x \in \mathbb{R}^{d}$. In this way, according to \cite{weakandstrong,COSCO2022127}, if the intensity of the noise is chosen as follows with a proper coupling constant $\beta>0$:
\begin{equation} \label{definition_beta_epsilon}
    \beta_{\varepsilon} := \beta \cdot\varepsilon^{\frac{d-2}{2}},
\end{equation}
then $\mathscr{U}_{\varepsilon}(x,t)$ has a non-trivial limit when the mollification is removed.
Although $\mathscr{U}_{\varepsilon}(x,t)$ depends on $\beta$, we ignore this dependence in the sequel for the sake of simplicity.

\subsection{Overview of the main results.} \label{section_overview_main_result}
Now, we give a brief clarification about the main theorems and their background. Before we proceed, let us first introduce some notations and terminologies. Recall that $U_{0}(x)$ is the initial datum of the mollified SHE. Throughout this paper, we adopt the following notation to denote the $N$-th moment of $\mathscr{U}_{\varepsilon}(x,t)$:
\begin{align} \label{definition_moment}
    \mathcal{Q}_{\varepsilon;t}^{\beta,N}U_{0}^{\otimes N}(\Vec{x}) := \mathbb{E}\sbigbkt{\mathscr{U}_{\varepsilon}(\Vec{x}(1),t)...\mathscr{U}_{\varepsilon}(\Vec{x}(N),t)} \quad \forall \Vec{x} \in \mathbb{R}^{N \times d}, \; t,\beta,\varepsilon >0, \; N\geq 2.
\end{align}
Here, $\Vec{x}(j) \in \mathbb{R}^{d}$ such that $\Vec{x} := (\Vec{x}(1),...,\Vec{x}(N))$, and we call it the $j$-th component of $\Vec{x}$. We call $\mathcal{Q}_{\varepsilon;t}^{\beta,N}U_{0}^{\otimes N}(\Vec{x})$ a higher moment if $N\geq 3$. Moreover, we say $\mathcal{Q}_{\varepsilon;t}^{\beta,N}U_{0}^{\otimes N}(\Vec{x})$ a \emph{mixed} $N$-th moment of $\mathscr{U}_{\varepsilon}(x,t)$ if $\Vec{x} \in \mathbb{R}^{N \times d} \setminus \Pi_{N}$. Here, $\Pi_{N}$ defined as follows denotes the collection of all $\Vec{x}\in \mathbb{R}^{N \times d}$ that has the same components:
\begin{align} \label{definition_big_pi}
    \Pi^{(N)} := \set{ \Vec{x} \in \mathbb{R}^{N \times d}: \Vec{x}(i) = \Vec{x}(i') \text{ for some } \forall 1\leq i \neq i' \leq N}.
\end{align}
Furthermore, the system is in the $L^{N}$-regime or at the $L^{N}$-criticality if the coupling constant $\beta$ of the mollified SHE (\ref{Mollified_SHE}) is in the interval $(0,\beta_{L^{N}})$ or $\beta = \beta_{L^{N}}$, respectively. Here, the critical coupling constant $\beta_{L^{\gamma}}$ is defined as follows, which denotes the largest positive number such that the random fields $(\mathscr{U}_{\varepsilon}(x,t))_{(x,t) \in \mathbb{R}^{d} \times \mathbb{R}_{+}}$ with $U_{0} \equiv 1$ are uniformly bounded in $L^{\gamma}$: 
\begin{align} \label{definition_beta_L_N}
    \beta_{L^{\gamma}} := \sup\set{\beta > 0: \sup_{\varepsilon > 0,\; (x,t) \in \mathbb{R}^{d} \times \mathbb{R}_{+}}\mathbb{E}[\mathscr{U}_{\varepsilon}(x,t)^{\gamma}] < \infty} \quad \forall \gamma\in (1,\infty).
\end{align}
%where the limit can be regarded as the Gaussian kernel $G_{t}U_{0}(x)$,by making use of
\paragraph{Part I: Higher moments and their critical regimes.} Studying the $\varepsilon \to 0^{+}$ asymptotics of the higher moments of the two-dimensional mollified SHE goes back to the works \cite{chen,Gu_2d,caravenna2019moments,Heat_flow}. See also \cite{caravenna2019moments,Heat_flow,DDP_regime,cosco2023moments,CCO} for related results about the higher moments of the partition function of the discrete directed polymer. Briefly speaking, the works \cite{chen,Gu_2d,caravenna2019moments} consider the asymptotics of the higher moments of the two-dimensional mollified SHE in order to describe its limiting behavior at the $L^{2}$-criticality. %spatially averaged limit of

Regarding the high-dimensional case, although the spatially averaged limit of $\mathscr{U}_{\varepsilon}(x,t)$ was proved in \cite{weakandstrong} and was extended in \cite{COSCO2022127} to the entire weak disorder regime, the asymptotics of the mixed higher moments for $\mathscr{U}_{\varepsilon}(x,t)$ are still unclear. Even though this problem loses its original application, we still believe that the asymptotics of the mixed higher moments are interesting in their own right. On the one hand, it is unknown whether the limit of the mixed $N$-th moment $\mathcal{Q}_{\varepsilon;t}^{\beta,N}U_{0}^{\otimes N}(\Vec{x})$ corresponds to the limit of $\mathscr{U}_{\varepsilon}(x,t)$, which can be regarded as the Gaussian kernel $G_{t}U_{0}(x)$ as the results in \cite{COSCO2022127}. In fact, it can be shown that the limit of the mixed second moment $\mathcal{Q}_{\varepsilon;t}^{\beta,2}(\Vec{x})$ does \emph{not} match the limit of $\mathscr{U}_{\varepsilon}(x,t)$ when $d = 3$ and $\beta = \beta_{L^{2}}$. This fact can be seen through (\ref{limit_two_particle}). On the other hand, the asymptotic of $N$-th moment $\mathcal{Q}_{\varepsilon;t}^{\beta,N}U_{0}^{\otimes N}(\Vec{x})$ provides information about the corresponding high-order critical regime. As a matter of fact, it is believed that the critical coupling constants $(\beta_{L^{N}})_{N\geq 2}$ are in strictly descending order.
\begin{conjecture} \label{Conjecture_order}
Suppose that $d\geq 3$. Then, the high-order critical regimes have the following property: 
\begin{align} \label{inequality_conjecture_order}
    ... < \beta_{L^{N+1}} < \beta_{L^{N}} < ... < \beta_{L^{3}} < \beta_{L^{2}}.
\end{align}
\end{conjecture}
According to our understanding, the above conjecture remains unexplored so far. By
contrast, for the two-dimensional discrete directed polymer, the conjecture about the high-order critical regime is well understood. See \cite{DDP_regime} for the precise result.

%See \cite{DDP_regime} for the high-order critical regimes of the two-dimensional discrete directed polymer.

With the above motivation, the first objective of our article is to compute the $\varepsilon \to 0^{+}$ asymptotics of the mixed higher moments of $\mathscr{U}_{\varepsilon}(x,t)$. To make the problem nontrivial, we aims to study the $\varepsilon\to 0^{+}$ asymptotics of the mixed higher moments when the system is within the $L^{2}$-regime or at the $L^{2}$-criticality. This reveals new difficulties as compared to the case of $\beta > \beta_{L^{2}}$. The primary reason is that the divergence of the limiting mixed higher moment directly follows from 
the divergence of the limiting mixed second moment when $\beta > \beta_{L^{2}}$, where the latter divergence is proved in our main result as well. However, if the system belongs to the $L^{2}$-regime or at the $L^{2}$-criticality, then the asymptotics of the mixed higher moments cannot be concluded in this way.

About the asymptotics of the mixed higher moments when the system is within the $L^{2}$-regime or at the $L^{2}$-criticality, it seems natural to conjecture that the limit of the mixed higher moment may be similar to the two-dimensional case as the results in \cite{chen,Gu_2d}. To be more precise, for a proper initial datum $U_{0}$ and for every $\Vec{x} \in \mathbb{R}^{N \times d} \setminus \Pi^{(N)}$, one may feel that $\mathcal{Q}_{\varepsilon;t}^{\beta,N}U_{0}^{\otimes N}(\Vec{x})$ tends to the Gaussian kernel $G^{(N)}_{t}U^{\otimes N}_{0}(\Vec{x})$ in the $L^{2}$-regime, and $\mathcal{Q}_{\varepsilon;t}^{\beta,N}U_{0}^{\otimes N}(\Vec{x})$ tends to a non-Gaussian kernel $\mathrm{Q}^{\beta,N}_{t}U^{\otimes N}_{0}(\Vec{x})$ at the $L^{2}$-criticality. Here, $U^{\otimes N}_{0}(\Vec{x}) := \bigotimes_{j=1}^{N} U_{0}(\Vec{x})$. This conjecture is convincing mainly because, thanks to the results in \cite[Chapter I.1]{albeverio2012solvable}, it can be shown that the above conjecture holds when $d = 3$ and $N = 2$. In this way, one may expect that the above conjecture holds for all $N\geq 3$ and when $d = 3$ by applying the same arguments used in the two-dimensional case \cite{chen}.

Although the above conjecture seems reasonable, our primary result aims to show that it does not hold. To be more precise, the first main theorem shows that for every $d\geq 3$ and for a large $N\geq 3$, \emph{there exists a nontrivial interval in the $L^{2}$-regime such that the $N$-moment $\mathcal{Q}_{\varepsilon;t}^{\beta,N}U_{0}^{\otimes N}(\Vec{x})$ is divergent for all $\Vec{x} \in \mathbb{R}^{N \times d}$}. As an application, we provide a partial result for Conjecture \ref{Conjecture_order} by showing the following estimate:
\begin{align*}
    \beta_{L^{N}}  \approx \frac{C_{d,\phi}}{\sqrt{N-1}},
\end{align*}
where $a\approx b$ denotes $c \cdot b \leq a \leq c'\cdot b$. See Section \ref{section_main_part_1} for the precise statement. In particular, by making use of the above estimate, we explore some properties of the continuous directed polymer. For the sake of simplicity, we postpone the elaboration of this until Section \ref{section_CDP}.

With the aim of proving the above divergence, the following two primary difficulties must be conquered. 
\begin{enumerate} [label=(\roman*)]
    \item The first barrier comes from the fact that the mixed higher moment is much smaller than the higher moment at the same space-time point. See the right-hand side of (\ref{definition_beta_L_N}) for the latter moment. Due to this reason, even though the condition $\beta > \beta_{L^{N}}$ provides us the unboundedness of the higher moment at the same space-time point, it is still nontrivial to conclude the divergence of the mixed higher moment by using this fact.
    \item The second obstacle arises because even if the mixed higher moment is divergent when the system has a certain coupling
    constant $\beta>\beta_{L^{N}}$, how to ensure that this coupling constant could be in the $L^{2}$-regime is still being determined.
\end{enumerate}
To overcome the above difficulties, we will re-express the mixed higher moment as a rescaled semigroup of a many-body quantum system. In this way, both of the above problems can be solved by studying the corresponding Hamiltonian through the theory of the unbounded operator. We refer to Section \ref{section_proofoutline} for a more explicit illustration.

%it is used to insure that the divergence in claim (A) occurs inside the $L^{2}$-regime. In fact, instead of proving $\beta_{L^{N}} < \beta_{L^{2}}$, we will compute the value of $\beta_{L^{N}}$ directly, where the result not only implies claim (B) but also show that the high-order critical regimes $(\beta_{L^{N}})_{N\geq 2}$ are in strictly descending order. Consequently, our main result will also answer the previous question about the high-order critical regimes of the $d$-dimensional continuous directed polymer since $\beta_{L^{N}}$ has the following representation that corresponds to these regimes:
%\begin{align} \label{identity_beta_L_N_2}
%    \beta_{L^{N}} = \sup\set{\beta > 0: \sup_{T > 0} \mathbb{E}[\mathcal{Z}^{\beta}_{T}(0;\xi)^{N}] < \infty},
%\end{align}
%where we have used (\ref{definition_beta_L_N}) and (\ref{connect_SHE_CDP}).

\paragraph{Part II: The sub-limiting higher moments for the 3D-SHE.}
Although the main theorem in the first part disproves the conjecture of the convergence of the higher moment, at this stage, one may still want to know a more detailed reason \emph{why the approach used in the two-dimensional case \cite{chen} does not work in high-dimensions}. To explain this problem more precisely, we consider the sub-limiting higher moment of the three-dimensional mollified SHE. Simply put, it is a natural three-dimensional counterpart of the limiting higher moment of the two-dimensional mollified SHE studied in \cite{chen}.

Here, we illustrate the above concept in detail. The author in \cite{chen} showed that the limiting mixed higher moment of the two-dimensional mollified SHE at the $L^{2}$-criticality is a series such that each term in the series is a Brownian path integral that consists of a sequence of nonconsecutive two-body contact interactions, where the interaction is described by the two-dimensional delta-Bose gas studied in \cite[Chapter I.5]{albeverio2012solvable}. See \cite[(2.19)]{chen} for an explicit expression. We also refer to \cite{Gu_2d,caravenna2019moments} for similar expansions of the higher moment. With this elaboration, to explain the problem mentioned in previous paragraph, we consider the same object as \cite[(2.19)]{chen}, but the two-body contact interaction is governed by the three-dimensional delta-Bose gas considered in \cite[Chapter I.1]{albeverio2012solvable}. In the sequel, we call this object \emph{the sub-limiting higher moment}. 

Therefore, the second main result of this article aims to prove \emph{the divergence of the sub-limiting higher moment for the three-dimensional mollified SHE with an initial datum that can be bounded from below by a Gaussian kernel}. We refer to Section \ref{section_main_part2_pseudo_limiting_moments} for the precise statement. In fact, from a quantum physics point of view, the sub-limiting higher moment can be regarded as a natural limit of an approximated semigroup of the many-body delta-Bose gas in three dimensions. From this perspective, our work is analogous to \cite{Teta}. To be more specific, \cite[Lemma 6.2]{Teta} proves the unboundedness of an expected limiting quadratic form. But, the unboundedness of an expected limiting semigroup is proved here. For the sake of simplicity, we postpone the clarification of this until Section \ref{section_related_literature}.

About the difficulty to obtain this divergence, the strategy mentioned in the first part, which uses the theory of the unbound operator, is no longer available. This is because the sub-limiting higher moment is even smaller than the limiting higher moment (see Proposition \ref{proposition_fatou}). Then, an obstacle arises because one has to develop \emph{a systematic strategy} to estimate the nonconsecutive contact interactions in the path integrals in the sub-limiting higher moment. This problem is similar to the difficulty of proving the first part in \cite[Theorem 2.3]{chen}, but we deal with a lower bound here. 

To overcome the above obstacle, taking advantage of the Gaussian lower bound for the initial condition, we can \emph{systematically compute} the path integrals in the sub-limiting three moment through an idea that re-expresses and combines the Gaussian kernels in the path integral. This strategy eventually leads to the following geometric-type lower bound with a ratio close to $1.008$:
\begin{align} \label{formal_lower_bound}
    \sum_{k\geq k_{0}}^{\infty} (1.008)^{k} \lesssim \textit{"The sub-limiting three moment".}
\end{align}
We refer to Section \ref{section_proofoutline} for a more detailed illustration about the above ideas.

As an application, the divergence of the sub-limiting higher moment gives a three-dimensional refinement of the main result in the first part. To be more precise, for the three-dimensional mollified SHE, \emph{the divergence of the $N$-th moment holds for every $N\geq 3$ inside the $L^{2}$-regime}. The other byproduct is a three-dimensional partial result for Conjecture \ref{Conjecture_order}, which shows that if $d = 3$, then \emph{the $L^{N}$-regime is a proper subset of the $L^{2}$-regime for all $N\geq 3$}. 

However, due to the connection between the three-dimensional mollified SHE and the three-dimensional many-body quantum system mentioned in the first part, the above applications are also consequences of a well-known result in quantum physics, which is called the Efimov effect. Loosely speaking, a three-body quantum  system in dimension three with suitable interactions has an infinite number of negative bound state energies if the Hamiltonians of the two-body subsystems do not have negative spectrum, and all of them are resonant. This property will be illustrated in detail in Section \ref{section_related_literature}. 

%our second main result proves the unboundedness of the sub-limiting moment when the system is at the $L^{2}$-criticality and $d = 3$. Roughly speaking, the sub-limiting moment is an object constructed in the same way as the limiting moment of the two-dimensional mollified SHE by using \cite{albeverio2012solvable}. 

\subsection{Organization of the paper.} \label{section_structure}
The remainder of this article is organized in the following way. In Section \ref{section_main}, we will first introduce all the main results of this paper and the connection between our work and various results in the literature. The proof ideas for the main results will be presented in Section \ref{section_proofoutline}. The rest of the sections are dedicated to the details of the proofs. See Section \ref{section_organization_proof} for more details.
\paragraph{List of notations.} 
\small
In the sequel, we adopt the following notations.
\begin{itemize}[label={--}]
    \item Throughout this article, $d \geq 3$ and $N\geq 2$ are positive integers, and $\gamma$ and $\beta$ are positive real numbers.
    \item Unless specified, all the functions in this article are complex-valued:
    \begin{align*}
        C_{c}^{\infty}(\mathbb{R}^{\mathbf{d}}):=C_{c}^{\infty}(\mathbb{R}^{\mathbf{d}};\mathbb{C}), \; L^{p}(\mathbb{R}^{\mathbf{d}}):=L^{p}(\mathbb{R}^{\mathbf{d}};\mathbb{C}), \; \mathscr{S}(\mathbb{R}^{\mathbf{d}}) := \mathscr{S}(\mathbb{R}^{\mathbf{d}};\mathbb{C}),
    \end{align*}
    where $\mathscr{S}(\mathbb{R}^{\mathbf{d}};\mathbb{C})$ denotes the Schwartz space on $\mathbb{R}^{\mathbf{d}}$. Also, all the functions in $C_{c}^{\infty}(\mathbb{R}^{\mathbf{d}};\mathbb{R}_{+})$ , $C_{b}(\mathbb{R}^{\mathbf{d}};\mathbb{R}_{+})$, and $\mathscr{S}(\mathbb{R}^{\mathbf{d}};\mathbb{R}_{+})$ are non-negative.
    \item We use $\Vec{x}$ to denote the vector in $\mathbb{R}^{N \times d}$, and use $\Vec{x}(j)$ to denote the $j$-th component in $\mathbb{R}^{d}$. In other words, $\Vec{x} = (\Vec{x}(1),...,\Vec{x}(N))$.
    \item The notations $\mathbf{P}^{\Vec{B}}_{\Vec{x}}$ and $\mathbf{P}^{B}_{x}$ denotes the distribution of $N \times d$-dimensional Brownian motion and $d$-dimensional Brownian motion, respectively, where $\mathbf{E}^{\Vec{B}}_{\Vec{x}}$ and $\mathbf{E}^{B}_{x}$ denote the corresponding expectations. 
    \item Let $G^{(N)}_{t}(\Vec{x})$ defined as follows be the transition density of $N \times d$-dimensional Brownian motion:
    \begin{align*}
        G^{(N)}_{t}(\Vec{x}):= \biggbkt{\frac{1}{2\pi t}}^{\frac{N\times d}{2}} \exp\biggbkt{-\frac{|\Vec{x}|^{2}}{2t}}, \quad \Vec{x} \in \mathbb{R}^{N \times d}.
    \end{align*}
    For the sake of simplicity, $G_{t}(x) := G^{(1)}_{t}(x)$. Also, $\mathscr{G}_{t}(k)$ defined in (\ref{formula_G_t}) denotes the Fourier transform of $G_{t}$.
    \item The measure $\mathbf{P}^{\Vec{B}}_{\Vec{x}}\sbkt{d(\Vec{w}(t))_{0\leq t\leq T}: \Vec{w}(T) = \Vec{z}}$ is the non-normalized distribution of the $(N\times d)$-dimensional Brownian bridge starting from $(\Vec{x},0)$ to $(\Vec{z},T)$. In other words, the following measure is the distribution of $(N\times d)$-dimensional Brownian bridge starting from $(\Vec{x},0)$ to $(\Vec{z},T)$:
    \begin{align*}
    \mathbf{P}^{\Vec{B}}_{\Vec{x}}\biggbkt{d(\Vec{w}(t))_{0\leq t\leq T}: \Vec{w}(T) = \Vec{z}}\biggl/G^{(N)}_{T}(\Vec{x} - \Vec{z}).
    \end{align*}
    \item For a Hilbert space $H$ and a bound operator $\mathbf{L}$ on it, the notations $\langle u,v \rangle_{H}$ and $||\mathbf{L}||_{H}$ denote the inner product and the operator norm, respectively. In particular, we set 
    \begin{align*}
        \langle f,g \rangle_{L^{2}(\mathbb{R}^{\mathbf{d}})} := \int_{\mathbb{R}^{\mathbf{d}}} dz  \cdot \overline{f}(z) g(z).
    \end{align*}
    \item For an unbounded operator $\mathbf{A}$, 
    $D(\mathbf{A})$ denotes its domain. In particular, if $\mathbf{A}$ is self-adjoint, all the quantities $\sup \mathbf{A}$, $\inf \mathbf{A}$, $\mathbf{P}_{A}(d\lambda)$, and $\exp(T \cdot \mathbf{A})$ are defined in Section \ref{section_op}.
    \item We use $\bigotimes_{j=1}^{N} f_{j}(\Vec{x}) := \prod_{j=1}^{N} f_{j}(\Vec{x}(j))$ to denote the tensor product of functions $f_{1},...,f_{N}$ on $\mathbb{R}^{d}$. For the sake of simplicity, $U_{0}^{\otimes N}(\Vec{x}) := \bigotimes_{j=1}^{N} U_{0} (\Vec{x})$.
    \item For each $h \in \mathscr{S}(\mathbb{R}^{\mathbf{d}})$, its Fourier transform is defined by 
    \begin{align}  \label{FT}
        \widehat{h}(k) := \int_{\mathbb{R}^{\mathbf{d}}} dx \exp(-2\pi i k \cdot x) h(x).
    \end{align}
    \item The notation $B(a,r)$ denotes the open ball in $\mathbb{R}^{d}$ with center $a \in \mathbb{R}^{d}$ and radius $r>0$.
    \item The operators $\mathbf{H}^{\beta,N}$, $\mathcal{H}_{\varepsilon}^{\beta_{\varepsilon},N}$, $\mathbf{T}_{\beta}^{\lambda}$, and $\mathbf{L}^{\beta}$ are given in (\ref{short_range_connection}), (\ref{definition_hamiltonian_epsilon}), (\ref{definition_Birman-Schwinger_operator}), and (\ref{definition_hamiltonian_L}), respectively.
    \item The semigroups $\mathbf{Q}^{\beta,N}_{T}$, $\mathcal{Q}_{\varepsilon;t}^{\beta,N}$, and $\mathscr{Q}_{t}^{N}$ are given in (\ref{definition_Q_T}), (\ref{definition_moment}), and (\ref{definition_sub_limit}), respectively.
    \item The path integrals $\mathbf{I}^{\beta,N}_{T,k}$ and $\mathscr{I}_{0;t}^{N;(\ell_{1},\ell'_{1},\mathfrak{i}_{1}),...,(\ell_{m},\ell'_{m},\mathfrak{i}_{m})}$ are defined in (\ref{definition_short_path_integral_order_k}) and (\ref{definition_sub_limiting_path_integral}), respectively. The two-body contact interaction $\mathscr{A}^{(\ell_{j},\ell_{j}',\mathfrak{i}_{j})}_{0;u_{j-1},s_{j},u_{j}}(\Vec{x}_{j-1},\Vec{x}_{j})$ is given in (\ref{definition_transition_kernel}).
    \item The functions $\phi_{\varepsilon}$, $R_{\varepsilon}$, $R$, $\gamma^{*}(\beta)$, $\mathcal{Z}_{T}^{\beta}(0;\xi)$, ($V_{\beta}$, $\mathfrak{V}_{\beta}$, $\mathcal{G}^{\lambda}$), $(\mathbf{v}_{j},h_{\mathbf{v}_{1}})$, $\mathscr{C}(z,z')$, $\theta$,  $L_{m}$, $\eta_{m}$, $\zeta_{m}$ are defined in section \ref{section_the_model}, (\ref{definition_R_epsilon}), (\ref{definition_R}), (\ref{definition_gamm_star}), (\ref{definition_partition_function_1}), (\ref{definition_V_and_Yukawa}), Lemma \ref{lemma_Hilbert_Schmidt_operator}, (\ref{definition_C_z_z'}), Proposition \ref{proposition_small_support}, (\ref{definition_Lm}), (\ref{definition_eta_m}), and (\ref{definition_zeta}), respectively.
    
    \item The coupling constants $\beta_{\varepsilon}$, $\beta_{L^{\gamma}}$, $\beta_{N,+}$, $\alpha_{N,+}$, and $\beta_{N,+}^{sym}$ are given in (\ref{definition_beta_epsilon}), (\ref{definition_beta_L_N}), (\ref{definition_beta_N_+}), (\ref{definition_alpha_N}), and (\ref{definition_beta_sym}), respectively.
    \item The positive constants $r_{\phi}$ and $\alpha_{\infty,+}$ are defined in Proposition \ref{proposition_small_support} and (\ref{definition_alpha_infinity}), respectively.
    \item The measures $(\mathbf{u}^{z,z'}_{\varepsilon;\sigma},\overline{\mathbf{u}}^{z,z'}_{\varepsilon;\sigma}(dt))$ and $\mathbb{P}$ are defined in Proposition \ref{proposition_critical_second_moment_point_to_point} and (\ref{definition_white}), respectively.
    
    \item The sets $\mathcal{E}^{(N)} := \{(i,i'): 1\leq i<i'\leq N\}$. Also, $\mathcal{S}(\mathbb{R}^{N \times d})$, $\Pi^{(N)}$, and  $\widetilde{\Pi}^{(N)}$ are given in the step 3 of Section \ref{section_estimate_betaLN_and_betaN+}, (\ref{definition_big_pi}), and (\ref{definition_tilde_pi}), respectively.
    \item The notations $(\Vec{x}^{(\ell,\ell')},\Vec{x}^{(\ell,\ell')^{*}},\Vec{x}^{(\ell,\ell')^{c}})$, $\mathbf{t}_{j}$, and $\overline{r}_{m}$ are defined in (\ref{notation_relative})
    (\ref{definition_tj}), and (\ref{definition_overline_rm}), respectively.
    \item $a_{\varepsilon} \overset{\varepsilon \to 0^{+}}{\approx} b_{\varepsilon}$ is just a formal notation, which means that we believe that $a_{\varepsilon}$ is closed to $b_{\varepsilon}$ as $\epsilon \to 0^{+}$, and will not be used in the proofs of this article.
    \item For notational simplicity, we will use $c$, $C$, $C'$ to denote generic constants, whose values may change from place to place.
    \item The notations $a\lesssim b$, $a\gtrsim b$, and $a\approx b$ are defined by $a \leq C b$, $C a \geq b$, and $C a \leq b \leq C'a$, respectively. Moreover, $a\oplus b$ is given as follows for each positive numbers $a,b$:
    \begin{align} \label{definition_oplus}
        \frac{1}{a\oplus b} := \frac{1}{a}+\frac{1}{b}.
    \end{align}
\end{itemize}
\normalsize
\section{Main results} \label{section_main}
Now, we are ready to introduce our main results. The remainder of this section is organized in the following way. To begin with, in Section \ref{section_main_part_1}, we will present the results which are related to the asymptotics of the higher moments of the mollified SHE and the estimates for the high-order coupling constants. Next, we will introduce the sub-limiting higher moments and establish their unboundedness in Section \ref{section_main_part2_pseudo_limiting_moments}. Third, the applications related to the continuous directed polymer will be presented in Section \ref{section_CDP}. Finally, we will discuss some related literature in Section \ref{section_related_literature}. In particular, we will explain our main results from a quantum physics perspective in this subsection.

Before we proceed, it has to be stressed that all the results in this article hold under the assumptions about $\phi(x)$ mentioned in Section \ref{section_the_model}. To be more specific, \emph{we assume that $\phi \in C_{c}^{\infty}(\mathbb{R}^{d};\mathbb{R}_{+})$ and $\phi$ is symmetric-decreasing.} This is because our proof holds when $\phi$ is invariant under rearrangement.

\subsection{Estimates for higher moments and their critical regimes.} \label{section_main_part_1}
Let us now introduce the results mentioned in the part I of Section \ref{section_overview_main_result}. To begin with, we present the divergence of the higher moment of the mollified SHE inside the $L^{2}$-regime. In particular, in order to show that the above divergence holds for all non-negative and bounded initial datum, including the case of the partition function of the continuous directed polymer, we suppose that initial datum of the mollified SHE (\ref{Mollified_SHE}) is non-negative and smooth, with compact support. In this way, the above divergence can be deduced immediately from the following estimates:
\begin{theorem} \label{Main_result_1}
Assume that $d\geq 3$. Recall that the $N$-th moment $\mathcal{Q}_{\varepsilon;t}^{\beta,N}U_{0}^{\otimes N}(\Vec{x})$ is defined in (\ref{definition_moment}), $U_{0}$ is the initial datum of the mollified SHE, and the critical coupling constant $\beta_{L^{N}}$ is defined in (\ref{definition_beta_L_N}). Then, there exist a positive integer $N_{0} \geq 3$ and a sequence of positive numbers $(\beta_{N,+})_{N\geq 2}$ such that 
\begin{align} \label{estimate_betaLN_less_betaL2}
    \beta_{L^{N}}\leq \beta_{N,+} < \beta_{L^{2}} \quad \forall N\geq N_{0}
\end{align}
and
\begin{align} \label{identity_beta2+_=_betaL2}
    \beta_{2,+} = \beta_{L^{2}}.
\end{align}
In particular, if $N\geq 2$ and $\beta > \beta_{N,+}$, then, given an arbitrary $t>0$, we have the following estimates. 
\begin{enumerate} [label=(\roman*)]
    \item Let $\Vec{x} \in \mathbb{R}^{N \times d}$ and $U_{0} \in C_{c}^{\infty}(\mathbb{R}^{d};\mathbb{R}_{+})$. Then
    \begin{align} \label{inequality_main_1_1}
        \mathcal{Q}_{\varepsilon;t}^{\beta,N}U_{0}^{\otimes N}(\Vec{x}) \geq C \cdot \set{
        \varepsilon^{N \times d - 4} \cdot \exp(a\cdot \varepsilon^{-2})}
        \quad\forall \varepsilon < \varepsilon_{0}.
    \end{align}
    \item Let $\varphi\in C_{c}^{\infty}(\mathbb{R}^{N \times d};\mathbb{R}_{+})$ and $U_{0} \in C_{c}^{\infty}(\mathbb{R}^{d};\mathbb{R}_{+})$. Then
    \begin{align}\label{inequality_main_1_2}
        \Bigl\langle \varphi, \mathcal{Q}_{\varepsilon;t}^{\beta,N}U_{0}^{\otimes N}\Bigr\rangle_{L^{2}(\mathbb{R}^{N \times d})}\geq C'\cdot \set{
        \varepsilon^{N\times d - 4} \cdot \exp(a\cdot \varepsilon^{-2}) }\quad\forall \varepsilon < \varepsilon_{0}.
    \end{align}
\end{enumerate}
Here, $C,C',a,$ and $\varepsilon_{0}$ are positive constants such that they all depend on $\beta,d,N,\phi,$ and $t$, where $\phi(x)$ comes from the approximation of the identity defined in Section \ref{section_the_model}. Moreover, $C$ and $C'$ also relate to $(\Vec{x},U_{0})$ and $(\varphi,U_{0})$, respectively.
\end{theorem}

We next introduce a more explicit partial result for Conjecture \ref{Conjecture_order} than (\ref{estimate_betaLN_less_betaL2}). Recalling the definition (\ref{definition_beta_L_N}), we see that $\beta_{L^{N+1}} \leq \beta_{L^{N}}$ for all $N\geq 2$. Although we are not able to prove $\beta_{L^{N+1}} < \beta_{L^{N}}$, we provide a proper estimate for $\beta_{L^{N}}$, which shows that $\beta_{L^{N_{2}}} < \beta_{L^{N_{1}}}$ if $N_{2}$ is much larger than $N_{1}$.

\begin{theorem} \label{Main_result_3}
Assume that $d\geq 3$. Recall that the critical coupling constant $\beta_{L^{N}}$ is defined in (\ref{definition_beta_L_N}). Then, $\beta_{L^{N}}$ satisfies the following estimate:
\begin{align} \label{estimate_main_for_betaLN}
    \frac{\beta_{L^{2}}}{\sqrt{N-1}} \leq \beta_{L^{N}} \leq \frac{\alpha_{\infty,+}}{\sqrt{N-1}} \quad \forall N\geq 3.
\end{align}
Here, $\alpha_{\infty,+}$ is a positive constant that only depends on $d$ and $\phi$, where $\phi(x)$ originates from the approximation of the identity defined in Section \ref{section_the_model}.
\end{theorem}

\begin{remark} \label{total_mass}
Here, we would like to add some comments about the above results and some related works.
\begin{enumerate} [label=(\roman*)]
    \item Theorem \ref{Main_result_1} may not imply the divergences for all the higher moments, but this problem can be corrected in three dimensions in the next subsection.
    \item We expect the divergence of the higher moment of the mollified SHE can be further extended to $\beta>\beta_{L^{N}}$. In this way, one may wonder that for the case of $N\geq 3$ and $d\geq 3$, what is the limit of $\mathcal{Q}_{\varepsilon;t}^{\beta,N}U_{0}^{\otimes N}(\Vec{x})$ when the system is at the $L^{N}$-criticality. We note that if $\Vec{x} = \Vec{0}$, then the problem may be connected with \cite[Corollary 2.12]{junk2024taildistributionfunctionpartition} and the discussion below this corollary. 
    \item Recall that $\phi$ comes from the approximation of the identity defined in Section \ref{section_the_model}. Although we suppose that $\int dx \phi(x) = 1$, it is clear from the proofs that all the results in this subsection hold without this assumption. As a consequence, Theorem \ref{Main_result_3} can be further extended to the higher fractional moments. More details will be introduced in Section \ref{section_CDP}.
\end{enumerate}
\end{remark}

\subsection{The sub-limiting higher moments at $L^{2}$-criticality for the 3D-SHE.}  \label{section_main_part2_pseudo_limiting_moments}
Let us now introduce the results explained in the part II of Section \ref{section_overview_main_result}. The overall goal of this section is to illustrate the unboundedness of the sub-limiting higher moments of the three-dimensional mollified SHE at $L^{2}$-criticality. To this aim, we first clarify the concept of the sub-limiting higher moments in Section \ref{section_definition_sub-limiting_moment}. Then, the main results will be presented in Section \ref{section_unboundedness_of_the_sub-limiting_moments}. Throughout this subsection, we assume that $d = 3$ and $\beta = \beta_{L^{2}}$. Then, for the sake of simplicity, the dependences of the notations on $\beta_{L^{2}}$ will be ignored in this subsection.

\subsubsection{Definition and interpretations.} \label{section_definition_sub-limiting_moment}
The prime objective of this part is to introduce the concept of the sub-limiting higher moment of the three-dimensional mollified SHE. Briefly speaking, it is a natural three-dimensional counterpart of the limiting higher moment of the two-dimensional mollified SHE. See \cite[(2.19)]{chen} for the two-dimensional case.

\paragraph{Series expansions of higher moments.} To explain the definition, the prime objective of this part is to derive a series expansion of the $N$-th moment $\mathcal{Q}_{\varepsilon;t}^{\beta,N}U_{0}^{\otimes N}(\Vec{x})$ with $U_{0} \in C_{b}(\mathbb{R}^{d};\mathbb{R}_{+})$, which is analogous to \cite[Lemma 4.1]{chen}. 

Due to the above reason, we first briefly interpret the series expansion \cite[(2.12)]{chen} in our framework. To do this, we note that the $N$-th moment $\mathcal{Q}_{\varepsilon;t}^{\beta,N}U_{0}^{\otimes N}(\Vec{x})$ admits the following Feynman-Kac representation, which can obtained by using the Feynman-Kac formula (\ref{representation_FK}) of $\mathscr{U}_{\varepsilon}(x,t)$:
\begin{align} \label{representation_FK_0}
    \mathcal{Q}_{\varepsilon;t}^{\beta,N}U_{0}^{\otimes N}(\Vec{x}) = \mathbf{E}_{\Vec{x}}^{\Vec{B}}\bigbkt{\exp\bkt{\beta_{\varepsilon}^{2}\int_{0}^{t} \sum_{1\leq i<i'\leq N} R_{\varepsilon}(\Vec{B}_{s}(i') - \Vec{B}_{s}(i)) ds }U_{0}^{\otimes N}(\Vec{B}(t))} \quad \forall \Vec{x} \in \mathbb{R}^{N \times d}.
\end{align}
Recalling the series expansion \cite[(2.6)]{chen},  $\mathcal{Q}_{\varepsilon;t}^{\beta;N}U_{0}^{\otimes N}(\Vec{x}_{0})$ can be decomposed as follows: 
\begin{align} \label{decomposition_first}
    &\mathcal{Q}_{\varepsilon;t}^{\beta;N}U_{0}^{\otimes N}(\Vec{x}_{0})
    = \mathbf{E}_{\Vec{x}_{0}}^{\Vec{B}}[U_{0}^{\otimes N}(\Vec{B}(t))]  \nonumber\\
    &+\sum_{m=1}^{\infty} \sum_{(\ell_{1},\ell'_{1}),...,(\ell_{m},\ell'_{m}) \in \mathcal{E}^{(N)}, \; (\ell_{1},\ell_{1}') \neq ...\neq (\ell_{m},\ell_{m}')} \int_{0<s_{1}<...<s_{m}<t} \prod_{j=1}^{m}ds_{j} \cdot \mathscr{L}^{N;(\ell_{1},\ell_{1}'),...,(\ell_{m},\ell_{m}')}_{\varepsilon;s_{1},s_{2},...,s_{m},t}U_{0}^{\otimes N}(\Vec{x}_{0}),
\end{align}
where the set $\mathcal{E}^{(N)}$ is defined in Section \ref{section_structure}. Here, for each $f\in C_{b}(\mathbb{R}^{N \times d};\mathbb{R}_{+})$, the path integral on the right-hand side of (\ref{decomposition_first}) is given by
\begin{align} \label{decomposition_first_path_integral}
    &\mathscr{L}^{N;(\ell_{1},\ell'_{1}),...,(\ell_{m},\ell'_{m})}_{\varepsilon;s_{1},s_{2},...,s_{m},t}f(\Vec{x}_{0})
     :=\nonumber\\
     &\mathbf{E}_{\Vec{x}_{0}}^{\Vec{B}}\biggl[ \prod_{j=1}^{m}\beta_{\varepsilon}^{2}R_{\varepsilon}(\Vec{B}_{s_{j}}(\ell_{j}')-\Vec{B}_{s_{j}}(\ell_{j})) \exp\bkt{\beta_{\varepsilon}^{2}\int_{s_{j}}^{s_{j+1}} R_{\varepsilon}(\Vec{B}_{u}(\ell_{j}')-\Vec{B}_{u}(\ell_{j}))du} \cdot    
     f(\Vec{B}(t))
     \biggr], \text{ where } s_{m+1} = t.
\end{align}
See Remark \ref{remak_series_expansion} for an analytic explanation of (\ref{decomposition_first}), and see \cite[Proposition 2.1]{chen} for the original probabilistic proof of (\ref{decomposition_first}). 

Next, we further decompose the path integral $\mathscr{L}^{N;(\ell_{1},\ell_{1}'),...,(\ell_{m},\ell_{m}')}_{\varepsilon;s_{1},s_{2},...,s_{m},t}f(\Vec{x}_{0})$. In view of the following fundamental theorem of calculus:
\begin{align*}
    &\exp\bkt{\beta_{\varepsilon}^{2}\int_{s_{j}}^{s_{j+1}} R_{\varepsilon}(\Vec{B}_{u}(\ell_{j}')-\Vec{B}_{u}(\ell_{j}))du}\\
    &= 1 + \int_{s_{j}}^{s_{j+1}} du_{j} \exp\bkt{\beta_{\varepsilon}^{2}\int_{s_{j}}^{u_{j}} R_{\varepsilon}(\Vec{B}_{v}(\ell_{j}')-\Vec{B}_{v}(\ell_{j}))dv} \beta_{\varepsilon}^{2}R_{\varepsilon}(\Vec{B}_{u_{j}}(\ell_{j}')-\Vec{B}_{u_{j}}(\ell_{j})),
\end{align*}
the path integral $\mathscr{L}^{N;(\ell_{1},\ell_{1}'),...,(\ell_{m},\ell_{m}')}_{\varepsilon;s_{1},s_{2},...,s_{m},t}f(\Vec{x}_{0})$ admits the following expansion:
\begin{align*}
    \mathscr{L}^{N;(\ell_{1},\ell_{1}'),...,(\ell_{m},\ell_{m}')}_{\varepsilon;s_{1},s_{2},...,s_{m},t}f(\Vec{x}_{0})
    = \sum_{\mathfrak{i}_{1},...,\mathfrak{i}_{m} \in \mathcal{E}^{(N)} \times \{0,1\}} \mathscr{L}^{N;(\ell_{1},\ell_{1}',\mathfrak{i}_{1}),...,(\ell_{m},\ell'_{m},\mathfrak{i}_{m})}_{\varepsilon;s_{1},s_{2},...,s_{m},t}f(\Vec{x}_{0}),
\end{align*}
where the path integral $\mathscr{L}^{N;(\ell_{1},\ell_{1}',\mathfrak{i}_{1}),...,(\ell_{m},\ell'_{m},\mathfrak{i}_{m})}_{\varepsilon;s_{1},s_{2},...,s_{m},t}f(\Vec{x}_{0})$ is given by
\begin{align} \label{definition_LN}
    &\mathscr{L}^{N;(\ell_{1},\ell_{1}',\mathfrak{i}_{1}),...,(\ell_{m},\ell'_{m},\mathfrak{i}_{m})}_{\varepsilon;s_{1},s_{2},...,s_{m},t}f(\Vec{x}_{0})
    := \mathbf{E}_{\Vec{x}_{0}}^{\Vec{B}}\biggl[\prod_{j=1}^{m}\biggbkt{\beta_{\varepsilon}^{2}R_{\varepsilon}(\Vec{B}_{s_{j}}(\ell'_{j})-\Vec{B}_{s_{j}}(\ell_{j}))\nonumber\\
    &\quad\quad\quad\cdot\biggbkt{\int_{s_{j}}^{s_{j+1}} du_{j} \exp\bkt{\beta_{\varepsilon}^{2}\int_{s_{j}}^{u_{j}} R_{\varepsilon}(\Vec{B}_{v}(\ell'_{j})-\Vec{B}_{v}(\ell_{j}))dv} \beta_{\varepsilon}^{2}R_{\varepsilon}(\Vec{B}_{u_{j}}(\ell'_{j})-\Vec{B}_{u_{j}}(\ell_{j}))}^{\mathfrak{i}_{j}}}\cdot f(\Vec{B}(t)) \biggr], \nonumber\\
    &\quad\quad\quad\quad\quad\quad\quad\quad\quad\quad\quad\quad\quad\quad\quad\quad\quad\quad\quad\quad\quad\quad\quad\quad\quad\quad\quad\quad\quad\quad\quad\quad\quad\text{where } s_{m+1} = t.
\end{align}
With the above notations, we conclude the following series expansion of the $N$-th moment of the mollified SHE:
\begin{align} \label{decomposition_second}
    \mathcal{Q}_{\varepsilon;t}^{\beta;N}U_{0}^{\otimes N}&(\Vec{x}_{0})
    = \mathbf{E}_{\Vec{x}_{0}}^{\Vec{B}}[U_{0}^{\otimes N}(\Vec{B}(t))]  \nonumber\\
    &+\sum_{m=1}^{\infty} \sum_{(\ell_{1},\ell'_{1},\mathfrak{i}_{1}),...,(\ell_{m},\ell'_{m},\mathfrak{i}_{m}) \in \mathcal{E}^{(N)} \times \{0,1\}, \; (\ell_{1},\ell_{1}') \neq ...\neq (\ell_{m},\ell_{m}')} \mathscr{I}_{\varepsilon;t}^{N;(\ell_{1},\ell'_{1},\mathfrak{i}_{1}),...,(\ell_{m},\ell'_{m},\mathfrak{i}_{m})}U^{\otimes N}_{0}(\Vec{x}_{0}).
\end{align}
Here, the path integral on the right-hand side of (\ref{decomposition_second}) is given by 
\begin{align}
    &\mathscr{I}_{\varepsilon;t}^{N;(\ell_{1},\ell'_{1},\mathfrak{i}_{1}),...,(\ell_{m},\ell'_{m},\mathfrak{i}_{m})}f(\Vec{x}_{0})\nonumber\\
    &\quad\quad\quad\quad\quad:= 
    \int_{0<s_{1}<...<s_{m}<t} \prod_{j=1}^{m}ds_{j} \cdot\mathscr{L}^{N;(\ell_{1},\ell_{1}',\mathfrak{i}_{1}),...,(\ell_{m},\ell_{m}',\mathfrak{i}_{m})}_{\varepsilon;s_{1},s_{2},...,s_{m},t}f(\Vec{x}_{0})\quad \forall f \in C_{b}(\mathbb{R}^{N \times d};\mathbb{R}_{+}),
\end{align}
where the path integral $\mathscr{L}^{N;(\ell_{1},\ell_{1}',\mathfrak{i}_{1}),...,(\ell_{m},\ell_{m}',\mathfrak{i}_{m})}_{\varepsilon;s_{1},s_{2},...,s_{m},t}f(\Vec{x}_{0})$ is defined in (\ref{definition_LN}).

Third, we derive an integral representation for the path integral $\mathscr{I}_{\varepsilon;t}^{N;(\ell_{1},\ell'_{1},\mathfrak{i}_{1}),...,(\ell_{m},\ell'_{m},\mathfrak{i}_{m})}f(\Vec{x}_{0})$. Due to this purpose, applying the Markov property of Brownian motion $\Vec{B}$ to (\ref{definition_LN}), we then obtain the following expression:
\begin{align} \label{identity_application_markov_property}
    &\mathscr{I}_{\varepsilon;t}^{N;(\ell_{1},\ell'_{1},\mathfrak{i}_{1}),...,(\ell_{m},\ell'_{m},\mathfrak{i}_{m})}f(\Vec{x}_{0}) = \int_{u_{0} = 0<s_{1}<u_{1}<s_{2}...<s_{m}<u_{m}<t} \prod_{j=1}^{m}ds_{j}du_{j} \cdot \int_{\mathbb{R}^{(m+1) \times N \times d}}  \prod_{j=1}^{m+1} d\Vec{x}_{j} \cdot \nonumber\\
    &\biggbkt{\prod_{j=1}^{m} \delta_{0}(s_{j}-u_{j})^{1-\mathfrak{i}_{j}}
    \cdot \biggbkt{\prod_{1\leq i\leq N,\; i\not \in \{\ell_{j},\ell'_{j}\}} G_{u_{j}-u_{j-1}}(\Vec{x}_{j}(i) - \Vec{x}_{j-1}(i))} 
    \cdot \mathscr{A}^{(\ell_{j},\ell_{j}',\mathfrak{i}_{j})}_{\varepsilon;u_{j-1},s_{j},u_{j}}(\Vec{x}_{j-1},\Vec{x}_{j}) 
    }\nonumber\\
    &\quad\quad\quad\quad\quad\quad\quad\quad\quad\quad\quad\quad\quad\quad\quad\quad\quad\quad\quad\quad\quad\quad\quad\quad\cdot G^{(N)}_{t-u_{m}}(\Vec{x}_{m+1}-\Vec{x}_{m}) f(\Vec{x}_{m+1}),
\end{align}
where the transition kernel
$\mathscr{A}^{(\ell_{j},\ell_{j}',\mathfrak{i}_{j})}_{\varepsilon;u_{j-1},s_{j},u_{j}}(\Vec{x}_{j-1},\Vec{x}_{j})$ is defined as follows:
\begin{align} \label{definition_epsilon_motion}
    &\mathscr{A}^{(\ell_{j},\ell_{j}',\mathfrak{i}_{j})}_{\varepsilon;u_{j-1},s_{j},u_{j}}(\Vec{x}_{j-1},\Vec{x}_{j})\nonumber\\
    &:= \int_{\mathbb{R}^{2 \times d}} dy_{j}dy_{j}' G_{s_{j}-u_{j-1}}(\Vec{x}_{j-1}(\ell_{j})-y_{j})G_{s_{j}-u_{j-1}}(\Vec{x}_{j-1}(\ell'_{j})-y_{j}') \beta_{\varepsilon}^{2} R_{\varepsilon}(y_{j}'-y_{j}) \cdot \nonumber\\
    &\mathbf{E}_{y_{j},y_{j}'}^{B,B'}\bigbkt{\exp\bkt{\beta_{\varepsilon}^{2}\int_{0}^{u_{j}-s_{j}}R_{\varepsilon}(B'(v) - B(v))dv}: (B(u_{j}-s_{j}),B'(u_{j}-s_{j})) = (\Vec{x}_{j}(\ell_{j}),\Vec{x}_{j}(\ell'_{j}))} \nonumber\\
    &\quad\quad\quad\quad\quad\quad\quad\quad\quad\quad\quad\quad\quad\quad\quad\quad\quad\quad\quad\quad\quad\quad\cdot\beta_{\varepsilon}^{2} R_{\varepsilon}(\Vec{x}_{j}(\ell_{j}')-\Vec{x}_{j}(\ell_{j})) \quad \text{if } \mathfrak{i}_{j} = 1;\nonumber\\
    &\mathscr{A}^{(\ell_{j},\ell_{j}',\mathfrak{i}_{j})}_{\varepsilon;u_{j-1},s_{j},u_{j}}(\Vec{x}_{j-1},\Vec{x}_{j})\nonumber\\
    &:= 
    G_{u_{j}-u_{j-1}}(\Vec{x}_{j-1}(\ell_{j})-\Vec{x}_{j}(\ell_{j}))G_{u_{j}-u_{j-1}}(\Vec{x}_{j-1}(\ell'_{j})-\Vec{x}_{j}(\ell'_{j})) \beta_{\varepsilon}^{2} R_{\varepsilon}(\Vec{x}_{j}(\ell'_{j})-\Vec{x}_{j}(\ell_{j})) \quad\text{if } \mathfrak{i}_{j} = 0.
\end{align}
Here, the measure $\mathbf{P}_{y,y'}^{B,B'}\sbkt{\cdot: (B(t),B'(t)) = (x,x')}$ is defined in Section \ref{section_structure}. In other words, $\mathscr{A}^{(\ell_{j},\ell_{j}',\mathfrak{i}_{j})}_{\varepsilon;u_{j-1},s_{j},u_{j}}(\Vec{x}_{j-1},\Vec{x}_{j})$ describes the motion of the $\varepsilon$-range interacting system within the time interval $(u_{j-1},u_{j})$, where the system consists of Brownian motions $\Vec{B}(\ell_{j})$ and $\Vec{B}(\ell'_{j})$. Simply put, if $\mathfrak{i}_{j} = 1$, then Brownian motions $\Vec{B}(\ell_{j})$ and $\Vec{B}(\ell'_{j})$ go through the $\varepsilon$-range interaction below within the time interval $(s_{j},u_{j})$:
\begin{align} \label{epsilon_interaction}
    &\beta_{\varepsilon}^{2} R_{\varepsilon}(y_{j}'-y_{j}) \cdot \mathbf{E}_{y_{j},y_{j}'}^{B,B'}\bigbkt{\exp\bkt{\beta_{\varepsilon}^{2}\int_{0}^{u_{j}-s_{j}}R_{\varepsilon}(B'(v) - B(v))dv}:\nonumber\\
    &\quad\quad\quad\quad\quad\quad\quad\quad(B(u_{j}-s_{j}),B'(u_{j}-s_{j})) = (\Vec{x}_{j}(\ell_{j}),\Vec{x}_{j}(\ell'_{j}))} \cdot \beta_{\varepsilon}^{2} R_{\varepsilon}(\Vec{x}_{j}(\ell'_{j})-\Vec{x}_{j}(\ell_{j})).
\end{align}
On the other hand, if $\mathfrak{i}_{j} = 0$, then Brownian motions $\Vec{B}(\ell_{j})$ and $\Vec{B}(\ell'_{j})$ only meet each other within an $\varepsilon$-range at time $u_{j}$.

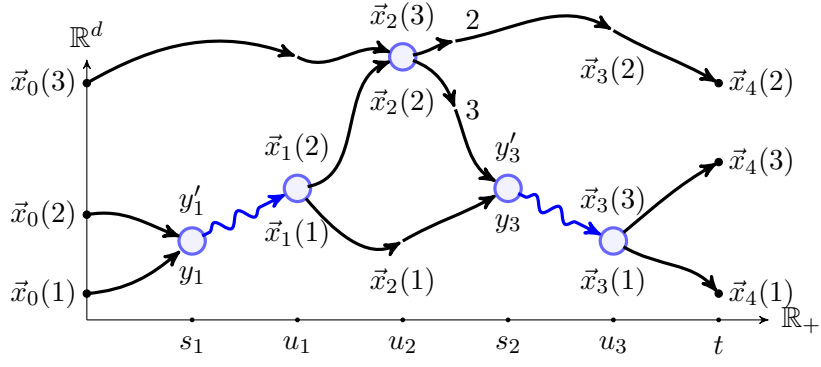
\begin{figure}
\centering
\begin{tikzpicture} [
        scale = 0.7,
        ->,    > = stealth',
       shorten > = 1pt,auto,
   node distance = 3cm,
      decoration = {snake,   % <-- added
                    pre length=3pt,post length=7pt,% <-- for better looking of arrow,
                    },
main node/.style = {circle,draw,fill=#1,
                    font=\sffamily\Large\bfseries},
                    ]

\path[black] (0,-0.12) node[anchor=north]{$s_{1}$};
\path[black] (2,-0.12) node[anchor=north]{$u_{1}$};
\path[black] (4,-0.12) node[anchor=north]{$u_{2}$};
\path[black] (6,-0.12) node[anchor=north]{$s_{2}$};
\path[black] (8,-0.12) node[anchor=north]{$u_{3}$};
\path[black] (10,-0.12) node[anchor=north]{$t$};

\filldraw[black] (0,0) circle (1pt);
\filldraw[black] (2,0) circle (1pt);
\filldraw[black] (4,0) circle (1pt);
\filldraw[black] (6,0) circle (1pt);
\filldraw[black] (8,0) circle (1pt);
\filldraw[black] (10,0) circle (1pt);

%\filldraw[black] (6,4) circle (8pt);

\draw[->] (-2,0) -- (11,0) node[anchor=west]{$\mathbb{R}_{+}$};
\draw[->] (-2,0) -- (-2,5) node[anchor=south]{$\mathbb{R}^{d}$};

\filldraw[black] (-2,0.5) circle (2pt) node[anchor=east]{$\Vec{x}_{0}(1)$};
\filldraw[black] (-2,2) circle (2pt)
node[anchor=east]{$\Vec{x}_{0}(2)$};
\filldraw[black] (-2,4.5) circle (2pt)
node[anchor=east]{$\Vec{x}_{0}(3)$};

\node (A) at (0,1.5) {};
\filldraw[color=blue!60, fill=blue!5, very thick](0,1.5) circle (0.25);
\node (B) at (2,2.5) {};
\filldraw[color=blue!60, fill=blue!5, very thick](2,2.5) circle (0.25);
\node (C) at (4,5) {};
\filldraw[color=blue!60, fill=blue!5, very thick](4,5) circle (0.25);
\node (D) at (6,2.5) {};
\filldraw[color=blue!60, fill=blue!5, very thick](6,2.5) circle (0.25);
\node (E) at (8,1.5) {};
\filldraw[color=blue!60, fill=blue!5, very thick](8,1.5) circle (0.25);

\draw[->] [very thick] (-2,0.5) to[out=5,in=-140] (A);
\draw[->] [very thick] (-2,2) to[out=10,in=160] (A);

%\draw[->] [very thick,blue] (A) to[out=10,in=-160] (B);
%\draw[->] [very thick,blue] (A) to[out=-5,in=-140] (B);
\path[draw=blue, very thick, decorate] (A) -- (B);
\path[draw=blue, very thick, decorate] (D) -- (E);

\draw[->] [very thick] (-2,4.5) to[out=30,in=160] (2,5);
\draw[->] [very thick] (2,5) to[out=-30,in=160] (C);
\draw[->] [very thick] (B) to[out=10,in=-160] (C);
\draw[->] [very thick] (B) to[out=-50,in=-150] (4,1.5);
\draw[->] [very thick] (4,1.5) to[out=30,in=-150] (D);
\draw[->] [very thick] (C) to[out=-20,in=110] (5,4);
\draw[->] [very thick] (5,4) to[out=-70,in=150] (D);
\draw[->] [very thick] (C) to[out=10,in=-170] (5,5.3);
\draw[->] [very thick] (5,5.3) to[out=10,in=150] (8,5.5);
\draw[->] [very thick] (8,5.5) to[out=-30,in=150] (10,4.5);
\draw[->] [very thick] (E) to[out=40,in=-150] (10,3);
\draw[->] [very thick] (E) to[out=-30,in=130] (10,0.5);

\filldraw[black] (10,0.5) circle (2pt) node[anchor=west]{$\Vec{x}_{4}(1)$};
\filldraw[black] (10,3) circle (2pt)
node[anchor=west]{$\Vec{x}_{4}(3)$};
\filldraw[black] (10,4.5) circle (2pt)
node[anchor=west]{$\Vec{x}_{4}(2)$};

\node at (0,1.2)[anchor=north]{$y_{1}$};
\node at (0,1.8)[anchor=south]{$y'_{1}$};

\node at (2,2.1)[anchor=north]{$\Vec{x}_{1}(1)$};
\node at (2,2.8)[anchor=south]{$\Vec{x}_{1}(2)$};

\node at (4,1.2)[anchor=north]{$\Vec{x}_{2}(1)$};
\node at (4,4.6)[anchor=north]{$\Vec{x}_{2}(2)$};
\node at (4,5.3)[anchor=south]{$\Vec{x}_{2}(3)$};

\node at (6,2.8)[anchor=south]{$y'_{3}$};
\node at (6,2.2)[anchor=north]{$y_{3}$};

\node at (5,5.7)[anchor=west]{$2$};
\node at (5,4)[anchor=west]{$3$};

\node at (8,1.2)[anchor=north]{$\Vec{x}_{3}(1)$};
\node at (8,1.8)[anchor=south]{$\Vec{x}_{3}(3)$};
\node at (8,5.2)[anchor=north]{$\Vec{x}_{3}(2)$};

\end{tikzpicture}
\caption{The above graph denotes the path integral $\mathscr{I}_{\varepsilon;t}^{N;(\ell_{1},\ell'_{1},\mathfrak{i}_{1}),...,(\ell_{m},\ell'_{m},\mathfrak{i}_{m})}U^{\otimes N}_{0}(\Vec{x}_{0})$ with $N = 3$, $m = 3$, $(\ell_{1},\ell'_{1},\mathfrak{i}_{1}) = (1,2,1)$, $(\ell_{2},\ell'_{2},\mathfrak{i}_{2}) = (2,3,0)$, and $(\ell_{3},\ell'_{3},\mathfrak{i}_{3}) = (1,3,1)$. Here, each of the above blue region denotes the $\varepsilon$-regions induced by $R_{\varepsilon}$, and each of the above blue wave line denotes that $\Vec{B}(\ell_{j}')$ and $\Vec{B}(\ell_{j})$ go through the $\varepsilon$-range interaction (\ref{epsilon_interaction}). Also, the variable $y'_{j}$ always relates to the larger index $\ell'_{j}$.}
\label{BD_model_3}
\end{figure}

Finally, we summarize the above results and notations by a graphical explanation. Briefly speaking, $\mathcal{Q}_{\varepsilon;t}^{\beta;N}U_{0}^{\otimes N}(\Vec{x}_{0})$ is a series of the form (\ref{decomposition_second}), where each term $\mathscr{I}_{\varepsilon;t}^{N;(\ell_{1},\ell'_{1},\mathfrak{i}_{1}),...,(\ell_{m},\ell'_{m},\mathfrak{i}_{m})}U^{\otimes N}_{0}(\Vec{x}_{0})$ in this series can be regarded as a graph, and Figure \ref{BD_model_3} is an example of this kind of graph.

\begin{remark} \label{remak_series_expansion}
Here, we explain the decomposition (\ref{decomposition_first}). Recalling the Feynman-Kac representation (\ref{representation_FK_0}). As the proof of Lemma \ref{lemma_series_expansion}, we apply the fundamental theorem of calculus (\ref{identity_fundamental}) to the exponential function below:
\begin{align*}
    \exp\bkt{\beta_{\varepsilon}^{2}\int_{0}^{t} \sum_{1\leq i<i'\leq N} R_{\varepsilon}(\Vec{B}_{s}(i') - \Vec{B}_{s}(i)) ds }.
\end{align*}
Then, one can conclude the following series expansion:
\begin{align} \label{decomposition_direct}
    &\mathcal{Q}_{\varepsilon;t}^{\beta;N}U_{0}^{\otimes N}(\Vec{x}_{0})
    = \mathbf{E}_{\Vec{x}_{0}}^{\Vec{B}}[U_{0}^{\otimes N}(\Vec{B}(t))] 
    + \sum_{m=1}^{\infty} \sum_{(\ell_{1},\ell'_{1}),...,(\ell_{m},\ell_{m}') \in \mathcal{E}^{(N)}} \int_{0<v_{1}<...<v_{m}<t} \prod_{j=1}^{m}dv_{j}  \nonumber \\
    &\mathbf{E}_{\Vec{x}_{0}}^{\Vec{B}}\bigbkt{\beta_{\varepsilon}^{2}R_{\varepsilon}(\Vec{B}_{v_{1}}(\ell_{1}')-\Vec{B}_{v_{1}}(\ell_{1}))
     ...\beta_{\varepsilon}^{2}R_{\varepsilon}(\Vec{B}_{v_{m}}(\ell_{m}')-\Vec{B}_{v_{m}}(\ell_{m})) 
     U_{0}^{\otimes N}(\Vec{B}(t))
     }.
\end{align}
To combine the above consecutive interactions, we set $s_{0} = 0$ and $s_{i} := v_{\pi(i)}$, where $\pi(i)$ is defined by
\begin{align*}
    \pi(i) := \min\set{j: j>\pi(i-1) \text{ and } (\ell_{j},\ell'_{j}) \neq (\ell_{\pi(i-1)},\ell'_{\pi(i-1)})} \quad \text{if } i\geq 2; \quad \pi(1) = 1.
\end{align*}
Note that $s_{i}$ denotes the starting time of the $i$-th consecutive interaction on the right-hand side of (\ref{decomposition_direct}).
Then, starting with the last interaction on the right-hand side of (\ref{decomposition_direct}), we can apply the following rule to combine the consecutive interactions:
\begin{align*}
    &\beta_{\varepsilon}^{2}R_{\varepsilon}(\Vec{B}_{s_{j}}(\ell'_{j})-\Vec{B}_{s_{j}}(\ell_{j}))\exp\bkt{\beta_{\varepsilon}^{2}\int_{s_{j}}^{s_{j+1}} R_{\varepsilon}(\Vec{B}_{u}(\ell'_{j})-\Vec{B}_{u}(\ell_{j}))du}\\
    &= \beta_{\varepsilon}^{2}R_{\varepsilon}(\Vec{B}_{s_{j}}(\ell'_{j})-\Vec{B}_{s_{j}}(\ell_{j})) \cdot \biggbkt{1+\sum_{m=1}^{\infty} \frac{1}{m!} \bkt{\beta_{\varepsilon}^{2}\int_{s_{j}}^{s_{j+1}} R_{\varepsilon}(\Vec{B}_{u}(\ell'_{j})-\Vec{B}_{u}(\ell_{j}))du}^{m} }\\
    &= \beta_{\varepsilon}^{2}R_{\varepsilon}(\Vec{B}_{s_{j}}(\ell'_{j})-\Vec{B}_{s_{j}}(\ell_{j}))\cdot \biggbkt{1+\\
    &\sum_{m = 1}^{\infty} \int_{s_{j} < v_{1} < ... < v_{m} < s_{j+1}} \prod_{j=1}^{m}dv_{j} \cdot\beta_{\varepsilon}^{2} R_{\varepsilon}(\Vec{B}_{v_{1}}(\ell'_{j})-\Vec{B}_{v_{1}}(\ell_{j}))...\beta_{\varepsilon}^{2} R_{\varepsilon}(\Vec{B}_{v_{m}}(\ell'_{j})-\Vec{B}_{v_{m}}(\ell_{j})) }.
\end{align*}
In this way, one can conclude (\ref{decomposition_first}).
\end{remark}

\paragraph{The sub-limiting higher moments.} With the above notations at hand, we are ready to introduce the sub-limiting higher moment of the mollified SHE. Briefly speaking, it is of the form (\ref{decomposition_second}), but the path integral $\mathscr{I}_{\varepsilon;t}^{N;(\ell_{1},\ell'_{1},\mathfrak{i}_{1}),...,(\ell_{m},\ell'_{m},\mathfrak{i}_{m})}U^{\otimes N}_{0}(\Vec{x}_{0})$ is replaced by its limit. To illustrate this definition, a formal derivation of this definition will be presented in the next part. 
\begin{definition}
Assume that $d = 3$. Then, given an arbitrary $t>0$,  for all $N\geq 3$ and for every $U_{0} \in C_{b}(\mathbb{R}^{d};\mathbb{R}_{+})$, the sub-limit of $N$-th mixed moment of the mollified SHE is defined by
\begin{align} \label{definition_sub_limit}
    \mathscr{Q}_{t}^{N}&U_{0}^{\otimes N}(\Vec{x}_{0})
    = \mathbf{E}_{\Vec{x}_{0}}^{\Vec{B}}[U_{0}^{\otimes N}(\Vec{B}(t))]  \nonumber\\
    &+\lim_{M\to\infty}\sum_{m=1}^{M} \sum_{(\ell_{1},\ell'_{1},\mathfrak{i}_{1}),...,(\ell_{m},\ell'_{m},\mathfrak{i}_{m}) \in \mathcal{E}^{(N)} \times \{0,1\}, \; (\ell_{1},\ell_{1}') \neq ...\neq (\ell_{m},\ell_{m}')} \mathscr{I}_{0;t}^{N;(\ell_{1},\ell'_{1},\mathfrak{i}_{1}),...,(\ell_{m},\ell'_{m},\mathfrak{i}_{m})} U_{0}^{\otimes N} (\Vec{x}_{0}) \nonumber\\
    &\quad\quad\quad\quad\quad\quad\quad\quad\quad\quad\quad\quad\quad\quad\quad\quad\quad\quad\quad\quad\quad\quad\quad\quad\quad\quad\quad\quad\forall\Vec{x}_{0} \in \mathbb{R}^{N \times d}\setminus \Pi^{(N)}.
\end{align}
where the sub-limiting path integral is given by
\begin{align} \label{definition_sub_limiting_path_integral}
    &\mathscr{I}_{0;t}^{N;(\ell_{1},\ell'_{1},\mathfrak{i}_{1}),...,(\ell_{m},\ell'_{m},\mathfrak{i}_{m})}f(\Vec{x}_{0}):= \int_{u_{0} = 0<s_{1}<u_{1}<s_{2}...<s_{m}<u_{m}<t} \prod_{j=1}^{m}ds_{j}du_{j} \cdot \int_{\mathbb{R}^{(m+1) \times N \times d}}  \prod_{j=1}^{m+1} d\Vec{x}_{j}\nonumber\\
    &\biggbkt{\prod_{j=1}^{m}
    \biggbkt{\prod_{1\leq i\leq N,\; i\not \in \{\ell_{j},\ell'_{j}\}} G_{u_{j}-u_{j-1}}(\Vec{x}_{j}(i) - \Vec{x}_{j-1}(i))} 
    \cdot \mathscr{A}^{(\ell_{j},\ell_{j}',\mathfrak{i}_{j})}_{0;u_{j-1},s_{j},u_{j}}(\Vec{x}_{j-1},\Vec{x}_{j})} \nonumber\\
    &\quad\quad\quad\quad\quad\quad\quad\quad\cdot G^{(N)}_{t-u_{m}}(\Vec{x}_{m}-\Vec{x}_{m+1}) f(\Vec{x}_{m+1}) \quad \forall\Vec{x}_{0} \in\mathbb{R}^{N \times d}\setminus \Pi^{(N)}, \; f \in C_{b}(\mathbb{R}^{N \times d};\mathbb{R}_{+}).
\end{align}
Here, the transition kernel $\mathscr{A}^{(\ell_{j},\ell_{j}',\mathfrak{i}_{j})}_{0;u_{j-1},s_{j},u_{j}}(\Vec{x}_{j-1},\Vec{x}_{j})$ describes the motion of the contact interacting two-particles system within the time interval $(u_{j-1},u_{j})$, and it is defined as follows:
\begin{align} \label{definition_transition_kernel}
    &\mathscr{A}^{(\ell_{j},\ell_{j}',\mathfrak{i}_{j})}_{0;u_{j-1},s_{j},u_{j}}(\Vec{x}_{j-1},\Vec{x}_{j}) :=
    \int_{\mathbb{R}^{2 \times d}} dy_{j}dy_{j}' G_{s_{j}-u_{j-1}}(\Vec{x}_{j-1}(\ell_{j})-y_{j})
    G_{s_{j}-u_{j-1}}(\Vec{x}_{j-1}(\ell'_{j})-y_{j}') \cdot \nonumber\\
    &\delta_{0}\bkt{\frac{y'_{j}-y_{j}}{\sqrt{2}}}\cdot\biggbkt{\sqrt{\frac{2\pi}{u_{j}-s_{j}}} G_{u_{j}-s_{j}}\bkt{\frac{y'_{j}+y_{j}}{\sqrt{2}}-\frac{\Vec{x}_{j}(\ell'_{j}) + \Vec{x}_{j}(\ell_{j})}{\sqrt{2}}} }\cdot \delta_{0}\bkt{\frac{\Vec{x}_{j}(\ell'_{j}) - \Vec{x}_{j}(\ell_{j})}{\sqrt{2}}}
    \quad \text{if $\mathfrak{i}_{j} = 1$;}\nonumber\\
    &\mathscr{A}^{(\ell_{j},\ell_{j}',\mathfrak{i}_{j})}_{0;u_{j-1},s_{j},u_{j}}(\Vec{x}_{j-1},\Vec{x}_{j}) := 0 \quad \text{if $\mathfrak{i}_{j} = 0$.}
\end{align}
Note that the above Dirac’s delta functions are defined in the sense of tempered distribution:
\begin{align*}
    \int_{\mathbb{R}^{d}} dz \delta_{0}\bkt{\frac{x-z}{\sqrt{2}}} f(z) = 2^{d/2} f(x) \quad \forall f \in \mathscr{S}(\mathbb{R}^{d}).
\end{align*}
\end{definition}

\begin{remark} \label{remark_after_definition} Here, we would like to add some comments concerning the well-definednesses of (\ref{definition_sub_limit}) and (\ref{definition_sub_limiting_path_integral}).
\begin{enumerate} [label=(\roman*)]
    \item 
    Note that the limit defined in (\ref{definition_sub_limit}) exists in line with the fact that both of the initial datum $U_{0}$ and $R_{\varepsilon}$ are non-negative.
    \item Here, we illustrate the primary reason why the assumption $\Vec{x}_{0} \in \mathbb{R}^{N\times d}\setminus \Pi^{(N)}$ is required. If $U_{0} = G_{\nu}$ for some $\nu > 0$, then one can show that
    \begin{align*}
        \mathscr{I}_{0;t}^{N;(\ell_{1},\ell'_{1},1)} U_{0}^{\otimes N} (\Vec{x}_{0}) \approx \prod_{1\leq i\leq N,\; i\not \in \{\ell_{1},\ell'_{1}\}}& G_{t+\nu}(\Vec{x}_{0}(i)) \cdot G_{t+\nu}\bkt{\frac{\Vec{x}_{0}(\ell'_{1}) + \Vec{x}_{0}(\ell_{1})}{\sqrt{2}}} \\
        &\cdot\int_{0<s_{1}<u_{1}<t} ds_{1} du_{1}  \sqrt{\frac{1}{u_{1}-s_{1}}} G_{s_{1}}\bkt{\frac{\Vec{x}_{0}(\ell'_{1}) - \Vec{x}_{0}(\ell_{1})}{\sqrt{2}}}.
    \end{align*}
    Consequently, if $\Vec{x}_{0}(\ell_{1}') = \Vec{x}_{0}(\ell_{1})$, then the latter integral is divergent.
    \item Due to the delta functions on the right-hand side of  (\ref{definition_sub_limiting_path_integral}), one may question about the finiteness of the sub-limiting path integral $\mathscr{I}_{0;t}^{N;(\ell_{1},\ell'_{1},\mathfrak{i}_{1}),...,(\ell_{m},\ell'_{m},\mathfrak{i}_{m})}f(\Vec{x}_{0})$. But, for $N = 3$ and for a proper non-negative function $f(\Vec{x})$, it can be shown that, for almost every $\Vec{x}_{0} \in \mathbb{R}^{N \times d}$, \begin{align*}
        \mathscr{I}_{0;t}^{N;(\ell_{1},\ell'_{1},\mathfrak{i}_{1}),...,(\ell_{m},\ell'_{m},\mathfrak{i}_{m})}f(\Vec{x}_{0})<\infty \quad \forall m\geq 0.
    \end{align*} 
\end{enumerate}
\end{remark}

\begin{figure}
\centering
\begin{tikzpicture} [
        scale = 0.7,
        ->,    > = stealth',
       shorten > = 1pt,auto,
   node distance = 3cm,
      decoration = {snake,   % <-- added
                    pre length=3pt,post length=7pt,% <-- for better looking of arrow,
                    },
main node/.style = {circle,draw,fill=#1,
                    font=\sffamily\Large\bfseries},
                    ]

\path[black] (0,-0.12) node[anchor=north]{$s_{1}$};
\path[black] (2,-0.12) node[anchor=north]{$u_{1}$};
\path[black] (4,-0.12) node[anchor=north]{$s_{2}$};
\path[black] (6,-0.12) node[anchor=north]{$u_{2}$};
\path[black] (8,-0.12) node[anchor=north]{$s_{3}$};
\path[black] (10,-0.12) node[anchor=north]{$u_{3}$};
\path[black] (12,-0.12) node[anchor=north]{$t$};

\filldraw[black] (0,0) circle (1pt);
\filldraw[black] (2,0) circle (1pt);
\filldraw[black] (4,0) circle (1pt);
\filldraw[black] (6,0) circle (1pt);
\filldraw[black] (8,0) circle (1pt);
\filldraw[black] (10,0) circle (1pt);
\filldraw[black] (12,0) circle (1pt);

%\filldraw[black] (6,4) circle (8pt);

\draw[->] (-2,0) -- (13,0) node[anchor=west]{$\mathbb{R}_{+}$};
\draw[->] (-2,0) -- (-2,5) node[anchor=south]{$\mathbb{R}^{d}$};

\filldraw[black] (-2,0.5) circle (2pt) node[anchor=east]{$\Vec{x}_{0}(1)$};
\filldraw[black] (-2,2) circle (2pt)
node[anchor=east]{$\Vec{x}_{0}(2)$};
\filldraw[black] (-2,4.5) circle (2pt)
node[anchor=east]{$\Vec{x}_{0}(3)$};

\node (A) at (0,1.5) {};
\filldraw[color=red!60, fill=red!5, very thick](0,1.5) circle (0.13);
\node (B) at (2,2.5) {};
\filldraw[color=red!60, fill=red!5, very thick](2,2.5) circle (0.13);
\node (C) at (4,5) {};
\filldraw[color=red!60, fill=red!5, very thick](4,5) circle (0.13);
\node (D) at (8,2.5) {};
\filldraw[color=red!60, fill=red!5, very thick](8,2.5) circle (0.13);
\node (E) at (10,1.5) {};
\filldraw[color=red!60, fill=red!5, very thick](10,1.5) circle (0.13);
\node (CC) at (6,5.5) {};
\filldraw[color=red!60, fill=red!5, very thick](6,5.5) circle (0.13);

\draw[->] [very thick] (-2,0.5) to[out=5,in=-140] (A);
\draw[->] [very thick] (-2,2) to[out=10,in=160] (A);

%\draw[->] [very thick,blue] (A) to[out=10,in=-160] (B);
%\draw[->] [very thick,blue] (A) to[out=-5,in=-140] (B);
\path[draw=red, very thick, decorate] (A) -- (B);
\path[draw=red, very thick, decorate] (C) -- (CC);
\path[draw=red, very thick, decorate] (D) -- (E);

\draw[->] [very thick] (-2,4.5) to[out=30,in=160] (2,5);
\draw[->] [very thick] (2,5) to[out=-30,in=160] (C);
\draw[->] [very thick] (B) to[out=10,in=-160] (C);
\draw[->] [very thick] (B) to[out=-50,in=-150] (6,1);
\draw[->] [very thick] (6,1) to[out=30,in=-150] (D);
\draw[->] [very thick] (CC) to[out=-10,in=100] (6.7,4);
\draw[->] [very thick] (6.7,4) to[out=-60,in=150] (D);
\draw[->] [very thick] (CC) to[out=25,in=-150] (7,6);
\draw[->] [very thick] (7,6) to[out=30,in=150] (10,6);
\draw[->] [very thick] (10,6) to[out=-30,in=130] (12,4.5);
\draw[->] [very thick] (E) to[out=40,in=-150] (12,3);
\draw[->] [very thick] (E) to[out=-30,in=130] (12,0.5);

\filldraw[black] (12,0.5) circle (2pt) node[anchor=west]{$\Vec{x}_{4}(1)$};
\filldraw[black] (12,3) circle (2pt)
node[anchor=west]{$\Vec{x}_{4}(3)$};
\filldraw[black] (12,4.5) circle (2pt)
node[anchor=west]{$\Vec{x}_{4}(2)$};

\node at (0,1.2)[anchor=north]{$y_{1}$};
\node at (0,1.8)[anchor=south]{$y'_{1}$};

\node at (2,2.1)[anchor=north]{$\Vec{x}_{1}(1)$};
\node at (2,2.8)[anchor=south]{$\Vec{x}_{1}(2)$};
\node at (2,5.3)[anchor=south]{$\Vec{x}_{1}(3)$};

\node at (4,4.6)[anchor=north]{$y_{2}$};
\node at (4,5.3)[anchor=south]{$y'_{2}$};

\node at (8,2.8)[anchor=south]{$y'_{3}$};
\node at (8,2.2)[anchor=north]{$y_{3}$};

\node at (10,1.2)[anchor=north]{$\Vec{x}_{3}(1)$};
\node at (10,1.8)[anchor=south]{$\Vec{x}_{3}(3)$};
\node at (6,1.3)[anchor=south]{$\Vec{x}_{2}(1)$};
\node at (6,5.8)[anchor=south]{$\Vec{x}_{2}(2)$};
\node at (5.8,5.2)[anchor=north]{$\Vec{x}_{2}(3)$};
\node at (10,5.7)[anchor=north]{$\Vec{x}_{3}(2)$};

\end{tikzpicture}
\caption{The above graph stands for the sub-limiting path integral $\mathscr{I}_{0;t}^{N;(\ell_{1},\ell'_{1},\mathfrak{i}_{1}),...,(\ell_{m},\ell'_{m},\mathfrak{i}_{m})}U^{\otimes N}_{0}(\Vec{x}_{0})$ with $N = 3$, $m = 3$, $(\ell_{1},\ell'_{1},\mathfrak{i}_{1}) = (1,2,1)$, $(\ell_{2},\ell'_{2},\mathfrak{i}_{2}) = (2,3,1)$, and $(\ell_{3},\ell'_{3},\mathfrak{i}_{3}) = (1,3,1)$. Here, each of the above red notation, consisting of two red circles and a red wave line, denotes the transition kernel $\mathscr{A}^{(\ell_{j},\ell_{j}',\mathfrak{i}_{j})}_{0;u_{j-1},s_{j},u_{j}}(\Vec{x}_{j-1},\Vec{x}_{j})$. Similarly, the variable $y'_{j}$ always relates to the larger index $\ell'_{j}$.
}
\label{BD_model_5}
\end{figure}
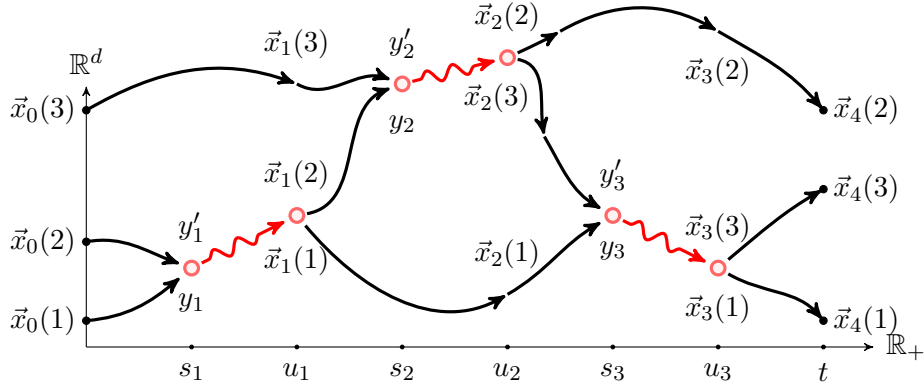

Finally, we summarize the above notations as follows. In view of (\ref{definition_sub_limit}), $\mathscr{Q}_{t}^{N}U_{0}^{\otimes N}(\Vec{x}_{0})$ has the same representation as  $\mathcal{Q}_{\varepsilon;t}^{\beta;N}U_{0}^{\otimes N}(\Vec{x}_{0})$. For the case that $\mathfrak{i}_{j} = 0$ for some $1\leq j\leq m$, the sub-limiting path integral $\mathscr{I}_{0;t}^{N;(\ell_{1},\ell'_{1},\mathfrak{i}_{1}),...,(\ell_{m},\ell'_{m},\mathfrak{i}_{m})}U_{0}^{\otimes N}(\Vec{x}_{0})$ is zero due to (\ref{definition_transition_kernel}). For the case that $\mathfrak{i}_{j} = 1$ for all $1\leq j\leq m$, $\mathscr{I}_{0;t}^{N;(\ell_{1},\ell'_{1},\mathfrak{i}_{1}),...,(\ell_{m},\ell'_{m},\mathfrak{i}_{m})}U_{0}^{\otimes N}(\Vec{x}_{0})$ has the same expression as $\mathscr{I}_{\varepsilon;t}^{N;(\ell_{1},\ell'_{1},\mathfrak{i}_{1}),...,(\ell_{m},\ell'_{m},\mathfrak{i}_{m})}U_{0}^{\otimes N}(\Vec{x}_{0})$, but the transition kernel $\mathscr{A}^{(\ell_{j},\ell_{j}',\mathfrak{i}_{j})}_{\varepsilon;u_{j-1},s_{j},u_{j}}(\Vec{x}_{j-1},\Vec{x}_{j})$ is replacing by the transition kernel $\mathscr{A}^{(\ell_{j},\ell_{j}',\mathfrak{i}_{j})}_{0;u_{j-1},s_{j},u_{j}}(\Vec{x}_{j-1},\Vec{x}_{j})$. See Figure \ref{BD_model_5} for a graphical clarification of $\mathscr{I}_{0;t}^{N;(\ell_{1},\ell'_{1},\mathfrak{i}_{1}),...,(\ell_{m},\ell'_{m},\mathfrak{i}_{m})}U_{0}^{\otimes N}(\Vec{x}_{0})$.

\paragraph{Heuristic for (\ref{definition_sub_limit}).}
In this part, let us elaborate on the reason why $\mathscr{Q}_{t}^{N}U_{0}^{\otimes N}(\Vec{x}_{0})$ is a formal limit of $\mathcal{Q}_{\varepsilon;t}^{\beta;N}U_{0}^{\otimes N}(\Vec{x}_{0})$ as $\varepsilon \to 0^{+}$. Briefly speaking, the formal limit means that \emph{we put the limit $\varepsilon \to 0^{+}$ inside the summation with respect to $m$ on the right-hand side of (\ref{decomposition_second})} as follows:
\begin{align} \label{inside_limit}
    &\mathbf{E}_{\Vec{x}_{0}}^{\Vec{B}}[U_{0}^{\otimes N}(\Vec{B}(t))]  \nonumber\\
    &+\sum_{m=1}^{\infty} \sum_{(\ell_{1},\ell'_{1},\mathfrak{i}_{1}),...,(\ell_{m},\ell'_{m},\mathfrak{i}_{m}) \in \mathcal{E}^{(N)} \times \{0,1\}, \; (\ell_{1},\ell_{1}') \neq ...\neq (\ell_{m},\ell_{m}')} \lim_{\varepsilon\to 0^{+}}\mathscr{I}_{\varepsilon;t}^{N;(\ell_{1},\ell'_{1},\mathfrak{i}_{1}),...,(\ell_{m},\ell'_{m},\mathfrak{i}_{m})}U^{\otimes N}_{0}(\Vec{x}_{0}).
\end{align}
Moreover, we will explain the reason why we have the formal approximation below:
\begin{align} \label{formal_limit}
    \mathscr{I}_{\varepsilon;t}^{N;(\ell_{1},\ell'_{1},\mathfrak{i}_{1}),...,(\ell_{m},\ell'_{m},\mathfrak{i}_{m})} U_{0}^{\otimes N} (\Vec{x}_{0}) \overset{\varepsilon \to 0^{+}}{\approx} \mathscr{I}_{0;t}^{N;(\ell_{1},\ell'_{1},\mathfrak{i}_{1}),...,(\ell_{m},\ell'_{m},\mathfrak{i}_{m})} U_{0}^{\otimes N} (\Vec{x}_{0}).
\end{align}
Consequently, (\ref{inside_limit}) and (\ref{formal_limit}) are the major reasons why $\mathscr{Q}_{t}^{N}U_{0}^{\otimes N}(\Vec{x}_{0})$ defined in (\ref{definition_sub_limit}) is a formal limit of $\mathcal{Q}_{\varepsilon;t}^{\beta;N}U_{0}^{\otimes N}(\Vec{x}_{0})$ as $\varepsilon \to 0^{+}$.

As a matter of fact, for $N = 3$ and for a proper initial datum $U_{0}$, (\ref{formal_limit}) can be proved rigorously. But, we feel that giving a heuristic for (\ref{formal_limit}) is enough instead of proving it since (\ref{formal_limit}) will not be used in this article. Due to this reason, it is appropriate to call (\ref{formal_limit}) a formal approximation.

The rest of this part concentrates on the heuristic for (\ref{formal_limit}). Before doing so, we first introduce the main idea, which is the following asymptotic of the rescaled Brownian exponential functional at $L^{2}$-criticality: 
\begin{align}\label{limit_two_particle}
    &\int_{\mathbb{R}^{2\times d}} dzdz' \beta_{\varepsilon}^{2} R(\sqrt{2}z) \mathbf{E}_{\varepsilon \cdot z}^{B}\bigbkt{\exp\bkt{\beta_{\varepsilon}^{2}\int_{0}^{t} R_{\varepsilon}(\sqrt{2}B(v)) dv}: B(t) = \varepsilon \cdot z'} \beta_{\varepsilon}^{2} R(\sqrt{2}z')\overset{\varepsilon\to 0^{+}}{\approx} \bkt{\frac{2\pi}{t}}^{1/2}.
\end{align}
The above approximation follows from the properties of three-dimensional delta-Bose gas studied in \cite[Chapter I.1]{albeverio2012solvable}. We refer to Section \ref{section_two_particles_system} for a rigorous proof of (\ref{limit_two_particle}).
\begin{remark}
%We want to emphasize that the above approximation relies on the condition $\beta = \beta_{L^{2}}$. This will be clarified in Section \ref{section_two_particles_system}. On the other hand, 
Note that if $\beta < \beta_{L^{2}}$, it can be shown that the left-hand side of (\ref{limit_two_particle}) vanishes as $\varepsilon\to 0^{+}$, and thus, the left-hand side of (\ref{formal_limit}) vanishes for all $m\geq 1$. In this way, the formal limit (\ref{inside_limit}) is nothing but the Gaussian kernel $\mathbf{E}_{\Vec{x}_{0}}^{\Vec{B}}[U_{0}^{\otimes N}(\Vec{B}(t))]$.
\end{remark}

Let us now elaborate on the heuristics for (\ref{formal_limit}). The core is to explain the reason why the transition kernel $\mathscr{A}^{(\ell_{j},\ell_{j}',\mathfrak{i}_{j})}_{\varepsilon;u_{j-1},s_{j},u_{j}}(\Vec{x}_{j-1},\Vec{x}_{j})$ can be replaced by the transition kernel $\mathscr{A}^{(\ell_{j},\ell_{j}',\mathfrak{i}_{j})}_{0;u_{j-1},s_{j},u_{j}}(\Vec{x}_{j-1},\Vec{x}_{j})$ by using (\ref{limit_two_particle}), where $\mathscr{A}^{(\ell_{j},\ell_{j}',\mathfrak{i}_{j})}_{\varepsilon;u_{j-1},s_{j},u_{j}}(\Vec{x}_{j-1},\Vec{x}_{j})$ and $\mathscr{A}^{(\ell_{j},\ell_{j}',\mathfrak{i}_{j})}_{0;u_{j-1},s_{j},u_{j}}(\Vec{x}_{j-1},\Vec{x}_{j})$ are defined in (\ref{definition_epsilon_motion}) and (\ref{definition_transition_kernel}), respectively.  It has to be stressed that the following argument is \emph{not} a rigorous proof. 
\begin{proof} [Formal derivation of (\ref{definition_sub_limit})]
In the sequel, the notation $a_{\varepsilon}\overset{\varepsilon \to 0^{+}}{\approx}b_{\varepsilon}$ only denotes a formal limit. In view of the above elaboration, it is enough to clarify the formal limit (\ref{formal_limit}). To compute the limit of the right-hand side of (\ref{identity_application_markov_property}) by using the limit (\ref{limit_two_particle}), we re-express the right-hand side of (\ref{identity_application_markov_property}) as the following form: 
\begin{align} \label{definition_F_epsilon}
    &\text{(R.H.S) of (\ref{identity_application_markov_property})} = \int_{u_{0} = 0<s_{1}<u_{1}<s_{2}...<s_{m}<u_{m}<t} \prod_{k=1}^{m}ds_{k}du_{k} \cdot \int_{\mathbb{R}^{N \times d}}  d\Vec{x}_{j-1} \cdot F^{1:j-1}_{\varepsilon}(\Vec{x}_{j-1})  \nonumber\\
    &\int_{\mathbb{R}^{(N-2)\times d}} \prod_{1\leq k\leq N,\;k \not \in \{\ell_{j},\ell'_{j}\}} d\Vec{x}_{j}(k) \cdot \biggbkt{\prod_{1\leq k\leq N,\; k \not \in \{\ell_{j},\ell'_{j}\}} G_{u_{k}-u_{j-1}}(\Vec{x}_{j}(k)-\Vec{x}_{j-1}(k))}  \cdot \delta_{0}(s_{j}-u_{j})^{1-\mathfrak{i}_{j}}\nonumber\\
    &\biggbkt{\int_{\mathbb{R}^{2 \times d}}  d\Vec{x}_{j}(\ell_{j})d\Vec{x}_{j}(\ell_{j}') \cdot \mathscr{A}^{(\ell_{j},\ell_{j}',\mathfrak{i}_{j})}_{\varepsilon;u_{j-1},s_{j},u_{j}}(\Vec{x}_{j-1},\Vec{x}_{j}) F^{j+1: m}_{\varepsilon}\Bigl(\Vec{x}_{j}(\ell_{j}),\Vec{x}_{j}(\ell'_{j}),(\Vec{x}_{j}(k))_{1\leq k\leq N,\; k \not \in \{\ell_{j},\ell'_{j}\}}\Bigr)}.
\end{align}
Now, we apply the following change of variables to $\mathscr{A}^{(\ell_{j},\ell_{j}',\mathfrak{i}_{j})}_{\varepsilon;u_{j-1},s_{j},u_{j}}(\Vec{x}_{j-1},\Vec{x}_{j})\cdot d\Vec{x}_{j}(\ell_{j})d\Vec{x}_{j}(\ell_{j}')$:
\begin{align*}
    &y_{j}^{(\ell_{j},\ell'_{j})} = \frac{y'_{j}-y_{j}}{\varepsilon\cdot\sqrt{2}}, \quad
    y_{j}^{(\ell_{j},\ell'_{j})^{*}} = \frac{y'_{j}+y_{j}}{\sqrt{2}}, \\
    &\quad\quad\quad\quad\quad\quad\Vec{x}_{j}^{(\ell_{j},\ell'_{j})} = \frac{\Vec{x}_{j}(\ell'_{j})-\Vec{x}_{j}(\ell_{j})}{\varepsilon \cdot \sqrt{2}}, \quad \text{and} \quad
    \Vec{x}_{j}^{(\ell_{j},\ell'_{j})^{*}} = \frac{\Vec{x}_{j}(\ell'_{j})+\Vec{x}_{j}(\ell_{j})}{\sqrt{2}} \quad \text{if }  \mathfrak{i}_{j} = 1;
\end{align*}
\begin{align*}
    \Vec{x}_{j}^{(\ell_{j},\ell'_{j})} = \frac{\Vec{x}_{j}(\ell'_{j})-\Vec{x}_{j}(\ell_{j})}{\varepsilon \cdot \sqrt{2}} \quad \text{and} \quad
    \Vec{x}_{j}^{(\ell_{j},\ell'_{j})^{*}} = \frac{\Vec{x}_{j}(\ell'_{j})+\Vec{x}_{j}(\ell_{j})}{\sqrt{2}} \quad \text{if }  \mathfrak{i}_{j} = 0.
\end{align*}
Then, we conclude the following identity for the last line of the right-hand side of (\ref{definition_F_epsilon}):
\begin{align} \label{identity_change_of_variable_formal_limit_1}
    &\int_{\mathbb{R}^{2\times d}} d\Vec{x}_{j}(\ell_{j})d\Vec{x}_{j}(\ell'_{j})\cdot F^{j+1: m}_{\varepsilon}\Bigl(\Vec{x}_{j}(\ell_{j}),\Vec{x}_{j}(\ell'_{j}),(\Vec{x}_{j}(k))_{1\leq k\leq N,\; k \not \in \{\ell_{j},\ell'_{j}\}}\Bigr)\mathscr{A}^{(\ell_{j},\ell_{j}',\mathfrak{i}_{j})}_{\varepsilon;u_{j-1},s_{j},u_{j}}(\Vec{x}_{j-1},\Vec{x}_{j})\nonumber\\
    &= \int_{\mathbb{R}^{2 \times d}} dy_{j}^{(\ell_{j},\ell'_{j})}dy_{j}^{(\ell_{j},\ell'_{j})^{*}} \int_{\mathbb{R}^{2\times d}} d\Vec{x}_{j}^{(\ell_{j},\ell'_{j})}d\Vec{x}_{j}^{(\ell_{j},\ell'_{j})^{*}}\cdot \nonumber\\
    &F^{j+1:m}_{\varepsilon}\bkt{\frac{\Vec{x}_{j}^{(\ell_{j},\ell'_{j})^{*}}-\varepsilon\cdot\Vec{x}_{j}^{(\ell_{j},\ell'_{j})}}{\sqrt{2}},\frac{\Vec{x}_{j}^{(\ell_{j},\ell'_{j})^{*}}+\varepsilon\cdot\Vec{x}_{j}^{(\ell_{j},\ell'_{j})}}{\sqrt{2}},(\Vec{x}_{j}(k))_{1\leq k\leq N,\; k \not \in \{\ell_{j},\ell'_{j}\}}} \cdot \nonumber\\ 
    &G_{s_{j}-u_{j-1}}\bkt{y_{j}^{(\ell_{j},\ell'_{j})^{*}} - \Vec{x}_{j-1}^{(\ell_{j},\ell'_{j})^{*}}}  G_{s_{j}-u_{j-1}}\bkt{\varepsilon \cdot y_{j}^{(\ell_{j},\ell_{j}')}-\Vec{x}_{j-1}^{(\ell_{j},\ell'_{j})}} \beta_{\varepsilon}^{2}R\bkt{\sqrt{2}y_{j}^{(\ell_{j},\ell_{j}')}}\cdot \nonumber\\
    &\biggbkt{\mathbf{E}_{\varepsilon \cdot y_{j}^{(\ell_{j},\ell_{j}')}}^{B}\bigbkt{\exp\bkt{\beta_{\varepsilon}^{2}\int_{0}^{u_{j}-s_{j}}R_{\varepsilon}(\sqrt{2}B(v))dv}: B(t) = \varepsilon \cdot \Vec{x}_{j}^{(\ell_{j},\ell'_{j})}} \nonumber\\
    &\quad\quad\quad\quad\quad\quad\quad\quad\quad\quad\quad\quad G_{u_{j}-s_{j}}\bkt{ \Vec{x}_{j}^{(\ell_{j},\ell'_{j})^{*}} - y_{j}^{(\ell_{j},\ell'_{j})^{*}}} }\cdot\beta_{\varepsilon}^{2} R\bkt{\sqrt{2} \Vec{x}_{j}^{(\ell_{j},\ell'_{j})}} \quad \text{if } \mathfrak{i}_{j} = 1;
\end{align}
\begin{align} \label{identity_change_of_variable_formal_limit_0}
    &\int_{\mathbb{R}^{2\times d}} d\Vec{x}_{j}(\ell_{j})d\Vec{x}_{j}(\ell'_{j})\cdot F^{j+1:m}_{\varepsilon}\Bigl(\Vec{x}_{j}(\ell_{j}),\Vec{x}_{j}(\ell'_{j}),(\Vec{x}_{j}(k))_{1\leq k\leq N,\; k \not \in \{\ell_{j},\ell'_{j}\}}\Bigr)\mathscr{A}^{(\ell_{j},\ell_{j}',\mathfrak{i}_{j})}_{\varepsilon;u_{j-1},s_{j},u_{j}}(\Vec{x}_{j-1},\Vec{x}_{j})\nonumber\\
    &= \int_{\mathbb{R}^{2\times d}} d\Vec{x}_{j}^{(\ell_{j},\ell'_{j})}d\Vec{x}_{j}^{(\ell_{j},\ell'_{j})^{*}}\cdot \nonumber\\
    &F^{j+1:m}_{\varepsilon}\bkt{\frac{\Vec{x}_{j}^{(\ell_{j},\ell'_{j})^{*}}-\varepsilon\cdot\Vec{x}_{j}^{(\ell_{j},\ell'_{j})}}{\sqrt{2}},\frac{\Vec{x}_{j}^{(\ell_{j},\ell'_{j})^{*}}+\varepsilon\cdot\Vec{x}_{j}^{(\ell_{j},\ell'_{j})}}{\sqrt{2}},(\Vec{x}_{j}(k))_{1\leq k\leq N,\; k \not \in \{\ell_{j},\ell'_{j}\}}} \cdot \nonumber\\ 
    &G_{u_{j}-u_{j-1}}\bkt{\Vec{x}_{j}^{(\ell_{j},\ell'_{j})^{*}} - \Vec{x}_{j-1}^{(\ell_{j},\ell'_{j})^{*}}}  G_{u_{j}-u_{j-1}}\bkt{\varepsilon\cdot\Vec{x}_{j}^{(\ell_{j},\ell'_{j})}-\Vec{x}_{j-1}^{(\ell_{j},\ell'_{j})}} \cdot \beta_{\varepsilon}^{2} R\bkt{\sqrt{2} \Vec{x}_{j}^{(\ell_{j},\ell'_{j})}} \quad \text{if } \mathfrak{i}_{j} = 0,
\end{align}
where $\Vec{x}_{j-1}^{(\ell_{j},\ell_{j}')}$ and $\Vec{x}_{j-1}^{(\ell_{j},\ell_{j}')^{*}}$ are functions with the following variables:
\begin{align*}
    \Vec{x}_{j-1}^{(\ell_{j-1},\ell'_{j-1})}, \quad 
    \Vec{x}_{j-1}^{(\ell_{j-1},\ell'_{j-1})^{*}}, \quad \text{and} \quad \Bigl\{\Vec{x}_{j-1}(i): 1\leq i\leq N,\; i \not \in \{\ell_{j-1},\ell'_{j-1}\}\Bigr\}.
\end{align*}
Here, for the case of $\mathfrak{i}_{j} = 1$, we have used the following facts in (\ref{identity_change_of_variable_formal_limit_1}) to re-express $\mathscr{A}^{(\ell_{j},\ell_{j}',\mathfrak{i}_{j})}_{\varepsilon;u_{j-1},s_{j},u_{j}}(\Vec{x}_{j-1},\Vec{x}_{j})$ on the right-hand side of (\ref{definition_F_epsilon}):
\begin{align*}
    G_{s_{j}-u_{j-1}}(y_{j}&-\Vec{x}_{j-1}(\ell_{j}))G_{s_{j}-u_{j-1}}(y_{j}'-\Vec{x}_{j-1}(\ell'_{j}))\\
    &= G_{s_{j}-u_{j-1}}\bkt{y_{j}^{(\ell_{j},\ell'_{j})^{*}} - \Vec{x}_{j-1}^{(\ell_{j},\ell'_{j})^{*}}} G_{s_{j}-u_{j-1}}\bkt{ y_{j}^{(\ell_{j},\ell_{j}')}-\Vec{x}_{j-1}^{(\ell_{j},\ell'_{j})}} 
\end{align*}
and 
\begin{align*}
    &\mathbf{E}_{y_{j},y_{j}'}^{B,B'}\bigbkt{\exp\bkt{\beta_{\varepsilon}^{2}\int_{0}^{u_{j}-s_{j}}R_{\varepsilon}(B'(v) - B(v))dv}: (B(u_{j}-s_{j}),B'(u_{j}-s_{j})) = (\Vec{x}_{j}(\ell_{j}),\Vec{x}_{j}(\ell'_{j}))}\\
    &= \mathbf{E}_{y_{j}^{(\ell_{j},\ell_{j}')}}^{B}\bigbkt{\exp\bkt{\beta_{\varepsilon}^{2}\int_{0}^{u_{j}-s_{j}}R_{\varepsilon}(\sqrt{2}B(v))dv}: B(t) = \Vec{x}_{j}^{(\ell_{j},\ell'_{j})}} G_{u_{j}-s_{j}}\bkt{ \Vec{x}_{j}^{(\ell_{j},\ell'_{j})^{*}} - y_{j}^{(\ell_{j},\ell'_{j})^{*}}}.
\end{align*}
Moreover, for the case of $\mathfrak{i}_{j} = 0$, the following identity in (\ref{identity_change_of_variable_formal_limit_0}) have been applied to re-express $\mathscr{A}^{(\ell_{j},\ell_{j}',\mathfrak{i}_{j})}_{\varepsilon;u_{j-1},s_{j},u_{j}}(\Vec{x}_{j-1},\Vec{x}_{j})$ on the right-hand side of (\ref{definition_F_epsilon}):
\begin{align*}
    G_{u_{j}-u_{j-1}}(\Vec{x}_{j}(\ell_{j})&-\Vec{x}_{j-1}(\ell_{j}))G_{u_{j}-u_{j-1}}(\Vec{x}_{j}(\ell_{j}')-\Vec{x}_{j-1}(\ell'_{j}))\\
    &= G_{u_{j}-u_{j-1}}\bkt{\Vec{x}_{j}^{(\ell_{j},\ell'_{j})^{*}} - \Vec{x}_{j-1}^{(\ell_{j},\ell'_{j})^{*}}} \cdot G_{u_{j}-u_{j-1}}\bkt{ \Vec{x}_{j}^{(\ell_{j},\ell_{j}')}-\Vec{x}_{j-1}^{(\ell_{j},\ell'_{j})}}. 
\end{align*}

Next, the following approximations follow by passing $\varepsilon \to 0^{+}$ for the right-hand side of (\ref{identity_change_of_variable_formal_limit_1}) and (\ref{identity_change_of_variable_formal_limit_0}):
\begin{align} \label{limit_derivation_1}
    &\text{(R.H.S) of (\ref{identity_change_of_variable_formal_limit_1})} \overset{\varepsilon\to 0^{+}}{\approx} \int_{\mathbb{R}^{d}} dy_{j}^{(\ell_{j},\ell'_{j})^{*}} \int_{\mathbb{R}^{d}} d\Vec{x}_{j}^{(\ell_{j},\ell'_{j})^{*}}\nonumber\\
    &F^{j+1:m}_{0}\bkt{\frac{\Vec{x}_{j}^{(\ell_{j},\ell'_{j})^{*}}-\Vec{x}_{j}^{(\ell_{j},\ell'_{j})}}{\sqrt{2}} \biggr|_{\Vec{x}_{j}^{(\ell_{j},\ell'_{j})} = 0},\frac{\Vec{x}_{j}^{(\ell_{j},\ell'_{j})^{*}}+\Vec{x}_{j}^{(\ell_{j},\ell'_{j})}}{\sqrt{2}}\biggr|_{\Vec{x}_{j}^{(\ell_{j},\ell'_{j})} = 0},(\Vec{x}_{j}(k))_{1\leq k\leq N,\; k \not \in \{\ell_{j},\ell'_{j}\}}} \cdot \nonumber\\
    &G_{s_{j}-u_{j-1}}\bkt{ \Vec{y}_{j}^{(\ell_{j},\ell'_{j})^{*}} - \Vec{x}_{j-1}^{(\ell_{j},\ell'_{j})^{*}} \biggl|_{\Vec{x}_{j-1}^{(\ell_{j-1},\ell'_{j-1})} = 0}} G_{s_{j}-u_{j-1}}\bkt{\Vec{y}_{j}^{(\ell_{j},\ell'_{j})}\biggr|_{\Vec{y}_{j}^{(\ell_{j},\ell'_{j})} = 0}-\Vec{x}_{j-1}^{(\ell_{j},\ell'_{j})}\biggr|_{\Vec{x}_{j-1}^{(\ell_{j-1},\ell'_{j-1})} = 0}} \cdot  \nonumber\\
    &\biggbkt{G_{u_{j}-s_{j}}\bkt{ \Vec{x}_{j}^{(\ell_{j},\ell'_{j})^{*}} - y_{j}^{(\ell_{j},\ell'_{j})^{*}}} \cdot \bkt{\frac{2\pi}{u_{j}-s_{j}}}^{1/2}}\quad \text{if $\mathfrak{i}_{j} = 1$ and $2\leq j\leq m$,}  
\end{align}
\begin{align} \label{limit_derivation_1_2}
    &\text{(R.H.S) of (\ref{identity_change_of_variable_formal_limit_1})} \overset{\varepsilon\to 0^{+}}{\approx} \int_{\mathbb{R}^{d}} dy_{j}^{(\ell_{j},\ell'_{j})^{*}} \int_{\mathbb{R}^{d}} d\Vec{x}_{j}^{(\ell_{j},\ell'_{j})^{*}}\nonumber\\
    &F^{j+1:m}_{0}\bkt{\frac{\Vec{x}_{j}^{(\ell_{j},\ell'_{j})^{*}}-\Vec{x}_{j}^{(\ell_{j},\ell'_{j})}}{\sqrt{2}} \biggr|_{\Vec{x}_{j}^{(\ell_{j},\ell'_{j})} = 0},\frac{\Vec{x}_{j}^{(\ell_{j},\ell'_{j})^{*}}+\Vec{x}_{j}^{(\ell_{j},\ell'_{j})}}{\sqrt{2}}\biggr|_{\Vec{x}_{j}^{(\ell_{j},\ell'_{j})} = 0},(\Vec{x}_{j}(k))_{1\leq k\leq N,\; k \not \in \{\ell_{j},\ell'_{j}\}}} \cdot \nonumber\\
    &G_{s_{j}-u_{j-1}}\bkt{ \Vec{y}_{j}^{(\ell_{j},\ell'_{j})^{*}} - \Vec{x}_{j-1}^{(\ell_{j},\ell'_{j})^{*}}}  G_{s_{j}-u_{j-1}}\bkt{\Vec{y}_{j}^{(\ell_{j},\ell'_{j})}\biggr|_{\Vec{y}_{j}^{(\ell_{j},\ell'_{j})} = 0}-\Vec{x}_{j-1}^{(\ell_{j},\ell'_{j})}} \cdot  \nonumber\\
    &\biggbkt{G_{u_{j}-s_{j}}\bkt{ \Vec{x}_{j}^{(\ell_{j},\ell'_{j})^{*}} - y_{j}^{(\ell_{j},\ell'_{j})^{*}}} \cdot\bkt{\frac{2\pi}{u_{j}-s_{j}}}^{1/2}}\quad\text{if $\mathfrak{i}_{j} = 1$ and $j = 1$,}  
\end{align}
and
\begin{align} \label{limit_derivation_2}
    \text{(R.H.S) of (\ref{identity_change_of_variable_formal_limit_0})} \overset{\varepsilon\to 0^{+}}{\approx} 0 \quad \text{if $\mathfrak{i}_{j} = 0$ and $1\leq j\leq m$,}
\end{align}
where $\Vec{x}_{j-1}^{(\ell_{j-1},\ell'_{j-1})^{*}} = 0$ on the right-hand side of (\ref{limit_derivation_1}) is because of $R_{\varepsilon}(\sqrt{2} \Vec{x}_{j-1}^{(\ell_{j-1},\ell'_{j-1})})$ in the previous transition kernel $\mathscr{A}^{(\ell_{j-1},\ell_{j-1}',\mathfrak{i}_{j-1})}_{\varepsilon;u_{j-2},s_{j-1},u_{j-1}}(\Vec{x}_{j-2},\Vec{x}_{j-1})$. Moreover, for the case of $\mathfrak{i}_{j} = 1$, we have used the approximation below which follows from (\ref{limit_two_particle}):
\begin{align*}
     &\int_{\mathbb{R}^{d}} dy_{j}^{(\ell_{j},\ell'_{j})} \int_{\mathbb{R}^{d}} d\Vec{x}_{j}^{(\ell_{j},\ell'_{j})}\beta_{\varepsilon}^{2}R\bkt{\sqrt{2}y_{j}^{(\ell_{j},\ell_{j}')}}\nonumber\\
    &\mathbf{E}_{\varepsilon \cdot y_{j}^{(\ell_{j},\ell_{j}')}}^{B}\bigbkt{\exp\bkt{\beta_{\varepsilon}^{2}\int_{0}^{u_{j}-s_{j}}R_{\varepsilon}(\sqrt{2}B(v))dv}: B(t) = \varepsilon \cdot \Vec{x}_{j}^{(\ell_{j},\ell'_{j})}} \beta_{\varepsilon}^{2} R\bkt{\sqrt{2} \Vec{x}_{j}^{(\ell_{j},\ell'_{j})}} \overset{\varepsilon \to 0^{+}}{\approx}\bkt{\frac{2\pi}{u_{j}-s_{j}}}^{1/2} 
\end{align*}
Furthermore, for the case of $\mathfrak{i}_{j} = 0$, we have used the limit below due to (\ref{definition_beta_epsilon}):
\begin{align*}
    \beta_{\varepsilon}^{2} R\bkt{\sqrt{2} \Vec{x}_{j}^{(\ell_{j},\ell'_{j})}} \overset{\varepsilon \to 0^{+}}{\approx} 0.
\end{align*}

Finally, we explain the reason why the the above results imply the formal limit (\ref{formal_limit}).
To this aim, we note that the following approximation follows from (\ref{limit_derivation_2}):
\begin{align*}
    \mathscr{I}_{\varepsilon;t}^{N;(\ell_{1},\ell'_{1},\mathfrak{i}_{1}),...,(\ell_{m},\ell'_{m},\mathfrak{i}_{m})} U_{0}^{\otimes N} (\Vec{x}_{0}) \overset{\varepsilon \to 0^{+}}{\approx} 0 \quad \text{if there exists $1\leq j\leq m$ such that $\mathfrak{i}_{j} = 0$,}
\end{align*}
which matches the definition of $\mathscr{I}_{0;t}^{N;(\ell_{1},\ell'_{1},\mathfrak{i}_{1}),...,(\ell_{m},\ell'_{m},\mathfrak{i}_{m})} U_{0}^{\otimes N} (\Vec{x}_{0})$. See the discussion below Remark \ref{remark_after_definition}. Moreover, if $\mathfrak{i}_{j} = 1$ for all $1\leq j\leq m$, then, thanks to (\ref{limit_derivation_1}) and (\ref{limit_derivation_1_2}), one can then replace $\mathscr{A}^{(\ell_{j},\ell_{j}',\mathfrak{i}_{j})}_{\varepsilon;u_{j-1},s_{j},u_{j}}(\Vec{x}_{j-1},\Vec{x}_{j})$ on the right-hand side of (\ref{definition_F_epsilon}) by $\mathscr{A}^{(\ell_{j},\ell_{j}',\mathfrak{i}_{j})}_{0;u_{j-1},s_{j},u_{j}}(\Vec{x}_{j-1},\Vec{x}_{j})$. To be more precise, the Dirac’s delta function $\delta_{0}((\Vec{x}_{j}(\ell'_{j}) - \Vec{x}_{j}(\ell_{j}))/\sqrt{2})$ in (\ref{definition_transition_kernel}) follows from the fact that $\Vec{x}_{j}^{(\ell_{j},\ell_{j}')} = 0$ on the right-hand sides of (\ref{limit_derivation_1}) and (\ref{limit_derivation_1_2}), and the Dirac’s delta function $\delta_{0}((y_{j}' - y_{j})/\sqrt{2})$ in (\ref{definition_transition_kernel}) is because of $y_{j}^{(\ell_{j},\ell_{j}')} = 0$ on the right-hand sides of (\ref{limit_derivation_1}) and (\ref{limit_derivation_1_2}). Therefore, we conclude
\begin{align*}
    \mathscr{I}_{\varepsilon;t}^{N;(\ell_{1},\ell'_{1},\mathfrak{i}_{1}),...,(\ell_{m},\ell'_{m},\mathfrak{i}_{m})} U_{0}^{\otimes N} (\Vec{x}_{0}) \overset{\varepsilon \to 0^{+}}{\approx} \mathscr{I}_{0;t}^{N;(\ell_{1},\ell'_{1},\mathfrak{i}_{1}),...,(\ell_{m},\ell'_{m},\mathfrak{i}_{m})} U_{0}^{\otimes N} (\Vec{x}_{0}) \text{ if $\mathfrak{i}_{j} = 1$ for all $1\leq j\leq m$.}
\end{align*}
Therefore, the elaboration is complete.
\end{proof}

\subsubsection{Unboundedness of the sub-limiting higher moments.} \label{section_unboundedness_of_the_sub-limiting_moments}
With the above terminologies in place, we are ready to introduce our second main result of this article. In particular, to show that the unboundedness of the sub-limiting higher moment has nothing to do with the initial condition, we consider all the initial conditions that can be bounded from below by a Gaussian initial datum. This includes the constant initial condition, which relates to the partition function of the continuous directed polymer. Then, the second main result is then the following.

\begin{theorem} \label{Main_result_2}
Assume that $d = 3$ and $N\geq 3$. Recall that $\mathscr{Q}_{t}^{N}U_{0}^{\otimes N}(\Vec{x}_{0})$ is the sub-limiting $N$-th moment is defined in (\ref{definition_sub_limit}) and $U_{0}$ is the initial datum of the mollified SHE. Suppose that $U_{0} \in C_{b}(\mathbb{R}^{d};\mathbb{R}_{+})$ so that $U_{0}(x) \geq C G_{\nu}(x)$ for some $C,\nu > 0$. Then, for each $t>0$, the following limits hold:
\begin{align*}
    \mathscr{Q}_{t}^{N}U_{0}^{\otimes N}(\Vec{x}_{0}) = \infty \quad \forall \Vec{x}_{0} \in \mathbb{R}^{N \times d} \setminus \Pi^{(N)}
\end{align*}
and
\begin{align} \label{unboundedness_sub_limit}
    \Bigl\langle \varphi,\mathscr{Q}_{t}^{N}U_{0}^{\otimes N}\Bigr\rangle_{L^{2}(\mathbb{R}^{N \times d})}= \infty \quad \forall \varphi \in C_{c}^{\infty}(\mathbb{R}^{N \times d}; \mathbb{R}_{+}),
\end{align}
where the above divergences are in the sense of the limit $M\to\infty$ on the right-hand side of (\ref{definition_sub_limit}).
\end{theorem}
\begin{remark}
An explanation of Theorem \ref{Main_result_2} from a quantum physics point of view will be presented in Section \ref{section_related_literature}.
\end{remark}

Recalling Theorem \ref{Main_result_1}, we cannot conclude the divergence for all the higher moments inside the $L^{2}$-regime. But, at least for the three-dimensional case, Theorem \ref{Main_result_2} allows us to conclude a refinement of the above divergence, and implies another partial result for Conjecture \ref{Conjecture_order}. Both of these can be deduced from the following proposition.

\begin{proposition} \label{Main_result_5}
Assume that $d = 3$. Recall that the positive integer $N_{0}$ defined in Theorem \ref{Main_result_1}. Then $N_{0} = 3$. In other words, the estimate (\ref{estimate_betaLN_less_betaL2}) holds for all $N\geq 3$.
\end{proposition}

%\begin{remark}Why $d = 3$? How about $d\geq 4$?\end{remark}

\begin{remark} \label{Efimov_effect_first_explanation}
Let us briefly explain the connection between Proposition \ref{Main_result_5} and a three-dimensional many-body quantum system. Loosely speaking, the major tool to conclude Proposition \ref{Main_result_5} is (\ref{positivity_beta_L2}), which follows from Theorem \ref{Main_result_2}. However, (\ref{positivity_beta_L2}) is also an application of a well-known result in quantum physics, which is called the Efimov effect. A more explicit elaboration of this will be presented in Section \ref{section_related_literature}.
\end{remark}

\subsection{Applications to continuous directed polymer.} \label{section_CDP}
Let us now introduce some applications to the continuous directed polymer based on previous results. The overall goal is to extend Theorem \ref{Main_result_3} to the fractional moment and prove an estimate for the critical exponent of the continuous directed polymer. Before doing this, we first introduce the model of continuous directed polymer and its connection with the mollified SHE. Recall the Feynman-Kac representation (\ref{representation_FK}). Then, if we set the initial datum $U_{0} \equiv 1$, then one has the following relation:
\begin{align} \label{connect_SHE_CDP}
    \mathscr{U}_{\varepsilon}(x,t) = \mathcal{Z}^{\beta}_{t \cdot \varepsilon^{-2}}(0;\xi^{(\varepsilon,x,t)}).
\end{align}
Here, $\mathcal{Z}^{\beta}_{T}(x;\xi)$ defined as follows is the partition function of the continuous directed polymer in white noise environment
at inverse temperature $\beta$:
\begin{align}
    \label{definition_partition_function_1}
    \mathcal{Z}^{\beta}_{T}(x;\xi)
    := \mathbf{E}^{B}_{x}\bigbkt{\exp\bkt{\beta \int_{0}^{T} ds \int_{\mathbb{R}^{d}} dy\phi(B(s) - y) \xi(y,s) - \frac{\beta^{2} R(0) T }{2}}},
\end{align}
where $\phi(x)$ comes from the approximation of the identity defined in Section \ref{section_the_model}, and 
\begin{align} \label{definition_R}
    R(x) := R_{1}(x) = \phi*\phi(x).
\end{align}
Moreover, $\xi^{(\varepsilon,x,t)}$ defined as follows is the diffusively rescaled, time-reversed space-time white noise:
\begin{equation} \label{definition_rescaled_space_time_noise}
    \xi^{(\varepsilon,x,t)}(y,s) := \varepsilon^{\frac{d+2}{2}}\xi(\varepsilon y + x,t - \varepsilon^{2}s).
\end{equation}

About the limit of the partition functions $(\mathcal{Z}^{\beta}_{T}(0;\xi))_{T\geq 0}$, a well-known result proved by the authors in \cite{weakandstrong} shows that there exists a critical inverse temperature $\beta_{c}>0$ such that $(\mathcal{Z}^{\beta}_{T}(0;\xi))_{T\geq 0}$ is uniformly integrable if and only if $\beta < \beta_{c}$. To be more specific, one has
\begin{equation} \label{limit_CDP}
    \lim_{T\to\infty}  \mathcal{Z}^{\beta}_{T}(0;\xi) = 
    \begin{cases}
        \mathcal{Z}^{\beta}_{\infty}(0;\xi), &\text{ if } \beta < \beta_{c}\\
        0, &\text{ if } \beta > \beta_{c},
    \end{cases}
\end{equation}
where $\mathcal{Z}^{\beta}_{\infty}(0;\xi)$ is a positive random variable. Consequently, the model is said to be in the weak disorder regime if $\beta \in (0,\beta_{c})$ and be in the strong disorder regime if $\beta \in (\beta_{c},\infty)$. In this subsection, we consider the $L^{\gamma}$-regime of the continuous directed polymer, where $\gamma$ is a real number such that $\gamma > 2$. Here, due to the connection (\ref{connect_SHE_CDP}), the $L^{\gamma}$-regime defined as follows is consistent with  the $L^{\gamma}$-regime of the mollified SHE defined in (\ref{definition_beta_L_N}):
\begin{align} \label{identity_beta_L_N_2}
    \beta_{L^{\gamma}} = \sup\set{\beta > 0: \sup_{T > 0} \mathbb{E}[\mathcal{Z}^{\beta}_{T}(0;\xi)^{\gamma}] < \infty}.
\end{align}

To the best of our knowledge, the distribution of limiting partition function is still unknown. Although it is hard to characterize the distribution of $\mathcal{Z}_{\infty}^{\beta}(0;\xi)$, it is possible to estimate the tail probability of $\mathcal{Z}_{\infty}^{\beta}(0;\xi)$. For the discrete directed polymer, the tail probability of the limiting partition function has been established in \cite{junk2024taildistributionfunctionpartition} when the system belongs to the weak disorder regime. To be more specific, they showed that the tail probability of limiting partition function has polynomial decay, where the exponent is the largest positive number $\overline{\gamma}$ such that the partition function is $L^{\overline{\gamma}}$-bounded, and it is called the critical exponent. 

We believe that the above result may have a continuous counterpart. Due to this reason, our primary concern is the critical exponent of the continuous directed polymer:
\begin{align} \label{definition_gamm_star}
    \gamma^{*}(\beta) := \sup\set{\gamma>1: \sup_{T>0} \mathbb{E}[\mathcal{Z}^{\beta}_{T}(0;\xi)^{\gamma}] < \infty}.
\end{align}
Then, the main theorem proves an estimate for the critical coupling constant $\beta_{L^{\gamma}}$ which relates to the higher fractional moment of the continuous directed polymer. Moreover, a byproduct of this property gives an estimate for the critical exponent $\gamma^{*}(\beta)$ of the polymer in the $L^{2}$-regime. Indeed, in view of (\ref{identity_beta_L_N_2}) and (\ref{definition_gamm_star}), one may notice that they are \emph{formally} inverse to each other (but, more information is required in order to conclude $\gamma^{*}(\beta_{L^{\gamma}}) = \gamma$).

\begin{theorem} \label{main_result_gamma}
Suppose that $d \geq 3$. Then, the critical coupling constant $\beta_{L^{\gamma}}$ and the critical exponent $\gamma^{*}(\beta)$ satisfy the following estimates:
\begin{align} \label{estimate_fractional_beta}
    \frac{\beta_{L^{2}}}{\sqrt{\gamma-1}}\leq \beta_{L^{\gamma}} \leq \frac{\alpha_{\infty,+}}{\sqrt{\gamma-1}} \quad \forall \gamma \in (2,\infty)
\end{align}
and 
\begin{align} \label{estimate_gamma}
    1+\biggbkt{\frac{\beta_{L^{2}}}{\beta}}^{2}\leq \gamma^{*}(\beta) \leq 1+\biggbkt{\frac{\alpha_{\infty,+}}{\beta}}^{2} \quad \forall \beta \in (0,\beta_{L^{2}}),
\end{align}
where $\alpha_{\infty,+}$ is a positive constant defined in Theorem \ref{Main_result_3}.
\end{theorem}

Here, we roughly clarify the proof ideas for Theorem \ref{main_result_gamma}. To do this, recalling the third part in Remark \ref{total_mass}, we know that the Theorem \ref{Main_result_3} holds whenever $\int dx \phi(x)$ is positive. As an application, scaling the function $\phi$ in (\ref{definition_partition_function_1}) properly and considering the corresponding partition function, we may have $\beta_{L^{\gamma}} \approx \widetilde{\beta}_{L^{N}}$. Here, $\widetilde{\beta}_{L^{N}}$ is a critical coupling constant for a scaled partition function. With this estimate, (\ref{estimate_fractional_beta}) then follows from Theorem \ref{Main_result_3}. See Section \ref{section_corollary_main_result_gamma} for the proof of Theorem \ref{main_result_gamma}.

\begin{remark}
Here, we would like to add some comments about our work and various results in the literature.
\begin{enumerate} [label=(\roman*)]
    \item Since we do not have the continuity of $\gamma^{*}(\beta)$, it has to be stressed that the above estimation cannot conclude $\gamma^{*}(\beta_{L^{2}}) = 2$, which is still an ongoing problem for the continuous directed polymer. But, we believe that this identity holds due to the similar work for the discrete directed polymer. We refer to \cite{betaclessthenbetaL2_3} for further details.
    \item Recently, several works \cite{COSCO2022127,GFDDP,Nikos3,me} have proven fluctuations in the entire $L^{2}$-regime of the SHE, KPZ equation, and the partition function and the free energy of the directed polymer in dimensions $d\geq 3$ while avoiding the use of the estimates for the higher moments. The upper bound in (\ref{estimate_gamma}) may indicate the necessity of doing this at least for the continuous case.
    \item For the discrete directed polymer, the critical exponent is also associated with the fluctuation of the partition function. See \cite{junk2024taildistributionfunctionpartition,junk_new_new} for more details. 
    \item The estimate (\ref{estimate_gamma}) could be used to show that the law of $\mathcal{Z}^{\beta}_{\infty}(0;\xi)$ is incompatible with some distributions. For example, combining the following characterization with (\ref{estimate_gamma}): 
    \begin{align} \label{characterization_gamma}
        \gamma^{*}(\beta) = \sup\set{\gamma>1: \mathbb{E}[\mathcal{Z}^{\beta}_{\infty}(0;\xi)^{\gamma}] < \infty} \quad \forall \gamma > 1, \; 0<\beta<\beta_{c},
    \end{align}
    we know that $\mathcal{Z}^{\beta}_{\infty}(0;\xi)$ is not log-normally distributed for every $\beta < \beta_{L^{2}}$. Here, (\ref{characterization_gamma}) is proved by using Fatou's lemma and conditional Jensen's inequality. 
\end{enumerate}

\end{remark}

\subsection{Related literature and additional comments.} \label{section_related_literature}
We now discuss the connection between our work and
various results in the literature. Since we have clarified the results related to the SHE and the directed polymer model in Section \ref{section_overview_main_result} and Section \ref{section_CDP}, respectively, we concentrate on the works in quantum physics in this subsection. Due this reason, we assume that $d = 3$ in the sequel.  
%In view of the Feynman-Kac representation (\ref{representation_FK_0}) of the higher moment

\subsubsection{Many-body delta-Bose gas in three dimensions.} The three-dimensional mollified SHE can be naturally connected with three-dimensional delta-Bose, which is a $N$-particles quantum system with a Hamiltonian that can be formally written by
\begin{align} \label{H_DBG}
    \mathscr{H}_{0}^{\beta,N} := -\frac{1}{2} \Delta_{\Vec{x}} - \beta^{2} \sum_{1\leq i<i'\leq N} \delta_{0}(\Vec{x}(i') - \Vec{x}(i)), \quad \Vec{x} \in \mathbb{R}^{N \times d}.
\end{align}
Here, the state space is $L^{2}(\mathbb{R}^{N \times d})$. To be more precise, the above Hamiltonian describes a system of $N$ quantum particles in dimension three, interacting through
attractive contact potentials. This model has a rich history in quantum physics. We refer to \cite{DBG_summary,ferretti2023contact} for a review of this model.

Now, we clarify the above relation in detail. From the mathematical point of view, the expression (\ref{H_DBG}) is ill-defined because the above Dirac’s delta function do not define an operator in $L^{2}(\mathbb{R}^{N \times d})$. To this aim, it is natural to consider the regularization scheme 
\begin{align} \label{definition_hamiltonian_epsilon}
    \mathscr{H}^{\beta_{\varepsilon},N}_{\varepsilon} := -\mathcal{H}_{\varepsilon}^{\beta_{\varepsilon},N},
    \quad\text{where}\quad\mathcal{H}_{\varepsilon}^{\beta_{\varepsilon},N} := \frac{1}{2} \Delta_{\Vec{x}} + \beta_{\varepsilon}^{2} \sum_{1\leq i<i'\leq N} R_{\varepsilon}(\Vec{x}(i')-\Vec{x}(i)).
\end{align}
Then, due to (\ref{representation_FK_0}), the moment of the mollified SHE can be connected with the semigroup of the approximated Hamiltonian as follows:
\begin{align} \label{moment_and_semigroup}
    \Bigl\langle \varphi, \mathcal{Q}_{\varepsilon;t}^{\beta,N}U_{0}^{\otimes N}\Bigr\rangle_{L^{2}(\mathbb{R}^{N \times d})} = \Bigl\langle \varphi, \exp\sbkt{-t \cdot \mathscr{H}^{\beta_{\varepsilon},N}_{\varepsilon}}U_{0}^{\otimes N}\Bigr\rangle_{L^{2}(\mathbb{R}^{N \times d})} \quad \forall U_{0} \in \mathscr{S}(\mathbb{R}^{d}), \; \varphi \in \mathscr{S}(\mathbb{R}^{N \times d}),
\end{align}
where the self-adjointness of $\mathcal{H}_{\varepsilon}^{\beta_{\varepsilon},N}$ is proved in Lemma \ref{lemma_self_adjoint}, the semigroup on the right-hand side of (\ref{moment_and_semigroup}) is defined in (\ref{definition_semigroup}). 

With the above notation, we now discuss approximated schemes of (\ref{H_DBG}) and their asymptotics. The case of $N = 2$ can be found in \cite[Chapter 1.1]{albeverio2012solvable}. Here, we focus on the case of $N\geq 3$.
\begin{enumerate} [label=(\roman*)]
    \item We begin by discussing the expected limits of the regularization schemes of (\ref{H_DBG}). This problem goes back to the work \cite{Teta}, where the authors considered a different regularization scheme of (\ref{H_DBG}). The unboundedness of an expected limiting quadratic form was proved in \cite[Lemma 6.2]{Teta}. Our consideration is similar to their work. For the approximated scheme $\mathscr{H}^{\beta_{\varepsilon},N}_{\varepsilon}$, recall that we have explained the reason why the left-hand side of (\ref{unboundedness_sub_limit}) is a formal limit of the right-hand side of (\ref{moment_and_semigroup}) in Section \ref{section_definition_sub-limiting_moment}. With the above interpretation, Theorem \ref{Main_result_2} equivalently proves the unboundedness of an expected limiting semigroup of the approximated Hamiltonian  $\mathscr{H}^{\beta_{\varepsilon},N}_{\varepsilon}$. 
    \item As $\varepsilon \to 0^{+}$, an instability is involved in the regularization scheme $\mathscr{H}^{\beta_{\varepsilon},N}_{\varepsilon}$. Simply put, with a specific coupling constant $\beta$, the ground state energy $\inf \mathscr{H}^{\beta_{\varepsilon},N}_{\varepsilon}$ runs towards $-\infty$ as $\varepsilon \to 0^{+}$. This instability is called the Thomas effect in the literature of quantum physics. This can be found in our proof of Lemma \ref{lemma_positivity_beta_L2}. See the third part of Section \ref{section_proof_Main_result_5} for more details. Also, we refer to \cite[Proposition 4.1]{griesemer2023weakness} and \cite[Section C]{griesemer2023weakness_thesis} for further illustration of this effect. 
    \item For the case of $N = 3$, in a recent work \cite{teta2} considering a different regularization scheme of (\ref{H_DBG}), the convergence and the stability were proved.
\end{enumerate}

\subsubsection{Many-body quantum system with short-range interactions.} Due to the scaling property of Brownian motion, the mollified SHE can also be related to a many-body quantum system with a Hamiltonian given by
\begin{align*}
    \mathrm{H}^{\beta,N} := -\mathbf{H}^{\beta,N},
\end{align*}
where $\mathbf{H}^{\beta,N}$ is defined in (\ref{short_range_hamiltonian}). We refer to (\ref{short_range_connection}) and Lemma \ref{short_range_semigroup} for a more explicit explanation of the above relation. Because of this connection, the mollified SHE and the above quantum system share some similar features.
\begin{enumerate} [label=(\roman*)]
    \item For the three-dimensional case, the above quantum system has a crucial property, which is called \emph{the Efimov effect} predicted in \cite{efimov2}. Simply put, this property indicates that, if a three-body system has a certain coupling constant $\beta_{E}$ such that the Hamiltonians of the two-body subsystems do not have negative spectrum, and all of them are resonant, then the Hamiltonian $\mathrm{H}^{\beta_{E},3}$ of this system has an infinite number of negative bound state  energies. In particular, this effect shows that the ground state energy $\inf \mathrm{H}^{\beta_{E},3}$ is negative. The first rigorous proof of this property was given in \cite{efimov0} but we also refer to \cite{TAMURA1991433}. The above phenomenon gives a quantum physics perspective why, for the three-dimensional continuous directed polymer, the $L^{3}$-regime is a proper subset of the $L^{2}$-regime. Roughly speaking, we will show that $\beta_{L^{2}}$ satisfies the above assumptions. Then, one has $\inf \mathrm{H}^{\beta_{L^{2}},3} < 0$, which implies (\ref{positivity_beta_L2}). In this way, $\beta_{L^{3}} < \beta_{L^{2}}$ can be proved by making use of (\ref{positivity_beta_L2}). We refer to Section \ref{section_heuristics_main_result_5} for further illustration.
    \item For the general case $d\geq 3$, $\beta_{L^{N}}$ shares a strong connection with $\beta_{N,+}$, where $\beta_{N,+}$ is given in (\ref{definition_beta_N_+}). This can be seen from (\ref{goal_Main_result_3}) and (\ref{identity_beta2+_=_betaL2}). In fact, we believe they are the same, although we fail to prove this. From quantum physics perspective, $\beta_{N,+}$ plays the role of the best constant of a many-body inequality with short-range interactions. See (\ref{definition_v_form}) for more details. This kind of consideration can also be found in \cite{hoffmann2008many,Guzu,frank2024hardy} and \cite[Section 3]{griesemer2023weakness}. Although the lower bound in the latter work can be extended to a lower bound for $\beta_{N,+}$, this is still \emph{not} enough to conclude (\ref{goal_Main_result_3}) since a lower bound for $\beta_{L^{N}}$ is what we need to prove Theorem \ref{Main_result_3}.  
\end{enumerate}

\section{Proof outline.} \label{section_proofoutline}
Before the proofs of our main results, let us now introduce the proof ideas and the primary tools arising in these proofs. An organization of these proofs will be presented in the last subsection.

\subsection{Heuristics for Theorem \ref{Main_result_1} and Theorem \ref{Main_result_3}.} \label{section_heuristics_main_result_1}
%Q:Where is the difficulty? Explain the difficulty first and then use the proof ideas to answer it.
The main ideas to conclude Theorem \ref{Main_result_1} and Theorem \ref{Main_result_3} rely on the connection between the mollified SHE and a many-body quantum system with short-range interactions. To interpret this, we observe that the following Feynman-Kac representation of the $N$-moments $\mathcal{Q}_{\varepsilon;t}^{\beta,N}U_{0}^{\otimes N}(\Vec{x})$ follows by applying the scaling property of the Brownian motion $\Vec{B}$ in (\ref{representation_FK_0}):
\begin{align} \label{short_range_connection}
    \mathcal{Q}_{\varepsilon;t}^{\beta,N}U_{0}^{\otimes N}(\Vec{x}) 
    = \int_{\mathbb{R}^{N \times d}} d\Vec{z} \mathbf{Q}_{Tt}^{\beta,N} \sbkt{\sqrt{T}\Vec{x},\sqrt{T}\Vec{z}}T^{\frac{N \cdot d}{2}}\cdot U^{\otimes N}_{0}(\Vec{z}),
\end{align}
where we have used the following change of variable for the mollification parameter $\varepsilon$:
\begin{align} \label{T_and_epsilon}
    T := \varepsilon^{-2}.
\end{align}
Here, $U^{\otimes N}_{0}(\Vec{z}):= \bigotimes_{j=1}^{N} U_{0}(\Vec{x})$ and the transition kernel $\mathbf{Q}_{Tt}^{\beta,N}(\Vec{x},\Vec{z})$ is given by
\begin{align} \label{definition_Q_T}
    \mathbf{Q}_{T}^{\beta,N} \sbkt{\Vec{x},\Vec{z}}
    := \mathbf{E}_{\Vec{x}}^{\Vec{B}}\bigbkt{\exp\bkt{\beta^{2}\int_{0}^{T} \sum_{1\leq i<i'\leq N} R(\Vec{B}_{s}(i') - \Vec{B}_{s}(i)) } : \Vec{B}_{T} = \Vec{z}},
\end{align}
where the measure $\mathbf{P}^{\Vec{B}}_{\Vec{x}}\sbkt{d(\Vec{w}(t))_{0\leq t\leq T}: \Vec{w}(T) = \Vec{z}}$ defined in Section \ref{section_structure} is the non-normalized distribution of the $N\times d$-dimensional Brownian bridge. Consequently, due to (\ref{definition_Q_T}) and the Feynman-Kac formula, one can show that the transition kernel $\mathbf{Q}_{T}^{\beta,N}$ is the semigroup of the unbounded operator $\mathbf{H}^{\beta,N}$. Here, $-\mathbf{H}^{\beta,N}$ defined below denotes the Hamiltonian of a quantum many-body system with short-range interactions:
\begin{align} \label{short_range_hamiltonian}
    \mathbf{H}^{\beta,N} := \frac{1}{2} \Delta_{\Vec{x}} + \beta^{2} \sum_{1\leq i<i'\leq N} R(\Vec{x}(i')-\Vec{x}(i)), \quad D(\mathbf{H}^{\beta,N}) = C_{c}^{\infty}(\mathbb{R}^{N \times d}),
\end{align}
Note that the function $R$ described the short-range interaction is defined in (\ref{definition_R}). With these notations, the above connection can be formulated as follows.
\begin{lemma} \label{short_range_semigroup}
Assume that $d\geq 3$ and $N\geq 2$. Then, for every $\beta > 0$ and for each $T>0$, the following identity holds:
\begin{align} 
    \Bigl \langle\varphi,\mathbf{Q}_{T}^{\beta,N} \varphi'\Bigr\rangle_{L^{2}(\mathbb{R}^{N \times d})} = \Bigl\langle\varphi,\exp\sbkt{T \cdot \mathbf{H}^{\beta,N}}  \varphi'\Bigr\rangle_{L^{2}(\mathbb{R}^{N \times d})} \quad \forall \varphi,\varphi' \in \mathscr{S}(\mathbb{R}^{N\times d})
\end{align}    
where the semigroup on the right-hand side is defined in (\ref{definition_semigroup}).
\end{lemma}

Note that $\mathbf{H}^{\beta,N}$ can be further extended to a self-adjoint operator. See Lemma \ref{lemma_self_adjoint} for a more explicit elaboration. For convenience, we still use the same notation $\mathbf{H}^{\beta,N}$ to denote its self-adjoint extension throughout this article.

\paragraph{Proof ideas for (\ref{inequality_main_1_1}) and (\ref{inequality_main_1_2}).} With the above connection, we now introduce the main ingredient to conclude (\ref{inequality_main_1_2}). Note that the idea to prove (\ref{inequality_main_1_1}) is similar. Due to Lemma \ref{short_range_semigroup} and (\ref{definition_semigroup}), it is natural to consider the following critical coupling constant $\beta_{N,+}$ in order to estimate the right-hand side of (\ref{short_range_connection}) from below by an exponential function:
\begin{align} \label{definition_beta_N_+}
    \beta_{N,+} := 
    \inf\set{\beta > 0: \sup \mathbf{H}^{\beta,N} > 0} \quad \forall N\geq 2.
\end{align}
Simply put, $\beta_{N,+}^{2}$ denotes the smallest positive constant such that the ground state energy of the quantum system described by $-\mathbf{H}^{\beta,N}$ is negative. In this way, the condition $\beta > \beta_{N,+}$ provides us an exponential lower bound for $\langle \theta,\mathbf{Q}_{T}^{\beta,N}\theta \rangle_{L^{2}(\mathbb{R}^{N \times d})}$. To be more accurate, by combining Lemma \ref{short_range_semigroup} and (\ref{definition_semigroup}), one can show that the condition $\beta > \beta_{N,+}$ implies the existence of the function $\theta\in C_{c}^{\infty}(\mathbb{R}^{N \times d})$ such that 
\begin{align} \label{estimate_proof_ideas_1}
    \Bigl\langle \theta,\mathbf{Q}_{T}^{\beta,N}\theta \Bigr\rangle_{L^{2}(\mathbb{R}^{N \times d})} \gtrsim \exp(a\cdot T) \quad \forall T>0.
\end{align}
This lower bound is \emph{the major ingredient} to conclude (\ref{inequality_main_1_2}).

However, although we have the above lower bound for $\langle \theta,\exp(T \cdot \mathbf{H}^{\beta,N})\theta \rangle_{L^{2}(\mathbb{R}^{N \times d})}$, there are still two obstacles arising in the proof of (\ref{inequality_main_1_2}). To simplify the following discussion, we assume that $t = 1$. Then, we have
\begin{align} \label{identity_proof_ideas_1}
    \Bigl\langle \varphi,\mathcal{Q}_{\varepsilon;t}^{\beta,N}U_{0}^{\otimes N} \Bigr\rangle_{L^{2}(\mathbb{R}^{N \times d})} = 
    \int_{\mathbb{R}^{N \times d \times 2}} d\Vec{x}d\Vec{z} \varphi(\Vec{x}) \cdot\mathbf{Q}_{T}^{\beta,N} \sbkt{\sqrt{T}\Vec{x},\sqrt{T}\Vec{z}}T^{\frac{N d}{2}} \cdot U_{0}^{\otimes N}(\Vec{z}).
\end{align}

The first difficulty to prove (\ref{inequality_main_1_2}) arises from the above $\sqrt{T}\Vec{x}$ and $\sqrt{T}\Vec{z}$ since 
the right-hand side (\ref{identity_proof_ideas_1}) is not of the form 
$\langle \varphi,\mathbf{Q}_{T}^{\beta\,N}U_{0}^{\otimes N} \rangle_{L^{2}(\mathbb{R}^{N \times d})}$. Moreover, the above $\sqrt{T}\Vec{x}$ and $\sqrt{T}\Vec{z}$ may cause the right-hand side (\ref{identity_proof_ideas_1}) to be more diminutive as $T\to \infty$. This decay obstructs the way to derive an exponential lower bound for the right-hand side of (\ref{identity_proof_ideas_1}). Indeed, since $\Pi^{(N)}$ is a measure zero subset of $\mathbb{R}^{N \times d}$ and
\begin{align} \label{identity_proof_ideas_2}
    &\mathbf{Q}_{T}^{\beta,N} \sbkt{\sqrt{T}\Vec{x},\sqrt{T}\Vec{z}} = \mathbf{E}_{\sqrt{T}\Vec{x}}^{\Vec{B}}\bigbkt{\exp\bkt{\beta^{2}\int_{0}^{T} \sum_{1\leq i<i'\leq N} R(\Vec{B}_{s}(i') - \Vec{B}_{s}(i)) } : \Vec{B}_{T} = \sqrt{T}\Vec{z}},
\end{align}
each pair of the components of Brownian motion $\Vec{B}$ are far apart as $T\to \infty$. However, the function $R$ described the interactions only has compact support. In this way, the expectation on the right-hand side of (\ref{identity_proof_ideas_2}) may become smaller when $T$ is large. 

The second barrier relates to the functions $\varphi(\Vec{x})$ and $U_{0}^{\otimes N}(\Vec{z})$ on the right-hand side of (\ref{identity_proof_ideas_1}). To be more specific, even if the $\sqrt{T}\Vec{x}$ and $\sqrt{T}\Vec{z}$ on the right-hand side of (\ref{identity_proof_ideas_1}) could be replaced by $\Vec{x}$ and $\Vec{z}$, we still do not know whether $\langle \varphi, \mathbf{Q}_{T}^{\beta,N} U_{0}^{\otimes N} \rangle_{L^{2}(\mathbb{R}^{N \times d})}$ has an exponential lower bound, since the condition $\beta>\beta_{N,+}$ \emph{only} gives an exponential lower bound for $\langle \theta, \mathbf{Q}_{T}^{\beta,N} \theta \rangle_{L^{2}(\mathbb{R}^{N \times d})}$ and the relations between $\theta$, $\varphi$, and $U_{0}^{\otimes N}$ are unknown.

To overcome the above issues, our primary tool is to construct a proper non-negative function $\theta(\Vec{x})$ such that \emph{the function $\theta$ has a small support while having the exponential lower bound (\ref{estimate_proof_ideas_1})}. 
\begin{proposition} \label{proposition_small_support}
Assume that $d\geq 3$, $N\geq 2$, and $\beta > \beta_{N,+}$. Recall that the short-range many-particles semigroup $\mathbf{Q}^{\beta,N}_{T}(\Vec{x},\Vec{z})$ is defined in (\ref{short_range_connection}). Then, there exists a function $\theta \in L^{2}(\mathbb{R}^{N \times d})$ and positive numbers $C,\alpha $ such that 
\begin{align} \label{property_bounded_and_compact_support}
    0\leq \theta(\Vec{x}) \leq 1 \quad \forall \Vec{x} \in \mathbb{R}^{N\times d}, \quad \textup{supp}(\theta) \subseteq \bigotimes_{j=1}^{N} B(0,r_{\phi}/2),
\end{align}
and
\begin{align} \label{inequality_C_exp_alpha_T}
    \Bigl\langle \theta, \mathbf{Q}_{T}^{\beta,N} \theta \Bigr\rangle_{L^{2}(\mathbb{R}^{N \times d})} \geq C \exp(\alpha \cdot T) \quad \forall T>0,
\end{align}
where $r_{\phi}$ is a positive number such that $\textup{supp}(R) = B(0,2r_{\phi})$. Here, $R(x)$ is defined in (\ref{definition_R}), where $\phi$ originates from the approximation of the identity defined in Section \ref{section_the_model}.
\end{proposition}
\begin{remark} \label{remark_R}
Notice that, since $\phi \in C_{c}^{\infty}(\mathbb{R}^{d};\mathbb{R}_{+})$ and $\phi$ is symmetric-decreasing, it follows that $\textup{supp}(\phi) = B(0,r_{\phi})$. As a result, we know that $R$ is symmetric-decreasing and
\begin{align} \label{support_R}
    \textup{supp}(R) = B(0,2r_{\phi}).
\end{align}
\end{remark}

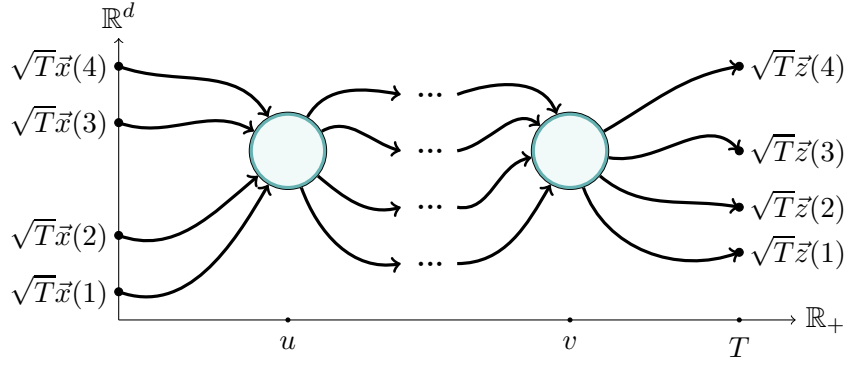
\begin{figure} 
\centering
\begin{tikzpicture} [scale=0.75]
\path[black] (1,-0.12) node[anchor=north]{$u$};
\path[black] (6,-0.12) node[anchor=north]{$v$};
\path[black] (9,-0.12) node[anchor=north]{$T$};

\filldraw[black] (1,0) circle (1pt);
\filldraw[black] (6,0) circle (1pt);
\filldraw[black] (9,0) circle (1pt);

%\filldraw[black] (1,3) circle (8pt);
\node (A) [circle,draw,inner sep=0pt,minimum width=1cm,very thick]  at (1,3) {};
\filldraw[color=teal!60, fill=teal!5, very thick](1,3) circle (0.65);
\node (B) [circle,draw,inner sep=0pt,minimum width=1cm,very thick]  at (6,3) {};
\filldraw[color=teal!60, fill=teal!5, very thick](6,3) circle (0.65);

%\filldraw[black] (6,4) circle (8pt);

\draw[->] (-2,0) -- (10,0) node[anchor=west]{$\mathbb{R}_{+}$};
\draw[->] (-2,0) -- (-2,5) node[anchor=south]{$\mathbb{R}^{d}$};

\filldraw[black] (-2,0.5) circle (2pt) node[anchor=east]{$\sqrt{T} \Vec{x}(1)$};
\filldraw[black] (-2,1.5) circle (2pt)
node[anchor=east]{$\sqrt{T} \Vec{x}(2)$};
\filldraw[black] (-2,3.5) circle (2pt)
node[anchor=east]{$\sqrt{T} \Vec{x}(3)$};
\filldraw[black] (-2,4.5) circle (2pt)
node[anchor=east]{$\sqrt{T} \Vec{x}(4)$};

%\draw[->] [black] plot [smooth] coordinates { (-2,0.5) (-1.3,1.35)  (-0.5,1.5) };
\draw[->] [very thick] (-2,0.5) to[out=-20,in=-120] (A);
\draw[->] [very thick] (-2,1.5) to[out=-20,in=-140] (A);
\draw[->] [very thick] (-2,3.5) to[out=-20,in=150] (A);
\draw[->] [very thick] (-2,4.5) to[out=-20,in=120] (A);
%\draw[->] [black] plot [smooth] coordinates { (-0.5,1.5) (0.75,2.2) (1,3) };

\draw[->] [very thick] (A) to[out=60,in=170] (3,4);
\draw[->] [very thick] (A) to[out=30,in=165] (3,3);
\draw[->] [very thick] (A) to[out=-40,in=-165] (3,2);
\draw[->] [very thick] (A) to[out=-70,in=-170] (3,1);

\node (C) at (3.1,1) [anchor=west]{\textbf{...}};
\node (D) at (3.1,2) [anchor=west]{\textbf{...}};
\node (E) at (3.1,3) [anchor=west]{\textbf{...}};
\node (F) at (3.1,4) [anchor=west]{\textbf{...}};

\draw[->] [very thick] (4,4) to[out=20,in=110] (B);
\draw[->] [very thick] (4,3) to[out=10,in=140] (B);
\draw[->] [very thick] (4,2) to[out=-10,in=-170] (B);
\draw[->] [very thick] (4,1) to[out=0,in=-120] (B);
%\draw[->] [very thick] (A) to[out=-20,in=-140] (2,5);
%\draw[->] [very thick] (A) to[out=-20,in=150] (2,4);
%\draw[->] [very thick] (A) to[out=-20,in=120] (2,3);

%\draw[->] [black] plot [smooth] coordinates { (-2,0.5) (-1.5,2.5) (-1,2) (-0.5,2.5) (0,1.5) } node[anchor=west]{$\textbf{u}(x,t)$};

\filldraw[black] (9,1.2) circle (2pt) node[anchor=west]{$\sqrt{T} \Vec{z}(1)$};
\filldraw[black] (9,2) circle (2pt)
node[anchor=west]{$\sqrt{T} \Vec{z}(2)$};
\filldraw[black] (9,3) circle (2pt)
node[anchor=west]{$\sqrt{T} \Vec{z}(3)$};
\filldraw[black] (9,4.5) circle (2pt)
node[anchor=west]{$\sqrt{T} \Vec{z}(4)$};

\draw[->] [very thick] (B) to[out=-70,in=-160] (9,1.2);
\draw[->] [very thick] (B) to[out=-40,in=170] (9,2);
\draw[->] [very thick] (B) to[out=-10,in=130] (9,3);
\draw[->] [very thick] (B) to[out=30,in=-170] (9,4.5);

%\draw [black] plot [smooth] (-0.75,2.5) -- (-0.25,3.5) -- (0.25,1.5) -- (0.75,2.5) -- (1.25,2) -- (1.75,2.5) -- (2.25,1.5) -- (2.75,2) -- (3.25,0.5) -- (3.75,3)-- (4.25,2) -- (4.75,2.5) -- (5.25,0.5) -- (5.75,2) -- (6.25,2.5) -- (6.75,1) -- (7.25,2)--(7.75,0.5);

\end{tikzpicture}
\caption{The above green region denotes the open ball $B(0,r_{\phi}/2)$, and the above $u$ and $v$ denote the times of the first and last interactions, respectively.}
\label{BD_model_4}
\end{figure}

With the above property, the main ideas are then the following. The key is to create the term $\langle \theta, \mathbf{Q}_{T}^{\beta,N} \theta \rangle_{L^{2}(\mathbb{R}^{N \times d})}$ in $\mathbf{Q}_{T}^{\beta,N} \sbkt{\sqrt{T}\Vec{x},\sqrt{T}\Vec{z}}$ on the right-hand side of (\ref{identity_proof_ideas_1}). Loosely speaking, by making use of the series expansion of $\mathbf{Q}_{T}^{\beta,N} \sbkt{\sqrt{T}\Vec{x},\sqrt{T}\Vec{z}}$ (see Lemma \ref{lemma_series_expansion}), we can bound the rescaled semigroup $\mathbf{Q}_{T}^{\beta,N} \sbkt{\sqrt{T}\Vec{x},\sqrt{T}\Vec{z}}$ on the right-hand side of (\ref{identity_proof_ideas_1}) from below by the sum of all the path integrals in the series expansion of $\mathbf{Q}_{T}^{\beta,N} \sbkt{\sqrt{T}\Vec{x},\sqrt{T}\Vec{z}}$ that force each of the components of Brownian motion $\Vec{B}$ to be in $B(0,r_{\phi}/2)$ at the first and the last interactions. See Figure \ref{BD_model_4} for a graphical clarification about the above idea.

In this way, we can estimate the first and the last interactions by using the properties in (\ref{property_bounded_and_compact_support}) and the follows estimate:
\begin{align*}
    \biggbkt{\sum_{1\leq i<i'\leq N} R(\Vec{x}(i')-\Vec{x}(i)) } \cdot \theta(\Vec{x}) \gtrsim \theta(\Vec{x}),
\end{align*}
where $\Vec{x}$ denotes the position of Brownian motion $\Vec{B}$ at the time of the first interaction or the last interaction. Here, we have used the fact that $R(\Vec{x}(i')-\Vec{x}(i)) \geq R(r_{\phi}) > 0$, where the estimate follows since $R$ is symmetric-decreasing. With the above estimate, by decomposing the rescaled semigroup $\mathbf{Q}_{T}^{\beta,N} \sbkt{\sqrt{T}\Vec{x},\sqrt{T}\Vec{z}}$ at the first and the last interactions, the sum of all the middle of the path integrals produces the term $\langle \theta, \mathbf{Q}_{T'}^{\beta,N} \theta \rangle_{L^{2}(\mathbb{R}^{N \times d})}$, where $T' < T$. Thanks to (\ref{inequality_C_exp_alpha_T}), this term creates an exponential lower bound for the right-hand side of (\ref{identity_proof_ideas_1}).

\paragraph{Proof ideas for (\ref{estimate_betaLN_less_betaL2}) and (\ref{estimate_main_for_betaLN}).} Let us now interpret the heuristics for the estimates (\ref{estimate_betaLN_less_betaL2}) and (\ref{estimate_main_for_betaLN}). Due to this purpose, the overall goal is to show the following estimate:
\begin{align} \label{goal_Main_result_3}
    \frac{\beta_{L^{2}}}{\sqrt{N-1}} \leq \beta_{L^{N}} \leq \beta_{N,+} \leq \frac{\alpha_{\infty,+}}{\sqrt{N-1}} \quad \forall N\geq 3, \quad \text{where } \alpha_{\infty,+} \in (0,\infty),
\end{align}
where $\beta_{N,+}$ is defined in (\ref{definition_beta_N_+}) and $\alpha_{\infty,+}$ only depends on $d$ and $R$. The first and the second inequalities follow from the hypercontractivity from \cite[Theorem 5.1]{GHS} and (\ref{inequality_main_1_1}), respectively. Now, we concentrate on the third inequality here. To do this, we estimate $\beta_{N,+}$ by making use of its variational form, which directly follows from its definition (\ref{definition_beta_N_+}):
\begin{align} \label{definition_v_form}
    \beta_{N,+}^{2} = \inf_{f \in D(\mathbf{H}^{\beta,N}), \; ||f||_{L^{2}(\mathbb{R}^{N \times d})} = 1} \mathcal{I}^{N}[f],
\end{align}
where $D(\mathbf{H}^{\beta,N})$ is given in (\ref{domain_laplace_all}), and the functional $\mathcal{I}^{N}[f]$ is defined by
\begin{align} \label{definition_functional}
    \mathcal{I}^{N}[f] := \biggbkt{\int_{\mathbb{R}^{N \times d}}  \frac{1}{2}|\nabla f(\Vec{x})|^{2} d\Vec{x} } \biggl/ \biggbkt{\sum_{1\leq i< i' \leq N} \int_{\mathbb{R}^{N \times d}}|f(\Vec{x})|^{2} R(\Vec{x}(i') - \Vec{x}(i)) d\Vec{x} }.
\end{align}
Therefore, further reducing the domain of the above variational form to the collection of symmetric functions gives the term $1/\sqrt{N-1}$ on the right-hand side of (\ref{goal_Main_result_3}). Indeed, for any proper symmetric function $f$, the following estimates hold:
\begin{align*}
    \int_{\mathbb{R}^{N \times d}}  \frac{1}{2}|\nabla f(\Vec{x})|^{2} d\Vec{x} \approx N \quad \text{and} \quad \sum_{1\leq i< i' \leq N} \int_{\mathbb{R}^{N \times d}}|f(\Vec{x})|^{2} R(\Vec{x}(i') - \Vec{x}(i)) d\Vec{x} \approx \frac{N\cdot (N-1)}{2}.
\end{align*}
Here, we call a function $f(\Vec{x})$ defined on $\mathbb{R}^{N \times d}$ a symmetric function if the following holds:
\begin{align} \label{symmetricity}
    f(\Vec{x}(1),...,\Vec{x}(N)) = f(\Vec{x}(\sigma(1)),...,\Vec{x}(\sigma(N))) \text{ for all the permutation $\sigma$ of the set $\{1,...,N\}$}.
\end{align}

\paragraph{Proof ideas for (\ref{identity_beta2+_=_betaL2}).}
Finally, we clarify the heuristics for the identity (\ref{identity_beta2+_=_betaL2}), where $\beta_{2,+}$ is defined in (\ref{definition_beta_N_+}). Here, we focus on the estimate $\beta_{2,+} \leq \beta_{L^{2}}$ since $\beta_{L^{2}} \leq \beta_{2,+}$ follows from (\ref{inequality_main_1_1}). The key to conclude the estimate $\beta_{2,+} \leq \beta_{L^{2}}$ is the following characterization of $\beta_{2,+}$:
\begin{align} \label{representation_beta_2_+}
    \beta_{2,+} = ||\mathbf{T}_{1}^{0}||_{L^{2}(\mathbb{R}^{d})}^{-\frac{1}{2}}.
\end{align}        
Here, the right-hand side denotes the operator norm of the Birman-Schwinger operator $\mathbf{T}_{1}^{0}$, where $\mathbf{T}_{1}^{0}$ is defined in (\ref{definition_Birman-Schwinger_operator}). With the above identity, the estimate $\beta_{2,+} \leq \beta_{L^{2}}$ can be deduced immediately from the following characterization and estimate:
\begin{align} \label{characterization_beta_L_2}
    \beta<\beta_{L^{2}} \quad \text{if and only if} \quad \int_{\mathbb{R}^{d}} dz_{0} \beta^{2} R(\sqrt{2}z_{0}) \mathbf{E}_{z_{0}}\bigbkt{\exp\bkt{\beta^{2}\int_{0}^{\infty} R(\sqrt{2}B(u)) du }} < \infty \quad \forall \beta>0.
\end{align}
and
\begin{align*}
    \int_{\mathbb{R}^{d}} dz_{0} \beta^{2}R(\sqrt{2}z_{0})\mathbf{E}_{z_{0}}\bigbkt{\exp\bkt{\beta^{2}\int_{0}^{\infty} R(\sqrt{2}B(u)) du }} 
    \lesssim \sum_{k=0}^{\infty} \beta^{2k}||\mathbf{T}_{1}^{0}||_{L^{2}(\mathbb{R}^{d})}^{k} \quad \forall \beta>0.
\end{align*}
Indeed, it is straightforward to see that $\beta < \beta_{L^{2}}$ for every $\beta < \beta_{2,+}$.

\subsection{Heuristics for Theorem \ref{Main_result_2}.} \label{section_heuristics_main_result_2}
Throughout this subsection, we assume that $d = 3$. To begin with, by taking advantage of (\ref{definition_sub_limiting_path_integral}) and our assumption about the initial datum in Theorem \ref{Main_result_2}, it is enough to prove the divergence of the sub-limiting three moment with Gaussian initial datum in order to conclude Theorem \ref{Main_result_2}. Indeed, (\ref{definition_sub_limiting_path_integral}) shows that for every $N\geq 4$, the sub-limiting $N$-th moment contains all the path integrals of the sub-limiting three moment. Consequently, \emph{we assume that $N = 3$ and $U_{0}(x) = G_{\nu}(x)$} in the sequel. Then, the main ideas of the proof \emph{heavily rely on these two assumptions}.

\paragraph{A short outline.} Before we proceed, we quickly present a short outline of the proof. The overall goal is to establish the following delicate lower bound for the sum of all the sub-limiting path integrals with length $m$ such that $\mathfrak{i}_{j} = 1$ for all $1\leq j\leq m$:
\begin{align*}
    \sum_{(\ell_{1},\ell'_{1},\mathfrak{i}_{1}),...,(\ell_{m},\ell'_{m},\mathfrak{i}_{m}) \in \mathcal{E}^{(N)} \times \{1\}, \; (\ell_{1},\ell_{1}') \neq ...\neq (\ell_{m},\ell_{m}')} \mathscr{I}_{0;t}^{N;(\ell_{1},\ell'_{1},\mathfrak{i}_{1}),...,(\ell_{m},\ell'_{m},\mathfrak{i}_{m})} U_{0}^{\otimes N} (\Vec{x}_{0})  \gtrsim \biggbkt{1.008}^{m}.
\end{align*}
Recalling that description under Remark \ref{remark_after_definition}, a sub-limiting path integral is equal to zero if and only if $\mathbf{i}_{j} = 0$ for some $j$. Hence, a sub-limiting path integral with this feature will not be considered in the proof. To this aim, in the following step 1 and step 2, we will first \emph{exactly evaluate} all the spatial integrals in the above sub-limiting path integrals by taking advantage of the condition $N = 3$ and the Gaussian initial datum. Eventually, the step 3 concentrates on estimating the temporal integrals in the above sub-limiting path integrals. 
%In this way, one can expect that the computation yields a hideous Gaussian kernel.

\paragraph{Step 1: The representation of relative motion.} 
The overall goal of the step 1 and the step 2 is to compute the spatial integrals in the right-hand side of (\ref{definition_sub_limiting_path_integral}). To do this, this step aims to drive a representation for $\mathscr{I}_{0;t}^{N;(\ell_{1},\ell'_{1},\mathfrak{i}_{1}),...,(\ell_{m},\ell'_{m},\mathfrak{i}_{m})} U_{0}^{\otimes N} (\Vec{x}_{0})$ by making use of the condition $N = 3$ and the nonconsecutive interactions. We call this expression \emph{the representation of relative motion}. In the sequel, we adopt the following notations. 
\begin{definition}
Recall that $d = 3$, $N = 3$, and $\mathcal{E}^{(N)}$ is defined in Section \ref{section_structure}. For each $(\ell,\ell')\in \mathcal{E}^{(N)}$ and for every $\Vec{x} \in \mathbb{R}^{N \times d}$, we define 
\begin{align} \label{notation_relative}
    \Vec{x}^{(\ell,\ell')} := \frac{\Vec{x}(\ell')-\Vec{x}(\ell)}{\sqrt{2}}, &\quad \Vec{x}^{(\ell,\ell')^{*}} := \frac{\Vec{x}(\ell')+\Vec{x}(\ell)}{\sqrt{2}}, \quad \text{and} \nonumber\\
    &\Vec{x}^{(\ell,\ell')^{c}} := \Vec{x}(\ell^{c}), \quad\text{where } \ell^{c} \in \{1,...,N\}\setminus\{\ell,\ell'\}.
\end{align}
\end{definition}

With this notation at hand, the sub-limiting path integral admits the following expression.
\begin{lemma} \label{lemma_relative_motion}
Assume that $d = 3$ and $N = 3$. Recall that $U_{0}^{\otimes N}(\Vec{x}) = G^{(N)}_{\nu}(\Vec{x})$ and the sub-limiting path integral $\mathscr{I}_{0;t}^{N;(\ell_{1},\ell'_{1},\mathfrak{i}_{1}),...,(\ell_{k},\ell'_{k},\mathfrak{i}_{k})}U_{0}^{\otimes N}(\Vec{x}_{0})$ is defined in (\ref{definition_sub_limit}). Then, given an arbitrary $t>0$, for every $m\geq 3$ and for each $\bigotimes_{j=1}^{m} (\ell_{j},\ell_{j}',\mathfrak{i}_{j}) \in \bigotimes_{j=1}^{m} (\mathcal{E}^{(N)} \times \{1\} )$ such that $(\ell_{j},\ell_{j}') \neq (\ell_{j-1},\ell_{j-1}')$ for all $2\leq j\leq m$, the sub-limiting path integral $\mathscr{I}_{0;t}^{N;(\ell_{1},\ell'_{1},\mathfrak{i}_{1}),...,(\ell_{k},\ell'_{k},\mathfrak{i}_{k})}U_{0}^{\otimes N}(\Vec{x}_{0})$ has the following representation:
\begin{align} \label{identity_lemma_relative_motion}
    &\mathscr{I}_{0;t}^{N;(\ell_{1},\ell'_{1},\mathfrak{i}_{1}),...,(\ell_{m},\ell'_{m},\mathfrak{i}_{m})}U_{0}^{\otimes N}(\Vec{x}_{0}) \nonumber\\
    &= \int_{v_{j},r_{j}>0, \; v_{0}+r_{1}+...+v_{m-1}+r_{m} <t} \prod_{j=1}^{m} dv_{j-1} dr_{j} \cdot \int_{\mathbb{R}^{(m-1)\times 2\times d}}  \prod_{j=1}^{m-1} d\Vec{z}_{j}^{(\ell_{j},\ell_{j}')^{*}} d\Vec{z}_{j}^{(\ell_{j},\ell_{j}')^{c}} 
    \nonumber\\
    &\biggbkt{G_{v_{0}}(\Vec{z}_{0}^{(\ell_{1},\ell_{1}')}) \cdot \sqrt{\frac{2\pi}{r_{1}}} \cdot G_{v_{0}+r_{1}}(\Vec{z}_{0}^{(\ell_{1},\ell_{1}')^{*}} - \Vec{z}_{1}^{(\ell_{1},\ell_{1}')^{*}}) G_{v_{0}+r_{1}}(\Vec{z}_{0}^{(\ell_{1},\ell_{1}')^{c}} - \Vec{z}_{1}^{(\ell_{1},\ell_{1}')^{c}})}\nonumber\\
    &\biggbkt{ \prod_{j=2}^{m-1} G_{v_{j-1}}\bkt{\frac{\Vec{z}_{j-1}^{(\ell_{j-1},\ell'_{j-1})^{c}}}{\sqrt{2}}- \frac{\Vec{z}_{j-1}^{(\ell_{j-1},\ell'_{j-1})^{*}}}{\sqrt{2}^{2}}} \cdot \sqrt{\frac{2\pi}{r_{j}}} \cdot \nonumber\\
    &\quad\quad G_{v_{j-1}+r_{j}}\bkt{\frac{\Vec{z}_{j-1}^{(\ell_{j-1},\ell'_{j-1})^{c}}}{\sqrt{2}}+ \frac{\Vec{z}_{j-1}^{(\ell_{j-1},\ell'_{j-1})^{*}}}{\sqrt{2}^{2}} - \Vec{z}_{j}^{(\ell_{j},\ell_{j}')^{*}}} G_{v_{j-1}+r_{j}}\bkt{\frac{\Vec{z}_{j-1}^{(\ell_{j-1},\ell'_{j-1})^{*}}}{\sqrt{2}} - \Vec{z}_{j}^{(\ell_{j},\ell_{j}')^{c}}} }\nonumber\\
    &\biggbkt{G_{v_{m-1}}\bkt{\frac{\Vec{z}_{m-1}^{(\ell_{m-1},\ell'_{m-1})^{c}}}{\sqrt{2}}- \frac{\Vec{z}_{m-1}^{(\ell_{m-1},\ell'_{m-1})^{*}}}{\sqrt{2}^{2}}} \cdot \sqrt{\frac{2\pi}{r_{m}}} \cdot \nonumber \\
    &G_{v_{m}+\nu}(0) G_{v_{m-1}+\overline{r}_{m}}\bkt{\frac{\Vec{z}_{m-1}^{(\ell_{m-1},\ell'_{m-1})^{c}}}{\sqrt{2}}+ \frac{\Vec{z}_{m-1}^{(\ell_{m-1},\ell'_{m-1})^{*}}}{\sqrt{2}^{2}}} G_{v_{m-1}+\overline{r}_{m}}\bkt{\frac{\Vec{z}_{m-1}^{(\ell_{m-1},\ell'_{m-1})^{*}}}{\sqrt{2}}}} \quad \forall \Vec{x}_{0} \in \mathbb{R}^{N \times d}\setminus \Pi^{(N)},
\end{align}
where $\Vec{z}_{0} = \Vec{x}_{0}$ and $v_{m} = t-(v_{0}+r_{1}+...+v_{m-1}+r_{m})$. Here,  $\Vec{z}_{0}^{(\ell_{1},\ell'_{1})}$, $\Vec{z}_{0}^{(\ell_{1},\ell'_{1})^{*}}$, and $\Vec{z}_{0}^{(\ell_{1},\ell'_{1})^{c}}$ are defined (\ref{notation_relative}), and $\overline{r}_{m}$ is given by
\begin{align}\label{definition_overline_rm}
    &\overline{r}_{m} := r_{m}+v_{m}+\nu.
\end{align}
\end{lemma}

Briefly speaking, the above representation is obtained by evaluating the Dirac’s delta functions on the right-hand side of (\ref{definition_sub_limiting_path_integral}). To this aim, the main idea is to apply the following change of variables:
\begin{align} \label{relative_map}
    (\Vec{x}_{j}(1),\Vec{x}_{j}(2),\Vec{x}_{j}(3),y_{j},y_{j}') \mapsto (\Vec{z}_{j}^{(\ell_{j},\ell_{j}')},\Vec{z}_{j}^{(\ell_{j},\ell_{j}')^{*}},\Vec{z}_{j}^{(\ell_{j},\ell_{j}')^{c}},w_{j},w'_{j}) \quad 1\leq j\leq m, 
\end{align}
where the right-hand side of the above map is defined as follows:
\begin{align} 
    \Vec{z}_{j}^{(\ell_{j},\ell_{j}')} = \Vec{x}_{j}^{(\ell_{j},\ell_{j}')}, \quad \Vec{z}_{j}^{(\ell_{j},\ell_{j}')^{*}} &= \Vec{x}_{j}^{(\ell_{j},\ell_{j}')^{*}}, \quad 
    \Vec{z}_{j}^{(\ell_{j},\ell_{j}')^{c}} = \Vec{x}_{j}^{(\ell_{j},\ell_{j}')^{c}}, \label{Heuristics_change_of_variable}\\
    w'_{j} &:= \frac{y'_{j}+y_{j}}{\sqrt{2}}, \quad\text{and} \quad
    w_{j} := \frac{y'_{j}-y_{j}}{\sqrt{2}} \quad \forall 1\leq j\leq m.\label{Heuristics_change_of_variable_1}
\end{align}
Before doing so, we first re-express the Gaussian kernels 
\begin{align*}
    G_{u_{j}-u_{j-1}}(\Vec{x}_{j}(i) - \Vec{x}_{j-1}(i)) \cdot \biggbkt{ G_{s_{j}-u_{j-1}}(\Vec{x}_{j-1}(\ell_{j})-y_{j})&
    G_{s_{j}-u_{j-1}}(\Vec{x}_{j-1}(\ell'_{j})-y_{j}')} \\
    &\cdot G_{u_{j}-s_{j}}\bkt{\frac{y'_{j}+y_{j}}{\sqrt{2}}-\frac{\Vec{x}_{j}(\ell'_{j}) + \Vec{x}_{j}(\ell_{j})}{\sqrt{2}}},
\end{align*}
which are on the right-hand side of (\ref{definition_sub_limiting_path_integral}), as follows: 
\begin{align*}
    G_{u_{j}-u_{j-1}}\bkt{\Vec{x}_{j}^{(\ell_{j},\ell_{j}')^{c}}-\Vec{x}_{j-1}^{(\ell_{j},\ell_{j}')^{c}}}\cdot\biggbkt{G_{s_{j}-u_{j-1}}\bkt{\Vec{x}_{j-1}^{(\ell_{j},\ell'_{j})^{*}} -\frac{y_{j}'+y_{j}}{\sqrt{2}}} &G_{s_{j}-u_{j-1}}\bkt{\Vec{x}_{j-1}^{(\ell_{j},\ell'_{j})} -\frac{y_{j}'-y_{j}}{\sqrt{2}}}} \\
    &\cdot G_{u_{j}-s_{j}}\bkt{\frac{y_{j}'+y_{j}}{\sqrt{2}}-\Vec{x}_{j}^{(\ell_{j},\ell'_{j})^{*}}},
\end{align*}
where $i \in \{1,...,3\}\setminus\{\ell_{j},\ell_{j}'\}$. Then, to apply the above change of variables, one has to express the above $\Vec{x}_{j-1}^{(\ell_{j},\ell_{j}')}$, $\Vec{x}_{j-1}^{(\ell_{j},\ell_{j}')^{*}}$, and $\Vec{x}_{j-1}^{(\ell_{j},\ell_{j}')^{c}}$ by using $\Vec{x}_{j-1}^{(\ell_{j-1},\ell_{j-1}')}$, $\Vec{x}_{j-1}^{(\ell_{j-1},\ell_{j-1}')^{*}}$, and $\Vec{x}_{j-1}^{(\ell_{j-1},\ell_{j-1}')^{c}}$. Thanks to $N = 3$ and the property of nonconsecutive interactions, the above $\Vec{x}_{j-1}^{(\ell_{j},\ell_{j}')}$, $\Vec{x}_{j-1}^{(\ell_{j},\ell_{j}')^{*}}$, and $\Vec{x}_{j-1}^{(\ell_{j},\ell_{j}')^{c}}$ admit simple expressions. See (\ref{expression_relative}), (\ref{expression_star}), and (\ref{expression_c}), respectively. In this way, the above change of variable allows us to compute all the Dirac’s delta functions on the right-hand side of (\ref{definition_sub_limiting_path_integral}).

\paragraph{Step 2: Spatial estimate for the sub-limiting path integral.} 
This step aims to estimate the spatial integrals on the right-hand side of (\ref{identity_lemma_relative_motion}).
%Based on the above representation, the spatial integrals in the path integral $\mathscr{I}_{0;t}^{N;(\ell_{1},\ell'_{1},\mathfrak{i}_{1}),...,(\ell_{m},\ell'_{m},\mathfrak{i}_{m})} U_{0}^{\otimes N} (\Vec{x}_{0})$ can be computed through Fourier transform. Although the computation yields hideous Guassian kernels due to our assumption on the initial datum of the mollified SHE (see the right-hand side of (\ref{integral_k_kc_kstar_fourier})), we will estimate these Gaussian kernels from below to conclude the following lower bound.
\begin{proposition} \label{proposition_partII_second_id}
Under the assumptions of Lemma \ref{lemma_relative_motion}, the sub-limiting path integral satisfies the following estimate:
\begin{align}\label{estimate_identity_fourier_3}
    &\mathscr{I}_{0;t}^{N;(\ell_{1},\ell'_{1},\mathfrak{i}_{1}),...,(\ell_{m},\ell'_{m},\mathfrak{i}_{m})}U_{0}^{\otimes N}(\Vec{x}_{0}) \geq C_{t,\nu} \nonumber\\
    &\int_{v_{j},r_{j}>0, \; v_{0}+r_{1}+...+v_{m-1}+r_{m} <t, \; v_{m} = t-(v_{0}+r_{1}+...+v_{m-1}+r_{m})} \biggbkt{\prod_{j=1}^{m} dv_{j-1} dr_{j} } \cdot\frac{1}{v_{0}^{3/2}} \exp\biggbkt{-\frac{|\Vec{x}_{0}|^{2}}{cv_{0}}}\nonumber\\
    &\biggbkt{\prod_{j=1}^{m}\sqrt{\frac{2\pi}{r_{j}}} }
    \cdot\sqrt{2}^{2 \times 3 \times m}  \cdot\biggbkt{\prod_{j=2}^{m-1} G_{v_{j-1}+(\mathbf{t}_{j+1}+3(v_{j-1}+r_{j}))}(0)}G_{v_{m-1}+3(v_{m-1}+\overline{r}_{m})}(0) G_{v_{m}+\nu}(0)\nonumber\\
    &\quad\quad\quad\quad\quad\quad\quad\quad\quad\quad\quad\quad\quad\quad\quad\quad\quad\quad\quad\quad\quad\quad\quad\quad\quad\quad\quad\quad \forall \Vec{x}_{0} \in \mathbb{R}^{N \times d}\setminus \Pi^{(N)}.
\end{align}
Here, $C_{t,\nu}$ and $c$ are positive constants. Moreover, $\overline{r}_{m}$ is defined in (\ref{definition_overline_rm}) and $\mathbf{t}_{j}$ are positive numbers defined as follows:
\begin{align} \label{definition_tj} 
    \mathbf{t}_{j} := 4\sbkt{v_{j-1}\oplus(\mathbf{t}_{j+1}+ 3(v_{j-1}+r_{j}))} \quad \forall 1\leq j\leq m-1 \quad \text{and} \quad\mathbf{t}_{m} := 4\sbkt{v_{m-1}\oplus 3(v_{m-1}+\overline{r}_{m})}.
\end{align}
Here, for each positive numbers $a,b$, $a\oplus b$ is the positive number defined in (\ref{definition_oplus}).
\end{proposition}

Roughly speaking, the proof of the lower bound (\ref{estimate_identity_fourier_3}) includes the following steps: (1) We will first evaluate the Fourier transform of the right-hand side of (\ref{identity_lemma_relative_motion}) with respect to the variables  $\Vec{z}_{0}^{(\ell_{1},\ell'_{1})}$, $\Vec{z}_{0}^{(\ell_{1},\ell'_{1})^{*}}$, and $\Vec{z}_{0}^{(\ell_{1},\ell'_{1})^{c}}$; (2) By re-expressing the above Fourier transform as a recursive expression, we will present an idea that computes the spatial integrals \emph{iteratively}; (3) The previous computation eventually leads to three hideous Gaussian kernels. Then, Proposition \ref{proposition_partII_second_id} follows by estimating these Gaussian kernels from below.

\paragraph{Step 3: Temporal estimate for the sub-limiting path integral.}
Finally, we gives a delicate lower bound for the temporal integrals on the right-hand side of (\ref{estimate_identity_fourier_3}).

\begin{proposition} \label{proposition_lower_bound_sum_of_path_integrals}
Assume that $d = 3$ and $N = 3$. Then, given an arbitrary $t>0$, for each $m \geq 3$, the following estimate holds:
\begin{align} \label{estimate_proposition_lower_bound_sum_of_path_integrals}
    &\sum_{(\ell_{1},\ell'_{1},\mathfrak{i}_{1}),...,(\ell_{m},\ell'_{m},\mathfrak{i}_{m}) \in \mathcal{E}^{(N)} \times \{1\}, \; (\ell_{1},\ell_{1}') \neq ...\neq (\ell_{m},\ell_{m}')} \text{(R.H.S) of (\ref{estimate_identity_fourier_3})}\nonumber\\
    &\quad\quad\quad\quad\quad\quad\quad\quad\quad\quad\quad\quad\quad\geq C \frac{1}{|\Vec{x}_{0}|} \cdot \biggbkt{1.008}^{m-1} \frac{1}{(m-1)(m-2)} \quad \forall\Vec{x}_{0} \in \mathbb{R}^{N \times d}\setminus \Pi^{(N)},
\end{align}
where $C$ is a positive constant that depends on $t$, $\nu$, and $d$.
\end{proposition}

Briefly speaking, the left-hand side of (\ref{estimate_proposition_lower_bound_sum_of_path_integrals}) can be estimated from below by the following:
\begin{align} \label{pattern}
    C_{1}^{m} \cdot \biggbkt{\int_{0<v_{j}<1,\; j = 1,2,...,m} \prod_{j=1}^{m}dv_{j} \cdot \prod_{j=2}^{m}  \frac{1}{v_{j}+v_{j-1}}}.
\end{align}
The trick is to make use of a property of the homogeneous iterated integral. Moreover, the above iterated integral can be estimated recursively from below by $C_{2}^{m}$. But, we need to be careful so that all the estimates are tight enough to ensure $C_{1} \cdot C_{2} \geq 1.008$. This is the main obstacle of the proof.

\begin{remark}
Iterated integrals analogous to the one in (\ref{pattern}) were also considered in two dimensions. See \cite[Lemma 4.15]{chen} and \cite[Lemma 5.4]{caravenna2019moments} for more details.
\end{remark}
% key lemma

\subsection{Heuristics for Proposition \ref{Main_result_5}.} \label{section_heuristics_main_result_5}
Let us now interpret the proof ideas for Proposition \ref{Main_result_5}, where all the details of the proof will be presented in Section \ref{section_proof_Main_result_5}. Due to (\ref{definition_beta_L_N}), we know that $\beta_{L^{N}} \leq \beta_{L^{3}}$ for all $N\geq 4$. Hence, to prove Proposition \ref{Main_result_5}, it is enough to consider the case of $N = 3$, and thus, the overall goal is to show
\begin{align} \label{estimate_L3_3+_L2}
    \beta_{L^{3}} \leq \beta_{3,+} < \beta_{L^{2}},
\end{align}
where $\beta_{3,+}$ is defined in (\ref{definition_beta_N_+}). 

In this way, our proof of Proposition \ref{Main_result_5} relies on two auxiliary properties. First, we will show that the Efimov effect occurred at $\beta = \beta_{L^{2}}$ in our framework, which implies the following.
\begin{lemma}\label{lemma_positivity_beta_L2}
Assume that $d = 3$, $\beta = \beta_{L^{2}}$, and $N = 3$. Recall the definition (\ref{definition_total_energy}). Then
\begin{align} \label{positivity_beta_L2}
    \sup \mathbf{H}^{\beta,N} > 0,
\end{align}
where $\mathbf{H}^{\beta,N}$ is defined in (\ref{short_range_hamiltonian}), and its self-adjointness follows from Lemma \ref{lemma_self_adjoint}. 
\end{lemma}
\begin{remark}
Moreover, we will give an another proof of Lemma \ref{lemma_positivity_beta_L2} without using the Efimov effect. The key is to apply Theorem \ref{Main_result_2}. In this proof, one can see a clearer connection between the moment of the mollified SHE and (\ref{positivity_beta_L2}).
\end{remark}

Next, the continuity of the following map will be proved:
\begin{align} \label{GS_map}
    \beta \in (0,\infty) \mapsto \sup \mathbf{H}^{\beta,N}.
\end{align}
\begin{lemma} \label{lemma_continuity}
Assume that $d \geq 3$ and $N\geq 2$. Recall that $R(x) = \phi*\phi(x)$ and $\phi(x)$ originates from the approximation of the identity defined in Section \ref{section_the_model}. Then, the map defined in (\ref{GS_map}) is continuous. 
\end{lemma}

Finally, let $d = 3$ and $N = 3$, then the following estimate follows from Lemma \ref{lemma_positivity_beta_L2} and Lemma \ref{lemma_continuity}:
\begin{align} \label{positivity_less_than_beta_L2}
    \sup \mathbf{H}^{\beta,N} > 0 \quad \text{for some} \quad \beta < \beta_{L^{2}}.
\end{align}
Therefore, the first inequality in (\ref{estimate_L3_3+_L2}) follows from (\ref{goal_Main_result_3}), and the second equality in (\ref{estimate_L3_3+_L2}) holds due to (\ref{positivity_less_than_beta_L2}) and (\ref{definition_beta_N_+}). This completes the proof of Proposition \ref{Main_result_5}.

\subsection{Organization of the proofs.} \label{section_organization_proof}
The rest of the paper is organized in the following way: 
\begin{enumerate} [label=(\roman*)]
    \item The proofs of Theorem \ref{Main_result_1} and Theorem \ref{Main_result_3} will be Section \ref{section_proof_Main_result_1_3}; 
    \item Theorem \ref{Main_result_2}, Proposition \ref{Main_result_5}, and Theorem \ref{main_result_gamma} will be proved in Section \ref{section_proof_Main_result_2}, Section \ref{section_proof_Main_result_5}, and Section \ref{section_corollary_main_result_gamma}, respectively.
\end{enumerate}

\section{Proofs of Theorem \ref{Main_result_1} and Theorem \ref{Main_result_3}} \label{section_proof_Main_result_1_3}
This section is dedicated to the proofs of Theorem \ref{Main_result_1} and Theorem \ref{Main_result_3}. For convenience, we will first introduce our primary tools, and then complete the proofs of Theorem \ref{Main_result_1} and Theorem \ref{Main_result_3}. The proofs of these tools will then be presented later. See the end of this section for an outline.

To conclude the estimates (\ref{inequality_main_1_1}) and (\ref{inequality_main_1_2}), we start with the following estimates for the rescaled short-range many-body semigroup.
\begin{proposition} \label{proposition_divergence_1}
Assume that $d\geq 3$, $N\geq 2$, and $\beta > \beta_{N,+}$. Recall that the critical coupling constant $\beta_{N,+}$ is defined in (\ref{definition_beta_N_+}) and the many-particles semigroup $\mathbf{Q}^{\beta,N}_{T}(\Vec{x},\Vec{z})$ is defined in (\ref{short_range_connection}). Then, given an arbitrary $t>0$, we have the following estimates: 
\begin{enumerate} [label=(\roman*)]
    \item Let $\Vec{x} \in \mathbb{R}^{N \times d}$ and $f\in C_{c}^{\infty}(\mathbb{R}^{N \times d};\mathbb{R}_{+})$. Then
    \begin{align} \label{inequality_proposition_divergence_1_1}
        \int_{\mathbb{R}^{N \times d}} d\Vec{z} \mathbf{Q}_{Tt}^{\beta,N} \sbkt{\sqrt{T}\Vec{x},\sqrt{T}\Vec{z}}T^{\frac{N \cdot d}{2}} \cdot f(\Vec{z}) \geq C \cdot \set{ T^{2-N \cdot d /2} \cdot \exp(a\cdot T)} \quad\forall T < T_{0}.
    \end{align}
    \item Let $\varphi,f\in C_{c}^{\infty}(\mathbb{R}^{N \times d};\mathbb{R}_{+})$. Then
    \begin{align} \label{inequality_proposition_divergence_1_2}
        \int_{\mathbb{R}^{N \times d \times 2}} d\Vec{x}d\Vec{z} \varphi(\Vec{x}) \cdot\mathbf{Q}_{Tt}^{\beta,N} \sbkt{\sqrt{T}\Vec{x},\sqrt{T}\Vec{z}}T^{\frac{N \cdot d}{2}} \cdot f(\Vec{z}) \geq C' \cdot \set{ T^{2-N \cdot d /2} \cdot \exp(a\cdot T)} \quad\forall T < T_{0}.
    \end{align}
\end{enumerate}
Here, $C,C',a$, and $ T_{0}$ are positive constants such that they all depend on $\beta,d,N,\phi$, and $t$, where $\phi(x)$ comes from the approximation of the identity defined in Section \ref{section_the_model}.  Moreover, $C$ and $C'$ also relate to $(\Vec{x},f)$ and $(\varphi,f)$, respectively.
\end{proposition}

The second primary tool to prove Theorem \ref{Main_result_1} and Theorem \ref{Main_result_3} is the following properties for the critical coupling constants $(\beta_{N,+})_{N\geq 2}$.

\begin{proposition} \label{proposition_betaLN_estimate_and_identity}
Assume that $d\geq 3$. Recall that the critical coupling constant $\beta_{N,+}$ is defined in (\ref{definition_beta_N_+}). Then, the sequence $(\beta_{N,+})_{N\geq 2}$ satisfies the estimate (\ref{goal_Main_result_3}) and the identity (\ref{identity_beta2+_=_betaL2}).
\end{proposition} 

Consequently, the proofs of Theorem \ref{Main_result_1} and Theorem \ref{Main_result_3} are then the following.

\begin{proof} [Proofs of Theorem \ref{Main_result_1} and Theorem \ref{Main_result_3}]
Notice that (\ref{estimate_betaLN_less_betaL2}), (\ref{identity_beta2+_=_betaL2}), and (\ref{estimate_main_for_betaLN}) can be deduced immediately from Proposition \ref{proposition_betaLN_estimate_and_identity}. Moreover, due to (\ref{short_range_connection}) and (\ref{T_and_epsilon}), (\ref{inequality_main_1_1}) and (\ref{inequality_main_1_2}) follows from Proposition \ref{proposition_divergence_1}. Then, the proofs of Theorem \ref{Main_result_1} and Theorem \ref{Main_result_3} are complete.
\end{proof}

\paragraph{Structure of the proofs.}
The proofs of the above propositions are arranged in the following way:
\begin{enumerate} [label=(\roman*)]
    \item Proposition \ref{proposition_divergence_1} will be proved in Section \ref{section_proposition_divergence_1} and Section \ref{section_proposition_small_support};
    \item For the proof of Proposition \ref{proposition_betaLN_estimate_and_identity}, the estimate (\ref{goal_Main_result_3}) and the identity (\ref{identity_beta2+_=_betaL2}) will be proved in Section \ref{section_estimate_betaLN_and_betaN+} and Section \ref{section_charaacter_betaL2}, respectively.
\end{enumerate}

\section{Lower bound for the rescaled many-particles semigroup} \label{section_proposition_divergence_1}
This section aims to prove Proposition \ref{proposition_divergence_1}. As the proof ideas explained in Section \ref{section_heuristics_main_result_1}, the major tools are Proposition \ref{proposition_small_support} and  the following series expansion for the many-body semigroup $\mathbf{Q}_{T}^{\beta,N}(\Vec{x},\Vec{y})$. For the sake of simplicity, we postpone the proof of Proposition \ref{proposition_small_support} until Section \ref{section_proposition_small_support}.

\begin{lemma} \label{lemma_series_expansion}
Suppose that $N \geq 2$ and $d\geq 3$. Recall that the short-range many-particles semigroup $\mathbf{Q}^{\beta,N}_{T}(\Vec{x},\Vec{z})$ is defined in (\ref{short_range_connection}). Then, for each $\Vec{x},\Vec{y}\in \mathbb{R}^{N \times d}$, we have the following decomposition:
\begin{align} \label{definition_simple_series_expansion}
    \mathbf{Q}_{T}^{\beta,N}(\Vec{x},\Vec{y}) = \sum_{k=0}^{\infty} \mathbf{I}^{\beta,N}_{T,k}(\Vec{x},\Vec{y}),
\end{align}
where $\mathbf{I}^{\beta,N}_{T,k}(\Vec{x},\Vec{y})$ is called the short-range interaction path integral of order $k$ from $\Vec{x}$ to $\Vec{y}$ within time $T$, and it is defined as follows:
\begin{align} \label{definition_short_path_integral_order_k}
    \mathbf{I}^{\beta,N}_{T,k}(\Vec{x},\Vec{y}) &:= G_{T}^{(N)}(\Vec{x},\Vec{y})  \quad \text{if } k = 0;\nonumber\\
    \mathbf{I}^{\beta,N}_{T,k}(\Vec{x},\Vec{y}) &:= \int_{0 = t_{0} < t_{1}<...<t_{k} < T, \; \Vec{z}_{0} = \Vec{x}}  \prod_{j=1}^{k} dt_{j} d\Vec{z}_{j} \cdot\nonumber\\
    &\biggbkt{\prod_{j=1}^{k} G^{(N)}_{t_{j}-t_{j-1}}(\Vec{z}_{j-1},\Vec{z}_{j}) \cdot \biggbkt{\beta^{2} \sum_{1\leq i<i'\leq N} R(\Vec{z}_{j}(i')-\Vec{z}_{j}(i)) }} \cdot G_{T-t_{k}}^{(N)}(\Vec{z}_{k},\Vec{y}) \quad \text{if } k\geq 1.
\end{align}
\end{lemma}
\begin{proof}
The series expansion (\ref{definition_simple_series_expansion}) can be deduced immediately from applying the following integration by parts to the exponential function on the right-hand side of (\ref{definition_Q_T}):
\begin{align} \label{identity_fundamental}
    \exp\biggbkt{\int_{0}^{t} f(s) ds} = 1 + \int_{0}^{t} f(s) \exp\biggbkt{\int_{s}^{t} f(u) du} ds.
\end{align}
\end{proof}

Now, we are ready to present the proof of Proposition \ref{proposition_divergence_1}.
\paragraph{Step 1: Simplification.} To prove (\ref{inequality_proposition_divergence_1_1}) and (\ref{inequality_proposition_divergence_1_2}), we may assume that $t = 1$. As a consequence, it is enough to show that
\begin{align} \label{inequality_proposition_divergence_1_1_proof_t_=_1_1}
    \int_{\mathbb{R}^{N \times d}} d\Vec{z}_{0}\mathbf{Q}_{T}^{\beta,N} \sbkt{\sqrt{T}\Vec{x}_{0},\sqrt{T}\Vec{z}_{0}}T^{\frac{N\cdot d}{2}} \cdot f(\Vec{z}_{0}) \geq C \cdot\set{ T^{2-\frac{N \cdot d}{2}} \cdot \exp(a \cdot T) }\quad\forall T < T_{0}
\end{align}
and    
\begin{align} \label{inequality_proposition_divergence_1_1_proof_t_=_1_2}
    \int_{\mathbb{R}^{N \times d \times 2}} d\Vec{x}_{0}d\Vec{z}_{0} \varphi(\Vec{x}_{0}) \cdot\mathbf{Q}_{T}^{\beta,N} \sbkt{\sqrt{T}\Vec{x}_{0},\sqrt{T}\Vec{z}_{0}}T^{\frac{N\cdot d}{2}} \cdot f(\Vec{z}_{0}) \geq C \cdot\set{ T^{2-\frac{N \cdot d}{2}} \cdot \exp(a \cdot T) } \quad\forall T < T_{0}.
\end{align}
Here, $C,C',a,$ and $T_{0}$ are positive constants such that they all depend on $\beta,d,N$, and $\phi$.  Moreover, $C$ and $C'$ also relate to $\{\Vec{x},f\}$ and $\{\varphi,f\}$, respectively. The reason why it is enough to prove (\ref{inequality_proposition_divergence_1_1_proof_t_=_1_1}) and (\ref{inequality_proposition_divergence_1_1_proof_t_=_1_2}) is as follows. Notice that the left-hand sides of (\ref{inequality_proposition_divergence_1_1}) and (\ref{inequality_proposition_divergence_1_2}) can be re-expressed as follows:
\begin{align*}
    \int_{\mathbb{R}^{N \times d}} d\Vec{z}_{0}&\mathbf{Q}_{Tt}^{\beta,N} \sbkt{\sqrt{T}\Vec{x}_{0},\sqrt{T}\Vec{z}_{0}}T^{\frac{N\cdot d}{2}} \cdot f(\Vec{z}_{0})= \int_{\mathbb{R}^{N \times d}} d\Vec{z}_{0}\mathbf{Q}_{Tt}^{\beta,N} \bkt{\sqrt{Tt}\cdot\frac{\Vec{x}_{0}}{\sqrt{t}},\sqrt{Tt}\cdot\Vec{z}_{0}}(Tt)^{\frac{N\cdot d}{2}} \cdot f(\sqrt{t}\cdot \Vec{z}_{0})
\end{align*}
and 
\begin{align*}
    \int_{\mathbb{R}^{N \times d \times 2}} d\Vec{x}_{0}d\Vec{z}_{0} &\varphi(\Vec{x}_{0}) \cdot\mathbf{Q}_{Tt}^{\beta,N} \sbkt{\sqrt{T}\Vec{x}_{0},\sqrt{T}\Vec{z}_{0}}T^{\frac{N\cdot d}{2}} \cdot f(\Vec{z}_{0})\\
    &= \int_{\mathbb{R}^{N \times d \times 2}} d\Vec{x}d\Vec{z}_{0} t^{\frac{N\cdot d}{2}}\varphi(\sqrt{t}\cdot\Vec{x}_{0})\mathbf{Q}_{Tt}^{\beta,N} \sbkt{\sqrt{Tt} \Vec{x}_{0},\sqrt{Tt}\Vec{z}_{0}}(Tt)^{\frac{N\cdot d}{2}} \cdot f(\sqrt{t}\cdot \Vec{z}_{0}).
\end{align*}
As a result, by replacing the $T,\Vec{x}_{0},\varphi(\Vec{x}_{0}),$ and $f(\Vec{z}_{0})$ in (\ref{inequality_proposition_divergence_1_1_proof_t_=_1_1}) and (\ref{inequality_proposition_divergence_1_1_proof_t_=_1_2}) by $Tt,\Vec{x}_{0}/\sqrt{t}, t^{\frac{N\cdot d}{2}} \varphi(\sqrt{t} \cdot \Vec{x}_{0}),$ and $f(\sqrt{t} \cdot \Vec{z}_{0})$, respectively, Proposition \ref{proposition_divergence_1} is an 
immediate consequence of
(\ref{inequality_proposition_divergence_1_1_proof_t_=_1_1}) and (\ref{inequality_proposition_divergence_1_1_proof_t_=_1_2}).

\paragraph{Step2: Decomposition of the path integral.} To prove (\ref{inequality_proposition_divergence_1_1_proof_t_=_1_1}) and (\ref{inequality_proposition_divergence_1_1_proof_t_=_1_2}), we first recall that the ideas mentioned in Section \ref{section_heuristics_main_result_1}. We will decompose each path integral $\mathbf{I}^{\beta,N}_{T,k}(\sqrt{T}\Vec{x}_{0},\sqrt{T}\Vec{z}_{0})$ in $\mathbf{Q}_{T}^{\beta,N} \sbkt{\sqrt{T}\Vec{x}_{0},\sqrt{T}\Vec{z}_{0}}$ into three separate parts by using Proposition \ref{proposition_small_support}. These parts consist of the path integral before the first interaction, the path integral between the first interaction and the last interaction, and the path integral after the last interaction.

To begin with, we bound the first and the last interactions. Notice that the series expansion (\ref{definition_simple_series_expansion}) implies the following estimate:
\begin{align} \label{proof_rescaled_semigroup_1}
    \mathbf{Q}_{T}^{\beta,N}\sbkt{\sqrt{T}\Vec{x}_{0},\sqrt{T}\Vec{z}_{0}} &\geq \int_{u+v < T} dudv \int_{\mathbb{R}^{N\times d \times 2}} d\Vec{x}_{1} d\Vec{z}_{1} \biggbkt{ G^{(N)}_{u}(\sqrt{T}\Vec{x}_{0},\Vec{x}_{1}) \biggbkt{\beta^{2} \sum_{1\leq i < i' \leq N}R(\Vec{x}_{1}(i') - \Vec{x}_{1}(i))}} \nonumber\\
    &\cdot\mathbf{Q}_{T-u-v}^{\beta,N} \sbkt{\Vec{x}_{1},\Vec{z}_{1}} \cdot\biggbkt{\biggbkt{\beta^{2} \sum_{1\leq \ell < \ell' \leq N}R(\Vec{z}_{1}(\ell') - \Vec{z}_{1}(\ell))}G^{(N)}_{v}(\Vec{z}_{1},\sqrt{T}\Vec{z}_{0})}.
\end{align}
Here, we have excluded the terms of $k = 0$ and $k = 1$ in the above lower bound. Moreover, due to the boundedness of $\theta$ in (\ref{property_bounded_and_compact_support}),  it holds that
\begin{align} \label{proof_rescaled_semigroup_2}
    &\text{(R.H.S) of (\ref{proof_rescaled_semigroup_1})}\geq \int_{u+v < T} dudv \int_{\mathbb{R}^{N\times d \times 2}} d\Vec{x}_{1} d\Vec{z}_{1} \nonumber\\
    &\quad\quad\quad\quad\quad\quad\quad\quad\biggbkt{G^{(N)}_{u}(\sqrt{T}\Vec{x}_{0},\Vec{x}_{1}) \biggbkt{\beta^{2} \sum_{1\leq i < i' \leq N}R(\Vec{x}_{1}(i') - \Vec{x}_{1}(i))} }
    \cdot\theta(\Vec{x}_{1})\mathbf{Q}_{T-u-v}^{\beta,N} \sbkt{\Vec{x}_{1},\Vec{z}_{1}} \theta(\Vec{z}_{1})\nonumber\\
    &\quad\quad\quad\quad\quad\quad\quad\quad\quad\quad\quad\cdot\biggbkt{\biggbkt{\beta^{2} \sum_{1\leq \ell < \ell' \leq N}R(\Vec{z}_{1}(\ell') - \Vec{z}_{1}(\ell))}  G^{(N)}_{v}(\Vec{z}_{1},\sqrt{T}\Vec{z}_{0})}.
\end{align}
Now, we interpret how to bound the interactions on the right-hand side of the above estimate. The reason why we need the support of $\theta$ to be so small is that, since all the components of Brownian motion $\Vec{B}$ are closed enough when $\Vec{B}$ in the support of $\theta$, the above interactions can be bounded from below by a constant. To be more specific, since $|\Vec{x}(i')-\Vec{x}(i)| \leq r_{\phi}$ for every $\Vec{x} \in \text{supp}(\theta)$, it follows that
\begin{align*}
    \biggbkt{\beta^{2} \sum_{1\leq i < i' \leq N}R(\Vec{x}(i') - \Vec{x}(i))} \geq \beta^{2} R(r_{\phi}) \frac{N(N-1)}{2} \quad \forall \Vec{x} \in \text{supp}(\theta).
\end{align*}
Here, we have used the fact that $R(r) = \phi*\phi(r)$ is symmetric-decreasing. Also, recalling (\ref{support_R}), we know that $R(r_{\phi}) > 0$. As a consequence, the above estimate implies
\begin{align} \label{proof_rescaled_semigroup_3}
    \text{(R.H.S) of (\ref{proof_rescaled_semigroup_2})}
    \geq C \int_{u+v < T} dudv &\int_{\mathbb{R}^{N\times d \times 2}} d\Vec{x}_{1} d\Vec{z}_{1} \biggbkt{G^{(N)}_{u}(\sqrt{T}\Vec{x}_{0},\Vec{x}_{1}) }\nonumber\\
    &\cdot 
    \theta(\Vec{x}_{1})\mathbf{Q}_{T-u-v}^{\beta,N}(\Vec{x}_{1},\Vec{z}_{1}) \theta(\Vec{z}_{1})\cdot \biggbkt{ G^{(N)}_{v}(\Vec{z}_{1},\sqrt{T}\Vec{z}_{0})}.
\end{align}
Moreover, applying change of variables to $u$ and $v$ to the right-hand side of (\ref{proof_rescaled_semigroup_3}), (\ref{proof_rescaled_semigroup_1}), (\ref{proof_rescaled_semigroup_2}), and (\ref{proof_rescaled_semigroup_3}) yield the following estimate:
\begin{align} \label{proof_rescaled_semigroup_4}
    \mathbf{Q}_{T}^{\beta,N}&\sbkt{\sqrt{T}\Vec{x}_{0},\sqrt{T}\Vec{z}_{0}} \geq C T^{2-N \times d}\int_{u+v < 1} dudv \int_{\mathbb{R}^{N\times d \times 2}} d\Vec{x}_{1} d\Vec{z}_{1} \nonumber\\
    &G^{(N)}_{u}\bkt{\Vec{x}_{0}-\frac{1}{\sqrt{T}}\cdot\Vec{x}_{1}} \cdot\biggbkt{\theta(\Vec{x}_{1}) 
    \mathbf{Q}_{T-T(u+v)}^{\beta,N}(\Vec{x}_{1},\Vec{z}_{1}) 
    \theta(\Vec{z}_{1}) }\cdot G^{(N)}_{v}\bkt{ \frac{1}{\sqrt{T}}\cdot \Vec{z}_{1}-\Vec{z}_{0}}.
\end{align}

Now, we break the spatial integrals on the right-hand side of (\ref{proof_rescaled_semigroup_4}) at time $u$ and $v$. The reason is that, due to Proposition \ref{proposition_small_support}, the middle term on the right-hand side of the above ineqaulity will produce $\langle \theta, \mathbf{Q}_{T-T(u+v)}^{\beta,N} \theta \rangle_{L^{2}(\mathbb{R}^{N \times d})}$, which will further give the exponential lower bounds in (\ref{inequality_proposition_divergence_1_1_proof_t_=_1_1}) and (\ref{inequality_proposition_divergence_1_1_proof_t_=_1_2}). As a result, we now separate the above Gaussian kernels from the above
middle term with respect to the spatial variables $\Vec{x}_{1}$ and $\Vec{z}_{1}$. To do this, we apply the following fact to the above Gaussian kernels:
\begin{align*}
    \abs{\Vec{x}_{0}-\frac{1}{\sqrt{T}}\cdot \Vec{x}_{1}}
    \leq |\Vec{x}_{0}|+1 \quad \forall \Vec{x}_{0} \in \text{supp}(\theta), \; T<T_{0}.
\end{align*}
where $T_{0}$ is a positive constant that makes the above estimate holds. Then, this implies 
\begin{align} \label{lowerbound_for_gaussain_kernel}
    G^{(N)}_{u}\biggbkt{\Vec{x}_{0}-\frac{1}{\sqrt{T}}\cdot\Vec{x}_{1}} \geq G^{(N)}_{u}(|\Vec{x}_{0}|+1).
\end{align}
Therefore, combining (\ref{proof_rescaled_semigroup_4}) and (\ref{lowerbound_for_gaussain_kernel}) implies
\begin{align*} 
    &\mathbf{Q}_{T}^{\beta,N} \sbkt{\sqrt{T}\Vec{x}_{0},\sqrt{T}\Vec{z}_{0}} \nonumber\\`
    &\geq C T^{2-N \times d} \int_{u+v < 1} dudv G^{(N)}_{u}(|\Vec{x}_{0}|+1) \cdot\Bigl\langle \theta, \mathbf{Q}_{T-T(u+v)}^{\beta,N} \theta \Bigr\rangle_{L^{2}(\mathbb{R}^{N \times d})}  \cdot G^{(N)}_{v}(|\Vec{z}_{0}|+1).
\end{align*}

Finally, to prove (\ref{inequality_proposition_divergence_1_1_proof_t_=_1_1}) and (\ref{inequality_proposition_divergence_1_1_proof_t_=_1_2}), it remains to break the above temporal integrals. Hence, we need to separate the middle term on the right-hand side of the above estimate from the above Gaussian kernels with respect to the temporal variables $u$ and $v$. To do this, we first note that
\begin{align*}
    &\mathbf{Q}_{T}^{\beta,N} \sbkt{\sqrt{T}\Vec{x}_{0},\sqrt{T}\Vec{z}_{0}}\\
    &\geq C T^{2-N \times d} \int_{0}^{\frac{1}{4}} du \int_{0}^{\frac{1}{4}}dv G^{(N)}_{u}(|\Vec{x}_{0}|+1) \cdot\Bigl\langle \theta, \mathbf{Q}_{T-T(u+v)}^{\beta,N} \theta \Bigr\rangle_{L^{2}(\mathbb{R}^{N \times d})}  \cdot G^{(N)}_{v}(|\Vec{z}_{0}|+1).
\end{align*}
Here, we have used the fact that $\theta$ is non-negative. As a result, by using the exponential lower bound (\ref{inequality_C_exp_alpha_T}) and the fact that $T-T\cdot(u+v) \geq T/2$, we conclude
\begin{align*}
    &\mathbf{Q}_{T}^{\beta,N} \sbkt{\sqrt{T}\Vec{x}_{0},\sqrt{T}\Vec{z}_{0}}\\
    &\geq C T^{2-N \times d} \int_{0}^{\frac{1}{4}} du \int_{0}^{\frac{1}{4}}dv G^{(N)}_{u}(|\Vec{x}_{0}|+1) \cdot \exp\biggbkt{(T-T(u+v))\cdot \alpha}  \cdot G^{(N)}_{v}(|\Vec{z}_{0}|+1)\\
    &\geq C
    T^{2-N \times d} \biggbkt{\int_{0}^{\frac{1}{4}} du G^{(N)}_{u}(|\Vec{x}_{0}|+1) }\cdot
    \biggbkt{\int_{0}^{\frac{1}{4}}dv G^{(N)}_{v}(|\Vec{z}_{0}|+1) }\cdot C\exp(\alpha \cdot T /2),
\end{align*}
which implies (\ref{inequality_proposition_divergence_1_1_proof_t_=_1_1}) and (\ref{inequality_proposition_divergence_1_1_proof_t_=_1_2}). Therefore, we have completed the proof.

\section{Construction of the small supported function.} \label{section_proposition_small_support}
This subsection is dedicated to the proof of Proposition \ref{proposition_small_support}. The detail will explain in Section \ref{section_trick_shrinking}. The rest of the subsections are the proofs of some auxiliary properties used in Section \ref{section_trick_shrinking}.

\subsection{The trick of shrinking the support.} \label{section_trick_shrinking}
The major difficulty to prove Proposition \ref{proposition_small_support} is to ensure that the function $\theta \in L^{2}(\mathbb{R}^{N \times d})$ has the exponential lower bound (\ref{inequality_C_exp_alpha_T}) while having a small support. Briefly, our strategy to solve this problem is to construct a sequence of non-negative bounded functions $(\theta_{k})_{k\geq 0}$ so that they all satisfy (\ref{inequality_C_exp_alpha_T}) with constants $(C_{k})_{k\geq 0}$ and their supports are shrinking as $k \to \infty$. Hence, there exists large $k\geq 0$ such that $\theta_{k}$ satisfies both of (\ref{property_bounded_and_compact_support}) and (\ref{inequality_C_exp_alpha_T}). 

To begin with, we note that applying the condition $\beta > \beta_{N,+}$ gives the following property. The proof is a straightforward consequence of \cite[Theorem 6.17]{lieb2001analysis} and the standard argument of mollification. To keep the following proof simple, we postpone the proof of Lemma \ref{lemma_psi} until Section \ref{section_proof_abs_value}.
\begin{lemma} \label{lemma_psi}
If $d\geq 3$, $N\geq 2$, and $\beta > \beta_{N,+}$, then there exists a positive number $\alpha$ and $\psi \in C_{c}^{\infty}(\mathbb{R}^{N \times d})$ such that 
\begin{align*}
    \Bigl\langle \psi,\mathbf{H}^{\beta,N}\psi \Bigr\rangle_{L^{2}(\mathbb{R}^{N \times d})} > \alpha, \quad ||\psi||_{L^{2}(\mathbb{R}^{N \times d})} = 1, \quad \text{and} \quad \psi \geq 0.
\end{align*}
\end{lemma}
Now, we are ready to present the proof of Proposition \ref{proposition_small_support}.
\paragraph{Step 1: Establishing the initial state.} The proof starts with showing that there exists a proper separable function $\theta_{0}$ that satisfies the exponential lower bound (\ref{inequality_C_exp_alpha_T}).
\begin{lemma} \label{lemma_theta_0}
If $d\geq 3$, $N\geq 2$, and $\beta > \beta_{N,+}$, then there exists a positive numbers $\ell_{0},C_{0},\alpha$ such that 
\begin{align*}
    \Bigl\langle \theta_{0},\mathbf{Q}_{T}^{\beta,N}\theta_{0} \Bigr\rangle_{L^{2}(\mathbb{R}^{N \times d})} \geq C_{0} \exp(T \cdot \alpha) \quad \forall T>0,
\end{align*}
where $\theta_{0}$ is given by
\begin{align*}
    \theta_{0}(\Vec{x}) := \bigotimes_{j=1}^{N} 1_{B(0,\ell_{0})}(\Vec{x}).
\end{align*}
\end{lemma}
\begin{proof}
To begin with, we notice that Lemma \ref{lemma_psi} gives the following estimate:
\begin{align} \label{lemma_theta_0_estimate_1}
    \Bigl\langle \psi, \exp(T\cdot \mathbf{H}^{\beta,N}) \psi \Bigr\rangle_{L^{2}(\mathbb{R}^{N \times d})} 
    &= \int_{\mathbb{R}} \exp(T \cdot \lambda) \mu_{\psi}(d\lambda)
    \geq \int_{\mathbb{R}} 1_{(\alpha,\infty)}(\lambda) \exp(T \cdot \lambda) \mu_{\psi}(d\lambda)\nonumber\\
    &\geq \mu_{\psi}((\alpha,\infty)) \exp(T \cdot \alpha) \quad \forall T>0,
\end{align}
where $\mu_{\psi}(\lambda) := \langle \psi, \mathbf{P}_{\mathbf{H}^{\beta,N}}(d\lambda)\psi \rangle_{L^{2}(\mathbb{R}^{N \times d})}$ is a probability measure on the real line, and $\psi \in C_{c}^{\infty}(\mathbb{R}^{d};\mathbb{R}_{+})$ is defined in Lemma \ref{lemma_psi}. Moreover, we observe that $\mu_{\psi}((\alpha,\infty)) > 0$. Indeed, if $\mu_{\psi}((\alpha,\infty)) = 0$, then
\begin{align*}
    \Bigl\langle\psi,\mathbf{H}^{\beta,N}\psi \Bigr\rangle_{L^{2}(\mathbb{R}^{N \times d})}= \int_{\mathbb{R}} \lambda \mu_{\psi}(d\lambda) \leq \alpha.
\end{align*}
On the other hand, since $\psi \in C_{c}^{\infty}(\mathbb{R}^{N \times d};\mathbb{R}_{+})$, there exists a positive constant $\ell_{0}$ such that
\begin{align} \label{lemma_theta_0_estimate_2}
    ||\psi||_{L^{\infty}(\mathbb{R}^{N \times d})} \cdot \bigotimes_{j=1}^{N} 1_{B(0,\ell_{0})}(\Vec{x}) \geq \psi(\Vec{x}) \quad \forall \Vec{x} \in \mathbb{R}^{N \times d}.
\end{align}
Therefore, by setting $\theta_{0}(\Vec{x}) := \bigotimes_{j=1}^{N} 1_{B(0,\ell_{0})}(\Vec{x})$ , and by applying Lemma \ref{short_range_semigroup}, (\ref{lemma_theta_0_estimate_1}) and (\ref{lemma_theta_0_estimate_2}) yield
\begin{align*}
    \Bigl\langle \theta_{0},\mathbf{Q}_{T}^{\beta,N}\theta_{0} \Bigr\rangle_{L^{2}(\mathbb{R}^{N \times d})}
    &\geq ||\psi||^{-2}_{L^{\infty}(\mathbb{R}^{N\times d})} \cdot\Bigl\langle \psi,\mathbf{Q}_{T}^{\beta,N}\psi \Bigr\rangle_{L^{2}(\mathbb{R}^{N \times d})} \\
    &= ||\psi||^{-2}_{L^{\infty}(\mathbb{R}^{N\times d})}\cdot
    \Bigl\langle \psi, \exp(T\cdot \mathbf{H}^{\beta,N}) \psi \Bigr\rangle_{L^{2}(\mathbb{R}^{N \times d})} 
    \geq C_{0}\exp(T \cdot \alpha) \quad \forall T>0,
\end{align*}
where we have used (\ref{definition_Q_T}) in the first inequality, and $C_{0}$ is a positive constant defined by $C_{0} := ||\psi||^{-2}_{L^{\infty}(\mathbb{R}^{N\times d})} \mu_{\psi}((\alpha,\infty))$. The proof is complete.
\end{proof}

\paragraph{Step 2: Rearrangement for the path integral.} 
Before we proceed, we recall the definition of the symmetric-decreasing rearrangement of a non-negative function from \cite[Section 3.3]{lieb2001analysis}. The rearrangement of a measurable set $\Gamma \subseteq \mathbb{R}^{d}$ is defined by an open ball $\Gamma^{*} := B(0,r)$, where $r$ is chosen such that $\Gamma$ and $\Gamma^{*}$ have the same Lebesgue measure. With this notation, the rearrangement of the following non-negative measurable function $u$ on $\mathbb{R}^{d}$ is defined by a non-negative radial symmetric function:
\begin{equation*}
    u^{*}(x) := \int_{0}^{\infty} 1_{\{u > \alpha\}^{*}}(x) d\alpha \quad \forall x\in \mathbb{R}^{d}.
\end{equation*}

As already explained at the beginning of this subsection, to prove Proposition \ref{proposition_small_support}, it remains to shrink the support of $\theta_{0}$. The key is a rearrangement inequality for the many-body semigroup $\mathbf{Q}^{\beta,N}_{T}$, which is a generalization of \cite[Theorem 3.8]{lieb2001analysis}. 
\begin{lemma} \label{lemma_GRI_path_integral}
Assume that $d\geq 3$, $N\geq 2$, and $\beta>0$. Let $u_{1},h_{1},...,u_{N},h_{N} \in L^{1}(\mathbb{R}^{d})$ be non-negative functions. Then, it holds that
\begin{align} \label{inequality_Qstar_Q}
    \Bigl\langle \bigotimes_{j=1}^{N} u^{*}_{j}, \mathbf{Q}^{\beta,N}_{T}\bigotimes_{j=1}^{N} h^{*}_{j} \Bigr\rangle_{L^{2}(\mathbb{R}^{N \times d})}
    \geq \Bigl\langle \bigotimes_{j=1}^{N} u_{j}, \mathbf{Q}^{\beta,N}_{T} \bigotimes_{j=1}^{N} h_{j} \Bigr\rangle_{L^{2}(\mathbb{R}^{N \times d})} \quad\forall T\geq 0.
\end{align}
where $\bigotimes_{j=1}^{N} u_{j}(\Vec{x})$ is defined in Section \ref{section_structure} and $\mathbf{Q}^{\beta,N}_{T}(\Vec{x},\Vec{y})$ is the semigroup defined in (\ref{definition_Q_T}).

\end{lemma}
\begin{proof}
We starts with a decomposition of the right-hand side of (\ref{inequality_Qstar_Q}). In view of Lemma \ref{lemma_series_expansion}, we have
\begin{align} \label{lemma_GRI_path_integral_0}
    \Bigl\langle \bigotimes_{j=1}^{N} u_{j}, \mathbf{Q}^{\beta,N}_{T} \bigotimes_{j=1}^{N} h_{j} \Bigr\rangle_{L^{2}(\mathbb{R}^{N \times d})} = \sum_{k=0}^{\infty} \Bigl\langle \bigotimes_{j=1}^{N} u_{j}, \mathbf{I}^{\beta,N}_{T,k} \bigotimes_{j=1}^{N} h_{j} \Bigr\rangle_{L^{2}(\mathbb{R}^{N \times d})}.
\end{align}
Moreover, each term in the above series can decomposed in the following way:
\begin{align} \label{identity_path_integral_decomposition}
   \Bigl\langle \bigotimes_{j=1}^{N} u_{j}&, \mathbf{I}^{\beta,N}_{T,k} \bigotimes_{j=1}^{N} h_{j} \Bigr\rangle_{L^{2}(\mathbb{R}^{N \times d})}\nonumber\\
   &= \sum_{(i_{1},i'_{1}) \in \mathcal{E}^{(N)}}... \sum_{(i_{N},i'_{N}) \in \mathcal{E}^{(N)}} \Bigl\langle \bigotimes_{j=1}^{N} u_{j},\mathbf{I}^{\beta,N;(i_{1},i'_{1}),...,(i_{N},i'_{N})}_{T,k}\bigotimes_{j=1}^{N} h_{j} \Bigr\rangle_{L^{2}(\mathbb{R}^{N \times d})} \quad \text{if } k\geq 1,
\end{align}
where $\mathcal{E}^{(N)}$ is defined in Section \ref{section_structure}. Here, for each $k\geq 1$, $\mathbf{I}^{\beta,N;(i_{1},i'_{1}),...,(i_{N},i'_{N})}_{T,k}(\Vec{x},\Vec{y})$ is the path integral defined as follows:
\begin{align*}
    &\mathbf{I}^{\beta,N;(i_{1},i'_{1}),...,(i_{N},i'_{N})}_{T,k}(\Vec{x},\Vec{y}) := \int_{0 = t_{0} < t_{1}<...<t_{k} < T}  \prod_{j=1}^{k} dt_{j} 
    \int_{\mathbb{R}^{k \times N \times d},\; \Vec{z}_{0} = \Vec{x}} \prod_{j=1}^{k} d\Vec{z}_{j}(1)...d\Vec{z}_{j}(N)
    \nonumber\\
    &\biggbkt{\prod_{j=1}^{k} \biggbkt{\prod_{\nu = 1}^{N} G_{t_{j}-t_{j-1}}(\Vec{z}_{j-1}(\nu)-\Vec{z}_{j}(\nu))} \cdot \beta^{2} R(\Vec{z}_{j}(i_{j}')-\Vec{z}_{j}(i_{j})) }\cdot \biggbkt{\prod_{\nu = 1}^{N} G_{T-t_{k}}(\Vec{z}_{k}(\nu)-\Vec{y}(\nu))}.
\end{align*}
The reason why we have the identity (\ref{identity_path_integral_decomposition}) is nothing but multiplying out all the interactions on the right-hand side of (\ref{definition_short_path_integral_order_k}).

To prove (\ref{inequality_Qstar_Q}), it is enough is to prove the following:
\begin{align} 
    &\Bigl\langle \bigotimes_{j=1}^{N} u^{*}_{j}, \mathbf{I}^{\beta,N;(i_{1},i'_{1}),...,(i_{N},i'_{N})}_{T,k}\bigotimes_{j=1}^{N} h^{*}_{j} \Bigr\rangle_{L^{2}(\mathbb{R}^{N \times d})}
    \geq \Bigl\langle \bigotimes_{j=1}^{N} u_{j}, \mathbf{I}^{\beta,N;(i_{1},i'_{1}),...,(i_{N},i'_{N})}_{T,k} \bigotimes_{j=1}^{N} h_{j} \Bigr\rangle_{L^{2}(\mathbb{R}^{N \times d})} \nonumber\\
    &\quad\quad\quad\quad\quad\quad\quad\quad\quad\quad\quad\quad\quad\quad\quad\quad\quad\quad\quad\quad \quad\quad\quad\quad\quad\quad\quad\quad\quad\quad\quad\quad\quad\text{if } k\geq 1;\label{representation_spatial_integral_3_1}\\
    &\Bigl\langle \bigotimes_{j=1}^{N} u^{*}_{j}, \mathbf{I}^{\beta,N}_{T,k}\bigotimes_{j=1}^{N} h^{*}_{j} \Bigr\rangle_{L^{2}(\mathbb{R}^{N \times d})}
    \geq \Bigl\langle \bigotimes_{j=1}^{N} u_{j}, \mathbf{I}^{\beta,N}_{T,k} \bigotimes_{j=1}^{N} h_{j} \Bigr\rangle_{L^{2}(\mathbb{R}^{N \times d})} \quad \text{if } k = 0. \label{representation_spatial_integral_3_2}
\end{align}
The key is to apply \cite[Theorem 3.8]{lieb2001analysis} to the right-hand sides of (\ref{representation_spatial_integral_3_1}) and (\ref{representation_spatial_integral_3_2}) . In the sequel, we will only concentrate on the proof of (\ref{representation_spatial_integral_3_1}) since (\ref{representation_spatial_integral_3_2}) can be obtained by \cite[Theorem 3.8]{lieb2001analysis} directly. Before applying \cite[Theorem 3.8]{lieb2001analysis} to the right-hand side of (\ref{representation_spatial_integral_3_1}), we show that it satisfies all the conditions in \cite[Theorem 3.8]{lieb2001analysis}. To this aim, we notice that the following functional is the spatial integrals in the quadratic form on the right-hand side of (\ref{representation_spatial_integral_3_1}):
\begin{align} \label{representation_spatial_integral}
    &\mathbf{A}_{t_{1},...,t_{k};k}^{\beta,N}:=\int_{\mathbb{R}^{2\times N\times d}} d\Vec{x}(1)... d\Vec{x}(N)\cdot d\Vec{y}(1)... d\Vec{y}(N) \cdot u_{1}(\Vec{x}(1))...u_{N}(\Vec{x}(N))\cdot h_{1}(\Vec{y}(1))...h_{N}(\Vec{y}(N))\nonumber\\
    &\int_{\mathbb{R}^{N \times k \times d}} \prod_{j=1}^{k} d\Vec{z}_{j}(1)...d\Vec{z}_{j}(N)\cdot
    \biggbkt{\prod_{j=1}^{k}  \biggbkt{\prod_{\nu = 1}^{N} G_{t_{j}-t_{j-1}}(\Vec{z}_{j-1}(\nu)-\Vec{z}_{j}(\nu))} \cdot \beta^{2} R(\Vec{z}_{j}(i_{j}')-\Vec{z}_{j}(i_{j})) }\nonumber\\
    &\cdot \biggbkt{\prod_{\nu = 1}^{N} G_{T-t_{k}}(\Vec{z}_{k}(\nu)-\Vec{y}(\nu))} 
\end{align}
In particular, $\mathbf{A}_{t_{1},...,t_{k};k}^{\beta,N}$ has the same expression as the functional defined in \cite[Theorem 3.8]{lieb2001analysis}. Indeed, the total number of variables and the total number of functions defined on $\mathbb{R}^{d}$ are equal to $K = 2N + kN$ and $M = 2N+k(N+1)+N$, respectively. Furthermore, all the functions $G_{t}(x)$, $R(x)$, $u_{j}(x)$, and $h_{j}(x)$ are non-negative functions in $L^{1}(\mathbb{R}^{d})$. 

Therefore, thanks to the above interpretation and \cite[Theorem 3.8]{lieb2001analysis}, $\mathbf{A}_{t_{1},...,t_{k};k}^{\beta,N}$ can be bounded from above by $(\mathbf{A}_{t_{1},...,t_{k};k}^{\beta,N})^{*}$ defined as follows:
\begin{align} \label{representation_spatial_integral_2}
    &(\mathbf{A}_{t_{1},...,t_{k};k}^{\beta,N})^{*}:=\int_{\mathbb{R}^{2\times N\times d}} d\Vec{x}(1)... d\Vec{x}(N)\cdot d\Vec{y}(1)... d\Vec{y}(N) \cdot u^{*}_{1}(\Vec{x}(1))...u^{*}_{N}(\Vec{x}(N))\cdot h^{*}_{1}(\Vec{y}(1))...h^{*}_{N}(\Vec{y}(N))\nonumber\\
    &\int_{\mathbb{R}^{N \times k \times d}} \prod_{j=1}^{k} d\Vec{z}_{j}(1)...d\Vec{z}_{j}(N)\cdot
    \biggbkt{\prod_{j=1}^{k}  \biggbkt{\prod_{\nu = 1}^{N} G_{t_{j}-t_{j-1}}(\Vec{z}_{j-1}(\nu)-\Vec{z}_{j}(\nu))} \cdot \beta^{2} R(\Vec{z}_{j}(i_{j}')-\Vec{z}_{j}(i_{j})) }\nonumber\\
    &\cdot \biggbkt{\prod_{\nu = 1}^{N} G_{T-t_{k}}(\Vec{z}_{k}(\nu)-\Vec{y}(\nu))}.
\end{align}
Here, we have used the facts that $G_{t}^{*}(x) = G_{t}(x)$ and that $R^{*}(x) = R(x)$, where the second identity follows from Remark \ref{remark_R}. Consequently, the above estimate implies (\ref{representation_spatial_integral_3_1}). Then, the proof is complete.
\end{proof}

\paragraph{Step 3: Construction of the shrinking sequence.} Finally, we present the main idea to shrink the support of $\theta_{0}$. In the following, we will construct a sequence of non-negative bounded functions $(\theta_{k})_{k\geq 0}$ which satisfies (\ref{inequality_C_exp_alpha_T}). In particular, the support of $\theta_{k}$ is shrinking with the rate $c \in (0,1)$ as $k \to \infty$. Consequently, by choosing a large $k \geq 1$ so that the support of $\theta_{k}$ is small enough, we complete the proof of Proposition \ref{proposition_small_support}.
\begin{lemma} \label{lemma_theta_j}
If $d\geq 3$, $N\geq 2$, and $\beta > \beta_{N,+}$, then there exists a sequence of functions $(\theta_{k})_{k\geq 0} \subseteq L^{2}(\mathbb{R}^{N \times d})$ such that for each $k\geq 0$, the following holds:
\begin{align} \label{property_bounded_and_compact_support_j_1}
    0\leq \theta_{k}(\Vec{x}) \leq 1 \quad \forall \Vec{x} \in \mathbb{R}^{N\times d}, \quad \textup{supp}(\theta_{k}) \subseteq B(0,\ell_{k}) \times ... \times B(0,\ell_{k}),
\end{align}
and
\begin{align} \label{inequality_C_exp_alpha_T_j_2}
    \Bigl\langle \theta_{k}, \mathbf{Q}_{T}^{\beta,N} \theta_{k} \Bigr\rangle_{L^{2}(\mathbb{R}^{N \times d})} \geq C_{k} \exp(\alpha \cdot T) \quad \forall T>0,
\end{align}
where $\theta_{0}, \ell_{0}, C_{0}$, and $\alpha$ are defined in Lemma \ref{lemma_theta_0}. Here, for every $k\geq 1$, $\ell_{k}$ and $C_{k}$ are defined as follows:
\begin{align*}
    \ell_{k+1} := c \cdot \ell_{k} \quad \text{and}\quad 
    C_{k+1} := 2^{-2N} C_{k}\quad
    \forall k \geq 0,
\end{align*}
where $c := (1-(1/2)^{d})^{1/d}$. 
\end{lemma}
\begin{proof}
Clearly, $\theta_{0}$ satisfies both of (\ref{property_bounded_and_compact_support_j_1}) and (\ref{inequality_C_exp_alpha_T_j_2}) by using Lemma \ref{lemma_theta_0}. Hence, it is enough to prove Lemma \ref{lemma_theta_j} by mathematical induction, where we assume that $\theta_{k}$ satisfies  (\ref{inequality_C_exp_alpha_T_j_2}) and has the following form as induction hypotheses:
\begin{align} \label{form_theta_j}
    \theta_{k}(\Vec{x}) = \bigotimes_{j=1}^{N}  1_{B(0,\ell_{k})}(\Vec{x}).
\end{align}

To do this, we now suppose that for some $k\geq 0$, $\theta_{k}$ satisfies  (\ref{inequality_C_exp_alpha_T_j_2}) and is of the form (\ref{form_theta_j}). Then, $\langle \theta_{k},\mathbf{Q}_{T}^{\beta,N}\theta_{k} \rangle_{L^{2}(\mathbb{R}^{N \times d})}$ can be decomposed as follows:
\begin{align}\label{lemma_theta_j_0}
    \Bigl\langle \theta_{k},\mathbf{Q}_{T}^{\beta,N}\theta_{k} \Bigr\rangle_{L^{2}(\mathbb{R}^{N \times d})}
    = \sum_{(i_{1},...,i_{N}) \in \{0,1\}^{N}} \sum_{(j_{1},...,j_{N}) \in \{0,1\}^{N}} \Bigl\langle 
    \bigotimes_{\nu=1}^{N} 1_{B^{(k)}_{\nu,i_{\nu}}}
    ,\mathbf{Q}_{T}^{\beta,N}\bigotimes_{\nu=1}^{N} 1_{B^{(k)}_{\nu,j_{\nu}}} \Bigr\rangle_{L^{2}(\mathbb{R}^{N \times d})},
\end{align}
where $B^{(k)}_{\nu,0} := B(0,\ell_{k}/2)$ and  $B^{(k)}_{\nu,1} := B(0,\ell_{k})\setminus B(0,\ell_{k}/2)$ for each $1\leq \nu\leq N$. Consequently, by applying Lemma \ref{lemma_GRI_path_integral} with $u_{\nu
}(x) = 1_{B^{(k)}_{\nu,i_{\nu}}}(x)$ and $h_{\nu
}(x) = 1_{B^{(k)}_{\nu,j_{\nu}}}(x)$, we obtain the following estimate:
\begin{align}\label{lemma_theta_j_1}
    &\Bigl\langle \theta_{k+1},\mathbf{Q}_{T}^{\beta,N}\theta_{k+1} \Bigr\rangle_{L^{2}(\mathbb{R}^{N \times d})} \geq \Bigl\langle \bigotimes_{\nu=1}^{N} 1_{B^{(k)}_{\nu,i_{\nu}}},\mathbf{Q}_{T}^{\beta,N}\bigotimes_{\nu=1}^{N} 1_{B^{(k)}_{\nu,j_{\nu}}} \Bigr\rangle_{L^{2}(\mathbb{R}^{N \times d})}\nonumber\\
    &\quad\quad\quad\quad\quad\quad\quad\quad\quad\quad\quad\quad\quad\quad\quad\quad \forall (i_{1},...,i_{N}) \in \{0,1\}^{N}, \; (j_{1},...,j_{N}) \in \{0,1\}^{N},
\end{align}
where we have used the following facts:
\begin{align*}
    (B^{(k)}_{\nu,1})^{*} = B(0,c \cdot \ell_{k}) \quad \text{and} \quad (B^{(k)}_{\nu,0})^{*} = B^{(k)}_{\nu,0} \subseteq B(0,c \cdot \ell_{k}) \quad \forall 1\leq \nu \leq N.
\end{align*}
Therefore, since $\theta_{k}$ satisfies (\ref{inequality_C_exp_alpha_T_j_2}), the following estimate follows from (\ref{lemma_theta_j_0}) and (\ref{lemma_theta_j_1}):
\begin{align*}
    \Bigl\langle \theta_{k+1},\mathbf{Q}_{T}^{\beta,N}\theta_{k+1} \Bigr\rangle_{L^{2}(\mathbb{R}^{N \times d})} &\geq 2^{-2N}\Bigl\langle \theta_{k},\mathbf{Q}_{T}^{\beta,N}\theta_{k} \Bigr\rangle_{L^{2}(\mathbb{R}^{N \times d})} \geq 2^{-2N} C_{k} \cdot \exp(\alpha \cdot T) \\
    &= C_{k+1} \cdot \exp(\alpha \cdot T)\quad \forall T>0.
\end{align*}
Then, the proof is complete.
\end{proof}

\subsection{Proof of Lemma \ref{lemma_psi}.} \label{section_proof_abs_value}
Now, we quickly sketch the proof of Lemma \ref{lemma_psi}. Let $d\geq 3$, $N\geq 2$, and $\beta > \beta_{N,+}$. Then, since $\sup \mathbf{H}^{\beta,N} > 0$, there exists a positive number $\alpha_{0}$ and $\varphi \in C_{c}^{\infty}(\mathbb{R}^{N \times d})$ such that 
\begin{align*}
    \Bigl\langle \varphi,\mathbf{H}^{\beta,N}\varphi \Bigr\rangle_{L^{2}(\mathbb{R}^{N \times d})} \geq \alpha_{0}.
\end{align*}
To prove Lemma \ref{lemma_psi}, it is natural to consider $\psi(\Vec{x}) := |\varphi(\Vec{x})|$. Due to \cite[Theorem 6.17]{lieb2001analysis}, we know that $\psi \in H^{1}(\mathbb{R}^{N \times d})$ and
\begin{align*}
    &-\frac{1}{2}\int_{\mathbb{R}^{N \times d}} d\Vec{x} |\nabla \psi(\Vec{x})|^{2} + \beta^{2} \sum_{1\leq i<j \leq N} \int_{\mathbb{R}^{N \times d}} d\Vec{x} |\psi(\Vec{x})|^{2} 
    \cdot R(\Vec{x}(j) - \Vec{x}(i))\geq \alpha_{0}.
\end{align*}
Consequently, the remaining part is nothing but applying mollification to $\psi(\Vec{x})$ in order to find a non-negative function $\theta \in C_{c}^{\infty}(\mathbb{R}^{N \times d};\mathbb{R}_{+})$ such that $\Bigl\langle \theta,\mathbf{H}^{\beta,N}\theta \Bigr\rangle_{L^{2}(\mathbb{R}^{N \times d})} \geq \alpha_{0}/2$. Since this is a standard argument, we skip the proof of this part. Therefore, we complete the proof of Lemma \ref{lemma_psi} with $\alpha = \alpha_{0}/2$.

\section{Estimates for $\beta_{L^{N}}$ and $\beta_{N,+}$} \label{section_estimate_betaLN_and_betaN+}
To conclude Proposition \ref{proposition_betaLN_estimate_and_identity}, the prime objective of this section is to prove the estimate (\ref{goal_Main_result_3}).

\paragraph{Step 1.}  To this aim, we begin by showing the first inequality in (\ref{goal_Main_result_3}). Briefly speaking, it can be deduced immediately from the hypercontractivity \cite[Theorem 5.1]{GHS}.
\begin{lemma} \label{lemma_HC_bound}
Assume that $d\geq 3$ and $\beta>0$. Let $\gamma \geq 2$ be a real number. Then, the following estimate holds: 
\begin{align*} 
    ||\mathcal{Z}^{\beta}_{T}(0;\xi)||_{L^{\gamma}} \leq ||\mathcal{Z}^{\beta \cdot \sqrt{\gamma-1}}_{T}(0;\xi)||_{L^{2}} \quad \forall T \geq 0.
\end{align*}
In particular, it holds that
\begin{align*}
    \frac{\beta_{L^{2}}}{\sqrt{\gamma-1}} \leq \beta_{L^{\gamma}} \quad \forall \gamma \geq 2.
\end{align*}
\end{lemma}
\begin{proof}
To prove the above moment bound, we first notice that $\mathcal{Z}^{\beta}_{T}(0;\xi)$ has the following Wiener chaos expansion:
\begin{align} \label{chaos}
    \mathcal{Z}^{\beta}_{T}(0;\xi) = 1+\sum_{k=1}^{\infty} \beta^{k} \int_{\mathbb{R}^{k \times d \times 2}}
    &\int_{0 = t_{0} < t_{1}<...<t_{u} < T} 
    \prod_{j=1}^{k} dt_{j} dy_{j} dq_{j} \nonumber\\
    &\cdot G_{t_{j} - t_{j-1}}(y_{j} - y_{j-1}) \phi(q_{j} - y_{j}) \xi(q_{j},t_{j}), \quad \text{where } y_{0} = 0,
\end{align}
where the above representation can be proved by appyling Itô's formula to $\mathcal{Z}^{\beta}_{T}(0;\xi)$. Therefore, applying \cite[Theorem 5.1]{GHS} with $A := (\gamma-1)^{-1/2} I$, $p := 2$, and $q := \gamma$ to $\mathcal{Z}^{\beta \cdot \sqrt{\gamma-1}}_{T}(0;\xi)$, where $I$ is the identity map on the Gaussian Hilbert space generated by the white noise, we conclude the following estimate:
\begin{align*}
    ||\mathcal{Z}^{\beta}_{T}(0;\xi)||_{L^{\gamma}}= ||\Gamma(A) \mathcal{Z}^{\beta \cdot \sqrt{\gamma-1}}_{T}(0;\xi)||_{L^{\gamma}} \leq ||\mathcal{Z}^{\beta \cdot \sqrt{\gamma-1}}_{T}(0;\xi)||_{L^{2}}.
\end{align*}
Then, it follows that if $\beta \cdot \sqrt{\gamma-1} < \beta_{L^{2}}$, then $\beta < \beta_{L^{\gamma}}$. The proof is complete.
\end{proof}

\paragraph{Step 2.} We will now conclude the second inequality in (\ref{goal_Main_result_3}). To this aim, we show a more general estimate as follows:
\begin{align} \label{ineqaulity_betaLN_betaN+}
    \beta_{L^{N}} \leq \beta_{N,+} \quad \forall N\geq 2.
\end{align}
We first observe that $\beta_{L^{N}}$ has the following characterization:
\begin{align}  \label{characterization_beta_LN}
    \beta_{L^{N}} = \sup\set{\beta > 0: \mathbf{E}_{\Vec{0}}^{\Vec{B}}\bigbkt{\exp\bkt{\beta^{2}\int_{0}^{\infty} \sum_{1\leq i<i'\leq N} R(\Vec{B}_{s}(i') - \Vec{B}_{s}(i)) ds }} < \infty} \quad \forall N\geq 2,
\end{align}
where (\ref{characterization_beta_LN}) is proved by showing that $\sup_{T\geq 0}\mathbb{E}[(\mathcal{Z}_{T}^{\beta})^{N}]$ is equal to the above Brownian exponential functional for every $\beta>0$. Here, we have used Fatou's lemma and the fact that $(\mathcal{Z}_{\cdot}^{\beta})^{N}$ is a submartingale. Then, (\ref{ineqaulity_betaLN_betaN+}) is a straightforward application of Proposition \ref{proposition_divergence_1} and the following fact:
\begin{align*}
    \int_{\mathbb{R}^{N \times d}} d\Vec{z} \mathbf{Q}_{T}^{\beta,N} \sbkt{\Vec{0}&,\sqrt{T}\Vec{z}}T^{\frac{N d}{2}} \cdot f(\Vec{z})\\
    &\leq ||f||_{L^{\infty}(\mathbb{R}^{N \times d})}\cdot \mathbf{E}_{\Vec{0}}^{\Vec{B}}\bigbkt{\exp\bkt{\beta^{2}\int_{0}^{\infty} \sum_{1\leq i<j\leq N} R(\Vec{B}_{s}(j) - \Vec{B}_{s}(i)) ds }} < \infty \\
    &\quad\quad\quad\quad\quad\quad\quad\quad\quad\quad\quad\quad\quad\quad\quad\forall f\in C_{c}^{\infty}(\mathbb{R}^{N \times d};\mathbb{R}_{+}), \; T>0, \; \beta < \beta_{L^{N}}.
\end{align*}
Indeed, Proposition \ref{proposition_divergence_1} implies that $\beta \leq \beta_{N,+}$ for every $\beta<\beta_{L^{N}}$. The proof is complete.

\paragraph{Step 3.}
To complete the proof of (\ref{goal_Main_result_3}), we prove the third inequality in (\ref{goal_Main_result_3}). As already explained in Section \ref{section_heuristics_main_result_1}, the proof relies on the representation (\ref{definition_v_form}) and the symmetricity (\ref{symmetricity}). In the following, $\mathcal{S}(\mathbb{R}^{N \times d})$ denotes the collection of all symmetric functions defined on $\mathbb{R}^{N \times d}$.

\paragraph{Step 3-1.} To begin with, we notice that the following estimate follows from the variational representation (\ref{definition_v_form}):
\begin{align} \label{estimate_symmtric}
    \beta_{N,+}^{2} \leq (\beta_{N,+}^{sym})^{2} = \frac{\alpha_{N,+}^{2}}{N-1},
\end{align}
where 
\begin{align} \label{definition_beta_sym}
    (\beta_{N,+}^{sym})^{2}:=\inf\set{ \mathcal{I}^{N}[f] :f \in \mathscr{S}(\mathbb{R}^{N\times d}) \cap \mathcal{S}(\mathbb{R}^{N \times d}), \quad ||f||_{L^{2}(\mathbb{R}^{N \times d})} = 1}, 
\end{align}
\begin{align} \label{definition_alpha_N}
    \alpha_{N,+}^{2} := \inf\set{ \mathcal{J}^{N}[f] :f \in \mathscr{S}(\mathbb{R}^{N\times d}) \cap \mathcal{S}(\mathbb{R}^{N \times d}), \quad ||f||_{L^{2}(\mathbb{R}^{N \times d})} = 1},
\end{align}
and
\begin{align*}
    \mathcal{J}^{N}[f] := \biggbkt{\int_{\mathbb{R}^{N \times d}} \frac{1}{2}|\nabla_{\Vec{x}(1)}f(\Vec{x})|^{2}+\frac{1}{2}|\nabla_{\Vec{x}(2)}f(\Vec{x})|^{2} d\Vec{x}} \biggl/ \biggbkt{\int_{\mathbb{R}^{N \times d}} |f(\Vec{x})|^{2}R(\Vec{x}(2) - \Vec{x}(1))}.
\end{align*}
Here, the functional $\mathcal{I}^{N}[f]$ in (\ref{definition_beta_sym}) is defined in (\ref{definition_functional}). The reason why we have the equality in (\ref{estimate_symmtric}) is the following:
\begin{align*}
    \int_{\mathbb{R}^{N \times d}} \frac{1}{2}|\nabla f(\Vec{x})|^{2}d\Vec{x}
    = \frac{N}{2}\cdot \biggbkt{\int_{\mathbb{R}^{N \times d}} \frac{1}{2}|\nabla_{\Vec{x}(1)}f(\Vec{x})|^{2}+\frac{1}{2}|\nabla_{\Vec{x}(2)}f(\Vec{x})|^{2} d\Vec{x}}
\end{align*}
and 
\begin{align*}
    \sum_{1\leq i<i'\leq N}\int_{\mathbb{R}^{N \times d}} |f(\Vec{x})|^{2}R(\Vec{x}(i') - \Vec{x}(i)) d\Vec{x} 
    = \frac{N(N-1)}{2} \cdot \biggbkt{\int_{\mathbb{R}^{N \times d}} |f(\Vec{x})|^{2}R(\Vec{x}(2) - \Vec{x}(1)) d\Vec{x}}.
\end{align*}

\paragraph{Step 3-2.}
Now, we explain the reason why the sequence $(\alpha_{N,+})_{N\geq 2}$ is increasing. To do this, we rewrite $\alpha_{N,+}$ as follows:
\begin{align*} 
    \alpha_{N,+} = &\inf\set{\alpha>0: -\bkt{\int_{\mathbb{R}^{N \times d}} \frac{1}{2}|\nabla_{\Vec{x}(1)}f(\Vec{x})|^{2}+\frac{1}{2}|\nabla_{\Vec{x}(2)}f(\Vec{x})|^{2} d\Vec{x}} \\
    &+ \alpha^{2} \bkt{\int_{\mathbb{R}^{N \times d}} |f(\Vec{x})|^{2}R(\Vec{x}(2) - \Vec{x}(1)) d\Vec{x} } > 0 \quad\text{for some } f \in \mathscr{S}(\mathbb{R}^{N\times d}) \cap \mathcal{S}(\mathbb{R}^{N \times d})}.
\end{align*}
As a result, let $\alpha > \alpha_{N+1,+}$, then there exists $f_{N+1} \in \mathscr{S}(\mathbb{R}^{(N+1)\times d}) \cap \mathcal{S}(\mathbb{R}^{(N+1) \times d})$ such that 
\begin{align*}
    \int_{\mathbb{R}^{d}}dz \cdot\biggbkt{ -\bkt{\int_{\mathbb{R}^{N \times d}} \frac{1}{2}|\nabla_{\Vec{x}(1)}f_{N+1}(\Vec{x},z)|^{2}&+\frac{1}{2}|\nabla_{\Vec{x}(2)}f_{N+1}(\Vec{x},z)|^{2} d\Vec{x}} \\
    &+ \alpha^{2} \bkt{\int_{\mathbb{R}^{N \times d}} |f_{N+1}(\Vec{x},z)|^{2}R(\Vec{x}(2) - \Vec{x}(1)) d\Vec{x} } }> 0.
\end{align*}
Hence, there exists $z \in \mathbb{R}^{d}$ such that the above integrand with respect to $dz$ is positive. In particular, since $f_{N+1}(\cdot,z) \in \mathscr{S}(\mathbb{R}^{N\times d}) \cap \mathcal{S}(\mathbb{R}^{N \times d})$, it follows that $\alpha > \alpha_{N,+}$. This implies that $\alpha_{N,+} \leq \alpha_{N+1,+}$.

\paragraph{Step 3-3.} Finally, in view of the above property of $(\alpha_{N,+})_{N\geq 2}$, (\ref{estimate_symmtric}) shows that
\begin{align}\label{definition_alpha_infinity}
    \beta_{N,+}^{2} \leq \frac{\alpha_{\infty,+}^{2}}{N-1}, \quad \text{where } \alpha_{\infty,+} := \lim_{N\to\infty} \alpha_{N,+}.
\end{align}
This conclude the third inequality in (\ref{goal_Main_result_3}). Hence, it remains to clarify the reason why $\alpha_{\infty,+}$ is finite. To this aim, we observe that $\alpha_{N,+}^{2} \leq \mathcal{J}^{N}[(G_{1}^{(N)})^{1/2}]$ for every $N\geq 2$. Furthermore, notice that $\mathcal{J}^{N}[(G_{1}^{(N)})^{1/2}]$ is independent of $N \geq 2$ since 
\begin{align*}
    \int_{\mathbb{R}^{(N-2) \times d}} d\Vec{x}(3)...d\Vec{x}(N) G_{1}^{(N-2)}(\Vec{x}(3),...,\Vec{x}(N)) = 1.
\end{align*}
Therefore, we know that $\alpha_{\infty,+}^{2}$ can be bounded from above by the finite number $\mathcal{J}^{2}[(G_{1}^{(2)})^{1/2}]$. Hence, the proof is complete.

\section{Variational characterization of $\beta_{L^{2}}$} \label{section_charaacter_betaL2}
To complete the proof of Proposition \ref{proposition_betaLN_estimate_and_identity}, this section is dedicated to the proof of the identity (\ref{identity_beta2+_=_betaL2}). Recall that $\beta_{2,+}$ is defined in (\ref{definition_beta_N_+}). Note that $\beta_{L^{2}} \leq \beta_{2,+}$ has been proved in (\ref{ineqaulity_betaLN_betaN+}). Hence, the goal of this section aims to show 
\begin{align} \label{estimate_beta2+_less_betaL2}
    \beta_{2,+} \leq \beta_{L^{2}}.
\end{align}
Before we proceed, let us first introduce our major tools. These auxiliary results will be proved in the later subsections.

To prove (\ref{estimate_beta2+_less_betaL2}),
the main idea is to connect both of $\beta_{L^{2}}$ and $\beta_{2,+}$ with the Birman-Schwinger operator. Then, the first step is to re-express $\beta_{L^{2}}$ as follows.
\begin{lemma} \label{lemma_characterization_beta_L_2}
Assume that $d\geq 3$. Recall that $R$ is defined in (\ref{definition_R}). Then, (\ref{characterization_beta_L_2}) holds.
\end{lemma}

Next, we develop an integral representation of the right-hand side of (\ref{characterization_beta_L_2}). Before doing so, we introduce the following bounded operators on $L^{2}(\mathbb{R}^{d})$, which are called the Birman-Schwinger operators at the level $\lambda$ with the potential function $V_{\beta}(x)$:
\begin{align} \label{definition_Birman-Schwinger_operator}
    &\mathbf{T}^{\lambda}_{\beta}h(x) := \mathfrak{V}_{\beta}(x) \int_{\mathbb{R}^{d}} \mathcal{G}^{\lambda}(x-y) \mathfrak{V}_{\beta}(y) h(y) dy \quad \forall \lambda \geq 0, \; h\in L^{2}(\mathbb{R}^{d}), \; \beta>0.
\end{align}
Here, the boundedness will be proved in Lemma \ref{lemma_properties_Birman-Schwinger_operators}, and we define
\begin{align} \label{definition_V_and_Yukawa}
    V_{\beta}(x) := \beta^{2} R(\sqrt{2}x),\quad \mathfrak{V}_{\beta}(x) := V_{\beta}(x)^{\frac{1}{2}}, \quad \text{and} \quad
    \mathcal{G}^{\lambda}(x) := \int_{0}^{\infty} \exp(-\lambda \cdot u) G_{u}(x) du,
\end{align}
where $\mathcal{G}^{\lambda}$ is the $d$-dimensional Yukawa potential and the function $R$ is defined in (\ref{definition_R}). With the above notations, the major step to prove (\ref{estimate_beta2+_less_betaL2}) is the following characterization of $\beta_{2,+}$.

\begin{lemma} \label{lemma_representation_beta_2_+}
Assume that $d\geq 3$. Then, $\mathbf{T}_{1}^{0}$ is a bounded operator on $L^{2}(\mathbb{R}^{d})$. In particular, $\beta_{2,+}$ admits the representation (\ref{representation_beta_2_+}).
\end{lemma}
%Birman-Schwinger principle For a brief explanation about the principle, we refer to the proof of \cite[Theorem 12.4]{lieb2001analysis}.

Now, we are ready to prove the identity (\ref{identity_beta2+_=_betaL2}).

\begin{proof}[Proof of (\ref{identity_beta2+_=_betaL2})]
Since $\beta_{L^{2}} \leq \beta_{2,+}$ has been proved in (\ref{ineqaulity_betaLN_betaN+}), it remains to show (\ref{estimate_beta2+_less_betaL2}). To this aim, applying the fundamental theorem of calculus (\ref{identity_fundamental}) to the exponential functional of Brownian motion on the right-hand side of (\ref{characterization_beta_L_2}), we have the following series expansion:
\begin{align} \label{representation_series_expansion_second_moment}
    &\int_{\mathbb{R}^{d}} dz_{0} \beta^{2}R(\sqrt{2}z_{0})\mathbf{E}_{z_{0}}\bigbkt{\exp\bkt{\beta^{2}\int_{0}^{\infty} R(\sqrt{2}B(u)) du }}
    = \int_{\mathbb{R}^{d}} dz_{0} \beta^{2}R(\sqrt{2}z_{0}) \nonumber\\
    &+\sum_{k=1}^{\infty} \int_{\mathbb{R}^{d \times k}} \prod_{j=1}^{k} dz_{k}\int_{0<t_{1}<...<t_{k} < \infty} \prod_{j=1}^{k} dt_{j} \cdot \beta^{2}R(\sqrt{2}z_{0}) \prod_{j=1}^{k} G_{t_{j}-t_{j-1}}(z_{j-1}-z_{j}) \beta^{2} R(\sqrt{2}z_{j})\nonumber\\
    &= \int_{\mathbb{R}^{d}} dz_{0} \beta^{2}R(\sqrt{2}z_{0})+\sum_{k=1}^{\infty} \int_{\mathbb{R}^{d \times k}} \prod_{j=1}^{k} dz_{k} \cdot \beta^{2}R(\sqrt{2}z_{0})  \prod_{j=1}^{k} \mathcal{G}^{0}(z_{j-1}-z_{j}) \beta^{2} R(\sqrt{2}z_{j}),
\end{align}
where we have used a change of variable $u_{j} = t_{j}-t_{j-1}$ in the last equality. Hence, we reformulate the right-hand side of (\ref{representation_series_expansion_second_moment}) by using the Birman-Schwinger operator $\mathbf{T}_{\beta}^{0}$ to conclude the following estimate:
\begin{align} \label{estimate_geometric_series}
    &\int_{\mathbb{R}^{d}} dz_{0} \beta^{2}R(\sqrt{2}z_{0})\mathbf{E}_{z_{0}}\bigbkt{\exp\bkt{\beta^{2}\int_{0}^{\infty} R(\sqrt{2}B(u)) du }}
    = \sum_{k=0}^{\infty} \bigl\langle \mathfrak{V}_{\beta}, (\mathbf{T}^{0}_{\beta})^{k}\mathfrak{V}_{\beta} \bigr\rangle_{L^{2}(\mathbb{R}^{d})}\nonumber\\
    &\lesssim \sum_{k=0}^{\infty} ||\mathbf{T}_{\beta}^{0}||_{L^{2}(\mathbb{R}^{d})}^{k} = \sum_{k=0}^{\infty} \beta^{2k}||\mathbf{T}_{1}^{0}||_{L^{2}(\mathbb{R}^{d})}^{k},
\end{align}
where we have used the fact that $||\mathbf{T}_{\beta}^{0}||_{L^{2}(\mathbb{R}^{d})} = \beta^{2} ||\mathbf{T}_{1}^{0}||_{L^{2}(\mathbb{R}^{d})}$. Therefore, by using Lemma \ref{lemma_representation_beta_2_+}, we know that if $\beta < \beta_{2,+}$, then the geometric series on the right-hand side of the above estimate converges. In this way, $\beta < \beta_{L^{2}}$ thanks to Lemma \ref{lemma_characterization_beta_L_2}. Consequently, we have $\beta_{2,+} \leq \beta_{L^{2}}$, and thus, the proof of (\ref{identity_beta2+_=_betaL2}) is complete.
\end{proof}

The remainder of this section is organized in the following way. We first will prove Lemma \ref{lemma_characterization_beta_L_2} in Section \ref{section_proof_of_characterization_beta_L_2}. Next, Lemma \ref{lemma_representation_beta_2_+} will be proved in Section \ref{section_proof_of_lemma_representation_beta_2_+} and Section \ref{section_Properties_of_the_Birman-Schwinger_operators}.

\subsection{Proof of Lemma \ref{lemma_characterization_beta_L_2}.} \label{section_proof_of_characterization_beta_L_2}
Recall the characterization (\ref{characterization_beta_LN}). To prove (\ref{characterization_beta_L_2}), it enough to show that the following for every $\beta>0$:
\begin{align} \label{estimation_second_moment_lower_and_upper}
    &C_{1}\mathbf{E}_{0}\bigbkt{\exp\bkt{\beta^{2}\int_{0}^{\infty} R(\sqrt{2}B(u)) du }} \leq\nonumber\\
    & 
    \int_{\mathbb{R}^{d}} dz_{0} R(\sqrt{2}z_{0})\mathbf{E}_{z_{0}}\bigbkt{\exp\bkt{\beta^{2}\int_{0}^{\infty} R(\sqrt{2}B(u)) du }} 
    \leq C_{2} \mathbf{E}_{0}\bigbkt{\exp\bkt{\beta^{2}\int_{0}^{\infty} R(\sqrt{2}B(u)) du }},
\end{align}
where $C_{j}$ is a positive constant which depends on $d$ and $R$.
\paragraph{Step 1.} To prove the lower bound in (\ref{estimation_second_moment_lower_and_upper}), we first observe that the following map is spherically symmetric:
\begin{align*}
    z\in \mathbb{R}^{d} \mapsto\mathbf{E}_{z}\bigbkt{\exp\bkt{\beta^{2}\int_{0}^{\infty} R(\sqrt{2}B(u)) du }}.
\end{align*}
Here, we have used the fact that $R(x) = \phi*\phi(x)$ and the spherical symmetricity of $\phi$ which comes from the approximation of the identity defined in Section \ref{section_the_model}. Then, by the Markov property of Brownian motion, one has the following identity:
\begin{align} \label{idenity_markov_0}
    &\mathbf{E}_{0}\bigbkt{\exp\bkt{\beta^{2}\int_{0}^{\infty} R(\sqrt{2}B(u)) du }}\nonumber\\
    &= \mathbf{E}_{0}\bigbkt{\exp\bkt{\beta^{2}\int_{0}^{\tau_{|z|}} R(\sqrt{2}B(u)) du } 
    } \cdot \mathbf{E}_{z}\bigbkt{\exp\bkt{\beta^{2}\int_{0}^{\infty} R(\sqrt{2}B(u)) du }}\quad \forall z \in \mathbb{R}^{d},
\end{align}
where $\tau_{r}$ is the first hitting time of $\partial B(0,r)$ by Brownian motion $B$. Therefore, it holds that
\begin{align*}
    \mathbf{E}_{0}\bigbkt{\exp\bkt{\beta^{2}\int_{0}^{\infty} R(\sqrt{2}B(u)) du }}
    \geq \mathbf{E}_{z}\bigbkt{\exp\bkt{\beta^{2}\int_{0}^{\infty} R(\sqrt{2}B(u)) du }} \quad \forall z \in \mathbb{R}^{d},
\end{align*}
Here, we have used the fact that $R$ is non-negative. As a consequence, the lower bound in (\ref{estimation_second_moment_lower_and_upper}) follows from the estimate below: 
\begin{align*}
    \int_{\mathbb{R}^{d}} dz_{0} R(\sqrt{2}z_{0}) \mathbf{E}_{z_{0}}\bigbkt{\exp\bkt{\beta^{2}\int_{0}^{\infty} R(\sqrt{2}B(u)) du }}\leq ||R(\sqrt{2}\cdot)||_{L^{1}(\mathbb{R}^{d})}  \mathbf{E}_{0}\bigbkt{\exp\bkt{\beta^{2}\int_{0}^{\infty} R(\sqrt{2}B(u)) du }}.
\end{align*}

\paragraph{Step 2.} Regarding the upper bound in (\ref{estimation_second_moment_lower_and_upper}), we set $\widetilde{\beta} := \beta_{L^{2}}+1$, and note that there exists a small $\widetilde{r} > 0$ such that
\begin{align*}
    \widetilde{C} :=\sup_{r\leq \widetilde{r}}\mathbf{E}_{0}\bigbkt{\exp\bkt{\widetilde{\beta}^{2}\int_{0}^{\tau_{r}} R(\sqrt{2}B(u)) du }} < \infty \quad\text{and}\quad B(0,\widetilde{r}) \subseteq \text{supp}(R(\sqrt{2}\cdot)) = B(0,\sqrt{2}r_{\phi}),
\end{align*}
where we have used (\ref{support_R}). Therefore, the estimate below follows from (\ref{idenity_markov_0}):
\begin{align*}
    ||R(\sqrt{2}\cdot)||_{L^{1}(B(0,\widetilde{r}))} &\cdot
    \mathbf{E}_{0}\bigbkt{\exp\bkt{\beta^{2}\int_{0}^{\infty} R(\sqrt{2}B(u)) du }} \\
    &\leq \widetilde{C}\cdot \int_{B(0,\widetilde{r})} dz_{0} R(\sqrt{2}z_{0}) \mathbf{E}_{z_{0}}\bigbkt{\exp\bkt{\beta^{2}\int_{0}^{\infty} R(\sqrt{2}B(u)) du }} \\
    &\leq \widetilde{C}\cdot \int_{\mathbb{R}^{d}} dz_{0} R(\sqrt{2}z_{0}) \mathbf{E}_{z_{0}}\bigbkt{\exp\bkt{\beta^{2}\int_{0}^{\infty} R(\sqrt{2}B(u)) du }}.
\end{align*}
Then, the above estimate implies the upper bound in (\ref{estimation_second_moment_lower_and_upper}). The proof is complete.

\subsection{Proof of Lemma \ref{lemma_representation_beta_2_+}.} \label{section_proof_of_lemma_representation_beta_2_+}
The prime objective of this subsection is to prove the variational representation (\ref{representation_beta_2_+}). Note that the boundedness of the operator $\mathbf{T}_{1}^{0}$ will be clarified in Lemma \ref{lemma_properties_Birman-Schwinger_operators}, which will be proved in Section \ref{section_Properties_of_the_Birman-Schwinger_operators}. To prove (\ref{representation_beta_2_+}), we first observe that
\begin{align} \label{estimate_first_lemma_representation_beta_2_+}
    ||\mathbf{T}_{1}^{0}||_{L^{2}(\mathbb{R}^{d})}^{-\frac{1}{2}} \leq \beta_{2,+}.
\end{align}
Indeed, if $\beta < ||\mathbf{T}_{1}^{0}||_{L^{2}(\mathbb{R}^{d})}^{-1/2}$, then $\beta< \beta_{L^{2}}$ follows from the estimate (\ref{estimate_geometric_series}) and the characterization (\ref{characterization_beta_L_2}). Hence, thanks to (\ref{ineqaulity_betaLN_betaN+}), one has $\beta<\beta_{2,+}$, and thus, (\ref{estimate_first_lemma_representation_beta_2_+}) holds. With the above clarification, it remains to prove the following estimate:
\begin{align} \label{estimate_second_lemma_representation_beta_2_+}
    \beta_{2,+} \leq ||\mathbf{T}_{1}^{0}||_{L^{2}(\mathbb{R}^{d})}^{-\frac{1}{2}}. 
\end{align}
The proof of above estimate is the core of the proof of Lemma \ref{lemma_representation_beta_2_+}.
In fact, this is an application of the Birman-Schwinger principle. See the proof of \cite[Theorem 12.4]{lieb2001analysis} for a more explicit explanation. Before the proof of (\ref{estimate_second_lemma_representation_beta_2_+}), we first introduce some auxiliary properties. 

To begin with, we derive a representation of $\beta_{2,+}$, which relates to an one-particle quantum system: 
\begin{lemma} \label{lemma_one_particle}
Assume that $d\geq 3$. Recall that $\beta_{2,+}$ is given in (\ref{definition_beta_N_+}). Then 
\begin{align} \label{representation_beta_2_+_one_particle}
    \beta_{2,+} = \inf\{\beta>0:\sup \mathbf{L}^{\beta} > 0\}.
\end{align}
Here, for each $\beta > 0$, the unbounded operator $\mathbf{L}^{\beta}$ is defined as follows:
\begin{align} \label{definition_hamiltonian_L}
    \mathbf{L}^{\beta} h(x) := \frac{1}{2}\Delta_{x} h(x) + V_{\beta}h(x), \quad D(\mathbf{L}^{\beta}) := C_{c}^{\infty}(\mathbb{R}^{d}),
\end{align}
where the potential function $V_{\beta}$ is defined in (\ref{definition_V_and_Yukawa}).
\end{lemma}
\begin{remark}
Recalling Lemma \ref{lemma_self_adjoint}, we know that $\mathbf{L}^{\beta}$ has a self-adjoint extension. For sake
of simplicity, we still use the same notation $\mathbf{L}^{\beta}$ to denote its self-adjoint extension.     
\end{remark}
\begin{proof}[Proof of Lemma \ref{lemma_one_particle}]
Let $\beta_{*}$ be the right-hand side of (\ref{representation_beta_2_+_one_particle}).
The key to prove (\ref{representation_beta_2_+_one_particle}) is the following identity: 
\begin{align} \label{identity_change_of_variable}
    &\int_{\mathbb{R}^{2 \times d}} -\frac{1}{2}|\nabla_{x_{1}} \rho(x_{1},x_{2})|^{2}-\frac{1}{2}|\nabla_{x_{2}} \rho(x_{1},x_{2})|^{2} + V_{\beta}\bkt{\frac{x_{2}-x_{1}}{\sqrt{2}}} |\rho(x_{1},x_{2})|^{2} dx_{1}dx_{2}\nonumber\\
    &= \int_{\mathbb{R}^{2 \times d}} -\frac{1}{2}|\nabla_{z} \widetilde{\rho}(z,z')|^{2}-\frac{1}{2}|\nabla_{z'} \widetilde{\rho}(z,z')|^{2} + V_{\beta}(z) |\widetilde{\rho}(z,z')|^{2} dzdz',
\end{align}
where the function $\widetilde{\rho}$ is defined as follows:
\begin{align*}
    \widetilde{\rho}(z,z') := \rho(x_{1},x_{2}), \quad \text{where} \quad (x_{1},x_{2}) = \bkt{\frac{z'-z}{\sqrt{2}},\frac{z'+z}{\sqrt{2}}}.
\end{align*}

Let $\beta > \beta_{2,+}$. Then, there exists $\rho \in C_{c}^{\infty}(\mathbb{R}^{2 \times d})$ such that the left-hand side of (\ref{identity_change_of_variable}) is positive. Hence, there exists $z'\in \mathbb{R}^{d}$ such that $h(\cdot) := \widetilde{\rho}(\cdot,z') \in C_{c}^{
\infty}(\mathbb{R}^{d})$ and
\begin{align} \label{energy_positive}
    \int_{\mathbb{R}^{d}} -\frac{1}{2}|\nabla_{z} h(z)|^{2}+ V_{\beta}(z) |h(z)|^{2} dz>0.
\end{align}
Consequently, we have $\beta > \beta_{*}$, and thus, $\beta_{2,+} \geq \beta_{*}$ 

On the other hand, let $\beta > \beta_{*}$, then there exists $h \in C_{c}^{\infty}(\mathbb{R}^{d})$ such that (\ref{energy_positive}) holds. Therefore, one can choose a large $u>0$ such that the right-hand side of (\ref{identity_change_of_variable}) is positive, where $\widetilde{\rho}(z,z') := h(z)G_{u}(z')^{1/2}$. Indeed, notice that 
\begin{align*}
    \int_{\mathbb{R}^{2 \times d}} -\frac{1}{2}|\nabla_{z'} \widetilde{\rho}(z,z')|^{2} \overset{u \to \infty}{\longrightarrow} 0.
\end{align*}
As a result, we obtain $\beta > \beta_{2,+}$, which implies $\beta_{*}\geq \beta_{2,+}$. The proof of (\ref{representation_beta_2_+_one_particle}) is complete. 
\end{proof}

Next, we introduce some properties which relates to the spectrum of the Birman-Schwinger opertaors $\mathbf{T}_{\beta}^{\lambda}$. For sake of simplicity, we postpone the proof of Lemma \ref{lemma_properties_Birman-Schwinger_operators} until Section \ref{section_Properties_of_the_Birman-Schwinger_operators}.
\begin{lemma} \label{lemma_properties_Birman-Schwinger_operators}
Assume that $d\geq 3$ and recall that $\mathbf{T}^{\lambda}_{\beta}$ is defined in (\ref{definition_Birman-Schwinger_operator}).
\begin{enumerate} [label=(\alph*)]
    \item For each $\lambda \geq 0$ and for every $\beta>0$, the Birman-Schwinger operators $\mathbf{T}_{\beta}^{\lambda}$ is a non-negative compact operator on $L^{2}(\mathbb{R}^{d})$.
    \item For any $\beta>0$, the largest eigenvalue $\mathscr{E}_{\beta}(\lambda):= ||\mathbf{T}_{\beta}^{\lambda}||_{L^{2}(\mathbb{R}^{d})}$ is continuous in $\lambda \geq 0$. Furthermore, $\mathscr{E}_{\beta}(\lambda)$ is decreasing and $\mathscr{E}_{\beta}(\lambda)$ vanishes as $\lambda \to \infty$.
\end{enumerate}
\end{lemma}
\begin{remark}
For the three-dimensional case, this lemma is much easier to prove because the Birman-Schwinger operators $\mathbf{T}_{\beta}^{\lambda}$ is a Hilbert Schmidt operator even when $\lambda = 0$. To be more specific, we have the estimate (\ref{Hilbert_Schmidt_operator_check}) in three dimensions. But, this estimate is false in general. Hence, a more careful argument is required to prove this lemma.
\end{remark}

Finally, we apply the Birman-Schwinger principle to conclude (\ref{estimate_second_lemma_representation_beta_2_+}). 
\begin{proof} [Proof of (\ref{estimate_second_lemma_representation_beta_2_+})]
To complete the proof of (\ref{estimate_second_lemma_representation_beta_2_+}), it suffices to show that if $\beta > ||\mathbf{T}_{1}^{0}||_{L^{2}(\mathbb{R}^{d})}^{-1/2}$, then $\beta> \beta_{2,+}$. Let $\beta > ||\mathbf{T}_{1}^{0}||_{L^{2}(\mathbb{R}^{d})}^{-1/2}$. Then, due to the representation (\ref{representation_beta_2_+_one_particle}), it is enough to show that $\mathbf{L}^{\beta}$ has a positive eigenvalue $\lambda$. To do this, the main idea is to apply the Birman-Schwinger principle.

We first note that $||\mathbf{T}_{\beta}^{0}||_{L^{2}(\mathbb{R}^{d})} = \beta^{2}||\mathbf{T}_{1}^{0}||_{L^{2}(\mathbb{R}^{d})}$. Then, $\mathscr{E}_{\beta}(0) > 1$ follows from our assumption on $\beta$. Hence, there exists $\lambda > 0$ such that $\mathscr{E}_{\beta}(\lambda) = 1$ due to the part (b) of Lemma \ref{lemma_properties_Birman-Schwinger_operators}. As a result, by applying the part (a) of Lemma \ref{lemma_properties_Birman-Schwinger_operators}, there exists $u \in L^{2}(\mathbb{R}^{d})$ such that 
\begin{align} \label{identity_eigenvalue_one}
    \mathbf{T}_{\beta}^{\lambda}u = u.
\end{align}
According to the Birman-Schwinger principle, we consider the following function:
\begin{align} \label{definition_h}
    h(x) := \int_{\mathbb{R}^{d}} \mathcal{G}^{\lambda}(x-y) \mathfrak{V}_{\beta}(y) u(y) dy,
\end{align}
where $\mathcal{G}^{\lambda}$ and $\mathfrak{V}_\beta$ are defined in (\ref{definition_V_and_Yukawa}). Then, the identity below follows from (\ref{identity_eigenvalue_one}):
\begin{align} \label{identity_u_and_h}
    u(x) = \mathfrak{V}_{\beta}(x) h(x).
\end{align}
Moreover, notice that $h \in D(\mathbf{L}^{\beta})$, where $D(\mathbf{L}^{\beta})$ is defined in (\ref{domain_laplace_one_particle}). Indeed, applying to the estimate below, we know that $|k|^{2} \widehat{h}(k) \in L^{2}(\mathbb{R}^{d})$: 
\begin{align*}
    \int_{\mathbb{R}^{d}} dk |k|^{4} \cdot|\widehat{h}(k)|^{2}
    &= \int_{\mathbb{R}^{d}} dk |\widehat{\mathfrak{V}_{\beta}u}(k)|^{2} \biggbkt{\frac{2|k|^{2}}{|2\pi k|^{2}+\lambda} }^{2} 
    \lesssim ||\mathfrak{V}_{\beta}u||_{L^{2}(\mathbb{R}^{d})}^{2} \lesssim ||R||^{2}_{L^{\infty}(\mathbb{R}^{d})} ||u||^{2}_{L^{2}(\mathbb{R}^{d})}.
\end{align*}
Here, we have used the following Fourier transform of the Yukawa potential $\mathcal{G}^{\lambda}$ in the first equality:
\begin{align} \label{Yukawa_potential_FT}
    \widehat{\mathcal{G}^{\lambda}}(k) = \frac{2}{|2\pi k|^{2} + \lambda} \quad \forall \lambda \geq 0,
\end{align}
where the case of $\lambda = 0$ is understood in the sense of \cite[Theorem 5.9]{lieb2001analysis} with $\alpha = 2$. As a consequence, due to the property of the Yukawa potential $\mathcal{G}^{\lambda}$ and (\ref{identity_u_and_h}), we get
\begin{align*}
    \biggbkt{\lambda-\frac{1}{2}\Delta}h(x) = \mathfrak{V}_{\beta}(x) u(x) = V_{\beta}(x) h(x).
\end{align*}
The above identity implies
\begin{align*}
    \mathbf{L}^{\beta}h(x) = \lambda h(x).
\end{align*}
Recalling the fact $h \in D(\mathbf{L}^{\beta})$, $\sup \mathbf{L}^{\beta} \geq \lambda > 0$ follows from the above identity. Thanks to (\ref{representation_beta_2_+_one_particle}), one has $\beta>\beta_{2,+}$. Therefore, we conclude that $\beta > \beta_{2,+}$ for every $\beta>||\mathbf{T}_{1}^{0}||_{L^{2}(\mathbb{R}^{d})}^{-1/2}$. The proof of (\ref{estimate_second_lemma_representation_beta_2_+}) is complete.

\end{proof}

\subsection{Properties of the general Birman-Schwinger operators.} \label{section_Properties_of_the_Birman-Schwinger_operators}
This part is dedicated to the proof of Lemma \ref{lemma_properties_Birman-Schwinger_operators}. Our proof consists of the following steps. A summarize of the proof will be presented in the last step.
\paragraph{Step 1.} As the first step to prove Lemma \ref{lemma_properties_Birman-Schwinger_operators}, we prove the boundedness of $\mathbf{T}^{\lambda}_{\beta}$ for all $\lambda \geq 0$ and every $\beta > 0$. We start with a simplification. Due to the symmetricity of $\mathbf{T}^{\lambda}_{\beta}$, we know that 
\begin{align} \label{estimate_lambda_less_zero}
    ||\mathbf{T}^{\lambda}_{\beta}||_{L^{2}(\mathbb{R}^{d})} = \sup_{||h||_{L^{2}(\mathbb{R}^{d})} \leq 1} &\biggl |\int_{\mathbb{R}^{d}} dx \overline{h(x)} \mathfrak{V}_{\beta}(x) 
    \int_{\mathbb{R}^{d}} dy \mathcal{G}^{\lambda}(x-y) \cdot \mathfrak{V}_{\beta}(y) h(y) \biggr|.
\end{align}
Hence, to prove the boundedness of $\mathbf{T}^{\lambda}_{\beta}$, it suffices to consider the case of $\lambda = 0$. Indeed, combining (\ref{estimate_lambda_less_zero}) and $\mathcal{G}^{\lambda}(x-y) \leq \mathcal{G}^{0}(x-y)$, we see that $||\mathbf{T}^{\lambda}_{\beta}||_{L^{2}(\mathbb{R}^{d})} \leq ||\mathbf{T}^{0}_{\beta}||_{L^{2}(\mathbb{R}^{d})}$.

Now, we focus one the case of $\lambda = 0$. In this way, one has the following estimate:
\begin{align*}
    ||\mathbf{T}^{0}_{\beta}||_{L^{2}(\mathbb{R}^{d})} \leq C_{d} \sup_{||h||_{L^{2}(\mathbb{R}^{d})} \leq 1}  \int_{\mathbb{R}^{d}} dx |h(x)| \mathfrak{V}_{\beta}(x) 
    \int_{\mathbb{R}^{d}} dy \frac{1}{|x-y|^{d-2}} \cdot \mathfrak{V}_{\beta}(y) |h(y)|,
\end{align*}
where $C_{d}$ is a positive constant such that $\mathcal{G}^{0}(x-y) = C_{d}|x-y|^{-(d-2)}$. To estimate the right-hand side of the above inequality, we decompose the above integral with respect to $dy$ as follows:
\begin{align*}
    \int_{\mathbb{R}^{d}} dx |h(x)| \mathfrak{V}_{\beta}(x) 
    \int_{\mathbb{R}^{d}} dy \frac{1}{|x-y|^{d-2}} \cdot \mathfrak{V}_{\beta}(y) |h(y)|
    \leq A_{\leq 1} + A_{\geq 1},
\end{align*}
where $A_{\bullet}$ is defined as follows:
\begin{align*}
    A_{\leq 1} := \int_{\mathbb{R}^{d}} dx |h(x)| \mathfrak{V}_{\beta}(x) 
    \int_{|x-y|\leq 1} dy \frac{1}{|x-y|^{d-2}} \cdot \mathfrak{V}_{\beta}(y) |h(y)|
\end{align*}
and 
\begin{align*}
    A_{\geq 1} := \int_{\mathbb{R}^{d}} dx |h(x)| \mathfrak{V}_{\beta}(x) 
    \int_{|x-y| \geq 1} dy \frac{1}{|x-y|^{d-2}} \cdot \mathfrak{V}_{\beta}(y) |h(y)|.
\end{align*}
With regard to the case of $A_{\geq 1}$, it is clear that $A_{\geq 1} \leq ||\mathfrak{V}_{\beta}||_{L^{\infty}(\mathbb{R}^{d})}^{2}$. Regarding the case of $A_{\leq 1}$, applying Young’s inequality implies
\begin{align*}
    A_{\leq 1} &\leq || \mathfrak{V}_{\beta} h ||_{L^{2}(\mathbb{R}^{d})} \cdot 
    \bignorm{\bkt{1_{B(0,1)}(\cdot) \frac{1}{|\cdot|^{d-2}}} * (\mathfrak{V}_{\beta} h) }_{L^{2}(\mathbb{R}^{d})}\\
    &\leq || \mathfrak{V}_{\beta} h ||_{L^{2}(\mathbb{R}^{d})} \cdot \bignorm{1_{B(0,1)}(\cdot) \frac{1}{|\cdot|^{d-2}} }_{L^{1}(\mathbb{R}^{d})} \cdot || \mathfrak{V}_{\beta} h ||_{L^{2}(\mathbb{R}^{d})} \lesssim ||\mathfrak{V}_{\beta}||_{L^{\infty}(\mathbb{R}^{d})}^{2},
\end{align*}
Therefore, the proof is complete.

\paragraph{Step 2.} Our goal of this step is to prove the compactness mentioned in the part (a) of Lemma \ref{lemma_properties_Birman-Schwinger_operators} without the case of $\lambda = 0$. To this aim, it is natural to prove the following properties:
\begin{enumerate} [label=(\alph*)]
    \item The operator $\mathbf{X}^{\lambda}_{\beta}$ defined as follows is a bounded operator from $L^{2}(\mathbb{R}^{d})$ to $H^{1}(\mathbb{R}^{d})$:
    \begin{align*}
        \mathbf{X}^{\lambda}_{\beta}u(x) := \int_{\mathbb{R}^{d}} \mathcal{G}^{\lambda}(x-y) \mathfrak{V}_{\beta}(y) u(y), \quad u\in L^{2}(\mathbb{R}^{d}).
    \end{align*}
    \item The following multiplicative operator is a compact operator from $H^{1}(\mathbb{R}^{d})$ to $L^{2}(\mathbb{R}^{d})$:
    \begin{align*}
        \mathbf{M}_{\beta}h(x) := \mathfrak{V}_{\beta}(x)h(x), \quad h\in H^{1}(\mathbb{R}^{d}).
    \end{align*}
\end{enumerate}
Here, $\mathcal{G}^{\lambda}$ and $\mathfrak{V}_{\beta}$ are defined in (\ref{definition_V_and_Yukawa}). Therefore, since $\mathbf{T}^{\lambda}_{\beta} = \mathbf{M}_{\beta} \mathbf{X}^{\lambda}_{\beta}$, the compactness of $\mathbf{T}^{\lambda}_{\beta}$ follows from \cite[Proposition 6.3]{brezis2011functional}, (a), and (b). 

Here, we briefly explain the proof ideas for the above properties. First, we observe that the property (a) follows from 
Plancherel’s theorem \cite[Theorem 5.3]{lieb2001analysis} and
the following Fourier transform of $\mathbf{X}^{\lambda}_{\beta}u$:
\begin{align*}
    \widehat{\mathbf{X}^{\lambda}_{\beta}u}(k) = \widehat{\mathfrak{V}_{\beta} u}(k) \cdot \widehat{\mathcal{G}^{\lambda}}(k) = \widehat{\mathfrak{V}_{\beta} u}(k) \cdot \frac{2}{|2\pi k|^{2} + \lambda}.
\end{align*}
Indeed, by using $\lambda > 0$ and the boundedness of $\mathfrak{V}_{\beta}$, $H^{1}$-norm of the right-hand side of the above identity can bounded from above by $C||u||_{L^{2}(\mathbb{R}^{d})}$. Next, to prove the property (b), the compactness is proved by using \cite[Theorem 3.18]{brezis2011functional} and the following argument. Let $(h_{j})_{j\geq 1}$ converges weakly to $h$ in $H^{1}(\mathbb{R}^{d})$. Recalling Remark \ref{remark_R}, we know that $\text{supp}(\mathfrak{V}_{\beta}) = B(0,\sqrt{2}r_{\phi})$ and $\mathfrak{V}_{\beta}$ is symmetric-decreasing. Then, it is straightforward to check that $(\mathbf{M}_{\beta}h_{j})_{j\geq 1}$ is bounded in $H^{1}(B(0,r))$, where $r \in (0, \sqrt{2}r_{\phi})$. Hence, if $r \in (0, \sqrt{2}r_{\phi})$, then applying Rellich-Kondrachov theorem \cite[Section 5.7]{evans2010partial} with $p = q = 2$ and $U = B(0,r)$ implies that there exists a subsequence  $(\mathbf{M}_{\beta}h_{j_{k}})_{k\geq 1}$ converges to $\mathbf{M}_{\beta}h$ in $L^{2}(B(0,r))$. Therefore, for each $\varepsilon>0$, there exists $r>0$ and $j_{*} \geq 1$ such that 
\begin{align} \label{RK_example}
    &\int_{\mathbb{R}^{d}} dx |\mathfrak{V}_{\beta}(x)h_{j_{*}}(x) - \mathfrak{V}_{\beta}(x)h(x)|^{2}\nonumber\\
    &\leq\int_{B(0,\sqrt{2}r_{\phi})\setminus B(0,r)} dx |\mathfrak{V}_{\beta}(x)h_{j_{*}}(x) - \mathfrak{V}_{\beta}(x)h(x)|^{2}+\int_{B(0,r)} dx |\mathfrak{V}_{\beta}(x)h_{j_{*}}(x) - \mathfrak{V}_{\beta}(x)h(x)|^{2} < \varepsilon.
\end{align}
Indeed, one can first bound the first term in the middle of (\ref{RK_example}) from above by the following through choosing a suitable $r \in (0,\sqrt{2}r_{\phi})$:
\begin{align*}
    ||\mathfrak{V}_{\beta}||_{L^{\infty}(B(0,\sqrt{2}r_{\phi})\setminus B(0,r))}\cdot \sup_{j\geq 1}||h_{j}||_{L^{2}(\mathbb{R}^{d})} < \varepsilon/2,
\end{align*}
and then the second term in the middle of (\ref{RK_example}) can be bounded from above by $\varepsilon/2$ through the convergence of $(\mathbf{M}_{\beta}h_{j_{k}})_{k\geq 1}$. Therefore, by using the above strategy iteratively, we see that there exists a subsequence of $(\mathbf{M}_{\beta}h_{j})_{j\geq 1}$ converging in $L^{2}(\mathbb{R}^{d})$. Hence, the proof of the property (b) is complete.

\paragraph{Step 3.} This step aims to prove the continuity of $\mathscr{E}_{\beta}(\lambda)$. To do this, it suffices to show that
\begin{align} \label{convergence_T_lambda}
    ||\mathbf{T}_{\beta}^{\lambda + \delta} - \mathbf{T}_{\beta}^{\lambda}||_{L^{2}(\mathbb{R}^{d})} = o(1) \quad \text{as } \delta \to 0 \quad \forall \lambda\geq 0.
\end{align}
Here, if $\lambda = 0$, then $\delta \to 0^{+}$. As already explained in Step 1,  the case of $\lambda = 0$ is the major case to be considered. Indeed, since  $\mathbf{T}_{\beta}^{\lambda}$ is bounded, the symmetricity of $\mathbf{T}_{\beta}^{\lambda}$ implies that
\begin{align} \label{identity_symmtric_bounded_opertaor}
    ||\mathbf{T}_{\beta}^{\lambda + \delta} - \mathbf{T}_{\beta}^{\lambda}||_{L^{2}(\mathbb{R}^{d})}
    = \sup_{||h||_{L^{2}(\mathbb{R}^{d})} \leq 1} &\biggl |\int_{\mathbb{R}^{d}} dx \overline{h(x)} \mathfrak{V}_{\beta}(x) 
    \cdot\int_{\mathbb{R}^{d}} dy (\mathcal{G}^{\lambda+\delta}(x-y) - \mathcal{G}^{\lambda}(x-y)) \cdot \mathfrak{V}_{\beta}(y) h(y) \biggr|,
\end{align}
and thus, it is straightforward to see that 
\begin{align*}
    ||\mathbf{T}_{\beta}^{\lambda + \delta} - \mathbf{T}_{\beta}^{\lambda}||_{L^{2}(\mathbb{R}^{d})}
    \leq ||\mathbf{T}_{\beta}^{0} - \mathbf{T}_{\beta}^{|\delta|}||_{L^{2}(\mathbb{R}^{d})}.
\end{align*}

With the above interpretation, it is enough to show that $||\mathbf{T}_{\beta}^{\delta} - \mathbf{T}_{\beta}^{0}||_{L^{2}(\mathbb{R}^{d})}$ vanishes as $\delta \to 0^{+}$. To this aim, we decompose the right-hand side of (\ref{identity_symmtric_bounded_opertaor}) as follows:
\begin{align} \label{decomposition_leq_geq}
    (\text{(R.H.S) of (\ref{identity_symmtric_bounded_opertaor}) with $\delta > 0$ and $\lambda = 0$})
    \leq A_{\leq 1} + A_{\geq 1}.
\end{align}
Here, $A_{\bullet}$ is given by
\begin{align*}
    A_{\bullet} := \sup_{||h||_{L^{2}(\mathbb{R}^{d})} \leq 1} \int_{\mathbb{R}^{d}} dx |h(x)| \mathfrak{V}_{\beta}(x) 
    \int_{\mathbb{R}^{d}} dy (\mathcal{G}_{\bullet}^{0}(x-y) - \mathcal{G}_{\bullet}^{\delta}(x-y)) \cdot \mathfrak{V}_{\beta}(y) |h(y)|,
\end{align*}
where the incomplete Yukawa potentials are defined as follows:
\begin{align*}
    \mathcal{G}_{\leq 1}^{\lambda}(z) := \int_{0}^{1} \exp(-\lambda \cdot u) G_{u}(z) du \quad\text{and} \quad\mathcal{G}_{\geq 1}^{\lambda}(z) := \int_{1}^{\infty} \exp(-\lambda \cdot u) G_{u}(z) du \quad \forall z\in \mathbb{R}^{d}\setminus\{0\}.
\end{align*}

To prove the right-hand side of (\ref{decomposition_leq_geq}) vanishes as $\delta \to 0^{+}$, we start with the case of $A_{\leq 1}$. To do this, notice that the following estimate follows from mean value theorem:
\begin{align*}
    \mathcal{G}_{\leq 1}^{0}(x-y) - \mathcal{G}_{\leq 1}^{\delta}(x-y) \leq \delta \cdot \int_{0}^{1} G_{u}(x-y) du \leq \delta \cdot \mathcal{G}^{0}(x-y),
\end{align*}
As a consequence, it holds that $A_{\leq 1} \leq \delta \cdot ||\mathbf{T}_{\beta}^{0}||_{L^{2}(\mathbb{R}^{d})} = o(1)$ as $\delta \to 0^{+}$,
where we have used the boundedness of $\mathbf{T}_{\beta}^{0}$ proved in the step 1 of this subsection. 
%Here, let us digress a bit on reason why we applied the above decomposition. The major reason is that applying mean value theorem produces a variable $u$ which causes $u G_{u}(x-y)$ is not integrable on $(0,\infty)$ when $d = 3,4$. As a result, the above estimates would be a problem to deal with the case of large $u$. 

Now, we turn the case of $A_{\geq 1}$. Clearly, the above strategy cannot work in this case. The key here is the use of
the integrability of $G_{u}(0)$ when $u$ is away from $0$. To be more specific, dominated convergence theorem implies that
\begin{align*}
    \mathcal{G}_{\geq 1}^{0}(x-y) - \mathcal{G}_{\geq 1}^{\delta}(x-y) &= 
    \int_{1}^{\infty} du \biggbkt{1-\exp(-\delta u)} G_{u}(x-y)
    \leq \int_{1}^{\infty} du \biggbkt{1-\exp(-\delta u)} G_{u}(0) \\
    &= o(1) \quad\text{as } \delta\to 0^{+} \text{ uniformly for all }x,y\in \mathbb{R}^{d}.
\end{align*}
Therefore, the estimate below follows from the above unifrom estimate:
\begin{align*}
    A_{\geq 1} \leq o(1) \cdot \sup_{||h||_{L^{2}(\mathbb{R}^{d})} \leq 1} \int_{\mathbb{R}^{d}} dx |h(x)| \mathfrak{V}_{\beta}(x) 
    \int_{\mathbb{R}^{d}} dy \mathfrak{V}_{\beta}(y) |h(y)|
    = o(1) \quad \text{as } \delta \to 0^{+},
\end{align*}
where we have used the fact that $\mathfrak{V}_{\beta},h \in L^{2}(\mathbb{R}^{d})$. In view of  $(\ref{decomposition_leq_geq})$ and the fact that both of $A_{\leq 1}$ and $A_{\geq 1}$ vanish as $\delta\to 0^{+}$, we see that $||\mathbf{T}_{\beta}^{0} - \mathbf{T}_{\beta}^{\delta}||_{L^{2}(\mathbb{R}^{d})}$ vanishes as $\delta \to 0^{+}$. The proof is complete.

%so that the claim follows at once
\paragraph{Step 4.} Finally, we summarize above results as follows in order to complete the proof of Lemma \ref{lemma_properties_Birman-Schwinger_operators}. We first note that the continuity of $\mathscr{E}_{\beta}(\lambda)$ has been proved in the step 3. Also, the compactness of $\mathbf{T}_{\beta}^{\lambda}$ for all $\lambda>0$ has been proved in the step 1. In this way, the compactness of $\mathbf{T}_{\beta}^{0}$ then follows from (\ref{convergence_T_lambda}) and \cite[Theorem 6.1]{brezis2011functional}. Furthermore, the following identities and estimate yield the claims that $\mathbf{T}_{\beta}^{\lambda}$ is non-negative and that $\mathscr{E}_{\beta}(\lambda)$ decrease to zero as $\lambda \to \infty$, where Parseval's theorem \cite[Theorem 5.3]{lieb2001analysis}, (\ref{Yukawa_potential_FT}), and (\ref{estimate_lambda_less_zero}) have been used to prove the following:
\begin{align*}
    \int_{\mathbb{R}^{d}} dx \overline{h(x)} \mathbf{T}_{\beta}^{\lambda}h(x)
    = \int_{\mathbb{R}^{d}} dk |\widehat{\mathfrak{V}_{\beta} h}(k)|^{2} \frac{2}{|2\pi k|^{2} + \lambda} \quad \forall h\in L^{2}(\mathbb{R}^{d}), \; \lambda \geq 0,
\end{align*}
\begin{align*} 
     \mathscr{E}_{\beta}(\lambda) = ||\mathbf{T}^{\lambda}_{\beta}||_{L^{2}(\mathbb{R}^{d})} = \sup_{||h||_{L^{2}(\mathbb{R}^{d})} \leq 1}  \int_{\mathbb{R}^{d}} dk |\widehat{\mathfrak{V}_{\beta} h}(k)|^{2} \frac{2}{|2\pi k|^{2} + \lambda} \quad \forall \lambda \geq 0,
\end{align*}
and
\begin{align*}
    \mathscr{E}_{\beta}(\lambda)
    \leq \frac{2}{\lambda} \sup_{||h||_{L^{2}(\mathbb{R}^{d})} \leq 1}||\mathfrak{V}_{\beta} h ||_{L^{2}(\mathbb{R}^{d})} \leq \frac{2}{\lambda} ||\mathfrak{V}_{\beta}||_{L^{\infty}(\mathbb{R}^{d})} \quad \forall \lambda > 0.
\end{align*}

\section{Proof of Theorem \ref{Main_result_2}} \label{section_proof_Main_result_2}
This section is dedicated to the proof Theorem \ref{Main_result_2}. As already explained at the beginning of Section \ref{section_heuristics_main_result_2},\emph{ we assume that $d = 3$, $N = 3$, and $U_{0}(x) = G_{\nu}(x)$ in the sequel}. With this assumption, Theorem \ref{Main_result_2} follows from Proposition \ref{proposition_partII_second_id} and Proposition \ref{proposition_lower_bound_sum_of_path_integrals}.
\begin{proof} [Proof of Theorem \ref{Main_result_2}.]
Let $m\geq 3$, $\Vec{x}_{0} \in \mathbb{R}^{N \times d}\setminus \Pi^{(N)}$, and $\varphi \in C_{c}^{\infty}(\mathbb{R}^{N \times d};\mathbb{R}_{+})$. Then, the following lower bounds about the sum of all the sub-limiting path integrals with length $m$ follow from Proposition \ref{proposition_partII_second_id} and Proposition \ref{proposition_lower_bound_sum_of_path_integrals}:
\begin{align*}
    \sum_{(\ell_{1},\ell'_{1},\mathfrak{i}_{1}),...,(\ell_{m},\ell'_{m},\mathfrak{i}_{m}) \in \mathcal{E}^{(N)} \times \{0,1\}, \; (\ell_{1},\ell_{1}') \neq ...\neq (\ell_{m},\ell_{m}')} &\mathscr{I}_{0;t}^{N;(\ell_{1},\ell'_{1},\mathfrak{i}_{1}),...,(\ell_{m},\ell'_{m},\mathfrak{i}_{m})} U_{0}^{\otimes N} (\Vec{x}_{0})\\
    &\geq C \frac{1}{|\Vec{x}_{0}|} \cdot\biggbkt{1.008}^{m-1} \frac{1}{(m-1)(m-2)}
\end{align*}
and
\begin{align*}
    \Biggl\langle\varphi,&\sum_{(\ell_{1},\ell'_{1},\mathfrak{i}_{1}),...,(\ell_{m},\ell'_{m},\mathfrak{i}_{m}) \in \mathcal{E}^{(N)} \times \{0,1\}, \; (\ell_{1},\ell_{1}') \neq ...\neq (\ell_{m},\ell_{m}')} \mathscr{I}_{0;t}^{N;(\ell_{1},\ell'_{1},\mathfrak{i}_{1}),...,(\ell_{m},\ell'_{m},\mathfrak{i}_{m})} U_{0}^{\otimes N}  \Biggr\rangle_{L^{2}(\mathbb{R}^{N \times d})}\\
    &\quad\quad\quad\quad\quad\quad\quad\quad\quad\quad\quad\quad\quad\geq C \biggbkt{\int_{\mathbb{R}^{N \times d}} d\Vec{x}_{0} \cdot \varphi(\Vec{x}_{0}) \frac{1}{|\Vec{x}_{0}|}} \cdot \biggbkt{1.008}^{m-1} \frac{1}{(m-1)(m-2)},
\end{align*}
where we have used the fact that the sub-limiting path integral is zero if $\mathfrak{i}_{j} = 0$ for some $j$. See the discussion under Remark \ref{remark_after_definition}. Consequently, in view of (\ref{definition_sub_limit}), the proof of Theorem \ref{Main_result_2} is complete.
\end{proof}

The remainder of this article is organized in the following way. In Section \ref{section_proof_lemma_relative_motion}, we will prove the representation of relative motion. Then, the proof of Proposition \ref{proposition_partII_second_id} will be presented in Section \ref{section_proposition_proposition_partII_second_id}. Finally, we will prove Proposition \ref{proposition_lower_bound_sum_of_path_integrals} in Section \ref{section_proof_proposition_lower_bound_sum_of_path_integrals}. 

\subsection{Derivation of the representation of relative motion.} \label{section_proof_lemma_relative_motion}
The prime objective of this subsection is to prove Lemma \ref{lemma_relative_motion}. Note that we have set $d = 3$, $N = 3$, and $U_{0}(x) = G_{\nu}(x)$ in the sequel. Before we proceed, we recall that the notations $\Vec{x}^{(\ell,\ell')}$, $\Vec{x}^{(\ell,\ell')^{*}}$, and $\Vec{x}^{(\ell,\ell')^{c}}$ are defined in (\ref{notation_relative}). With these notations, the sub-limiting path integral is then the following:
\begin{align} \label{recall_definition_sub-limiting_path_integral}
    &\mathscr{I}_{0;t}^{N;(\ell_{1},\ell'_{1},\mathfrak{i}_{1}),...,(\ell_{m},\ell'_{m},\mathfrak{i}_{m})}U_{0}^{\otimes N}(\Vec{x}_{0})\nonumber\\
    &= \int_{u_{0} = 0<s_{1}<u_{1}<s_{2}...<s_{m}<u_{m}<t} \prod_{j=1}^{m}ds_{j}du_{j} \int_{\mathbb{R}^{(m+1) \times N \times d}}  \prod_{j=1}^{m+1} d\Vec{x}_{j}\cdot\biggbkt{\prod_{j=1}^{m}
    G_{u_{j}-u_{j-1}}\bkt{\Vec{x}_{j}^{(\ell_{j},\ell'_{j})^{c}} - \Vec{x}_{j-1}^{(\ell_{j},\ell'_{j})^{c}}}
    \cdot \nonumber\\
    &\biggbkt{
    \int_{\mathbb{R}^{2 \times d}} dy_{j}dy_{j}' G_{s_{j}-u_{j-1}}(\Vec{x}_{j-1}(\ell_{j})-y_{j})
    G_{s_{j}-u_{j-1}}(\Vec{x}_{j-1}(\ell'_{j})-y_{j}') \delta_{0}\bkt{\frac{y'_{j}-y_{j}}{\sqrt{2}}} \cdot  \nonumber\\
    &\biggbkt{\sqrt{\frac{2\pi}{u_{j}-s_{j}}} G_{u_{j}-s_{j}}\bkt{\frac{y'_{j}+y_{j}}{\sqrt{2}}-\Vec{x}_{j}^{(\ell_{j},\ell_{j}')^{*}}} }\cdot\delta_{0}\bkt{\Vec{x}_{j}^{(\ell_{j},\ell_{j}')}}
    }
    } \cdot G^{(N)}_{t-u_{m}}(\Vec{x}_{m}-\Vec{x}_{m+1}) U_{0}^{\otimes N}(\Vec{x}_{m+1}).
\end{align}
As already explained in the step 1 of Section \ref{section_heuristics_main_result_2}, the main idea to conclude Lemma \ref{lemma_relative_motion} is to compute all the above Dirac’s delta functions by using the change of variable (\ref{relative_map}). 

\paragraph{Step 1.} To this aim, the first step is to apply the change of variables (\ref{Heuristics_change_of_variable_1}) to the spatial integral with respect to $dy_{j}dy'_{j}$ on the right-hand side of (\ref{recall_definition_sub-limiting_path_integral}). Before doing so, we simplify the Gaussian kernels on the third and the fourth lines of (\ref{recall_definition_sub-limiting_path_integral}). Observe that Gaussian kernel $G_{t}$ admits the following identity:
\begin{align} \label{identity_key_gaussian}
    G_{t}(z_{1})\cdot G_{t}(z_{2}) = G_{t}\bkt{\frac{z_{2}-z_{1}}{\sqrt{2}}}\cdot G_{t}\bkt{\frac{z_{2}+z_{1}}{\sqrt{2}}}.
\end{align}
Then, one can re-express the Gaussian kernels 
\begin{align*}
    \biggbkt{G_{s_{j}-u_{j-1}}(\Vec{x}_{j-1}(\ell_{j})-y_{j})
    G_{s_{j}-u_{j-1}}(\Vec{x}_{j-1}(\ell'_{j})-y_{j}')} \cdot G_{u_{j}-s_{j}}\biggbkt{\frac{y'_{j}+y_{j}}{\sqrt{2}}-\Vec{x}_{j}^{(\ell_{j},\ell_{j}')^{*}}}
\end{align*}
as follows:
\begin{align*}
    \biggbkt{G_{s_{j}-u_{j-1}}\bkt{\Vec{x}_{j-1}^{(\ell_{j},\ell'_{j})^{*}} -w'_{j}} G_{s_{j}-u_{j-1}}\bkt{\Vec{x}_{j-1}^{(\ell_{j},\ell'_{j})} -w_{j}}} \cdot G_{u_{j}-s_{j}}\bkt{w'_{j}-\Vec{x}_{j}^{(\ell_{j},\ell'_{j})^{*}}}
\end{align*}
Consequently, the following representation of the sub-limiting path integral then follows by evaluating the Dirac’s delta function $\delta_{0}(w_{j})$ and by using the Chapman–Kolmogorov equation with respect to $dw'_{j}$
\begin{align} \label{identity_first_change_variable}
    &\mathscr{I}_{0;t}^{N;(\ell_{1},\ell'_{1},\mathfrak{i}_{1}),...,(\ell_{m},\ell'_{m},\mathfrak{i}_{m})}U_{0}^{\otimes N}(\Vec{x}_{0})= \int_{u_{0} = 0<s_{1}<u_{1}<s_{2}...<s_{m}<u_{m}<t} \prod_{j=1}^{m}ds_{j}du_{j} \int_{\mathbb{R}^{m \times N \times d}}  \prod_{j=1}^{m+1} d\Vec{x}_{j}\nonumber\\
    &\biggbkt{\prod_{j=1}^{m}
    G_{u_{j}-u_{j-1}}\bkt{\Vec{x}_{j}^{(\ell_{j},\ell'_{j})^{c}} - \Vec{x}_{j-1}^{(\ell_{j},\ell'_{j})^{c}}} 
    G_{u_{j}-u_{j-1}}\bkt{\Vec{x}_{j}^{(\ell_{j},\ell'_{j})^{*}} - \Vec{x}_{j-1}^{(\ell_{j},\ell'_{j})^{*}}} G_{s_{j}-u_{j-1}}\bkt{\Vec{x}_{j-1}^{(\ell_{j},\ell'_{j})}}  \nonumber\\
    &\quad\quad\quad\quad \cdot \sqrt{\frac{2\pi}{u_{j}-s_{j}}}\cdot \delta_{0}\bkt{\Vec{x}_{j}^{(\ell_{j},\ell_{j}')}} }\cdot G^{(N)}_{t-u_{m}}(\Vec{x}_{m}-\Vec{x}_{m+1}) U_{0}^{\otimes N}(\Vec{x}_{m+1}).
\end{align}

\paragraph{Step 2.} To compute the rest of the Dirac’s delta functions on the right-hand side of (\ref{identity_first_change_variable}), the main idea is to apply the change of variables (\ref{Heuristics_change_of_variable}). But, before doing so, one has to express the above 
\begin{align*}
    \Vec{x}_{j-1}^{(\ell_{j},\ell_{j}')}, \;\Vec{x}_{j-1}^{(\ell_{j},\ell_{j}')^{*}}, \;\text{and}\;\Vec{x}_{j-1}^{(\ell_{j},\ell_{j}')^{c}} \quad\text{by using} \quad\Vec{x}_{j-1}^{(\ell_{j-1},\ell_{j-1}')}, \;\Vec{x}_{j-1}^{(\ell_{j-1},\ell_{j-1}')^{*}},\; \text{and}\; \Vec{x}_{j-1}^{(\ell_{j-1},\ell_{j-1}')^{c}},
\end{align*}
where all of them are defined in (\ref{notation_relative}). The derivation of these expressions heavily relies on the condition $N=3$ and the nonconsecutive interactions. 
\paragraph{Step 2-1.} Due to the above purpose, we first observe that
\begin{align} \label{expression_c}
    \Vec{x}_{j-1}^{(\ell_{j},\ell'_{j})^{c}}
    = \frac{\Vec{x}_{j-1}^{(\ell_{j-1},\ell'_{j-1})^{*}}\pm\Vec{x}_{j-1}^{(\ell_{j-1},\ell'_{j-1})}}{\sqrt{2}}.
\end{align}
Indeed, if $c \not \in \{\ell_{j},\ell_{j}'\}$, then it must be in $\{\ell_{j-1},\ell'_{j-1}\}$ thanks to $N = 3$ and $(\ell_{j},\ell_{j}') \neq (\ell_{j-1},\ell_{j-1}')$. In other words, $\Vec{x}_{j-1}^{(\ell_{j},\ell'_{j})^{c}}$ is either $\Vec{x}_{j-1}(\ell_{j-1})$ or $\Vec{x}_{j-1}(\ell'_{j-1})$, and thus, the above sign $\pm$ depends on which $c = \ell_{j-1}$ or $c = \ell_{j-1}'$ holds. As a consequence of (\ref{expression_c}), one then has
\begin{align} \label{identity_first_gaussian}
    G_{u_{j}-u_{j-1}}\bkt{\Vec{x}_{j}^{(\ell_{j},\ell'_{j})^{c}} - \Vec{x}_{j-1}^{(\ell_{j},\ell'_{j})^{c}}} 
    = G_{u_{j}-u_{j-1}}\bkt{\Vec{x}_{j}^{(\ell_{j},\ell'_{j})^{c}} - \frac{\Vec{x}_{j-1}^{(\ell_{j-1},\ell'_{j-1})^{*}}\pm\Vec{x}_{j-1}^{(\ell_{j-1},\ell'_{j-1})}}{\sqrt{2}}} \quad \forall 2\leq j\leq m.
\end{align}
Note that the above sign $\pm$ is actually not important because the variable $\Vec{x}_{j-1}^{(\ell_{j-1},\ell'_{j-1})}$ on the right-hand side of the above identity will eventually equal to zero thanks to the Dirac’s delta functions on the right-hand side of (\ref{identity_first_change_variable}).

\paragraph{Step 2-2.} Second, we notice that, in veiw of the conditions $N = 3$ and $(\ell_{j},\ell_{j}') \neq (\ell_{j-1},\ell_{j-1}')$, it is well-defined to denote by $b_{j-1}$ the common element that is in both of the sets $\{\ell_{j},\ell_{j}'\}$ and $\{\ell_{j-1},\ell_{j-1}'\}$. Also, we denote by $a_{j-1}$ the other element in the set $\{\ell_{j},\ell_{j}'\}$, and thus, the above reason shows that $a_{j-1}$ is not in $\{\ell_{j-1},\ell_{j-1}'\}$. With these notations, one has the following expression:
\begin{align} \label{expression_star}
    \Vec{x}_{j-1}^{(\ell_{j},\ell'_{j})^{*}}
    =
    \frac{\Vec{x}_{j-1}(a_{j-1})+\Vec{x}_{j-1}(b_{j-1})}{\sqrt{2}}
    =
    \frac{\Vec{x}_{j-1}^{(\ell_{j-1},\ell_{j-1}')^{c}}}{\sqrt{2}} + \frac{1}{\sqrt{2}}\bkt{\frac{\Vec{x}_{j-1}^{(\ell_{j-1},\ell'_{j-1})^{*}}\pm\Vec{x}_{j-1}^{(\ell_{j-1},\ell'_{j-1})}}{\sqrt{2}}}.
\end{align}
Again, the above sign $\pm$ depends on which $b_{j-1} = \ell_{j-1}$ or $b_{j-1} = \ell'_{j-1}$ holds. In this way, we conclude
\begin{align} \label{identity_second_gaussian}
    &G_{u_{j}-u_{j-1}}\bkt{\Vec{x}_{j}^{(\ell_{j},\ell'_{j})^{*}} - \Vec{x}_{j-1}^{(\ell_{j},\ell'_{j})^{*}}}\nonumber\\
    &= G_{u_{j}-u_{j-1}}\bkt{\Vec{x}_{j}^{(\ell_{j},\ell'_{j})^{*}}-\frac{\Vec{x}_{j-1}^{(\ell_{j-1},\ell_{j-1}')^{c}}}{\sqrt{2}} - \frac{\Vec{x}_{j-1}^{(\ell_{j-1},\ell'_{j-1})^{*}}\pm\Vec{x}_{j-1}^{(\ell_{j-1},\ell'_{j-1})}}{\sqrt{2}^{2}}} \quad \forall 2\leq j\leq m.
\end{align}
\paragraph{Step 2-3.} Third, with the above notations, one has
\begin{align} \label{expression_relative}
    \Vec{x}_{j-1}^{(\ell_{j},\ell'_{j})}
    =
    \pm\bkt{\frac{\Vec{x}_{j-1}(b_{j-1})-\Vec{x}_{j-1}(a_{j-1})}{\sqrt{2}}}
    =
    \pm\bkt{\frac{\Vec{x}_{j-1}^{(\ell_{j-1},\ell_{j-1}')^{c}}}{\sqrt{2}} - \frac{\Vec{x}_{j-1}^{(\ell_{j-1},\ell'_{j-1})^{*}}\pm\Vec{x}_{j-1}^{(\ell_{j-1},\ell'_{j-1})}}{\sqrt{2}^{2}}},
\end{align}
where the above signs depend on which of the numbers $b_{j-1}$ and $a_{j-1}$ is larger. Then, we conclude
\begin{align} \label{identity_third_gaussian}
    G_{s_{j}-u_{j-1}}(\Vec{x}_{j-1}^{(\ell_{j},\ell_{j}')}) = G_{s_{j}-u_{j-1}}\bkt{\frac{\Vec{x}_{j-1}^{(\ell_{j-1},\ell_{j-1}')^{c}}}{\sqrt{2}} - \frac{\Vec{x}_{j-1}^{(\ell_{j-1},\ell'_{j-1})^{*}}\pm\Vec{x}_{j-1}^{(\ell_{j-1},\ell'_{j-1})}}{\sqrt{2}^{2}}} \quad \forall 2\leq j\leq m.
\end{align}

\paragraph{Step 3.} Now, we are ready to apply the change of variables (\ref{Heuristics_change_of_variable}) to (\ref{identity_first_change_variable}). We also set
\begin{align*}
    &v_{j-1} = s_{j}-u_{j-1} \quad \text{and} \quad r_{j} = u_{j} - s_{j} \quad \forall 1\leq j\leq m.
\end{align*}
Then, the following expression for the sub-limiting path integral follows by applying the identities (\ref{identity_first_gaussian}), (\ref{identity_second_gaussian}), and (\ref{identity_third_gaussian}), and by evaluating the Dirac’s delta functions on the right-hand side of (\ref{identity_first_change_variable}):
\begin{align} \label{identity_last_lemma_relative_motion}
    &\mathscr{I}_{0;t}^{N;(\ell_{1},\ell'_{1},\mathfrak{i}_{1}),...,(\ell_{m},\ell'_{m},\mathfrak{i}_{m})}U_{0}^{\otimes N}(\Vec{x}_{0}) = \int_{v_{j},r_{j}>0, \; v_{0}+r_{1}+...+v_{m-1}+r_{m} <t} \prod_{j=1}^{m} dv_{j-1} dr_{j} \nonumber\\
    &\int_{\mathbb{R}^{m\times 2\times d + 3 \times d},\; \Vec{z}_{0} = \Vec{x}_{0}}  \bkt{\prod_{j=1}^{m} d\Vec{z}_{j}^{(\ell_{j},\ell_{j}')^{*}} d\Vec{z}_{j}^{(\ell_{j},\ell_{j}')^{c}} } \cdot \bkt{d\Vec{z}_{m+1}^{(\ell_{m},\ell_{m}')} d\Vec{z}_{m+1}^{(\ell_{m},\ell_{m}')^{*}} d\Vec{z}_{m+1}^{(\ell_{m},\ell_{m}')^{c}}}
    \nonumber\\
    &\biggbkt{G_{v_{0}}\bkt{\Vec{z}_{0}^{(\ell_{1},\ell_{1}')}} \cdot \sqrt{\frac{2\pi}{r_{1}}} \cdot G_{v_{0}+r_{1}}\bkt{\Vec{z}_{0}^{(\ell_{1},\ell_{1}')^{*}} - \Vec{z}_{1}^{(\ell_{1},\ell_{1}')^{*}}} G_{v_{0}+r_{1}}\bkt{\Vec{z}_{0}^{(\ell_{1},\ell_{1}')^{c}} - \Vec{z}_{1}^{(\ell_{1},\ell_{1}')^{c}}}}\nonumber\\
    &\biggbkt{ \prod_{j=2}^{m} G_{v_{j-1}}\bkt{\frac{\Vec{z}_{j-1}^{(\ell_{j-1},\ell'_{j-1})^{c}}}{\sqrt{2}}- \frac{\Vec{z}_{j-1}^{(\ell_{j-1},\ell'_{j-1})^{*}}}{\sqrt{2}^{2}}} \cdot \sqrt{\frac{2\pi}{r_{j}}} \cdot \nonumber\\
    &\quad\quad\quad\quad\quad\quad G_{v_{j-1}+r_{j}}\bkt{\frac{\Vec{z}_{j-1}^{(\ell_{j-1},\ell'_{j-1})^{c}}}{\sqrt{2}}+ \frac{\Vec{z}_{j-1}^{(\ell_{j-1},\ell'_{j-1})^{*}}}{\sqrt{2}^{2}} - \Vec{z}_{j}^{(\ell_{j},\ell_{j}')^{*}}} G_{v_{j-1}+r_{j}}\bkt{\frac{\Vec{z}_{j-1}^{(\ell_{j-1},\ell'_{j-1})^{*}}}{\sqrt{2}} - \Vec{z}_{j}^{(\ell_{j},\ell_{j}')^{c}}}}\nonumber\\
    &\biggbkt{G_{v_{m}}\bkt{0 - \Vec{z}_{m+1}^{(\ell_{m},\ell_{m}')}}G_{v_{m}}\bkt{\Vec{z}_{m}^{(\ell_{m},\ell_{m}')^{*}} - \Vec{z}_{m+1}^{(\ell_{m},\ell_{m}')^{*}}}G_{v_{m}}\bkt{\Vec{z}_{m}^{(\ell_{m},\ell_{m}')^{c}} - \Vec{z}_{m+1}^{(\ell_{m},\ell_{m}')^{c}}} }\nonumber\\
    &\biggbkt{G_{\nu}\bkt{\Vec{z}_{m+1}^{(\ell_{m},\ell_{m}')}}G_{\nu}\bkt{ \Vec{z}_{m+1}^{(\ell_{m},\ell_{m}')^{*}}}G_{\nu}\bkt{ \Vec{z}_{m+1}^{(\ell_{m},\ell_{m}')^{c}}}},
\end{align}
where, we have applied the following identities to the Gaussian kernel $G^{(N)}_{t-u_{m}}(\Vec{x}_{m}-\Vec{x}_{m+1})$ and the initial datum $U_{0}^{\otimes N}(\Vec{x}_{m+1})$ on the right-hand side of (\ref{identity_first_change_variable}):
\begin{align*}
    &G^{(N)}_{t-u_{m}}(\Vec{x}_{m}-\Vec{x}_{m+1}) = G_{v_{m}}\bkt{\Vec{z}_{m}^{(\ell_{m},\ell_{m}')} - \Vec{z}_{m+1}^{(\ell_{m},\ell_{m}')}}G_{v_{m}}\bkt{\Vec{z}_{m}^{(\ell_{m},\ell_{m}')^{*}} - \Vec{z}_{m+1}^{(\ell_{m},\ell_{m}')^{*}}}G_{v_{m}}\bkt{\Vec{z}_{m}^{(\ell_{m},\ell_{m}')^{c}} - \Vec{z}_{m+1}^{(\ell_{m},\ell_{m}')^{c}}}
\end{align*}
and 
\begin{align*}
    U_{0}^{\otimes N}(\Vec{x}_{m+1}) = G_{\nu}\bkt{\Vec{z}_{m+1}^{(\ell_{m},\ell_{m}')}}G_{\nu}\bkt{ \Vec{z}_{m+1}^{(\ell_{m},\ell_{m}')^{*}}}G_{\nu}\bkt{ \Vec{z}_{m+1}^{(\ell_{m},\ell_{m}')^{c}}}.
\end{align*}
Here, the $0$ in $G_{v_{m}}\bkt{0 - \Vec{z}_{m+1}^{(\ell_{m},\ell_{m}')}}$ on the right-hand side of (\ref{identity_last_lemma_relative_motion}) follows from the last Dirac’s delta function $\delta_{0}(\Vec{z}_{m}^{(\ell_{m},\ell_{m}')})$. 

\paragraph{Step 4.} Finally, we compute the integrals with respect to the variables $d\Vec{z}_{m+1}^{(\ell_{m},\ell_{m}')} d\Vec{z}_{m+1}^{(\ell_{m},\ell_{m}')^{*}} d\Vec{z}_{m+1}^{(\ell_{m},\ell_{m}')^{c}}$ and $d\Vec{z}_{m}^{(\ell_{m},\ell_{m}')^{*}} d\Vec{z}_{m}^{(\ell_{m},\ell_{m}')^{c}}$ on the right-hand side of (\ref{identity_last_lemma_relative_motion}). This can be done by applying the Chapman-Kolmogorov equation in the following order:
\begin{align*}
    G_{v_{m}+\nu}(0) = \int_{\mathbb{R}^{d}} d\Vec{z}_{m+1}^{(\ell_{m},\ell_{m}')} G_{v_{m}}\bkt{\Vec{z}_{m+1}^{(\ell_{m},\ell_{m}')}} G_{\nu}\bkt{\Vec{z}_{m+1}^{(\ell_{m},\ell_{m}')}},
\end{align*}
\begin{align*}
    G_{v_{m}+\nu}\bkt{\Vec{z}_{m}^{(\ell_{m},\ell_{m}')^{*}}} = \int_{\mathbb{R}^{d}} d\Vec{z}_{m+1}^{(\ell_{m},\ell_{m}')^{*}} G_{v_{m}}\bkt{\Vec{z}_{m}^{(\ell_{m},\ell_{m}')^{*}} - \Vec{z}_{m+1}^{(\ell_{m},\ell_{m}')^{*}}}G_{\nu}\bkt{\Vec{z}_{m+1}^{(\ell_{m},\ell_{m}')^{*}}},
\end{align*}
\begin{align*}
    G_{v_{m}+\nu}\bkt{\Vec{z}_{m}^{(\ell_{m},\ell_{m}')^{c}}} = \int_{\mathbb{R}^{d}} d\Vec{z}_{m+1}^{(\ell_{m},\ell_{m}')^{*}} G_{v_{m}}\bkt{\Vec{z}_{m}^{(\ell_{m},\ell_{m}')^{c}} - \Vec{z}_{m+1}^{(\ell_{m},\ell_{m}')^{c}}}G_{\nu}\bkt{\Vec{z}_{m+1}^{(\ell_{m},\ell_{m}')^{c}}},
\end{align*}
\begin{align*}
     &G_{v_{m-1}+\overline{r}_{m}}\bkt{\frac{\Vec{z}_{m-1}^{(\ell_{m-1},\ell'_{m-1})^{c}}}{\sqrt{2}}+ \frac{\Vec{z}_{m-1}^{(\ell_{m-1},\ell'_{m-1})^{*}}}{\sqrt{2}^{2}}}\\
     &= \int_{\mathbb{R}^{d}} d\Vec{z}_{m}^{(\ell_{m},\ell_{m}')^{*}}  G_{v_{m-1}+r_{m}}\bkt{\frac{\Vec{z}_{m-1}^{(\ell_{m-1},\ell'_{m-1})^{c}}}{\sqrt{2}}+ \frac{\Vec{z}_{m-1}^{(\ell_{m-1},\ell'_{m-1})^{*}}}{\sqrt{2}^{2}} - \Vec{z}_{m}^{(\ell_{m},\ell_{m}')^{*}}} G_{v_{m}+\nu}\bkt{\Vec{z}_{m}^{(\ell_{m},\ell_{m}')^{*}}},
\end{align*}
and
\begin{align*}
     &G_{v_{m-1}+\overline{r}_{m}}\bkt{ \frac{\Vec{z}_{m-1}^{(\ell_{m-1},\ell'_{m-1})^{*}}}{\sqrt{2}}}= \int_{\mathbb{R}^{d}} d\Vec{z}_{m}^{(\ell_{m},\ell_{m}')^{c}}  
     G_{v_{m-1}+r_{m}}\bkt{ \frac{\Vec{z}_{m-1}^{(\ell_{m-1},\ell'_{m-1})^{c}}}{\sqrt{2}}-\Vec{z}_{m}^{(\ell_{m},\ell_{m}')^{c}}} G_{v_{m}+\nu}\bkt{\Vec{z}_{m}^{(\ell_{m},\ell_{m}')^{c}}},
\end{align*}
where $\overline{r}_{m}$ is given in (\ref{definition_overline_rm}). Therefore, we complete the proof of Lemma \ref{lemma_relative_motion}.

\subsection{Spatial estimates for the sub-limiting path integrals.} \label{section_proposition_proposition_partII_second_id}
Thanks to Lemma \ref{lemma_relative_motion}, this subsection aims to prove the lower bound in Proposition \ref{proposition_partII_second_id}. The proof is based on a careful analysis of the Fourier transform of the sub-limiting path integral. Note that we have set $d = 3$, $N = 3$, and $U_{0}(x) = G_{\nu}(x)$. 

To evaluate the Fourier transform of the sub-limiting path integral, we regard it as a function 
\begin{align} \label{function_new_variable}
    \mathscr{I}_{0;t}^{N;(\ell_{1},\ell'_{1},\mathfrak{i}_{1}),...,(\ell_{k},\ell'_{k},\mathfrak{i}_{k})}U_{0}^{\otimes N}: \bkt{\Vec{z}_{0}^{(\ell_{1},\ell'_{1})},\Vec{z}_{0}^{(\ell_{1},\ell'_{1})^{*}},\Vec{z}_{0}^{(\ell_{1},\ell'_{1})^{c}}} \in \mathbb{R}^{N \times d}\setminus \widetilde{\Pi}^{(N)} \mapsto \mathbb{R}_{+},
\end{align}
where the above function is defined by the right-hand side of (\ref{identity_lemma_relative_motion}), and 
\begin{align} \label{definition_tilde_pi}
    \widetilde{\Pi}^{(N)} := \set{(z,z^{*},z^{c}) \in \mathbb{R}^{N \times d}: z \neq 0}.
\end{align}
Here, we have used the following equivalence:
\begin{align*}
    \Vec{x} \in \mathbb{R}^{N \times d} \setminus\Pi^{(N)}\quad\text{if and only if}\quad \bkt{\Vec{x}^{(\ell,\ell')},\Vec{x}^{(\ell,\ell')^{*}},\Vec{x}^{(\ell,\ell')^{c}}} \in \mathbb{R}^{N \times d} \setminus\widetilde{\Pi}^{(N)} \quad \forall (\ell,\ell') \in \mathcal{E}^{(N)},
\end{align*}
where $\Pi^{(N)}$ is defined in (\ref{definition_big_pi}), $\mathcal{E}^{(N)}$ is given in Section \ref{section_structure}, and the vector on the right-hand side above are defined in (\ref{notation_relative}). Moreover, instead of (\ref{FT}), we adopt the following notation to denote the Fourier transform of a function $f(z,z^{*},z^{c})$ defined on $\mathbb{R}^{N \times d}\setminus \widetilde{\Pi}^{(N)}$: 
\begin{align*}
    \mathscr{F}^{(N)}f(k,k^{*},k^{c}) := \int_{\mathbb{R}^{N \times d}} dzdz^{*}dz^{c} \exp\biggbkt{-2\pi i (k\cdot z +k^{*} \cdot z^{*} + k^{c} \cdot z^{c}) } f(z,z^{*},z^{c}).
\end{align*}

With the above notations, Proposition \ref{proposition_partII_second_id} will be proved according to the three steps mentioned at the end of the step 2 in Section \ref{section_heuristics_main_result_2}.

\subsubsection{Fourier transform of the sub-limiting path integral.} This part aims to evaluate the Fourier transform of the function defined in (\ref{function_new_variable}). Before we proceed, recalling \cite[Theorem 5.2]{lieb2001analysis}, we know that
the Fourier transform of the Gaussian kernel $G_{t}$ is as follows:
\begin{align} \label{formula_G_t}
    \mathscr{G}_{t}(k) := \widehat{G_{t}}(k) =\exp(-2\pi^{2}t|k|^{2}).
\end{align}

\begin{lemma} \label{lemma_identity_fourier_1}
Under the assumptions of Lemma \ref{lemma_relative_motion}, we have the following representation:
\begin{align} \label{identity_fourier}
    &\mathscr{F}^{(N)}\mathscr{I}_{0;t}^{N;(\ell_{1},\ell'_{1},\mathfrak{i}_{1}),...,(\ell_{m},\ell'_{m},\mathfrak{i}_{m})}U_{0}^{\otimes N}(k,k^{*},k^{c}) \nonumber\\
    &=\int_{v_{j},r_{j}>0, \; v_{0}+r_{1}+...+v_{m-1}+r_{m} <t, \; v_{m} = t-(v_{0}+r_{1}+...+v_{m-1}+r_{m})} \prod_{j=1}^{m} dv_{j-1} dr_{j} \nonumber\\
    &\biggbkt{\mathscr{G}_{v_{0}}(k)  \mathscr{G}_{v_{0}+r_{1}}(k^{*}) \mathscr{G}_{v_{0}+r_{1}}(k^{c})}\biggbkt{\biggbkt{\prod_{j=1}^{m}\sqrt{\frac{2\pi}{r_{j}}} }
    \cdot\sqrt{2}^{2 \times 3 \times (m-1)}  \cdot G_{v_{m}+\nu}(0)} \cdot \mathscr{B}_{v_{1},r_{2},...,v_{m-1},r_{m},v_{m}}(k^{*},k^{c}).
\end{align}
Here, we set
\begin{align} \label{definition_B}
    &\mathscr{B}_{v_{1},r_{2},...,v_{m-1},r_{m},v_{m}}(k^{*},k^{c}) := \int_{\mathbb{R}^{(m-1)\times d}} \prod_{j=2}^{m} dh_{j} \nonumber\\
    &\quad\quad\biggbkt{\mathscr{G}_{v_{1}+r_{2}}(\sqrt{2}h_{2})\mathscr{G}_{v_{1}}(-k^{*}+k^{c}/\sqrt{2}+h_{2})\mathscr{G}_{v_{1}+r_{2}}(k^{*}+k^{c}/\sqrt{2}-h_{2})}\nonumber\\
    &\quad\quad\biggbkt{\prod_{j=3}^{m-1} \mathscr{G}_{v_{j-1}+r_{j}}(\sqrt{2}h_{j})\mathscr{G}_{v_{j-1}}(-k^{*}-k^{c}/\sqrt{2}+2h_{j-1}+h_{j})\mathscr{G}_{v_{j-1}+r_{j}}(k^{*}+k^{c}/\sqrt{2}-h_{j})}\nonumber\\
    &\quad\quad\biggbkt{\mathscr{G}_{v_{m-1}+\overline{r}_{m}}(\sqrt{2}h_{m})\mathscr{G}_{v_{m-1}}(-k^{*}-k^{c}/\sqrt{2}+2h_{m-1}+h_{m})\mathscr{G}_{v_{m-1}+\overline{r}_{m}}(k^{*}+k^{c}/\sqrt{2}-h_{m})}
\end{align}
where the above product on the right-hand side about the Gaussian kernels $\mathscr{G}(\cdot)$ is equal to $1$ if $m = 3$, and $\overline{r}_{m}$ is defined in (\ref{definition_overline_rm}).
\end{lemma}

We will carry out the proof of Lemma \ref{lemma_identity_fourier_1} in the following steps.
\paragraph{Step 1.} To begin with, we re-express the right-hand side of (\ref{identity_lemma_relative_motion}) as the following form in order to apply Fourier transform:
\begin{align}
    &\mathscr{I}_{0;t}^{N;(\ell_{1},\ell'_{1},\mathfrak{i}_{1}),...,(\ell_{m},\ell'_{m},\mathfrak{i}_{m})}U_{0}^{\otimes N}\bkt{\Vec{z}_{0}^{(\ell_{1},\ell'_{1})},\Vec{z}_{0}^{(\ell_{1},\ell'_{1})^{*}},\Vec{z}_{0}^{(\ell_{1},\ell'_{1})^{c}}} \nonumber\\
    &= \int_{v_{j},r_{j}>0, \; v_{0}+r_{1}+...+v_{m-1}+r_{m} <t, \; v_{m} = t-(v_{0}+r_{1}+...+v_{m-1}+r_{m})} \biggbkt{\prod_{j=1}^{m} dv_{j-1} dr_{j}} \cdot\int_{\mathbb{R}^{2\times d}}   d\Vec{z}_{1}^{(\ell_{1},\ell_{1}')^{*}} d\Vec{z}_{1}^{(\ell_{1},\ell_{1}')^{c}} \nonumber\\
    &\cdot\biggbkt{G_{v_{0}}\bkt{\Vec{z}_{0}^{(\ell_{1},\ell_{1}')}}\cdot \sqrt{\frac{2\pi}{r_{1}}} \cdot G_{v_{0}+r_{1}}\bkt{\Vec{z}_{0}^{(\ell_{1},\ell_{1}')^{*}} - \Vec{z}_{1}^{(\ell_{1},\ell_{1}')^{*}}} G_{v_{0}+r_{1}}\bkt{\Vec{z}_{0}^{(\ell_{1},\ell_{1}')^{c}} - \Vec{z}_{1}^{(\ell_{1},\ell_{1}')^{c}}}} \nonumber\\
    &\quad\quad\quad\quad\quad\quad\quad\quad\quad\quad\quad\quad\quad\quad\quad\quad\quad\quad\cdot\mathscr{A}_{v_{1},r_{2},...,v_{m-1},r_{m},v_{m}}^{(\ell_{2},\ell'_{2},\mathfrak{i}_{2}),...,(\ell_{m},\ell'_{m},\mathfrak{i}_{m})}U_{0}^{\otimes N}(\Vec{z}_{1}^{(\ell_{1},\ell_{1}')^{*}},\Vec{z}_{1}^{(\ell_{1},\ell_{1}')^{c}}),
\end{align}
where the last term above is given by
\begin{align} \label{definition_A_recursion}
    &\mathscr{A}_{v_{n-1},r_{n},...,v_{m-1},r_{m},v_{m}}^{(\ell_{n},\ell'_{n},\mathfrak{i}_{n}),...,(\ell_{m},\ell'_{m},\mathfrak{i}_{m})}U_{0}^{\otimes N}(\Vec{z}_{n-1}^{(\ell_{n-1},\ell_{n-1}')^{*}},\Vec{z}_{n-1}^{(\ell_{n-1},\ell_{n-1}')^{c}}):= \int_{\mathbb{R}^{(m-n)\times 2\times d}}  \biggbkt{\prod_{j=n}^{m-1} d\Vec{z}_{j}^{(\ell_{j},\ell_{j}')^{*}} d\Vec{z}_{j}^{(\ell_{j},\ell_{j}')^{c}} }
    \cdot\nonumber\\
    &\biggbkt{ \prod_{j=n}^{m-1} G_{v_{j-1}}\bkt{\frac{\Vec{z}_{j-1}^{(\ell_{j-1},\ell'_{j-1})^{c}}}{\sqrt{2}}- \frac{\Vec{z}_{j-1}^{(\ell_{j-1},\ell'_{j-1})^{*}}}{\sqrt{2}^{2}}} \cdot \sqrt{\frac{2\pi}{r_{j}}} \cdot \nonumber\\
    &\quad\quad\quad G_{v_{j-1}+r_{j}}\bkt{\frac{\Vec{z}_{j-1}^{(\ell_{j-1},\ell'_{j-1})^{c}}}{\sqrt{2}}+ \frac{\Vec{z}_{j-1}^{(\ell_{j-1},\ell'_{j-1})^{*}}}{\sqrt{2}^{2}} - \Vec{z}_{j}^{(\ell_{j},\ell_{j}')^{*}}} G_{v_{j-1}+r_{j}}\bkt{\frac{\Vec{z}_{j-1}^{(\ell_{j-1},\ell'_{j-1})^{*}}}{\sqrt{2}} - \Vec{z}_{j}^{(\ell_{j},\ell_{j}')^{c}}}}\nonumber\\
    &\biggbkt{G_{v_{m-1}}\bkt{\frac{\Vec{z}_{m-1}^{(\ell_{m-1},\ell'_{m-1})^{c}}}{\sqrt{2}}- \frac{\Vec{z}_{m-1}^{(\ell_{m-1},\ell'_{m-1})^{*}}}{\sqrt{2}^{2}}} \cdot \sqrt{\frac{2\pi}{r_{m}}} \cdot \nonumber\\
    &\quad\quad\quad G_{v_{m}+\nu}(0) G_{v_{m-1}+\overline{r}_{m}}\bkt{\frac{\Vec{z}_{m-1}^{(\ell_{m-1},\ell'_{m-1})^{c}}}{\sqrt{2}}+ \frac{\Vec{z}_{m-1}^{(\ell_{m-1},\ell'_{m-1})^{*}}}{\sqrt{2}^{2}}} G_{v_{m-1}+\overline{r}_{m}}\bkt{\frac{\Vec{z}_{m-1}^{(\ell_{m-1},\ell'_{m-1})^{*}}}{\sqrt{2}}}}\nonumber\\
    &\quad\quad\quad\quad\quad\quad\quad\quad\quad\quad\quad\quad\quad\quad\quad\quad\quad\quad\quad\quad\quad\quad\quad\quad\quad\quad\quad\text{if } 2\leq n\leq m-1;
\end{align}
\begin{align} \label{definition_A_recursion_n=m}
    &\mathscr{A}_{v_{n-1},r_{n},...,v_{m-1},r_{m},v_{m}}^{(\ell_{n},\ell'_{n},\mathfrak{i}_{n}),...,(\ell_{m},\ell'_{m},\mathfrak{i}_{m})}U_{0}^{\otimes N}(\Vec{z}_{n-1}^{(\ell_{n-1},\ell_{n-1}')^{*}},\Vec{z}_{n-1}^{(\ell_{n-1},\ell_{n-1}')^{c}}):= G_{v_{m-1}}\bkt{\frac{\Vec{z}_{m-1}^{(\ell_{m-1},\ell'_{m-1})^{c}}}{\sqrt{2}}- \frac{\Vec{z}_{m-1}^{(\ell_{m-1},\ell'_{m-1})^{*}}}{\sqrt{2}^{2}}} \cdot \sqrt{\frac{2\pi}{r_{m}}} \nonumber \\
    &\quad \cdot G_{v_{m}+\nu}(0) G_{v_{m-1}+\overline{r}_{m}}\bkt{\frac{\Vec{z}_{m-1}^{(\ell_{m-1},\ell'_{m-1})^{c}}}{\sqrt{2}}+ \frac{\Vec{z}_{m-1}^{(\ell_{m-1},\ell'_{m-1})^{*}}}{\sqrt{2}^{2}}} G_{v_{m-1}+\overline{r}_{m}}\bkt{\frac{\Vec{z}_{m-1}^{(\ell_{m-1},\ell'_{m-1})^{*}}}{\sqrt{2}}}\quad\text{if } n = m.
\end{align}
Then, by \cite[Theorem 5.8]{lieb2001analysis}, it holds that
\begin{align} \label{identity_initial_recursion}
    &\mathscr{F}^{(N)}\mathscr{I}_{0;t}^{N;(\ell_{1},\ell'_{1},\mathfrak{i}_{1}),...,(\ell_{m},\ell'_{m},\mathfrak{i}_{m})}U_{0}^{\otimes N}(k,k^{*},k^{c}) \nonumber\\
    &= \int_{v_{j},r_{j}>0, \; v_{0}+r_{1}+...+v_{m-1}+r_{m} <t, \; v_{m} = t-(v_{0}+r_{1}+...+v_{m-1}+r_{m})} \biggbkt{\prod_{j=1}^{m} dv_{j-1} dr_{j}} \cdot\nonumber\\
    &\quad\quad\quad\biggbkt{\mathscr{G}_{v_{0}}(k)\cdot \sqrt{\frac{2\pi}{r_{1}}}\cdot \mathscr{G}_{v_{0}+r_{1}}(k^{*}) \mathscr{G}_{v_{0}+r_{1}}(k^{c})} \cdot \mathscr{F}^{(N-1)}\mathscr{A}_{v_{1},r_{2},...,v_{m-1},r_{m},v_{m}}^{(\ell_{2},\ell'_{2},\mathfrak{i}_{2}),...,(\ell_{m},\ell'_{m},\mathfrak{i}_{m})}U_{0}^{\otimes N}(k^{*},k^{c}),
\end{align}
where the last term above is given by
\begin{align*}
    \mathscr{F}^{(N-1)}f(k^{*},k^{c}) := \int_{\mathbb{R}^{(N-1) \times d}} dz^{*}dz^{c} \exp\biggbkt{-2\pi i (k^{*} \cdot z^{*} + k^{c} \cdot z^{c}) } f(z^{*},z^{c}).
\end{align*}

\paragraph{Step 2.}
To prove Lemma \ref{lemma_identity_fourier_1} , it remains to compute the last term on the right-hand side of (\ref{identity_initial_recursion}). To this aim, we show in this step that it can be computed recursively as follows:
\begin{align} \label{identity_first_recursion}
    &\mathscr{F}^{(N-1)}\mathscr{A}_{v_{1},r_{2},...,v_{m-1},r_{m},v_{m}}^{(\ell_{2},\ell'_{2},\mathfrak{i}_{2}),...,(\ell_{m},\ell'_{m},\mathfrak{i}_{m})}U_{0}^{\otimes N}(k^{*},k^{c}) \nonumber\\
    &= \int_{\mathbb{R}^{d}} dh_{2} \cdot \biggbkt{\sqrt{\frac{2\pi}{r_{2}}} \cdot \sqrt{2}^{2 \times 3} \cdot \mathscr{G}_{v_{1}+r_{2}}(\sqrt{2}h_{2})\mathscr{G}_{v_{1}}(-k^{*}+k^{c}/\sqrt{2}+h_{2})\mathscr{G}_{v_{1}+r_{2}}(k^{*}+k^{c}/\sqrt{2}-h_{2})} \cdot\nonumber\\
    &\mathscr{F}^{(N-1)}\mathscr{A}_{v_{2},r_{3},...,v_{m-1},r_{m},v_{m}}^{(\ell_{3},\ell'_{3},\mathfrak{i}_{3}),...,(\ell_{m},\ell'_{m},\mathfrak{i}_{m})}U_{0}^{\otimes N}(k^{*}+k^{c}/\sqrt{2}-h_{2},\sqrt{2}h_{2}).
\end{align}
In view of (\ref{definition_A_recursion}), we notice that
\begin{align*}
    &\mathscr{F}^{(N-1)}\mathscr{A}_{v_{1},r_{2},...,v_{m-1},r_{m},v_{m}}^{(\ell_{2},\ell'_{2},\mathfrak{i}_{2}),...,(\ell_{m},\ell'_{m},\mathfrak{i}_{m})}U_{0}^{\otimes N}(k^{*},k^{c}) =  \int_{\mathbb{R}^{2\times d}}  d\Vec{z}_{2}^{(\ell_{2},\ell_{2}')^{*}} d\Vec{z}_{2}^{(\ell_{2},\ell_{2}')^{c}} \\&\mathscr{A}_{v_{2},r_{3},...,v_{m-1},r_{m},v_{m}}^{(\ell_{3},\ell'_{3},\mathfrak{i}_{3}),...,(\ell_{m},\ell'_{m},\mathfrak{i}_{m})}U_{0}^{\otimes N}(\Vec{z}_{2}^{(\ell_{1},\ell_{1}')^{*}},\Vec{z}_{2}^{(\ell_{1},\ell_{2}')^{c}}) \cdot\biggbkt{\int_{\mathbb{R}^{2\times d}} 
    d\Vec{z}_{1}^{(\ell_{2},\ell_{2}')^{*}} d\Vec{z}_{1}^{(\ell_{2},\ell_{2}')^{c}}\\
    &\exp\biggbkt{-2\pi i \Vec{z}_{1}^{(\ell_{2},\ell_{2}')^{*}} \cdot k^{*}}
    \exp\biggbkt{-2\pi i \Vec{z}_{1}^{(\ell_{2},\ell_{2}')^{c}} \cdot k^{c}}\cdot\biggbkt{G_{v_{1}}\bkt{\frac{\Vec{z}_{1}^{(\ell_{1},\ell'_{1})^{c}}}{\sqrt{2}}- \frac{\Vec{z}_{1}^{(\ell_{1},\ell'_{1})^{*}}}{\sqrt{2}^{2}}} \cdot \sqrt{\frac{2\pi}{r_{2}}} \cdot \\
    &\quad\quad\quad\quad\quad\quad\quad\quad\quad G_{v_{1}+r_{2}}\bkt{\frac{\Vec{z}_{1}^{(\ell_{1},\ell'_{1})^{c}}}{\sqrt{2}}+ \frac{\Vec{z}_{1}^{(\ell_{1},\ell'_{1})^{*}}}{\sqrt{2}^{2}} - \Vec{z}_{2}^{(\ell_{2},\ell_{2}')^{*}}} G_{v_{1}+r_{2}}\bkt{\frac{\Vec{z}_{1}^{(\ell_{1},\ell'_{1})^{*}}}{\sqrt{2}} - \Vec{z}_{2}^{(\ell_{2},\ell_{2}')^{c}}} }}.
\end{align*}
Consequently, to prove (\ref{identity_first_recursion}), it is enough to prove the following identity:
\begin{align} \label{identity_recursive_start}
    &\int_{\mathbb{R}^{2\times d}} 
    d\Vec{z}_{1}^{(\ell_{2},\ell_{2}')^{*}} d\Vec{z}_{1}^{(\ell_{2},\ell_{2}')^{c}}
    \exp\biggbkt{-2\pi i \Vec{z}_{1}^{(\ell_{2},\ell_{2}')^{*}} \cdot k^{*}}
    \exp\biggbkt{-2\pi i \Vec{z}_{1}^{(\ell_{2},\ell_{2}')^{c}} \cdot k^{c}} \cdot\nonumber\\
    &\biggbkt{G_{v_{1}}\bkt{\frac{\Vec{z}_{1}^{(\ell_{1},\ell'_{1})^{c}}}{\sqrt{2}}- \frac{\Vec{z}_{1}^{(\ell_{1},\ell'_{1})^{*}}}{\sqrt{2}^{2}}} \cdot \sqrt{\frac{2\pi}{r_{2}}} \cdot \nonumber\\
    &G_{v_{1}+r_{2}}\bkt{\frac{\Vec{z}_{1}^{(\ell_{1},\ell'_{1})^{c}}}{\sqrt{2}}+ \frac{\Vec{z}_{1}^{(\ell_{1},\ell'_{1})^{*}}}{\sqrt{2}^{2}} - \Vec{z}_{2}^{(\ell_{2},\ell_{2}')^{*}}} G_{v_{1}+r_{2}}\bkt{\frac{\Vec{z}_{1}^{(\ell_{1},\ell'_{1})^{*}}}{\sqrt{2}} - \Vec{z}_{2}^{(\ell_{2},\ell_{2}')^{c}}}}\nonumber\\
    &= \int_{\mathbb{R}^{d}} dh_{2} \cdot \biggbkt{
    \sqrt{\frac{2\pi}{r_{2}}}\cdot \sqrt{2}^{2 \times 3} \cdot
    \mathscr{G}_{v_{1}+r_{2}}(\sqrt{2}h_{2}) 
    \mathscr{G}_{v_{1}}(-k^{*}+k^{c}/\sqrt{2}+h_{2})
    \mathscr{G}_{v_{1}+r_{2}}(k^{*}+k^{c}/\sqrt{2}-h_{2})}\nonumber\\
    &\exp\biggbkt{-2\pi i \Vec{z}_{2}^{(\ell_{2},\ell'_{2})^{*}} \cdot (k^{*}+k^{c}/\sqrt{2}-h_{2}) } \exp\biggbkt{-2\pi i \Vec{z}_{2}^{(\ell_{2},\ell'_{2})^{c}} \cdot \sqrt{2}  h_{2}}.
\end{align}

To prove (\ref{identity_recursive_start}), we first observe that 
\begin{align} \label{indentity_fourier_N-1}
    &\text{(L.H.S) of (\ref{identity_recursive_start})}= \int_{\mathbb{R}^{2\times d}} 
    dx_{1}dx_{1}'
    \exp\biggbkt{-2\pi i x_{1} \cdot (-k^{*}+k^{c}/\sqrt{2}) }
    \exp\biggbkt{-2\pi i x_{1}' \cdot (k^{*}+k^{c}/\sqrt{2}) }\nonumber\\
    &\biggbkt{\sqrt{\frac{2\pi}{r_{2}}} \cdot\sqrt{2}^{3}\cdot G_{v_{1}}(x_{1}) G_{v_{1}+r_{2}}\bkt{x_{1}'-\Vec{z}_{2}^{(\ell_{2},\ell'_{2})^{*}}} G_{v_{1}+r_{2}}\bkt{\frac{x_{1}'-x_{1}}{\sqrt{2}}-\Vec{z}_{2}^{(\ell_{2},\ell'_{2})^{c}}}},
\end{align}
where we have applied the following change of variables:
\begin{align*}
    x_{1} := \frac{\Vec{z}_{1}^{(\ell_{1},\ell'_{1})^{c}}}{\sqrt{2}}- \frac{\Vec{z}_{1}^{(\ell_{1},\ell'_{1})^{*}}}{\sqrt{2}^{2}}\quad\text{and} \quad x_{1}' := \frac{\Vec{z}_{1}^{(\ell_{1},\ell'_{1})^{c}}}{\sqrt{2}}+ \frac{\Vec{z}_{1}^{(\ell_{1},\ell'_{1})^{*}}}{\sqrt{2}^{2}},
\end{align*}
and the number $\sqrt{2}^{3}$ on the right-hand side of (\ref{indentity_fourier_N-1}) arises because of the Jacobian of the above change of variables. Next, recalling \cite[Theorem 5.8]{lieb2001analysis}, we know that
\begin{align*}
    \widehat{f \cdot g}(w) = \widehat{f} * \widehat{g}(w).
\end{align*}
Consequently, the integral with respect to $dx_{1}'$ on the right-hand side of $(\ref{indentity_fourier_N-1})$ can be compute as follows:
\begin{align*}
    &\int_{\mathbb{R}^{d}} 
    dx_{1}'
    \exp\biggbkt{-2\pi i x_{1}' \cdot (k^{*}+k^{c}/\sqrt{2}) } G_{v_{1}+r_{2}}\bkt{x_{1}'-\Vec{z}_{2}^{(\ell_{2},\ell'_{2})^{*}}} G_{v_{1}+r_{2}}\bkt{\frac{x_{1}'-x_{1}}{\sqrt{2}}-\Vec{z}_{2}^{(\ell_{2},\ell'_{2})^{c}}}\\
    &= \int_{\mathbb{R}^{d}} dh_{2} \biggbkt{ \exp\biggbkt{-2\pi i \Vec{z}_{2}^{(\ell_{2},\ell'_{2})^{*}} \cdot (k^{*}+k^{c}/\sqrt{2}-h_{2}) }\cdot \mathscr{G}_{v_{1}+r_{2}}(k^{*}+k^{c}/\sqrt{2}-h_{2})} \cdot\\
    &\biggbkt{\exp\biggbkt{-2\pi i (\sqrt{2}\Vec{z}_{2}^{(\ell_{2},\ell'_{2})^{c}}+x_{1}) \cdot  h_{2}}\cdot \mathscr{G}_{v_{1}+r_{2}}(\sqrt{2}h_{2}) \sqrt{2}^{3}}.
\end{align*}
Then, one has
\begin{align*}
    &\text{((L.H.S) of (\ref{identity_recursive_start}))} = \int_{\mathbb{R}^{d}} dh_{2}
    \biggbkt{\sqrt{\frac{2\pi}{r_{2}}} \sqrt{2}^{2 \times 3} \cdot\mathscr{G}_{v_{1}+r_{2}}(k^{*}+k^{c}/\sqrt{2}-h_{2})
    \mathscr{G}_{v_{1}+r_{2}}(\sqrt{2}h_{2})}\cdot\\
    &\biggbkt{\exp\biggbkt{-2\pi i \Vec{z}_{2}^{(\ell_{2},\ell'_{2})^{*}} \cdot (k^{*}+k^{c}/\sqrt{2}-h_{2}) } \exp\biggbkt{-2\pi i \Vec{z}_{2}^{(\ell_{2},\ell'_{2})^{c}} \cdot \sqrt{2}  h_{2}}}\cdot\\
    &\int_{\mathbb{R}^{d}} dx_{1} \exp\biggbkt{-2\pi i x_{1} \cdot (-k^{*}+k^{c}/\sqrt{2}+h_{2}) } G_{v_{1}}(x_{1}).
\end{align*}
Finally, the above integral with respect to $dx_{1}$ can be computed as follows by using (\ref{formula_G_t}):
\begin{align*}
    \int_{\mathbb{R}^{d}} dx_{1} \exp\biggbkt{-2\pi i x_{1} \cdot (-k^{*}+k^{c}/\sqrt{2}+h_{2}) } G_{v_{1}}(x_{1}) = \mathscr{G}_{v_{1}}(-k^{*}+k^{c}/\sqrt{2}+h_{2}).
\end{align*}
Therefore, we conclude (\ref{identity_recursive_start}), which implies (\ref{identity_first_recursion}).

\paragraph{Step 3.}
To complete the proof of Lemma \ref{lemma_identity_fourier_1}, it remains to compute the last line on the right-hand side of (\ref{identity_first_recursion}). By using the same computation to conclude (\ref{identity_first_recursion}), where $k^{*}$ and $k^{c}$ in (\ref{identity_first_recursion}) are replaced by $k^{*}+k^{c}/\sqrt{2}-h_{n-1}$ and $\sqrt{2}h_{n-1}$, respectively, one has the following recursive relation:
\begin{align} \label{identity_second_recursion}
    &\mathscr{F}^{(N-1)}\mathscr{A}_{v_{n-1},r_{n},...,v_{m-1},r_{m},v_{m}}^{(\ell_{n},\ell'_{n},\mathfrak{i}_{n}),...,(\ell_{m},\ell'_{m},\mathfrak{i}_{m})}U_{0}^{\otimes N}(k^{*}+k^{c}/\sqrt{2}-h_{n-1},\sqrt{2}h_{n-1}) = \int_{\mathbb{R}^{d}} dh_{n} \nonumber\\
    &\biggbkt{\sqrt{\frac{2\pi}{r_{n}}} \cdot\sqrt{2}^{2 \times 3} \cdot\mathscr{G}_{v_{n-1}+r_{n}}(\sqrt{2}h_{n})\mathscr{G}_{v_{n-1}}(-k^{*}-k^{c}/\sqrt{2}+2h_{n-1}+h_{n})\mathscr{G}_{v_{n-1}+r_{n}}(k^{*}+k^{c}/\sqrt{2}-h_{n})}\nonumber\\
    &\mathscr{F}^{(N-1)}\mathscr{A}_{v_{n},r_{n+1},...,v_{m-1},r_{m},v_{m}}^{(\ell_{n+1},\ell'_{n+1},\mathfrak{i}_{n+1}),...,(\ell_{m},\ell'_{m},\mathfrak{i}_{m})}U_{0}^{\otimes N}(k^{*}+k^{c}/\sqrt{2}-h_{n},\sqrt{2}h_{n}) \quad \forall 3\leq n\leq m-1.
\end{align}
As a consequence, it remains to compute $\mathscr{F}^{(N-1)}\mathscr{A}_{v_{m-1},r_{m},v_{m}}^{(\ell_{m},\ell'_{m},\mathfrak{i}_{m})}U_{0}^{\otimes N}(k^{*}+k^{c}/\sqrt{2}-h_{m-1},\sqrt{2}h_{m-1})$. Then, by applying (\ref{identity_recursive_start}), where $k^{*}$ and $k^{c}$ are replaced by $k^{*}+k^{c}/\sqrt{2}-h_{m-1}$ and $\sqrt{2}h_{m-1}$, respectively, and also, both of $\Vec{z}_{2}^{(\ell_{2},\ell_{2}')^{*}}$ and $\Vec{z}_{2}^{(\ell_{2},\ell_{2}')^{c}}$ are replaced by zero vectors, we know that
\begin{align} \label{identity_third_recursion}
    &\mathscr{F}^{(N-1)}\mathscr{A}_{v_{m-1},r_{m},v_{m}}^{(\ell_{m},\ell'_{m},\mathfrak{i}_{m})}U_{0}^{\otimes N}(k^{*}+k^{c}/\sqrt{2}-h_{m-1},\sqrt{2}h_{m-1})= \int_{\mathbb{R}^{d}} dh_{m} \biggbkt{\sqrt{\frac{2\pi}{r_{m}}} \cdot\sqrt{2}^{2 \times 3}\cdot G_{v_{m}+\nu}(0)\cdot \nonumber\\
    &\mathscr{G}_{v_{m-1}+\overline{r}_{m}}(\sqrt{2}h_{m})\mathscr{G}_{v_{m-1}}(-k^{*}-k^{c}/\sqrt{2}+2h_{m-1}+h_{m})\mathscr{G}_{v_{m-1}+\overline{r}_{m}}(k^{*}+k^{c}/\sqrt{2}-h_{m})}.
\end{align}

\begin{proof} [Proof of Lemma \ref{lemma_identity_fourier_1}.]
Finally, by combining (\ref{identity_initial_recursion}), (\ref{identity_first_recursion}), (\ref{identity_second_recursion}), and (\ref{identity_third_recursion}), we complete the proof of (\ref{identity_fourier}).
    
\end{proof}

\subsubsection{Algorithm for computing the spatial integrals.}
The overall goal of this part is to evaluate the spatial integrals on the right-hand side of (\ref{identity_fourier}).
\begin{lemma} \label{lemma_identity_fourier_2}
Under the assumptions of Lemma \ref{lemma_relative_motion}, the following identity holds:
\begin{align} \label{identity_lemma_identity_fourier_2}
    &\mathscr{F}^{(N)}\mathscr{I}_{0;t}^{N;(\ell_{1},\ell'_{1},\mathfrak{i}_{1}),...,(\ell_{m},\ell'_{m},\mathfrak{i}_{m})}U_{0}^{\otimes N}(k,k^{*},k^{c}) \nonumber\\
    &= \int_{v_{j},r_{j}>0, \; v_{0}+r_{1}+...+v_{m-1}+r_{m} <t, \; v_{m} = t-(v_{0}+r_{1}+...+v_{m-1}+r_{m})} \prod_{j=1}^{m} dv_{j-1} dr_{j} \nonumber\\
    &\biggbkt{\mathscr{G}_{v_{0}}(k)  \mathscr{G}_{v_{0}+r_{1}}(k^{*}) \mathscr{G}_{v_{0}+r_{1}}(k^{c})}\biggbkt{\biggbkt{\prod_{j=1}^{m}\sqrt{\frac{2\pi}{r_{j}}} }
    \sqrt{2}^{2 \times 3 \times (m-1)}  G_{v_{m}+\nu}(0)}\cdot\biggbkt{\mathscr{G}_{\mathbf{t}_{1}}\bkt{\frac{1}{3}k^{*}-\frac{2}{3}\frac{k^{c}}{\sqrt{2}}}\nonumber\\
    &\mathscr{G}_{(v_{1}+r_{2}+...+v_{m-1}+\overline{r}_{m})\cdot 2/3} (k^{*}+k^{c}/\sqrt{2}) 
    \biggbkt{\prod_{j=2}^{m-1} G_{v_{j-1}+(\mathbf{t}_{j+1}+3(v_{j-1}+r_{j}))}(0)}G_{v_{m-1}+3(v_{m-1}+\overline{r}_{m})}(0)},
\end{align}
where $\overline{r}_{m}$ is defined in (\ref{definition_overline_rm}), and $\mathbf{t}_{j}$ is given in (\ref{definition_tj}).
\end{lemma}

To prove it, thanks to the expression on the right-hand side of (\ref{identity_fourier}), it remains to evaluate $\mathscr{B}_{v_{1},r_{2},...,v_{m-1},r_{m},v_{m}}(k^{*},k^{c})$ given in (\ref{definition_B}). To this aim, the key is an algorithm to merge the Gaussian kernels on the right-hand side of (\ref{definition_B}) properly so that we can compute every integral with respect to $dh_{j}$ iteratively.

\paragraph{Step 1.} To prove Lemma \ref{lemma_identity_fourier_2}, it is enough to prove the following identity thanks to Lemma \ref{lemma_identity_fourier_1}: 
\begin{align} \label{integral_h1...h_m}
    &\mathscr{B}_{v_{1},r_{2},...,v_{m-1},r_{m},v_{m}}(k^{*},k^{c}) = \mathscr{G}_{\mathbf{t}_{1}}\bkt{\frac{1}{3}k^{*}-\frac{2}{3\sqrt{2}}k^{c}}\nonumber\\
    &\mathscr{G}_{(v_{1}+r_{2}+...+v_{m-1}+\overline{r}_{m})\cdot 2/3} (k^{*}+k^{c}/\sqrt{2}) 
    \biggbkt{\prod_{j=2}^{m-1} G_{v_{j-1}+(\mathbf{t}_{j+1}+3(v_{j-1}+r_{j}))}(0)}G_{v_{m-1}+3(v_{m-1}+\overline{r}_{m})}(0).
\end{align}
To this aim, we re-express the integrals on the right-hand side of (\ref{definition_B}) as follows:
\begin{align} \label{definition_A_2_kstar_kc}
    \mathscr{B}_{v_{1},r_{2},...,v_{m-1},r_{m},v_{m}}(k^{*},k^{c}) &:= \int_{\mathbb{R}^{d}} dh_{2} \mathscr{G}_{v_{1}+r_{2}} (\sqrt{2}h_{2})\nonumber\\
    &\mathscr{G}_{v_{1}}(-k^{*}+k^{c}/\sqrt{2}+h_{2})\mathscr{G}_{v_{1}+r_{2}}(k^{*}+k^{c}/\sqrt{2}-h_{2}) A_{3}(k^{*},k^{c},h_{2}),
\end{align}
\begin{align} \label{definition_A_j_kstar_kc}
    A_{n}(k^{*},k^{c},h_{n-1}) := \int_{\mathbb{R}^{d}} dh_{n}  &\mathscr{G}_{v_{n-1}+r_{n}}(\sqrt{2}h_{n})\mathscr{G}_{v_{n-1}}(-k^{*}-k^{c}/\sqrt{2}+2h_{n-1}+h_{n})\nonumber\\
    &\mathscr{G}_{v_{n-1}+r_{n}}(k^{*}+k^{c}/\sqrt{2}-h_{n}) A_{n+1}(k^{*},k^{c},h_{j}) \quad \forall 3\leq n\leq m-1, 
\end{align}
and
\begin{align} \label{Am_last}
    A_{m}(k^{*},k^{c},h_{m-1}) := \int_{\mathbb{R}^{d}} &dh_{m} 
    \mathscr{G}_{v_{m-1}+\overline{r}_{m}}(\sqrt{2}h_{m})\nonumber\\
    &\mathscr{G}_{v_{m-1}}(-k^{*}-k^{c}/\sqrt{2}+2h_{m-1}+h_{m})\mathscr{G}_{v_{m-1}+\overline{r}_{m}}(k^{*}+k^{c}/\sqrt{2}-h_{m}).
\end{align}
Then, the remainder of the proof is dedicated to compute $A_{n}$ iteratively. 

\paragraph{Step 2.}
To evaluate $\mathscr{B}_{v_{1},r_{2},...,v_{m-1},r_{m},v_{m}}(k^{*},k^{c})$, let us first compute $A_{m}$. To compute an integral consisting of three Gaussian kernels, one must merge these Gaussian kernels before the integration. Recalling the formula (\ref{formula_G_t}) of $\mathscr{G}_{t}(k)$, we combine the first and the third Gaussian kernels on the right-hand side of (\ref{Am_last}) as follows:
\begin{align} \label{decompistion_Am}
    \mathscr{G}_{v_{m-1}+\overline{r}_{m}}(\sqrt{2}h_{m})&\mathscr{G}_{v_{m-1}+\overline{r}_{m}}(k^{*}+k^{c}/\sqrt{2}-h_{m})\nonumber\\
    &= \mathscr{G}_{3\cdot (v_{m-1}+\overline{r}_{m})}\bkt{h_{m} - \frac{1}{3}(k^{*}+k^{c}/\sqrt{2})} \cdot \mathscr{G}_{(v_{m-1}+\overline{r}_{m})\cdot 2/3} (k^{*}+k^{c}/\sqrt{2}).
\end{align}
Moreover, observe that the following identity follows from Plancherel's theorem \cite[Theorem 5.3]{lieb2001analysis}:
\begin{align} \label{identity_mathscr_G_combine}
    \int_{\mathbb{R}^{d}} dw \mathscr{G}_{s_{1}}(w-k_{1})\mathscr{G}_{s_{2}}(w-k_{2}) = \mathscr{G}_{s_{1}\oplus s_{2}}(k_{1}-k_{2}) G_{s_{1}+s_{2}}(0),
\end{align}
where $a\oplus b$ is defined in (\ref{definition_oplus}). Then, thanks to (\ref{identity_mathscr_G_combine}), the following identity follows by combining (\ref{Am_last}) and (\ref{decompistion_Am}):
\begin{align*}
    &A_{m}(k^{*},k^{c},h_{m-1}) 
    = \mathscr{G}_{(v_{m-1}+\overline{r}_{m})\cdot 2/3} (k^{*}+k^{c}/\sqrt{2}) G_{v_{m-1}+3(v_{m-1}+\overline{r}_{m})}(0) \mathscr{G}_{\mathbf{t}_{m}}\bkt{h_{m-1} - \frac{1}{3}(k^{*}+k^{c}/\sqrt{2})},
\end{align*}
where $\mathbf{t}_{m}$ is given in (\ref{definition_tj}).

\paragraph{Step 3.} To compute $\mathscr{B}_{v_{1},r_{2},...,v_{m-1},r_{m},v_{m}}(k^{*},k^{c})$, the primary step is to show that the following expression solves the recursive equation (\ref{definition_A_j_kstar_kc}):
\begin{align} \label{identity_A_2}
    &A_{n}(k^{*},k^{c},h_{n-1})
    =\nonumber\\ &\mathscr{G}_{\mathbf{t}_{n}}\bkt{h_{n-1}-\frac{1}{3}(k^{*}+k^{c}/\sqrt{2})} 
    \biggbkt{\prod_{j = n}^{m-1} G_{v_{j-1}+(\mathbf{t}_{j+1}+3(v_{j-1}+r_{j}))}(0) \cdot \mathscr{G}_{(v_{j-1}+r_{j})\cdot 2/3}(k^{*}+k^{c}/\sqrt{2}) } \nonumber\\
    &G_{v_{m-1}+3(v_{m-1}+\overline{r}_{m})}(0) \mathscr{G}_{(v_{m-1}+\overline{r}_{m})\cdot 2/3} (k^{*}+k^{c}/\sqrt{2}) \quad \forall 3\leq n\leq m-1,
\end{align}
where $\mathbf{t}_{j+1}$ is defined in (\ref{definition_tj}). This is done by using mathematical induction. Observe that $A_{m}$ is also of the form (\ref{identity_A_2}), but the above product is $1$. Now, we assume that $A_{n}$ satisfies the above formula, where $4\leq n\leq m$. In the following, we show that $A_{n-1}$ has the above representation. Then, the identity below follows from (\ref{definition_A_j_kstar_kc}) and the assumption of mathematical induction:
\begin{align} 
    &A_{n-1}(k^{*},k^{c},h_{n-2})
    =  
    \biggbkt{\prod_{j = n}^{m-1} G_{v_{j-1}+(\mathbf{t}_{j+1}+3(v_{j-1}+r_{j}))}(0) \cdot \mathscr{G}_{(v_{j-1}+r_{j})\cdot 2/3}(k^{*}+k^{c}/\sqrt{2}) } \nonumber\\
    &\biggbkt{G_{v_{m-1}+3(v_{m-1}+\overline{r}_{m})}(0) \mathscr{G}_{(v_{m-1}+\overline{r}_{m})\cdot 2/3} (k^{*}+k^{c}/\sqrt{2})} \cdot I_{n-1}(k^{*},k^{c},h_{n-1}),\label{identity_A_j-1} \\
    &\text{where} \quad I_{n-1}(k^{*},k^{c},h_{n-1}):=\biggbkt{\int_{\mathbb{R}^{d}} dh_{n-1}  \mathscr{G}_{v_{n-2}+r_{n-1}}(\sqrt{2}h_{n-1})\mathscr{G}_{v_{n-2}}(-k^{*}-k^{c}/\sqrt{2}+2h_{n-2}+h_{n-1})\nonumber\\
    &\mathscr{G}_{v_{n-2}+r_{n-1}}(k^{*}+k^{c}/\sqrt{2}-h_{n-1}) \mathscr{G}_{\mathbf{t}_{n}}\bkt{h_{n-1}-\frac{1}{3}(k^{*}+k^{c}/\sqrt{2})}}.\label{identity_I_j-1}
\end{align}
Now, the key is to evaluate $I_{n-1}(k^{*},k^{c},h_{n-1})$. Similarly, to compute an integral consisting of four Gaussian kernels, one must first merge these Gaussian kernels. To this aim, we combine the first and the third Gaussian kernels in $I_{n-1}(k^{*},k^{c},h_{n-1})$:
\begin{align} \label{identity_A_j-1_0}
    \mathscr{G}_{v_{n-2}+r_{n-1}}(\sqrt{2}&h_{n-1}) \mathscr{G}_{v_{n-2}+r_{n-1}}(k^{*}+k^{c}/\sqrt{2}-h_{n-1}) \nonumber\\
    &= \mathscr{G}_{3\cdot (v_{n-2}+r_{n-1})}\bkt{h_{n-1}-\frac{1}{3}(k^{*}+k^{c}/\sqrt{2})} \mathscr{G}_{(v_{n-2}+r_{n-1})\cdot 2/3}(k^{*}+k^{c}/\sqrt{2}),
\end{align}
Moreover, the forth Gaussian kernel in $I_{n-1}(k^{*},k^{c},h_{n-1})$ and the first Gaussian kernel on the right-hand side of (\ref{identity_A_j-1_0}) can be merged as follows:
\begin{align} \label{identity_A_j-1_1}
    \mathscr{G}_{\mathbf{t}_{n}}\bkt{h_{n-1}-\frac{1}{3}(k^{*}+k^{c}/\sqrt{2})} &\mathscr{G}_{3\cdot (v_{n-2}+r_{n-1})}\bkt{h_{n-1}-\frac{1}{3}(k^{*}+k^{c}/\sqrt{2})}\nonumber\\
    &=\mathscr{G}_{\mathbf{t}_{n}+3\cdot (v_{n-2}+r_{n-1})}\bkt{h_{n-1}-\frac{1}{3}(k^{*}+k^{c}/\sqrt{2})}.
\end{align}
Then, thanks to (\ref{identity_A_j-1_0}) and (\ref{identity_A_j-1_1}), one has
\begin{align} \label{identity_A_j-1_2}
    &I_{n-1}(k^{*},k^{c},h_{n-1})=\mathscr{G}_{(v_{n-2}+r_{n-1})\cdot 2/3}(k^{*}+k^{c}/\sqrt{2}) \nonumber\\
    &\quad\quad\quad\int_{\mathbb{R}^{d}} dh_{n-1}  \mathscr{G}_{v_{n-2}}(-k^{*}-k^{c}/\sqrt{2}+2h_{n-2}+h_{n-1})\mathscr{G}_{\mathbf{t}_{n}+3\cdot (v_{n-2}+r_{n-1})}\bkt{h_{n-1}-\frac{1}{3}(k^{*}+k^{c}/\sqrt{2})}\nonumber\\
    &= \mathscr{G}_{(v_{n-2}+r_{n-1})\cdot 2/3}(k^{*}+k^{c}/\sqrt{2}) G_{v_{n-2}+(\mathbf{t}_{n}+3\cdot (v_{n-2}+r_{n-1}))}(0) \mathscr{G}_{\mathbf{t}_{n-1}}\bkt{h_{n-2}-\frac{1}{3}(k^{*}+k^{c}/\sqrt{2})},
\end{align}
where we have used (\ref{identity_mathscr_G_combine}) in the second equality. Then, (\ref{identity_A_2}) follows by combining (\ref{identity_A_j-1}) and (\ref{identity_A_j-1_2}).

\paragraph{Step 4.}
Finally, we evaluate $\mathscr{B}_{v_{1},r_{2},...,v_{m-1},r_{m},v_{m}}(k^{*},k^{c})$ by computing $A_{3}(k^{*},k^{c})$. Recalling (\ref{definition_A_2_kstar_kc}) and applying (\ref{identity_A_2}), we conclude
\begin{align} 
    &\mathscr{B}_{v_{1},r_{2},...,v_{m-1},r_{m},v_{m}}(k^{*},k^{c}) = 
    \biggbkt{\prod_{j = 3}^{m-1} G_{v_{j-1}+(\mathbf{t}_{j+1}+3(v_{j-1}+r_{j}))}(0) \cdot \mathscr{G}_{(v_{j-1}+r_{j})\cdot 2/3}(k^{*}+k^{c}/\sqrt{2}) } \nonumber\\
    &\biggbkt{G_{v_{m-1}+3(v_{m-1}+\overline{r}_{m})}(0) \mathscr{G}_{(v_{m-1}+\overline{r}_{m})\cdot 2/3} (k^{*}+k^{c}/\sqrt{2})} \cdot I_{2}(k^{*},k^{c}),\label{identity_A_2_by_A_3}\\
    &\text{where} \quad I_{2}(k^{*},k^{c}):=\biggbkt{\int_{\mathbb{R}^{d}} dh_{2} \mathscr{G}_{v_{1}+r_{2}} (\sqrt{2}h_{2})\nonumber\\
    &\mathscr{G}_{v_{1}}(-k^{*}+k^{c}/\sqrt{2}+h_{2})\mathscr{G}_{v_{1}+r_{2}}(k^{*}+k^{c}/\sqrt{2}-h_{2})\mathscr{G}_{\mathbf{t}_{3}}\bkt{h_{2}-\frac{1}{3}(k^{*}+k^{c}/\sqrt{2})}}.\label{identity_I_2}
\end{align}
Due to (\ref{identity_I_j-1}) and (\ref{identity_I_2}), we observe that $I_{2}(k^{*},k^{c})$ can be computed by using (\ref{identity_A_j-1_2}) with $h_{1} = 0$:
\begin{align} \label{identity_A_2_by_A_3_1}
    &I_{2}(k^{*},k^{c})
    = \mathscr{G}_{(v_{1}+r_{2}) \cdot 2/3}(k^{*}+k^{c}/\sqrt{2}) G_{v_{1}+(\mathbf{t}_{3}+3\cdot (v_{1}+r_{2}))}(0) \mathscr{G}_{\mathbf{t}_{1}}\bkt{\frac{1}{3}k^{*}-\frac{2}{3}\frac{k^{c}}{\sqrt{2}}},
\end{align}
Therefore, thanks to (\ref{identity_A_2_by_A_3}) and (\ref{identity_A_2_by_A_3_1}), one has
\begin{align*}
    \mathscr{B}_{v_{1},r_{2},...,v_{m-1},r_{m},v_{m}}(k^{*},k^{c})
    &= \mathscr{G}_{\mathbf{t}_{1}}\bkt{\frac{1}{3}k^{*}-\frac{2}{3}\frac{k^{c}}{\sqrt{2}}}\mathscr{G}_{(v_{1}+r_{2}+...+v_{m-1}+\overline{r}_{m})\cdot 2/3} (k^{*}+k^{c}/\sqrt{2}) \\&\biggbkt{\prod_{j=2}^{m-1} G_{v_{j-1}+(\mathbf{t}_{j+1}+3(v_{j-1}+r_{j}))}(0)} G_{v_{m-1}+3(v_{m-1}+\overline{r}_{m})}(0),
\end{align*}
where we have used the fact that
\begin{align*}
    &\mathscr{G}_{(v_{1}+r_{2}+...+v_{m-1}+\overline{r}_{m})\cdot 2/3} (k^{*}+k^{c}/\sqrt{2})= \biggbkt{\prod_{j = 2}^{m-1} \mathscr{G}_{(v_{j-1}+r_{j})\cdot 2/3}(k^{*}+k^{c}/\sqrt{2}) } \mathscr{G}_{(v_{m-1}+\overline{r}_{m})\cdot 2/3} (k^{*}+k^{c}/\sqrt{2}).
\end{align*}
Hence, we obtain (\ref{integral_h1...h_m}). Then, the proof of Lemma \ref{lemma_identity_fourier_2} is complete.

\subsubsection{Estimates for the final Gaussian kernels.} 
Now, we are ready to prove the lower bound in Proposition \ref{proposition_partII_second_id}. By making use of Lemma \ref{lemma_identity_fourier_2}, it remains to evaluate the inverse Fourier transform of the following functions:
\begin{align} \label{fourier_integrand}
    \mathscr{G}_{v_{0}}(k)  \mathscr{G}_{v_{0}+r_{1}}(k^{*}) \mathscr{G}_{v_{0}+r_{1}}(k^{c})\mathscr{G}_{\mathbf{t}_{1}}\bkt{\frac{1}{3}k^{*}-\frac{2}{3}\frac{k^{c}}{\sqrt{2}}}
    \mathscr{G}_{(v_{1}+r_{2}+...+v_{m-1}+\overline{r}_{m})\cdot 2/3} (k^{*}+k^{c}/\sqrt{2}),
\end{align}
where (\ref{fourier_integrand}) consists of all the functions on the right-hand side of (\ref{identity_lemma_identity_fourier_2}) with variables $k$, $k^{*}$, and $k^{c}$. Then, it is enough to compute the following Fourier integral.
\begin{lemma} \label{lemma_fourier_integral}
There exist functions $u_{1}$, $u_{2}$, and $u_{3}$ with variables $v_{0}$,$r_{1}$,...,$v_{m-1}$,$r_{m}$,$t$, and $\nu$ such that 
\begin{align} \label{integral_k_kc_kstar_fourier}
    &\int_{\mathbb{R}^{3 \times d}} dkdk^{*}dk^{c} \exp\biggbkt{2\pi i k \cdot \Vec{z}_{0}^{(\ell_{1},\ell_{1}')}} \exp\biggbkt{2\pi i k^{*} \cdot \Vec{z}_{0}^{(\ell_{1},\ell_{1}')^{*}}}
    \exp\biggbkt{2\pi i k^{c} \cdot \Vec{z}_{0}^{(\ell_{1},\ell_{1}')^{c}}}\nonumber\\
    &\mathscr{G}_{v_{0}}(k)  \mathscr{G}_{v_{0}+r_{1}}(k^{*}) \mathscr{G}_{v_{0}+r_{1}}(k^{c})\mathscr{G}_{\mathbf{t}_{1}}\bkt{\frac{1}{3}k^{*}-\frac{2}{3}\frac{k^{c}}{\sqrt{2}}}
    \mathscr{G}_{(v_{1}+r_{2}+...+v_{m-1}+\overline{r}_{m})\cdot 2/3} (k^{*}+k^{c}/\sqrt{2})\nonumber\\
    &= G_{v_{0}}\bkt{\Vec{z}_{0}^{(\ell_{1},\ell_{1}')}} G_{u_{1}}\bkt{\Vec{z}_{0}^{(\ell_{1},\ell_{1}')^{c}}} G_{u_{2}}\bkt{\Vec{z}_{0}^{(\ell_{1},\ell_{1}')^{*}}+u_{3}\Vec{z}_{0}^{(\ell_{1},\ell_{1}')^{c}}}.
\end{align}
Here, $u_{1}$, $u_{2}$, and $u_{3}$ satisfy the following estimate:
\begin{align} \label{estimate_u1u2u3} 
    c_{1,1}(t+\nu) \leq u_{1} \leq c_{1,2}(t+\nu), \quad c_{2,1}v_{0} \leq u_{2} \leq c_{2,2} (t+\nu), \quad\text{and}\quad |u_{3}| \leq c_{3},
\end{align}
where $c_{1,1}$, $c_{1,2}$, $c_{2,1}$, $c_{2,2}$ and $c_{3}$ are positive constants, and $\overline{r}_{m}$ and $\mathbf{t}_{1}$ are defined in (\ref{definition_overline_rm}) and (\ref{definition_tj}), respectively.
\end{lemma}
Let us first explain the reason why Lemma \ref{lemma_fourier_integral} allows us to conclude Proposition \ref{proposition_partII_second_id}.
\begin{proof} [Proof of Proposition \ref{proposition_partII_second_id}]
Recall that (\ref{fourier_integrand}) consists of all the functions on the right-hand side of (\ref{identity_lemma_identity_fourier_2}) with variables $k$, $k^{*}$, and $k^{c}$. Also, their inverse Fourier transforms are computed in (\ref{integral_k_kc_kstar_fourier}). Moreover, since $|x+y|^{2} \leq 2|x|^{2}+2|y|^{2}$ for all $x,y \in \mathbb{R}^{d}$, the following estimates follows from (\ref{estimate_u1u2u3}):
\begin{align*}
    &G_{v_{0}}\bkt{\Vec{z}_{0}^{(\ell_{1},\ell_{1}')}} \geq C \frac{1}{|v_{0}|^{3/2}} \exp\bkt{-\frac{|\Vec{z}_{0}^{(\ell_{1},\ell_{1}')}|^{2}}{2v_{0}}}, \quad G_{u_{1}}\bkt{\Vec{z}_{0}^{(\ell_{1},\ell_{1}')^{c}}} \geq C_{t,\nu} \exp\bkt{-\frac{|\Vec{z}_{0}^{(\ell_{1},\ell_{1}')^{c}}|^{2}}{2c_{1,1}v_{0}}}, \\
    &\text{and}\quad G_{u_{2}}\bkt{\Vec{z}_{0}^{(\ell_{1},\ell_{1}')^{*}}+u_{3}\Vec{z}_{0}^{(\ell_{1},\ell_{1}')^{c}}} \geq C_{t,\nu} \exp\bkt{-\frac{|\Vec{z}_{0}^{(\ell_{1},\ell_{1}')^{*}}|^{2}}{c_{2,1}v_{0}}} \exp\bkt{-\frac{c_{3}^{2} |\Vec{z}_{0}^{(\ell_{1},\ell_{1}')^{c}}|^{2}}{c_{2,1}v_{0}}}.
\end{align*}
Recall that $\Vec{z}_{0} = \Vec{x}_{0}$ from Lemma \ref{lemma_relative_motion}. Then, by Lemma \ref{lemma_identity_fourier_2}, Proposition \ref{proposition_partII_second_id} follows from the above lower bounds.
\end{proof}
To complete the proof of Proposition \ref{proposition_partII_second_id}, it remains to prove Lemma \ref{lemma_fourier_integral}.
\begin{proof} [Proof of Lemma \ref{lemma_fourier_integral}]
To begin with, thanks to (\ref{formula_G_t}), one has
\begin{align*}
    &\int_{\mathbb{R}^{d}} dk 
    \exp\biggbkt{2\pi i k \cdot \Vec{z}_{0}^{(\ell_{1},\ell_{1}')}}
    \mathscr{G}_{v_{0}}(k)
    = G_{v_{0}}\bkt{\Vec{z}_{0}^{(\ell_{1},\ell_{1}')}}.
\end{align*}
Then, to prove (\ref{integral_k_kc_kstar_fourier}), it remains to show that 
\begin{align} \label{identity_gaussian_kc_kstar}
    &\int_{\mathbb{R}^{2 \times d}} dk^{*}dk^{c} \exp\biggbkt{2\pi i k^{*} \cdot \Vec{z}_{0}^{(\ell_{1},\ell_{1}')^{*}}}
    \exp\biggbkt{2\pi i k^{c} \cdot \Vec{z}_{0}^{(\ell_{1},\ell_{1}')^{c}}}\mathscr{G}_{v_{0}+r_{1}}(k^{*}) \mathscr{G}_{v_{0}+r_{1}}(k^{c})\mathscr{G}_{\mathbf{t}_{1}}\bkt{\frac{1}{3}k^{*}-\frac{2}{3}\frac{k^{c}}{\sqrt{2}}}\nonumber\\
    &\mathscr{G}_{(v_{1}+r_{2}+...+v_{m-1}+\overline{r}_{m})\cdot 2/3} (k^{*}+k^{c}/\sqrt{2}) = G_{u_{1}}\bkt{\Vec{z}_{0}^{(\ell_{1},\ell_{1}')^{c}}} G_{u_{2}}\bkt{\Vec{z}_{0}^{(\ell_{1},\ell_{1}')^{*}}+u_{3}\Vec{z}_{0}^{(\ell_{1},\ell_{1}')^{c}}},
\end{align}
and $u_{1},u_{2},$ and $u_{3}$ satisfy the conditions (\ref{estimate_u1u2u3}).
To this aim, we first simplify the integrands on the left-hand side of (\ref{identity_gaussian_kc_kstar}). Notice that the Gaussian kernels on the left-hand side of (\ref{identity_gaussian_kc_kstar}) can be re-expressed as follows:
\begin{align*}
    \exp(-2\pi^{2}u_{1}|k^{c} + u_{3}k^{*}|^{2}) \exp(-2\pi^{2}u_{2}|k^{*}|^{2}),
\end{align*}
where $u_{j}$ is defined as follows:
\begin{align*}
    &u_{1}:= v_{0}+r_{1}+\frac{2}{9}\mathbf{t}_{1} + \frac{1}{3}(v_{1}+r_{2}+...+v_{m-1}+\overline{r}_{m}), \\
    &u_{3} := \frac{1}{u_{1}} \biggbkt{-\mathbf{t}_{1}\frac{\sqrt{2}}{9}+\frac{\sqrt{2}}{3} (v_{1}+r_{2}+...+v_{m-1}+\overline{r}_{m})},\\
    &u_{2} := u_{4}-u_{1}u_{3}^{2},\quad \text{and} \quad u_{4}:= v_{0}+r_{1}+\frac{1}{9}\mathbf{t}_{1} + \frac{2}{3}(v_{1}+r_{2}+...+v_{m-1}+\overline{r}_{m}).
\end{align*}
Then, the integrals on the left-hand side of (\ref{identity_gaussian_kc_kstar}) can be computed in the following order:
\begin{align*}
    \int_{\mathbb{R}^{d}} dk^{c} 
    \exp\biggbkt{2\pi i k^{c} \cdot \Vec{z}_{0}^{(\ell_{1},\ell_{1}')^{c}}}\exp(-2\pi^{2}u_{1}|k^{c} - u_{3}k^{*}|^{2})
    = \exp\biggbkt{2\pi i k^{*} \cdot u_{3} \Vec{z}_{0}^{(\ell_{1},\ell_{1}')^{c}}} G_{u_{1}}\bkt{\Vec{z}_{0}^{(\ell_{1},\ell_{1}')^{c}}}
\end{align*}
and
\begin{align*}
    &\int_{\mathbb{R}^{d}} dk^{*} 
    \exp\biggbkt{2\pi i k^{*} \cdot (\Vec{z}_{0}^{(\ell_{1},\ell_{1}')^{*}}+u_{3} \Vec{z}_{0}^{(\ell_{1},\ell_{1}')^{c}})} \exp(-2\pi^{2}u_{2}|k^{*}|^{2})= G_{u_{2}}\bkt{\Vec{z}_{0}^{(\ell_{1},\ell_{1}')^{*}}+u_{3} \Vec{z}_{0}^{(\ell_{1},\ell_{1}')^{c}}}.
\end{align*}
Therefore, we conclude (\ref{identity_gaussian_kc_kstar}).

Finally, to complete the proof of Lemma \ref{lemma_fourier_integral}, it remains to show that $u_{1}$, $u_{2}$, and $u_{3}$ satisfy the conditions (\ref{estimate_u1u2u3}). To do this, we first recall that 
\begin{align*}
    \mathbf{t}_{1} = 4(v_{0}\oplus(\mathbf{t}_{2}+ 3(v_{0}+r_{1}))) \quad \text{and} \quad \overline{r}_{m} = r_{m}+v_{m}+\nu.
\end{align*}
Due to (\ref{definition_oplus}), it is clear that $\mathbf{t}_{1} \approx v_{0}$. Consequently, since $v_{0}+r_{1}+...+v_{m-1}+r_{m}+v_{m} = t$, where $v_{m}$ is given in Lemma \ref{lemma_relative_motion}, one has $u_{1} \approx (t+\nu)$. Moreover, $|u_{3}|\approx 1$ follows from 
\begin{align*}
    -1\lesssim \frac{1}{u_{1}} \biggbkt{-\mathbf{t}_{1}\frac{\sqrt{2}}{9}}\leq u_{3} \leq \frac{1}{u_{1}} \biggbkt{\frac{\sqrt{2}}{3} (v_{1}+r_{2}+...+v_{m-1}+\overline{r}_{m})} \lesssim 1.
\end{align*}
Now, it remains to estimate $u_{2}$. About the upper bound of $u_{2}$, it is clear that $u_{2} \leq u_{4} \lesssim (t+\nu)$. Regarding the lower bound of $u_{2}$, we note that 
\begin{align*}
    &u_{2} = \frac{1}{u_{1}} \biggbkt{\bkt{v_{0}+r_{1}+\frac{1}{9}\mathbf{t}_{1} + \frac{2}{3}(v_{1}+r_{2}+...+v_{m-1}+\overline{r}_{m})} \\
    &\bkt{v_{0}+r_{1}+\frac{2}{9}\mathbf{t}_{1} + \frac{1}{3}(v_{1}+r_{2}+...+v_{m-1}+\overline{r}_{m})}-\bkt{-\mathbf{t}_{1}\frac{\sqrt{2}}{9}+\frac{\sqrt{2}}{3} (v_{1}+r_{2}+...+v_{m-1}+\overline{r}_{m})}^{2}}\\
    &= \frac{1}{u_{1}} \biggbkt{ (v_{0}+r_{1}) \frac{3}{9} \mathbf{t}_{1} + \bkt{\frac{1}{9}\mathbf{t}_{1} + \frac{2}{3}(v_{1}+r_{2}+...+v_{m-1}+\overline{r}_{m})}\bkt{\frac{2}{9}\mathbf{t}_{1} + \frac{1}{3}(v_{1}+r_{2}+...+v_{m-1}+\overline{r}_{m})}\\
    &-\bkt{-\mathbf{t}_{1}\frac{\sqrt{2}}{9}+\frac{\sqrt{2}}{3} (v_{1}+r_{2}+...+v_{m-1}+\overline{r}_{m})}^{2}}\gtrsim  \frac{t+\nu}{(t+\nu)} \mathbf{t}_{1} \gtrsim v_{0}.
\end{align*}
Then, the proof of Lemma \ref{lemma_fourier_integral} is complete.
\end{proof}

\subsection{Temporal estimates for the sub-limiting path integrals.} \label{section_proof_proposition_lower_bound_sum_of_path_integrals}
This subsection is dedicated to the proof of Proposition \ref{proposition_lower_bound_sum_of_path_integrals}. The overall goal is to prove the following delicate temporal estimate. 
\begin{lemma} \label{lemma_simplex_integral_lower_bound}
For every $t>0$, there exists a positive constant $C_{t}$ such that the following estimate holds:
\begin{align} \label{estimate_simplex_integral_lower_bound}
    &\biggbkt{\frac{1}{\pi}}^{m} \cdot \int_{v_{j},r_{j}>0, \;v_{1}+r_{2}+v_{2}+...+v_{m-1}+r_{m}+v_{m} <t}\biggbkt{\prod_{j=1}^{m} dv_{j} \prod_{j=2}^{m} dr_{j}}  \cdot \biggbkt{\prod_{j=2}^{m} \frac{1}{(r_{j})^{1/2}} \frac{1}{(v_{j-1}+v_{j}+3r_{j}/4)^{3/2}}}\nonumber\\&\geq C_{t} \cdot  \biggbkt{1.008}^{m} \frac{1}{m(m-1)} \quad \forall m\geq 3.
\end{align}
\end{lemma}
For the sake of simplicity, we postpone the proof of Lemma \ref{lemma_simplex_integral_lower_bound} in Section \ref{section_proof_lemma_simplex_integral_lower_bound}. Let us now present the proof of Proposition \ref{proposition_lower_bound_sum_of_path_integrals}. Note that we use $C$ to denote a positive constant that depends on $t$, $\nu$, and $d$ throughout the proof.
\begin{proof} [Proof of Proposition \ref{proposition_lower_bound_sum_of_path_integrals}]
To begin with, recalling that $\overline{r}_{m} = r_{m}+v_{m}+\nu$ and $v_{0}+r_{1}+...+r_{m} +v_{m} = t$, we have
\begin{align} \label{estimate_time_integral_2}
    &\text{((R.H.S) of (\ref{estimate_identity_fourier_3}))}\geq  C\int_{v_{j},r_{j}>0, \; v_{0}+r_{1}+...+v_{m-1}+r_{m} <t, \; v_{m} = t-(v_{0}+r_{1}+...+v_{m-1}+r_{m})} \prod_{j=1}^{m} dv_{j-1} dr_{j} \nonumber\\
    &\frac{1}{v_{0}^{3/2}} \exp\biggbkt{-\frac{|\Vec{x}_{0}|^{2}}{cv_{0}}}\cdot\biggbkt{ \prod_{j=1}^{m}\sqrt{\frac{2\pi}{r_{j}}} }
    \cdot\biggbkt{ \prod_{j=2}^{m-1} G_{(v_{j-1}+(\mathbf{t}_{j+1}+3(v_{j-1}+r_{j})))/4}(0)},
\end{align}
where we have bounded both of $G_{v_{m-1}+3(v_{m-1}+\overline{r}_{m})}(0)$ and $G_{v_{m}+\nu}(0)$ from below by $G_{t+\nu}(0)$ in the first inequality, and $\sqrt{2}^{2 \cdot 3}\cdot G_{t}(0) = G_{t/4}(0)$ has been used. To obtain the expression on the left-hand side of (\ref{estimate_simplex_integral_lower_bound}), we now estimate $\mathbf{t}_{j+1}$ from above in order to obtain a lower bound for $G_{(v_{j-1}+(\mathbf{t}_{j+1}+3(v_{j-1}+r_{j})))/4}(0)$. Recalling (\ref{definition_tj}), we have the following estimate:
\begin{align*}
    \frac{1}{\mathbf{t}_{j+1}/4} \geq \frac{1}{v_{j}} \quad \forall 2\leq j\leq m-2 \quad \text{and} \quad\frac{1}{\mathbf{t}_{m}/4} \geq \frac{1}{v_{m-1}}.
\end{align*}
Then, it holds that
\begin{align} \label{estimate_time_integral_3}
    &\text{(R.H.S) of (\ref{estimate_time_integral_2})}
    \geq C\int_{v_{j},r_{j}>0, \; v_{0}+r_{1}+...+v_{m-1}+r_{m} <t, \; v_{m} = t-(v_{0}+r_{1}+...+v_{m-1}+r_{m})} \prod_{j=1}^{m} dv_{j-1} dr_{j} \nonumber\\
    &\frac{1}{v_{0}^{3/2}} \exp\biggbkt{-\frac{|\Vec{x}_{0}|^{2}}{cv_{0}}} \cdot \biggbkt{\prod_{j=1}^{m}\sqrt{\frac{2\pi}{r_{j}}} }
    \cdot \biggbkt{ \prod_{j=2}^{m-1} G_{v_{j-1}+v_{j}+3r_{j}/4}(0)}.
\end{align}
Furthermore, we separate the integrals with respect to the variables $v_{0},r_{1},r_{m}$ from the rest of the variables on the right-hand side of the above estimate:
\begin{align} \label{estimate_time_integral_4}
    &\text{((R.H.S) of (\ref{estimate_time_integral_3}))}\nonumber \geq C\int_{v_{j},r_{j}>0, \; v_{1}+r_{2}+...+r_{m-1}+v_{m-1} <t/2} \int_{0<v_{0}+r_{1}+r_{m}<t/2}\prod_{j=1}^{m} dv_{j-1} dr_{j} \nonumber\\
    &\frac{1}{v_{0}^{3/2}} \exp\biggbkt{-\frac{|\Vec{x}_{0}|^{2}}{cv_{0}}} \cdot\biggbkt{\prod_{j=1}^{m}\sqrt{\frac{2\pi}{r_{j}}} }
    \cdot \biggbkt{\prod_{j=2}^{m-1} G_{v_{j-1}+v_{j}+3r_{j}/4}(0)}\nonumber\\
    &\geq C \biggbkt{\frac{1}{2}}^{m-2} \cdot\frac{1}{|\Vec{x}_{0}|} \cdot\biggbkt{ \frac{1}{\pi^{m-2}}\int_{v_{j},r_{j}>0, \; v_{1}+r_{2}+...+r_{m-1}+v_{m-1}<t/2} \prod_{j=1}^{m-1} dv_{j} \prod_{j=2}^{m-1} dr_{j} \nonumber\\
    &\prod_{j=2}^{m-1} \frac{1}{r_{j}^{1/2}} \cdot 
    \prod_{j=2}^{m-1} \frac{1}{(v_{j-1}+v_{j}+3r_{j}/4)^{3/2}}},
\end{align}
where following estimates have been used:
\begin{align*}
    \int_{0}^{t/6} \frac{1}{v_{0}^{3/2}} \exp\biggbkt{-\frac{|\Vec{x}_{0}|^{2}}{cv_{0}}} dv_{0} = C \frac{1}{|\Vec{x}_{0}|} \quad \text{and} \quad \int_{0}^{t/6} \frac{1}{r_{j}^{1/2}} dr_{j} = \biggbkt{\frac{t}{6}}^{1/2} \quad \forall j = 1,m.
\end{align*}
Then, thanks to Lemma \ref{lemma_simplex_integral_lower_bound}, we conclude
\begin{align} \label{estimate_time_integral_5}
    \text{(R.H.S) of (\ref{estimate_time_integral_4})}
    \geq C\frac{1}{|\Vec{x}_{0}|} \cdot\frac{1}{2^{m-2}} \biggbkt{1.008}^{m-1} \frac{1}{(m-1)(m-2)}.
\end{align}

Finally, we complete the proof of Proposition \ref{proposition_lower_bound_sum_of_path_integrals} by using the above results. By collecting the above estimates  (\ref{estimate_time_integral_2}), (\ref{estimate_time_integral_3}),  (\ref{estimate_time_integral_4}), and 
(\ref{estimate_time_integral_5}), we know that the right-hand side of (\ref{estimate_time_integral_5}) is a lower bound for the sub-limiting path integral with length $m$ uniformly for each $\bigotimes_{j=1}^{m} (\ell_{j},\ell_{j}',\mathfrak{i}_{j}) \in \bigotimes_{j=1}^{m} (\mathcal{E}^{(N)} \times \{1\} )$ such that $(\ell_{j},\ell_{j}') \neq (\ell_{j-1},\ell_{j-1}')$ for all $2\leq j\leq m$. Then, one has
\begin{align*}
    &\sum_{(\ell_{1},\ell'_{1},\mathfrak{i}_{1}),...,(\ell_{m},\ell'_{m},\mathfrak{i}_{m}) \in \mathcal{E}^{(N)} \times \{1\}, \; (\ell_{1},\ell_{1}') \neq ...\neq (\ell_{k},\ell_{m}')} \mathscr{I}_{0;t}^{N;(\ell_{1},\ell'_{1},\mathfrak{i}_{1}),...,(\ell_{m},\ell'_{m},\mathfrak{i}_{m})} U_{0}^{\otimes N}(\Vec{x}_{0})\\
    &\geq \biggbkt{3\cdot 2^{m-1} } \cdot \biggbkt{C\frac{1}{|\Vec{x}_{0}|}  \cdot\frac{1}{2^{m-2}} \biggbkt{1.008}^{m-1} \frac{1}{(m-1)(m-2)}} \geq C \frac{1}{|\Vec{x}_{0}|} \cdot \biggbkt{1.008}^{m-1} \frac{1}{(m-1)(m-2)}.
\end{align*}
Therefore, the proof of Proposition \ref{proposition_lower_bound_sum_of_path_integrals} is complete.
\end{proof}

\subsubsection{Singular iterated integral.} \label{section_proof_lemma_simplex_integral_lower_bound}
To complete the proof of Proposition \ref{proposition_lower_bound_sum_of_path_integrals}, it remains to prove Lemma \ref{lemma_simplex_integral_lower_bound}. We will carry out the proof in the following steps.

\paragraph{Step1: Lower bound for the temporal integrals.} 
The starting point of our proof is to compute all the integrals with respect to $dr_{j}$ on the left-hand side of (\ref{estimate_simplex_integral_lower_bound}). The primary reason is that the computation eventually reduces the left-hand side of (\ref{estimate_simplex_integral_lower_bound}) to an iterated integral of single variable, which is much easier to control. 

\begin{lemma} \label{lemma_goal_1_lemma_simplex_integral_lower_bound}
There exists a positive constant $\delta$ depending on $t$ such that
\begin{align} \label{goal_1_lemma_simplex_integral_lower_bound}
    &\text{(R.H.S) of (\ref{estimate_simplex_integral_lower_bound})}\geq\frac{1}{2} \bkt{\frac{3}{4}}^{1/2} \cdot \biggbkt{\bkt{\frac{4}{3}}^{1/2} \frac{2}{\pi} \cdot 0.99}^{m} \cdot \frac{\delta}{m(m-1)} \cdot L_{m} \quad \forall m\geq 3.
\end{align}    
Here, $L_{m}$ is the singular iterated integral defined as follows:
\begin{align} \label{definition_Lm}
    L_{m+1} := \int_{0}^{1} \eta_{m}(v) dv \quad \forall m\geq 0.
\end{align}
where the function $\eta_{m}(v):[0,1] \mapsto \mathbb{R}_{+}$ is given by
\begin{align} \label{definition_eta_m}
    \eta_{m}(v) := \int_{0}^{1} \frac{1}{v+s} \eta_{m-1}(s) ds \quad \forall m\geq 1 \quad \text{and} \quad \eta_{0}(v) := 1.
\end{align}
\end{lemma}

The key of the proof is to make use of a property of \emph{the homogeneous iterated integrals}. For example, thanks to a change of variable, the growth rate of the left-hand side of (\ref{estimate_simplex_integral_lower_bound}) with respect to $m$ is independent of the size of the domain.

\begin{proof} [Proof of Lemma \ref{lemma_goal_1_lemma_simplex_integral_lower_bound}.]
To begin with, we analyze the integral with respect to $dr_{j}$. By applying a change of variable, one has
\begin{align} \label{estimate_proof_lemma_simplex_integral_lower_bound_1}
    &\text{(L.H.S) of (\ref{goal_1_lemma_simplex_integral_lower_bound})} \geq \bkt{\frac{3}{4}}^{1-1/2}\biggbkt{\bkt{\frac{3}{4}}^{-1+1/2} \frac{1}{\pi}}^{m} \cdot \nonumber\\
    &\int_{v_{j},r_{j}>0, \;v_{1}+4r_{2}/3+v_{2}+...+v_{m-1}+4r_{m}/3+v_{m} <t} \biggbkt{\prod_{j=1}^{m} dv_{j} \prod_{j=2}^{m} dr_{j} } \cdot \biggbkt{\prod_{j=2}^{m} \frac{1}{(r_{j})^{1/2}} \frac{1}{(v_{j-1}+v_{j}+r_{j})^{3/2}}}\nonumber\\
    &\geq \bkt{\frac{3}{4}}^{1/2} \biggbkt{\bkt{\frac{4}{3}}^{1/2} \frac{1}{\pi}}^{m} \cdot \nonumber\\
    &\int_{v_{j},r_{j}>0, \;v_{1}+r_{2}+v_{2}+...+v_{m-1}+r_{m}+v_{m} < 3t/4} \biggbkt{\prod_{j=1}^{m} dv_{j} \prod_{j=2}^{m} dr_{j} } \cdot \biggbkt{\prod_{j=2}^{m} \frac{1}{(r_{j})^{1/2}} \frac{1}{(v_{j-1}+v_{j}+r_{j})^{3/2}}},
\end{align}
where we have decreased the domain of the integral in the second inequality.  Now, we compute the integral with respect to $dr_{j}$ by using 
\begin{align} \label{integral_au}
    \int_{0}^{u} \frac{1}{r^{1/2}(a+r)^{3/2}} dr = \frac{2}{a} \biggbkt{\frac{u}{a+u}}^{1/2} \quad \forall u, a > 0.
\end{align}
Before doing so, we decrease the domain as follows:
\begin{align}\label{estimate_proof_lemma_simplex_integral_lower_bound_2}
    &\text{(R.H.S) of (\ref{estimate_proof_lemma_simplex_integral_lower_bound_1})} \nonumber\\
    &\geq \bkt{\frac{3}{4}}^{1/2} \biggbkt{\bkt{\frac{4}{3}}^{1/2} \frac{1}{\pi}}^{m} \cdot \int_{v_{j},r_{j}>0, \;v_{1}+v_{2}+...+v_{m-1}+v_{m} < 3t/8,\; r_{2}<(3t/8)/(m-1),\;...,\;r_{m}<(3t/8)/(m-1)} \nonumber\\
    &\prod_{j=1}^{m} dv_{j} \prod_{j=2}^{m} dr_{j}  \cdot \prod_{j=2}^{m} \frac{1}{(r_{j})^{1/2}} \frac{1}{(v_{j-1}+v_{j}+r_{j})^{3/2}}. 
\end{align}
Then, applying (\ref{integral_au}) shows that
\begin{align} \label{estimate_proof_lemma_simplex_integral_lower_bound_3}
    &\text{(R.H.S) of (\ref{estimate_proof_lemma_simplex_integral_lower_bound_2})}
    \geq \frac{1}{2} \bkt{\frac{3}{4}}^{1/2} \biggbkt{\bkt{\frac{4}{3}}^{1/2} \frac{2}{\pi}}^{m} \cdot \nonumber\\
    &\int_{v_{j}>0, \;v_{1}+v_{2}+...+v_{m-1}+v_{m} < 3t/8}\prod_{j=1}^{m} dv_{j} \cdot \biggbkt{\prod_{j=2}^{m} \frac{1}{v_{j-1}+v_{j}} \biggbkt{\frac{\frac{1}{m-1} (3t/8)}{v_{j-1}+v_{j}+\frac{1}{m-1} (3t/8)}}^{1/2}}\nonumber\\
    &= \frac{1}{2} \bkt{\frac{3}{4}}^{1/2}  \biggbkt{\bkt{\frac{4}{3}}^{1/2} \frac{2}{\pi}}^{m} \cdot (m-1)^{(m-1)-m}\int_{v_{j}>0, \;v_{1}+v_{2}+...+v_{m-1}+v_{m} < 3(m-1)t/8} \nonumber\\
    &\prod_{j=1}^{m} dv_{j} \cdot \biggbkt{\prod_{j=2}^{m} \frac{1}{v_{j-1}+v_{j}} \biggbkt{\frac{(3t/8)}{v_{j-1}+v_{j}+(3t/8)}}^{1/2}},
\end{align}
where we have used change of variables in the first equality. 

To obtain $L_{m}$, it remains to remove the term on the right of the above product. Also, notice that the above domain depends on $m$. Then, the main idea to bypass these two problems is to further decrease the domain of the above integral so that it is independent of $m$ and
\begin{align*}
    \biggbkt{\frac{(3t/8)}{v_{j-1}+v_{j}+(3t/8)}}^{1/2} \longrightarrow 1 \quad \text{when domain is small.}
\end{align*}
Thanks to the homogeneity, we will see that shrinking the domain does not affect the growth rate with respect to $m$. Choose a small positive number $\delta$ depending on $t$ such that
\begin{align}
    3t\cdot2/8 > \delta > 0 \quad \text{and}\quad \biggbkt{\frac{(3t/8)}{\delta+(3t/8)}}^{1/2} \geq 0.99.
\end{align}
Then, thanks to the above properties of $\delta$ and the fact that $m\geq 3$, we know that
\begin{align} \label{estimate_proof_lemma_simplex_integral_lower_bound_4}
    &\text{(R.H.S) of (\ref{estimate_proof_lemma_simplex_integral_lower_bound_3})}\nonumber\\
    &\geq \frac{1}{2} \bkt{\frac{3}{4}}^{1/2} \biggbkt{\bkt{\frac{4}{3}}^{1/2} \frac{2}{\pi} \cdot 0.99}^{m}  \cdot \frac{1}{m-1}\int_{v_{j}>0, \;v_{1}+v_{2}+...+v_{m-1}+v_{m} < \delta}  \prod_{j=1}^{m} dv_{j} \cdot \prod_{j=2}^{m} \frac{1}{v_{j-1}+v_{j}},
\end{align}
where we have used the facts that $\delta < (m-1) \cdot 3t/8$ for all $m\geq 3$ and that $v_{j-1}+v_{j} < \delta$ for all $2\leq j\leq m$. 
Consequently, by using the homogeneity of the iterated integral on the right-hand side of (\ref{estimate_proof_lemma_simplex_integral_lower_bound_4}), the following estimates follows from change of variables:
\begin{align}\label{estimate_proof_lemma_simplex_integral_lower_bound_5}
    &\text{(R.H.S) of (\ref{estimate_proof_lemma_simplex_integral_lower_bound_4})}  \nonumber\\
    &= \frac{1}{2} \bkt{\frac{3}{4}}^{1/2} \biggbkt{\bkt{\frac{4}{3}}^{1/2} \frac{2}{\pi} \cdot 0.99}^{m}  \cdot \frac{\delta}{m-1}  \int_{v_{j}>0, \;v_{1}+v_{2}+...+v_{m-1}+v_{m} < 1}  \prod_{j=1}^{m} dv_{j} \cdot \prod_{j=2}^{m} \frac{1}{v_{j-1}+v_{j}}\nonumber\\
    &\geq \frac{1}{2} \bkt{\frac{3}{4}}^{1/2} \biggbkt{\bkt{\frac{4}{3}}^{1/2} \frac{2}{\pi} \cdot 0.99}^{m}  \cdot \frac{\delta}{m-1}  \int_{0<v_{j}<1/m, \; j = 1,2,...,m}  \prod_{j=1}^{m} dv_{j} \cdot \prod_{j=2}^{m} \frac{1}{v_{j-1}+v_{j}}\nonumber\\
    &= \frac{1}{2} \bkt{\frac{3}{4}}^{1/2}\biggbkt{\bkt{\frac{4}{3}}^{1/2} \frac{2}{\pi} \cdot 0.99}^{m} \cdot \frac{\delta}{m(m-1)} \cdot L_{m}.
\end{align}
Therefore, by collecting (\ref{estimate_proof_lemma_simplex_integral_lower_bound_1}), (\ref{estimate_proof_lemma_simplex_integral_lower_bound_2}),  (\ref{estimate_proof_lemma_simplex_integral_lower_bound_3}), 
(\ref{estimate_proof_lemma_simplex_integral_lower_bound_4}), and
(\ref{estimate_proof_lemma_simplex_integral_lower_bound_5}), we complete the proof of the estimate (\ref{goal_1_lemma_simplex_integral_lower_bound}).
\end{proof}
%Not large enough
\paragraph{Step 2: Estimate for the singular iterated integral.}
To complete the proof of Lemma \ref{lemma_simplex_integral_lower_bound}, it remains to estimate the singular iterated integral $L_{m}$. 
\begin{lemma} \label{lemma_estimate_Lm}
Recall that $L_{m}$ is defined in (\ref{definition_Lm}). Then, the following estimate holds:
\begin{align*}
    L_{m} \geq (2\ln(2))^{m-1} \quad \forall m\geq 1.
\end{align*}
\end{lemma}

To prove the above lemma, we consider an auxiliary function which relates to $\theta_{m}(v)$. 
\begin{lemma}
Let $\zeta_{m}(v):[0,1] \mapsto \mathbb{R}_{+}$ be the function defined as follows:
\begin{align} \label{definition_zeta}
    \zeta_{m}(v) := \int_{v}^{1} \frac{1}{v+s} \zeta_{m-1}(s) ds \quad \forall m\geq 1 \quad \text{and} \quad \zeta_{0}(v) := 1.
\end{align}
Then, we have the following properties:
\begin{align}\label{eta_1_and_zeta_1}
    \eta_{1}(v) = \ln(2) + \zeta_{1}(v) \quad \text{and} \quad \zeta'_{m}(v) \leq 0 \quad \forall m\geq 0,
\end{align}
where $\eta_{m}(v)$ is defined in (\ref{definition_eta_m}).
\end{lemma}
\begin{proof}
The first property in (\ref{eta_1_and_zeta_1}) follows from a simple calculation. The second property in (\ref{eta_1_and_zeta_1}) follows from the fact that $\zeta_{m} \geq 0$ and the identity below:
\begin{align*}
    \zeta_{m}'(v) = -\frac{1}{2v} \zeta_{m-1}(v) + \int_{v}^{1} \frac{-1}{(v+s)^{2}} \zeta_{m-1}(s) ds.
\end{align*}
\end{proof}
With the above notation, the main estimate to prove Lemma \ref{lemma_estimate_Lm} is the following estimate for $\eta_{m}(v)$.
\begin{lemma}\label{lemma_estimate_eta_m}
Recall that $\eta_{m}(v)$ and $\zeta_{m}(v)$ are defined in (\ref{definition_eta_m}) and (\ref{definition_zeta}), respectively. Then, the following estimate holds:
\begin{align}\label{estimate_lemma_estimate_eta_m}
    \eta_{m}(v) \geq \sum_{k = 0}^{m} \binom{m}{k} \zeta_{k}(v) \cdot (\ln 2)^{m-k} \quad \forall m\geq 1, \; v \in [0,1].
\end{align}
\end{lemma}
\begin{proof}
We prove (\ref{estimate_lemma_estimate_eta_m}) by mathematical induction. Note that the case of $m = 1$ follows from (\ref{eta_1_and_zeta_1}). Assuming that $\eta_{m}(v)$ satisfy the estimate (\ref{estimate_lemma_estimate_eta_m}), we now prove the estimate for $\eta_{m+1}(v)$. Due to (\ref{definition_eta_m}) and the assumption of the mathematical induction, we know that
\begin{align} \label{estimate_lemma_estimate_eta_m_1}
    &\eta_{m+1}(v) \geq \sum_{k = 0}^{m} \binom{m}{k} \int_{0}^{1}\zeta_{k}(s) ds \cdot (\ln 2)^{m-k}\nonumber\\
    &= \sum_{k = 0}^{m} \binom{m}{k} \int_{0}^{v}\frac{1}{v+s} \zeta_{k}(s) ds \cdot (\ln 2)^{m-k} + \sum_{k = 0}^{m} \binom{m}{k} \int_{v}^{1}\frac{1}{v+s} \zeta_{k}(s) ds \cdot (\ln 2)^{m-k}\nonumber\\
    &\geq \sum_{k = 0}^{m} \binom{m}{k} \zeta_{k}(v)\ln(2) \cdot (\ln 2)^{m-k} + \sum_{k = 0}^{m} \binom{m}{k} \zeta_{k+1}(v) \cdot (\ln 2)^{m-k},
\end{align}
where we have used the second property in (\ref{eta_1_and_zeta_1}) and (\ref{definition_zeta}) in the last inequality. It is clear that
\begin{align*}
    \text{(R.H.S) of (\ref{estimate_lemma_estimate_eta_m_1})}
    &= \sum_{k = 0}^{m} \binom{m}{k} \zeta_{k}(v) \cdot (\ln 2)^{m+1-k} + \sum_{k = 1}^{m+1} \binom{m}{k-1} \zeta_{k}(v) \cdot (\ln 2)^{m+1-k}\\
    &= \sum_{k = 0}^{m+1} \binom{m+1}{k} \zeta_{k}(v) \cdot (\ln 2)^{m+1-k}.
\end{align*}
Then, the proof of Lemma \ref{lemma_estimate_eta_m} is complete.
\end{proof}
Now, we are ready to prove Lemma \ref{lemma_estimate_Lm}. 
\begin{proof} [Proof of Lemma \ref{lemma_estimate_Lm}]
Let $m\geq 1$. Due to Lemma \ref{lemma_estimate_eta_m}, we have the following estimate: 
\begin{align} \label{lower_bound_Lm}
    L_{m+1} = \int_{0}^{1}\eta_{m}(v)dv \geq \sum_{k = 0}^{m} \binom{m}{k} \int_{0}^{1}\zeta_{k}(v) dv \cdot (\ln 2)^{m-k}.
\end{align}
Since $\int_{0}^{1} dv \zeta_{0}(v) = 1$ and 
\begin{align*}
    \int_{0}^{1}\zeta_{k}(v) dv =
    \int_{0}^{1} ds \zeta_{k-1}(s) \int_{0}^{s} dv\frac{1}{v+s} = \ln 2 \cdot \int_{0}^{1} ds \zeta_{k-1}(s) \quad \forall k\geq 1,
\end{align*}
it follows that
\begin{align} \label{integral_incomplete_eta}
    \int_{0}^{1}\zeta_{k}(v) dv = (\ln 2)^{k} \quad \forall k\geq 0.
\end{align}
Consequently, combining (\ref{lower_bound_Lm}) and (\ref{integral_incomplete_eta}) gives $L_{m+1} \geq (2 \ln 2)^{m}$. Therefore, the proof of Lemma \ref{lemma_estimate_Lm} is complete.
\end{proof}

\paragraph{Step 3.} Finally, we will now conclude Lemma \ref{lemma_simplex_integral_lower_bound} by applying the above results.
\begin{proof} [Proof of Lemma \ref{lemma_simplex_integral_lower_bound}]
In view of Lemma \ref{lemma_goal_1_lemma_simplex_integral_lower_bound} and Lemma \ref{lemma_estimate_Lm}, it holds that
\begin{align*}
    \text{(R.H.S) of (\ref{estimate_simplex_integral_lower_bound})}\geq C\biggbkt{\bkt{\frac{4}{3}}^{1/2} \frac{2}{\pi} \cdot 0.99 \cdot 2\ln (2)}^{m} \frac{1}{m(m-1)},
\end{align*}
where $C$ is a positive constant that only depends on $t$. Then, Lemma \ref{lemma_simplex_integral_lower_bound} follows from a simple computation $(4/3)^{1/2} \cdot 4\ln(2) \cdot 0.99/\pi \geq 1.008$.
\end{proof}

\section{Proof of Proposition \ref{Main_result_5}} \label{section_proof_Main_result_5}
The prime objective of this section is to prove Proposition \ref{Main_result_5}. As already explained in Section \ref{section_heuristics_main_result_5}, it is enough to prove Lemma \ref{lemma_positivity_beta_L2} and Lemma \ref{lemma_continuity}. 

\paragraph{Step 1.} We start with the latter lemma.
\begin{proof} [Proof of Lemma \ref{lemma_continuity}]
The continuity of (\ref{GS_map}) is proved by making use of the convexity of the mapping $A_{N}(\lambda): \mathbb{R}_{+} \mapsto \mathbb{R}$ given by
\begin{align*}
    &A_{N}(\lambda) := \sup \mathbf{H}^{\sqrt{\lambda},N} = \\
    &\sup_{\varphi \in D(\mathbf{H}^{\sqrt{\lambda},N}),\; ||\varphi||_{L^{2}(\mathbb{R}^{N \times d})} = 1} -\frac{1}{2}\int_{\mathbb{R}^{N \times d}} d\Vec{x} |\nabla \varphi(\Vec{x})|^{2} + \lambda \sum_{1\leq i<j \leq N} \int_{\mathbb{R}^{N \times d}} d\Vec{x} |\varphi(\Vec{x})|^{2} 
    \cdot R(\Vec{x}(j) - \Vec{x}(i))
\end{align*}
Indeed, due to the above expression, it is straightforward to show the convexity of $A_{N}$. Hence, to prove the continuity, it remains to prove the finiteness of $A_{N}(\lambda)$. About the upper bound, due to the boundedness of $R$, we know that $A_{N}(\lambda) \lesssim ||R||_{L^{\infty}(\mathbb{R}^{d})} \lambda$. Finally, the lower bound $A_{N}(\lambda) \geq 0$ follows from the positivity of $R$ and the fact that $\sup \mathbf{H}^{\sqrt{\lambda},N} \geq \sup \mathbf{H}^{0,N} = 0$. As a result, the proof of Lemma \ref{lemma_continuity} is complete.
\end{proof}

\paragraph{Step 2.} To prove Lemma \ref{lemma_positivity_beta_L2}, the key is to apply \cite[Theorem 5.1]{TAMURA1991433} to the Hamiltonian $-\mathbf{H}^{\beta_{L^{2}},N}$.

\begin{proof} [Proof of Lemma \ref{lemma_positivity_beta_L2}]
Recall that $d = 3$ and $N = 3$. To prove Lemma \ref{lemma_positivity_beta_L2}, it is enough to show that $\inf \mathbf{A} \geq 0$ and $\mathbf{A}$ has a zero resonance energy (see \cite[Section 3]{TAMURA1991433}), where $\mathbf{A}$ is the following Hamiltonian:
\begin{align*}
    \mathbf{A} := -\Delta_{x} - \beta_{L^{2}}^{2}R(x).
\end{align*}
Indeed, since the above statements corresponds to the conditions (A.2) and (A.3) in \cite[Theorem 5.1]{TAMURA1991433} with $m_{j} = 1$, respectively, $\inf -\mathbf{H}^{\beta_{L^{2}},N} < 0$ follows from \cite[Theorem 5.1]{TAMURA1991433}, where we have used the fact that the ground state energy is invariant under any change of basis. This is equivalent to (\ref{positivity_beta_L2}).

We will now conclude the above two conditions. First, notice that $\langle h,\mathbf{A}h \rangle_{L^{2}(\mathbb{R}^{d})} = -2^{d/2}\cdot \langle \widetilde{h},\mathbf{L}_{\beta}\widetilde{h} \rangle_{L^{2}(\mathbb{R}^{d})}$, where $h \in C_{c}^{\infty}(\mathbb{R}^{d})$, $\widetilde{h}(x) := h(\sqrt{2}x)$, and $\mathbf{L}_{\beta}$ is defined in (\ref{definition_hamiltonian_L}). Then, $\inf \mathbf{A} = 0$ follows from the fact that $\sup \mathbf{L}_{\beta} = 0$ at $\beta = \beta_{L^{2}}$. Indeed, since one can show that $\beta \mapsto \sup \mathbf{L}^{\beta}$ is continuous by using the proof of Lemma \ref{lemma_continuity}, the above fact follows from (\ref{identity_beta2+_=_betaL2}).

Now, it remains to show there exists a zero resonance, where its definition is given in \cite[Section 3]{TAMURA1991433}. The proof of this property requires some properties of the Birman-Schwinger operator $\mathbf{T}^{\lambda}_{\beta}$ in Lemma \ref{lemma_Hilbert_Schmidt_operator}. For the sake of simplicity, we postpone the proof of these properties until Section \ref{section_two_particles_system}. In the following, we will explain the reason why $u(x) := h_{\mathbf{v}_{1}}(x/\sqrt{2})$ is a zero resonance of $\mathbf{A}$, where $h_{\mathbf{v}_{1}}(x)$ is defined in (\ref{definition_h_again}). By the part (ii) of Lemma \ref{lemma_Hilbert_Schmidt_operator}, we know that $u \not \in L^{2}(\mathbb{R}^{d})$, $\mathbf{A}u = 0$ in the sense of distribution, and  $u \in L^{2}_{-s}(\mathbb{R}^{d})$ for any $s>-1/2$. In particular, by (\ref{equation_zero_ro reson}) and the fact that $\mathcal{G}^{0}(x) = 1/(2\pi |x|)$ if $d = 3$, where $\mathcal{G}^{0}(x)$ is defined in (\ref{definition_V_and_Yukawa}), we know that
\begin{align} 
    u(x)  = \frac{1}{2^{3/2}}\int_{\mathbb{R}^{d}} \mathcal{G}^{0}(x/\sqrt{2}-y/\sqrt{2}) V_{\beta}(y/\sqrt{2}) h_{\mathbf{v}_{1}}(y/\sqrt{2}) dy= \int_{\mathbb{R}^{d}} \frac{-1}{4\pi|x-y|} (-\beta_{L^{2}}^{2} R(y)) u(y) dy.
\end{align}
Therefore, thanks to the above properties, $u$ is a zero resonance of $\mathbf{A}$ in the sense of the definition given in \cite[Section 3]{TAMURA1991433}. Then, the proof is complete.
\end{proof}

\paragraph{Alternative approach.} Before we close this section, we provide another way to conclude Lemma \ref{lemma_positivity_beta_L2}. In the following proof, one can see an explicit relation between (\ref{positivity_beta_L2}) and the higher moments of the mollified SHE.

\begin{proof} [Another proof of Lemma \ref{lemma_positivity_beta_L2}]
First, notice that $\mathbf{H}^{\beta,N}$ can be connected with $\mathcal{H}_{\varepsilon}^{\beta_{\varepsilon},N}$ by using a change of variable:
\begin{align} \label{sup_H_and_H_epsilon}
    \varepsilon^{-2} \cdot \sup \mathbf{H}^{\beta,N} = \sup \mathcal{H}_{\varepsilon}^{\beta_{\varepsilon},N} \quad \forall \beta>0,
\end{align}
where $\mathcal{H}_{\varepsilon}^{\beta_{\varepsilon},N}$ is defined in (\ref{definition_hamiltonian_epsilon}), and it is self-adjoint thanks to Lemma \ref{lemma_self_adjoint}. Indeed, if $\varphi \in C_{c}^{\infty}(\mathbb{R}^{N \times d})$, then one has the following relation:
\begin{align*}
    \varepsilon^{-2}\cdot \Bigl\langle \varphi^{(\varepsilon)},\mathbf{H}^{\beta,N}\varphi^{(\varepsilon)} \Bigr\rangle_{L^{2}(\mathbb{R}^{N \times d})} = \Bigl\langle \varphi,\mathcal{H}_{\varepsilon}^{\beta_{\varepsilon},N}\varphi \Bigr\rangle_{L^{2}(\mathbb{R}^{N \times d})}, \quad \text{where}\quad \varphi^{(\varepsilon)}(\Vec{x}) := \varepsilon^{\frac{N\cdot d}{2}} \varphi(\varepsilon \cdot \Vec{x}).
\end{align*}
Consequently, thanks to (\ref{sup_H_and_H_epsilon}), (\ref{positivity_beta_L2}) is equivalent to the following limit:
\begin{align} \label{limit_sup_H_epsilon}
    \lim_{\varepsilon \to 0^{+}}\sup \mathcal{H}_{\varepsilon}^{\beta_{\varepsilon},N} = \infty \quad \text{when} \quad \beta = \beta_{L^{2}}.
\end{align}
Next, observe that the instability (\ref{limit_sup_H_epsilon}) follows from the divergence of the three moment at $L^{2}$-criticality:
\begin{align} \label{divergence_at_beta_L2}
    \lim_{\varepsilon \to 0^{+}}\Bigl\langle \varphi, \mathcal{Q}_{\varepsilon;t}^{\beta,N}U_{0}^{\otimes N} \Bigr\rangle_{L^{2}(\mathbb{R}^{N \times d})} = \infty \quad \text{when} \quad \beta = \beta_{L^{2}},
\end{align}
where we have set $\varphi(\Vec{x}) = G^{(N)}_{\nu}(\Vec{x})$ and $U_{0}(x) = G_{\nu}(x)$. Indeed, in view of (\ref{upperbound_for_semigroup}) and (\ref{moment_and_semigroup}), one can conclude the following estimate:
\begin{align} \label{estimate_exp_sup}
    \exp\biggbkt{\sup \mathcal{H}_{\varepsilon}^{\beta_{\varepsilon},N}} \cdot ||G_{\nu}^{(N)}||_{L^{2}(\mathbb{R}^{N \times d})}^{2} \geq \Bigl\langle \varphi, \mathcal{Q}_{\varepsilon;t}^{\beta,N}U_{0}^{\otimes N} \Bigr\rangle_{L^{2}(\mathbb{R}^{N \times d})}.
\end{align}
Then, combining (\ref{estimate_exp_sup}) and (\ref{divergence_at_beta_L2}) implies (\ref{limit_sup_H_epsilon}). 

With the above explanation, to prove Lemma \ref{lemma_positivity_beta_L2}, it remains to show (\ref{divergence_at_beta_L2}), which is the core of our proof of (\ref{positivity_beta_L2}). To this aim, we need to establish a proper lower bound for the three moment of the mollified SHE. The following trivial Gaussian lower bound is helpless to conclude (\ref{divergence_at_beta_L2}):
\begin{align*}
    \Bigl\langle\varphi, \mathcal{Q}_{\varepsilon;t}^{\beta,N}U_{0}^{\otimes N} \Bigr\rangle_{L^{2}(\mathbb{R}^{N \times d})} \geq \Bigl\langle\varphi, G^{(N)}U_{0}^{\otimes N} \Bigr\rangle_{L^{2}(\mathbb{R}^{N \times d})},
\end{align*}
However, the condition $\beta = \beta_{L^{2}}$ allows us to conclude the following non-Gaussian lower bound for the limiting three moment of the mollified SHE.
\begin{proposition} \label{proposition_fatou}
Assume that $d = 3$, $\beta = \beta_{L^{2}}$, and $N = 3$.  Recall that the $N$-th moment $\mathcal{Q}_{\varepsilon;t}^{\beta,N}U_{0}^{\otimes N}(\Vec{x})$ is given in (\ref{definition_moment}), $U_{0}$ is the initial datum of the mollified SHE, $\mathscr{Q}_{t}^{N}U_{0}^{\otimes N}(\Vec{x}_{0})$ is the sub-limiting $N$-th moment defined in (\ref{definition_sub_limit}). Let $U_{0}(x) = G_{\nu}(x)$. Then, given an arbitrary $t>0$, for every $\varphi \in \mathscr{S}(\mathbb{R}^{N \times d};\mathbb{R}_{+})$, the following estimate holds:
\begin{align} \label{non_Gaussian_lower_bound}
    \liminf_{\varepsilon \to 0^{+}} 
    \Bigl\langle \varphi, \mathcal{Q}_{\varepsilon;t}^{\beta,N}U_{0}^{\otimes N} \Bigr\rangle_{L^{2}(\mathbb{R}^{N \times d})}
    \geq \Bigl\langle \varphi,\mathscr{Q}_{t}^{N}U_{0}^{\otimes N} \Bigr\rangle_{L^{2}(\mathbb{R}^{N \times d})}.
\end{align}
\end{proposition}
\begin{proof}[Sketch of the proof]
Here, we only elaborate on the ideas to prove Proposition \ref{proposition_fatou} 
while skipping the details. Simply put, the proof follows from a similar line of thought as mentioned in the third part of Section \ref{section_definition_sub-limiting_moment} (the formal derivation of the sub-limiting higher moment). 

To begin with, by applying Fatou's lemma to put the limit on the left-hand side of (\ref{non_Gaussian_lower_bound}) inside the summation as (\ref{inside_limit}), it remains to show the following lower bound to conclude (\ref{non_Gaussian_lower_bound}):
\begin{align} \label{lower_bound_nonGaussian_path_integral}
    \liminf_{\varepsilon\to 0^{+}}\Bigr\langle \varphi,\mathscr{I}_{\varepsilon;t}^{N;(\ell_{1},\ell'_{1},\mathfrak{i}_{1}),...,(\ell_{m},\ell'_{m},\mathfrak{i}_{m})} U_{0}^{\otimes N} \Bigr\rangle_{L^{2}(\mathbb{R}^{N \times d})}
    \geq \Bigr\langle \varphi,\mathscr{I}_{0;t}^{N;(\ell_{1},\ell'_{1},\mathfrak{i}_{1}),...,(\ell_{m},\ell'_{m},\mathfrak{i}_{m})} U_{0}^{\otimes N}  \Bigr\rangle_{L^{2}(\mathbb{R}^{N \times d})},
\end{align}
where $\mathfrak{i}_{j} = 1$ for all $1\leq j\leq m$. Indeed, since we can bound the left-hand side of (\ref{lower_bound_nonGaussian_path_integral}) from below by zero for any other case in (\ref{inside_limit}) (see also the description of $\mathscr{Q}_{t}^{N}U_{0}^{\otimes N}(\Vec{x}_{0})$ under Remark \ref{remark_after_definition}), (\ref{non_Gaussian_lower_bound}) follows from (\ref{lower_bound_nonGaussian_path_integral}) and (\ref{definition_sub_limit}). Moreover, by apply Fatou's lemma again to put the above limit inside all the spatial integrals on the left-hand side of (\ref{lower_bound_nonGaussian_path_integral}), the integrands of these spatial integrals have the following form (see the right-hand side of (\ref{identity_change_of_variable_formal_limit_1})):
\begin{align} \label{lower_bound_nonGaussian_time_integral}
    \liminf_{\varepsilon \to 0^{+}} \int_{v_{j},r_{j}>0, \;v_{1}+r_{2}+v_{2}+...+v_{m-1}+r_{m}+v_{m} <t}\biggbkt{\prod_{j=1}^{m} dv_{j} \prod_{j=2}^{m} \mathbf{u}^{y_{j}^{(\ell_{j},\ell'_{j})},\Vec{x}_{j}^{(\ell_{j},\ell'_{j})}}_{\varepsilon;t}(dr_{j})} \cdot F_{\varepsilon}(v_{1},r_{2}...,r_{m},v_{m}),
\end{align}
where we have used a change of variable $v_{j} = s_{j}-u_{j-1}$ and $r_{j+1} = u_{j+1} - s_{j+1}$. Here, $F_{\varepsilon}(v_{1},r_{2}...,r_{m},v_{m})$ is a product of Gaussian kernels, which converges pointwise, and $\mathbf{u}^{z,z'}_{\varepsilon;t}(dr_{j})$ is defined in (\ref{measure_epsilon_z_z'}). In this way,
the key observation here is to re-express 
$dv_{j}$ and
$\mathbf{u}^{z,z'}_{\varepsilon;t}(dr_{j})$ as $t \cdot dv_{j}/t$ and $\mathbf{u}^{z,z'}_{\varepsilon;t}([0,t]) \cdot \overline{\mathbf{u}}^{z,z'}_{\varepsilon;t}(dr_{j})$, respectively, where $\overline{\mathbf{u}}^{z,z'}_{\varepsilon;t}(dr_{j})$ is defined in (\ref{normalized_measure_epsilon_z_z'}). Note that the weak convergence of $\overline{\mathbf{u}}^{z,z'}_{\varepsilon;t}(dr_{j})$ and the convergence of the total mass $\mathbf{u}^{z,z'}_{\varepsilon;t}([0,t])$ will be proved in Proposition \ref{proposition_critical_second_moment_point_to_point}. Then, the limit in (\ref{lower_bound_nonGaussian_time_integral}) can be put inside the temporal integral by using Fatou’s lemma and Skorokhod’s representation theorem. Finally, by combining the identity (\ref{identity_one}) and an argument that is analogous to the ideas in the third part of Section \ref{section_definition_sub-limiting_moment} to simplify the limiting integrand, we conclude the right-hand side of (\ref{lower_bound_nonGaussian_path_integral}).

%(\ref{limit_two_particle}) the third part of Section \ref{section_definition_sub-limiting_moment}:normalized_measure_epsilon_z_z'
\end{proof}

Finally, (\ref{divergence_at_beta_L2}) follows from Proposition \ref{proposition_fatou} and Theorem \ref{Main_result_2}. Then, the proof of (\ref{positivity_beta_L2}) is complete.    
\end{proof}

\section{Proof of Theorem \ref{main_result_gamma}} \label{section_corollary_main_result_gamma}
This section is dedicated to the proof of Theorem \ref{main_result_gamma}. We begin by showing (\ref{estimate_gamma}) by (\ref{estimate_fractional_beta}). 
\begin{proof} [Proof of (\ref{estimate_gamma})]
Fix $\beta < \beta_{L^{2}}$ and a small positive number $\varepsilon > 0$. Then, one has the following estimates for $\beta$:
\begin{align} \label{estimate_gamma2_beta_gamma1}
    \frac{\alpha_{\infty,+}}{\sqrt{\gamma_{\varepsilon;2}(\beta)-1}} < \beta < \frac{\beta_{L^{2}}}{\sqrt{\gamma_{\varepsilon;1}(\beta)-1}},
\end{align}
where $\gamma_{\varepsilon;1}(\beta) := 1+\sbkt{\frac{\beta_{L^{2}}}{\beta}}^{2} -\varepsilon$ and $\gamma_{\varepsilon;2}(\beta):= 1+\sbkt{\frac{\alpha_{\infty,+}}{\beta}}^{2}+\varepsilon$. Let $\varepsilon$ be small enough such that $\gamma_{\varepsilon;1} > 2$. Also, since $\beta < \beta_{L^{2}} \leq \alpha_{\infty,+}$ follows from Theorem \ref{Main_result_3}, we know that $\gamma_{\varepsilon;2} > 2$. Consequently, the \emph{strictly inequality} below follows from (\ref{estimate_fractional_beta}) and the estimate (\ref{estimate_gamma2_beta_gamma1}):
\begin{align} \label{estimate_gamma2_beta_gamma1_2}
    \beta_{L^{\gamma_{\varepsilon;2}(\beta)}} < \beta  < \beta_{L^{\gamma_{\varepsilon;1}(\beta)}}.
\end{align}

To prove (\ref{estimate_gamma}), it is enough to prove the following estimate for the critical exponent since taking $\varepsilon \to 0^{+}$ yields (\ref{estimate_gamma}):
\begin{align} \label{inequality_epsilon_12}
    \gamma_{\varepsilon;1}(\beta)\leq \gamma^{*}(\beta) \leq \gamma_{\varepsilon;2}(\beta).
\end{align}
Now, we show that (\ref{inequality_epsilon_12}) follows from (\ref{estimate_gamma2_beta_gamma1_2}). To be more specific, on the one hand, the upper bound in (\ref{inequality_epsilon_12}) follows from (\ref{identity_beta_L_N_2}) and the lower bound in (\ref{estimate_gamma2_beta_gamma1_2}). On the other hand, the lower bound in (\ref{inequality_epsilon_12}) can concluded by using (\ref{identity_beta_L_N_2}) and the upper bound in (\ref{estimate_gamma2_beta_gamma1_2}), where we have used the fact that for each real number $\gamma > 1$, the following holds: 
\begin{align} \label{uniform_moment_bound}
    \sup_{T>0}\mathbb{E}[(\mathcal{Z}^{\alpha}_{T})^{\gamma}] < \infty \quad \forall \alpha < \beta_{L^{\gamma}}.
\end{align}
Here, we clarify the reason why we have the above property. For the case of $\gamma \in \mathbb{N}$, (\ref{uniform_moment_bound}) follows from (\ref{representation_FK_0}). But, for the general case, the following property is required to prove (\ref{uniform_moment_bound}):
\begin{align*}
    \alpha \in (0,\infty) \mapsto \sup_{T>0}\mathbb{E}[(\mathcal{Z}^{\alpha}_{T})^{\gamma}] \text{ is increasing}.
\end{align*}
Here, the above property can be proved by using the following estimate:
\begin{align*}
    ||\mathcal{Z}^{\alpha}_{T}(0;\xi)||_{L^{\gamma}} = ||\Gamma(A) \mathcal{Z}^{\alpha'}_{T}(0;\xi)||_{L^{\gamma}} \leq ||\mathcal{Z}^{\alpha'}_{T}(0;\xi)||_{L^{\gamma}} \quad \forall \alpha < \alpha', \; \gamma>1,
\end{align*}
where we have applied \cite[Theorem 5.1]{GHS} to (\ref{chaos}) with $A := \alpha/\alpha' \cdot I$ and $p = q := \gamma$. 
\end{proof}

Let us now prove (\ref{estimate_fractional_beta}). As the proof ideas mentioned under Theorem \ref{main_result_gamma}, (\ref{estimate_fractional_beta}) can be proved by scaling the function $\phi$ in (\ref{definition_partition_function_1}).
\begin{proof} [Proof of (\ref{estimate_fractional_beta})]
Let $N\leq \gamma \leq N+1$ for some $N\geq 2$. Recall the third part of Remark \ref{total_mass} and (\ref{definition_partition_function_1}). Simply put, the main idea to prove (\ref{estimate_fractional_beta}) is scaling the function $\phi$ in order to match $\beta_{L^{\gamma}}$ with $\widetilde{\beta}_{L^{N+1}}$, where $\widetilde{\beta}_{L^{N+1}}$ is a critical coupling constant for the scaled partition function $\mathcal{W}_{T}^{\beta}(0;\xi)$. Since we have an estimate for $\widetilde{\beta}_{L^{N+1}}$ due to Theorem \ref{Main_result_3}, we can use this estimate to bound $\beta_{L^{\gamma}}$.

To prove (\ref{estimate_fractional_beta}), we first note that the lower bound in (\ref{estimate_fractional_beta}) has been proved in Lemma \ref{lemma_HC_bound}. Hence, it remains to prove the upper bound in (\ref{estimate_fractional_beta}). As the proof of Lemma \ref{lemma_HC_bound}, by applying the hypercontractivity from \cite[Theorem 5.1]{GHS} to the partition function $\mathcal{Z}_{T}^{\beta}(0;\xi)$, where $A := (\frac{\gamma-1}{N})^{1/2}I$, $p := \gamma$, and $q := N+1$, one has the following estimate:
\begin{align}  \label{estimate_N+1_gamma}
    \Bigl|\Bigl|\mathcal{Z}^{\beta \cdot \sqrt{(\gamma-1)/N}}_{T}(0;\xi)\Bigr|\Bigr|_{L^{N+1}} \leq ||\mathcal{Z}^{\beta}_{T}(0;\xi)||_{L^{\gamma}} \quad \forall T \geq 0, \; \beta>0.
\end{align}
Now, we rewrite the partition function $\mathcal{Z}^{\beta \cdot \sqrt{(\gamma-1)/N}}_{T}(0;\xi)$ by the following partition function $\mathcal{W}^{\beta}_{T}(x;\xi)$, which has the same form as (\ref{definition_partition_function_1}) but the function $\phi$ is replaced by $\widetilde{\phi}(x)$:
\begin{align*}
    \mathcal{W}^{\beta}_{T}(x;\xi)
    := \mathbf{E}^{B}_{x}\bigbkt{\exp\bkt{\beta \int_{0}^{T} ds \int_{\mathbb{R}^{d}} dy\widetilde{\phi}(B(s) - y) \xi(y,s) - \frac{\beta^{2} \widetilde{R}(0) T }{2}}},
\end{align*}
\begin{align*}
    \widetilde{\phi}(x) := \biggbkt{\frac{\gamma-1}{N}}^{1/2}\phi(x), \quad \text{and} \quad \widetilde{R}(x) := \widetilde{\phi}*\widetilde{\phi}(x).
\end{align*}
Then, the estimate (\ref{estimate_N+1_gamma}) shows that if $\beta < \beta_{L^{\gamma}}$, then $\beta<\widetilde{\beta}_{L^{N+1}}$. Here, $\widetilde{\beta}_{L^{N+1}}$ defined as (\ref{identity_beta_L_N_2}) denotes the  critical coupling constant of $(N+1)$-th moment for the partition function $\mathcal{W}_{T}^{\beta}(0;\xi)$. Consequently, one has the following estimate:
\begin{align} \label{estimate_gamma_widetilde_N+1}
    \beta_{L^{\gamma}} \leq \widetilde{\beta}_{L^{N+1}}.
\end{align}
Consequently, recalling the third part in Remark \ref{total_mass}, the following estimate follows from Theorem \ref{Main_result_3}:
\begin{align} \label{representation_widetilde_beta_L_N+1}
    \widetilde{\beta}_{L^{N+1}} \leq \frac{\widetilde{\alpha}_{\infty,+}}{\sqrt{N}}.
\end{align}
Here, $\widetilde{\alpha}_{\infty,+}:=\lim_{M\to\infty}\widetilde{\alpha}_{M,+}$, where $\widetilde{\alpha}_{M,+}$ is defined in (\ref{definition_alpha_N}), but the function $R$ is replaced by $\widetilde{R}$. Then, thanks to (\ref{definition_alpha_N}) and (\ref{definition_alpha_infinity}), the following relation holds:
\begin{align} \label{representation_widetilde_beta_L_2}
    \widetilde{\alpha}_{\infty,+} = \biggbkt{\frac{N}{\gamma-1}}^{1/2} \alpha_{\infty,+}.
\end{align}
Therefore, by combining (\ref{estimate_gamma_widetilde_N+1}), (\ref{representation_widetilde_beta_L_N+1}), and (\ref{representation_widetilde_beta_L_2}), the upper bound in (\ref{estimate_fractional_beta}) is proved. Then, the proof is complete.
    
\end{proof}

\appendix

\section{Self-adjointness and semigroup} \label{section_op}
This section aims to present some properties of the unbounded operators used in the previous proofs. Before doing so, we illustrate some notations. For a self-adjoint operator $\mathbf{A}$ defined on $D(\mathbf{A}) \subseteq H$, where $H$ is a Hilbert space, we define
\begin{align} \label{definition_total_energy}
    \sup \mathbf{A} := \sup_{u \in D(\mathbf{A}), \; ||u||_{H} = 1} \langle u,\mathbf{A} u \rangle_{H} \quad \text{and} \quad \inf \mathbf{A} := \inf_{u \in D(\mathbf{A}), \; ||u||_{H} = 1} \langle u,\mathbf{A} u \rangle_{H}
\end{align}
Throughout this article, we use $\mathbf{P}_{\mathbf{A}}(d\lambda)$ to denote the projection-valued measure defined in \cite[Theorem VIII.6]{reed1981functional}. 
Then, as \cite[(VIII.5)]{reed1981functional}, we define the semigroup of $\mathbf{A}$ as follows:
\begin{align} \label{definition_semigroup}
    \exp(T\cdot \mathbf{A}) := \int_{\mathbb{R}} \exp(T \cdot \lambda) \mathbf{P}_{\mathbf{A}}(d\lambda), 
\end{align}
where the domain is defined by
\begin{align*}
    D(\exp(T\cdot \mathbf{A})) := \biggl\{u\in H: \int_{\mathbb{R}}\exp(T\cdot \lambda) \langle u,\mathbf{P}_{\mathbf{A}}(d\lambda)u \rangle_{H} <\infty\biggr\}.
\end{align*}

Recall that the unbounded operator $\mathcal{H}_{\varepsilon}^{\beta_{\varepsilon},N}$ is given in (\ref{definition_hamiltonian_epsilon}) with $D(\mathcal{H}_{\varepsilon}^{\beta_{\varepsilon},N}) := C_{c}^{\infty}(\mathbb{R}^{N \times d})$, and $\mathbf{H}^{\beta,N}$ and $\mathbf{L}^{\beta}$ are defined in  (\ref{short_range_hamiltonian}) and (\ref{definition_hamiltonian_L}).
Here, the corresponding Hilbert spaces are $L^{2}(\mathbb{R}^{N \times d})$ for the first two operators and $L^{2}(\mathbb{R}^{d})$ for the last opertaor, respectively. Now, we present some properties of the above operators. 

\begin{lemma} \label{lemma_self_adjoint}
Assume that $d\geq 3$, $N\geq 2$, $\varepsilon>0$, and $\beta>0$. Then, unbounded operators $\mathcal{H}_{\varepsilon}^{\beta_{\varepsilon},N}$, $\mathbf{H}^{\beta,N}$, and $\mathbf{L}^{\beta}$ defined on the following domains are self-adjoint:
\begin{align} \label{domain_laplace_all}
    D(\mathcal{H}_{\varepsilon}^{\beta_{\varepsilon},N}) = D(\mathbf{H}^{\beta,N}) := \set{f\in L^{2}(\mathbb{R}^{N \times d}): |k|^{2}\widehat{f}(k) \in L^{2}(\mathbb{R}^{N \times d})}
\end{align}
and
\begin{align} \label{domain_laplace_one_particle}
    D(\mathbf{L}^{\beta}) := \set{h\in L^{2}(\mathbb{R}^{d}): |k|^{2}\widehat{h}(k) \in L^{2}(\mathbb{R}^{d})},
\end{align}
where Fourier transform defined in (\ref{FT}).
\end{lemma}
\begin{proof}
We first note that $\Delta_{\Vec{x}}$ and $\Delta_{x}$ are self-adjoint operators on the domains on the right-hand sides of (\ref{domain_laplace_all}) and (\ref{domain_laplace_one_particle}), respectively. See for \cite[Proposition 9.34]{hall2013quantum} more details. Moreover, thanks to the boundedness of $R$, the self-adjointness of all the above operators follow from Kato-Rellich theorem \cite[Theorem X.12]{reed1975ii} with the corresponding relative bound $a = 0$. Here, the relative bound is defined in \cite[(X.19a)]{reed1975ii}.
\end{proof}

\begin{lemma} \label{lemma_semigroup}
Suppose that $d\geq 3$, $N\geq 2$, $\varepsilon>0$, $\beta>0$, and $t,T > 0$. Then, both of $\mathcal{H}_{\varepsilon}^{\beta_{\varepsilon},N}$ and $\mathbf{H}^{\beta,N}$ are bounded from above. In particular, their semigroup $\exp(T\cdot\mathbf{H}^{\beta,N})$ and $\exp(t\cdot\mathcal{H}_{\varepsilon}^{\beta,N})$ are bounded operators on $L^{2}(\mathbb{R}^{N \times d})$ such that
\begin{align} \label{upperbound_for_semigroup}
    ||\exp(\mathbf{t}\cdot \mathbf{A})||_{L^{2}(\mathbb{R}^{N \times d})} \leq \exp\sbkt{\mathbf{t} \cdot \sup \mathbf{A}} \quad\forall \mathbf{t}>0.
\end{align}
\end{lemma}
\begin{proof}
Lemma \ref{lemma_semigroup} follows directly by using the fact that
\begin{align*}
    \mu_{f}((\sup \mathbf{A},\infty)) = 0, \quad \text{where}\quad \mu_{f}(d\lambda) := \langle f,\mathbf{P}_{\mathbf{A}}f\rangle_{L^{2}(\mathbb{R}^{N \times d})}.
\end{align*}
\end{proof}

\section{Asymptotics of the $\varepsilon$-range two-body interactions.} \label{section_two_particles_system}
The prime objective of this subsection is to prove a convergence analogous to (\ref{limit_two_particle}). The main idea is the following asymptotic of the Laplace transform of the rescaled Brownian exponential functional, which is an application of \cite[Lemma 1.2.4]{albeverio2012solvable}. 
\begin{lemma} \label{lemma_critical_second_moment_point_to_point}
Assume that $d = 3$ and $\beta = \beta_{L^{2}}$. Recall that $\mathfrak{V}_{\beta}(x)$ and $\mathcal{G}^{0}(x)$ are defined in (\ref{definition_V_and_Yukawa}), and $\mathbf{v}_{1}(x)$ is defined in Lemma \ref{lemma_Hilbert_Schmidt_operator}. Then, we have the following properties. 
\begin{enumerate} [label=(\roman*)]
    \item Given an arbitrary $\Lambda > 0$, for each $z\neq z'$, it holds that
    \begin{align} \label{part_i_lemma_critical_second_moment_point_to_point}
        \int_{0}^{\infty} dt \exp(-\Lambda \cdot t)
        \beta_{\varepsilon}^{2}  \mathbf{E}^{B}_{\varepsilon \cdot z}\bigbkt{\exp\bkt{\beta_{\varepsilon}^{2} \int_{0}^{t} R_{\varepsilon}(\sqrt{2}B(r)) dr} &: B(t) = \varepsilon \cdot z'} \beta_{\varepsilon}^{2}\overset{\varepsilon \to 0^{+}}{\longrightarrow} \frac{2\pi}{\sqrt{2\Lambda}} \cdot \mathscr{C}(z,z'),
    \end{align}
    \begin{align} \label{definition_C_z_z'}
            \text{where}\quad\mathscr{C}(z,z') := 
            \biggbkt{\int_{\mathbb{R}^{d}} dx\beta^{2} \mathcal{G}^{0}(z&-x) \mathfrak{V}_{\beta}(x) \mathbf{v}_{1}(x) } \cdot \Bigl|\langle \mathbf{v}_{1}, \mathfrak{V}_{\beta} \rangle_{L^{2}(\mathbb{R}^{d})}\Bigr|^{-2}  \cdot \nonumber\\
            &\biggbkt{\int_{\mathbb{R}^{d}} dx'\mathbf{v}_{1}(x')\mathfrak{V}_{\beta}(x') \mathcal{G}^{0}(x'-z') \beta^{2}} \quad \forall z,z'\in \mathbb{R}^{d}.
    \end{align}
    \item The first eigenvector $\mathbf{v}_{1}(x)$ defined in Lemma \ref{lemma_Hilbert_Schmidt_operator} can be chosen to be a non-negative function in $L^{2}(\mathbb{R}^{d})$.
\end{enumerate}

\end{lemma}
For the sake of simplicity, we postpone the
proof of Lemma \ref{lemma_critical_second_moment_point_to_point} until Section \ref{section_proof_lemma_critical_second_moment_point_to_point}. With the above properties, we now prove (\ref{limit_two_particle}) in the weak sense:
\begin{align} \label{weak_convergence_integral_h}
    &\int_{0}^{\sigma} dt h(t)\int_{\mathbb{R}^{2\times d}} dzdz'  \beta_{\varepsilon}^{2} R(\sqrt{2}z) \mathbf{E}_{\varepsilon \cdot z}^{B}\bigbkt{\exp\bkt{\beta_{\varepsilon}^{2}\int_{0}^{t} R_{\varepsilon}(\sqrt{2}B(v)) dv}: B(t) = \varepsilon \cdot z'} \beta_{\varepsilon}^{2} R(\sqrt{2}z') \nonumber\\
    &\overset{\varepsilon \to 0^{+}}{\longrightarrow}  \int_{0}^{\sigma} dt  h(t) \biggbkt{\frac{2\pi}{t}}^{1/2} \quad \forall h\in C_{b}((0,\sigma);\mathbb{R}).
\end{align}
\begin{proposition} \label{proposition_critical_second_moment_point_to_point}
Assume that $d = 3$ and $\beta = \beta_{L^{2}}$. Then, $\mathscr{C}(z,z') > 0$ for all $z,z'\in \mathbb{R}^{d}$, where $\mathscr{C}(z,z')$ is defined in (\ref{definition_C_z_z'}). Moreover, let $\sigma>0$ and we define $\mathbf{u}^{z,z'}_{\varepsilon;\sigma}(dt)$ by the Borel measure on $(0,\sigma)$ that has the following expression:
\begin{align} \label{measure_epsilon_z_z'}
    \mathbf{u}^{z,z'}_{\varepsilon;\sigma}(dt) := \beta_{\varepsilon}^{2} \cdot \mathbf{E}_{\varepsilon \cdot z}\bigbkt{\exp\bkt{\beta_{\varepsilon}^{2} \int_{0}^{t} R_{\varepsilon}(\sqrt{2}B(r)) dr} : B(t) = \varepsilon \cdot z'} \cdot\beta_{\varepsilon}^{2} dt \quad \text{if } \varepsilon>0;
\end{align}
and
\begin{align}
    \mathbf{u}^{z,z'}_{0;\sigma}(dt) := \mathscr{C}(z,z') \biggbkt{\frac{2\pi}{t}}^{1/2} dt \quad \text{if } \varepsilon = 0,
\end{align}
then one has the following asymptotics:
\begin{align} \label{weak_convergence_probability_measure}
    \overline{\mathbf{u}}^{z,z'}_{\varepsilon;\sigma}(dt) 
    \Longrightarrow \overline{\mathbf{u}}^{z,z'}_{0,\sigma}(dt) 
    \quad \text{and} \quad\mathbf{u}^{z,z'}_{\varepsilon;\sigma}((0,\sigma)) \longrightarrow \mathbf{u}^{z,z'}_{0;\sigma}((0,\sigma))
    \quad \text{as} \quad \varepsilon\to 0^{+},
\end{align}
where $\overline{\mathbf{u}}^{z,z'}_{\varepsilon;\sigma}(dt)$ is a probability measure given by 
\begin{align} \label{normalized_measure_epsilon_z_z'}
    \overline{\mathbf{u}}^{z,z'}_{\varepsilon;\sigma}(dt) := \mathbf{u}^{z,z'}_{\varepsilon;\sigma}(dt)\bigl/\mathbf{u}^{z,z'}_{\varepsilon;\sigma}((0,\sigma)) \quad \forall \varepsilon\geq 0.
\end{align}
\end{proposition}
\begin{remark} \label{remark_explain_delta-Bose_gas}
Our proof of Proposition \ref{proposition_fatou} only needs the following pointwise version of (\ref{weak_convergence_integral_h}), which is an application of (\ref{weak_convergence_probability_measure}):
\begin{align*} 
    \int_{0}^{\sigma} dt h(t)  \beta_{\varepsilon}^{2}  \mathbf{E}_{\varepsilon \cdot z}^{B}\bigbkt{\exp\bkt{\beta_{\varepsilon}^{2}\int_{0}^{t} R_{\varepsilon}(\sqrt{2}B(v)) dv}&: B(t) = \varepsilon \cdot z'} \beta_{\varepsilon}^{2}  \overset{\varepsilon \to 0^{+}}{\longrightarrow}  \int_{0}^{\sigma} dt h(t) \mathscr{C}(z,z') \biggbkt{\frac{2\pi}{t}}^{1/2}.
\end{align*}
But, (\ref{weak_convergence_integral_h}) 
follows from a similar line of thought as the argument to conclude Proposition \ref{proposition_critical_second_moment_point_to_point} and the following identity:
\begin{align} \label{identity_one}
    \int_{\mathbb{R}^{2\times d}} R(\sqrt{2}z) \mathscr{C}(z,z') R(\sqrt{2}z') = \langle \mathfrak{V}_{\beta},\mathbf{v}_{1} \rangle_{L^{2}(\mathbb{R}^{d})} \cdot \Bigl|\langle \mathbf{v}_{1}, \mathfrak{V}_{\beta} \rangle_{L^{2}(\mathbb{R}^{d})}\Bigr|^{-2} \cdot  \langle \mathfrak{V}_{\beta},\mathbf{v}_{1} \rangle_{L^{2}(\mathbb{R}^{d})} = 1.
\end{align}
Here, we have used the facts that
$\mathbf{T}^{0}_{\beta}\mathbf{v}_{1} = \mathbf{v}_{1}$ and that $\mathbf{v}_{1}$ is a real-valued function, which will be proved in Lemma \ref{lemma_Hilbert_Schmidt_operator}.
\end{remark}

\begin{proof}[Proof of Proposition \ref{proposition_critical_second_moment_point_to_point}]
We first note that $\langle \mathbf{v}_{1}, \mathfrak{V}_{\beta} \rangle_{L^{2}(\mathbb{R}^{d})} > 0$ follows from the non-negativity of $\mathbf{v}_{1}$ in Lemma \ref{lemma_critical_second_moment_point_to_point} and the part (ii) in Lemma \ref{lemma_Hilbert_Schmidt_operator}. Consequently, there exists a set of positive measure $A \subseteq \mathbb{R}^{d}$ such that $\mathfrak{V}_{\beta}(x) \mathbf{v}_{1}(x) > 0$ for all $x \in A$. This clearly gives $\mathscr{C}(z,z') > 0$. 

Now, it remains to prove the convergences in (\ref{weak_convergence_probability_measure}). Fix $\sigma > 0$ and $z,z' \in \mathbb{R}^{d}$ such that $z\neq z'$. To prove (\ref{weak_convergence_probability_measure}), we first introduce some auxiliary convergences. Let $\Lambda > 0$. Then, we consider the following Borel measure on $\mathbb{R}_{+}$:
\begin{align}
    \mathbf{m}^{z,z';\Lambda}_{\varepsilon}(dt) := \exp(-\Lambda \cdot t)\cdot\mathbf{u}^{z,z';\Lambda}_{\varepsilon}(dt) \quad \forall \varepsilon\geq0;
\end{align}
Notice that the following convergence of the total mass follows from Lemma \ref{lemma_critical_second_moment_point_to_point}:
\begin{align} \label{weak_convergence_total_mass_1}
    \lim_{\varepsilon\to 0^{+}}\mathbf{m}^{z,z';\Lambda}_{\varepsilon}(\mathbb{R}_{+}) = \mathscr{C}(z,z') \frac{2\pi}{\sqrt{2\Lambda}} = \mathbf{m}^{z,z';\Lambda}_{0}(\mathbb{R}_{+}).
\end{align}
Also, we observe that the following convergence of the moment-generating functions follows from Lemma \ref{lemma_critical_second_moment_point_to_point} and (\ref{weak_convergence_total_mass_1}):
\begin{align*}
    \lim_{\varepsilon\to 0^{+}}\int_{0}^{\infty} \exp(\xi \cdot t) \overline{\mathbf{m}}^{z,z';\Lambda}_{\varepsilon}(dt) = \biggbkt{\frac{\Lambda}{\Lambda-\xi}}^{1/2} = \int_{0}^{\infty} \exp(\xi \cdot t) \overline{\mathbf{m}}^{z,z';\Lambda}_{0}(dt)\quad \forall 0<\xi < \Lambda.
\end{align*}
Then, it holds that
\begin{align} \label{weak_convergence_probability_measure_1}
    \overline{\mathbf{m}}^{z,z';\Lambda}_{\varepsilon}(dt)\Longrightarrow\overline{\mathbf{m}}^{z,z';\Lambda}_{0}(dt) \quad\text{as} \quad\varepsilon \to 0^{+}. 
\end{align}

Now, we prove the convergence of total mass of $\mathbf{u}^{z,z'}_{\varepsilon;\sigma}(dt)$. To this aim, we consider $h(t) := \exp(\Lambda \cdot t) \cdot 1_{(0,\sigma)}(t)$. Then, the total mass of $\mathbf{u}^{z,z'}_{\varepsilon;\sigma}(dt)$ can be connected with $\mathbf{m}^{z,z'}_{\varepsilon;\sigma}(dt)$ as follows:
\begin{align*}
    \mathbf{u}^{z,z'}_{\varepsilon;\sigma}((0,\sigma)) = \int_{0}^{\infty} dt \cdot h(t) \mathbf{m}^{z,z';\Lambda}_{\varepsilon}(dt) \quad \forall \varepsilon \geq 0.
\end{align*}
To prove the second convergence in (\ref{weak_convergence_probability_measure}), it is enough to show that 
\begin{align} \label{weak_convergence_m_epsilon}
    \lim_{\varepsilon\to 0^{+}}\int_{0}^{\infty} dt \cdot h(t) \mathbf{m}^{z,z';\Lambda}_{\varepsilon}(dt) = \int_{0}^{\infty} dt \cdot h(t) \mathbf{m}^{z,z';\Lambda}_{0}(dt).
\end{align}
To this aim, we notice that there exists functions $(h_{\delta;-}(t))_{\delta>0}$ and $(h_{\delta;+}(t))_{\delta>0}$ in $C_{b}(\mathbb{R}_{+};\mathbb{R}_{+})$ such that $0\leq h_{\delta;-}(t) \leq h(t) \leq h_{\delta;+}(t)$ and $\lim_{\delta\to 0^{+}} h_{\delta;+}(t) = \lim_{\delta\to 0^{+}} h_{\delta;-}(t)= h(t)$. Then, for every $\delta > 0$, the following estimate holds:
\begin{align} \label{estimate_liminf_limsup}
    \int_{0}^{\infty} dt \cdot h_{\delta;-}(t) &\mathbf{m}^{z,z';\Lambda}_{0}(dt) \leq \liminf_{\varepsilon\to 0^{+}}\int_{0}^{\infty} dt \cdot h(t) \mathbf{m}^{z,z';\Lambda}_{\varepsilon}(dt)\nonumber\\
    &\leq \limsup_{\varepsilon\to 0^{+}}\int_{0}^{\infty} dt \cdot h(t) \mathbf{m}^{z,z';\Lambda}_{\varepsilon}(dt)\leq \int_{0}^{\infty} dt \cdot h_{\delta;+}(t) \mathbf{m}^{z,z';\Lambda}_{0}(dt),
\end{align}
where (\ref{weak_convergence_total_mass_1}) and (\ref{weak_convergence_probability_measure_1}) have been used in the first and the second inequalities. Therefore, by taking $\delta \to 0^{+}$, (\ref{estimate_liminf_limsup}) gives (\ref{weak_convergence_m_epsilon}), and thus, we conclude the second convergence in (\ref{weak_convergence_probability_measure}).

Finally, we prove the first convergence in (\ref{weak_convergence_probability_measure}). To do this, it suffices to show that 
\begin{align} \label{weak_convergence_widetilde_h}
    \lim_{\varepsilon\to 0^{+}}\int_{0}^{\infty} dt \cdot \widetilde{h}(t) \mathbf{u}^{z,z';\Lambda}_{\varepsilon}(dt) = \int_{0}^{\infty} dt \cdot \widetilde{h}(t) \mathbf{u}^{z,z';\Lambda}_{0}(dt) \quad \widetilde{h} \in C_{b}((0,\sigma);\mathbb{R}_{+}).
\end{align}
This can be done by considering $h(t) := \widetilde{h}(t) \exp(\Lambda \cdot t) \cdot 1_{(0,\sigma)}(t)$ a similar argument as above.
Thus, combining (\ref{weak_convergence_widetilde_h}) and the second convergence in (\ref{weak_convergence_probability_measure}) implies the first convergence in (\ref{weak_convergence_probability_measure}). Therefore, the proof of Proposition \ref{proposition_critical_second_moment_point_to_point} is complete.

\end{proof}

\section{Properties of the 3D Birman-Schwinger operators.} 
Before the proof of Lemma \ref{lemma_critical_second_moment_point_to_point}, we first introduce some preliminary properties of the Birman-Schwinger operator $\mathbf{T}^{\lambda}_{\beta}$ at $L^{2}$-criticality.
\begin{lemma} \label{lemma_Hilbert_Schmidt_operator}
Assume that $d = 3$ and $\beta = \beta_{L^{2}}$. Recall the Birman-Schwinger operator $\mathbf{T}^{\lambda}_{\beta}$ defined in \ref{definition_Birman-Schwinger_operator}, and  $\mathfrak{V}_{\beta}(x)$ is defined in (\ref{definition_V_and_Yukawa}) Then, we have the following properties.
\begin{enumerate} [label=(\roman*)]
    \item The Birman-Schwinger operator $\mathbf{T}^{0}_{\beta}$ is a non-negative Hilbert–Schmidt operator with eigenvalues $\lambda_{1} \geq \lambda_{2}\geq ...\geq 0$ and eigenvectors $\{\mathbf{v}_{j}:j = 1,2,...\} \subseteq L^{2}(\mathbb{R}^{d})$, where $\lambda_{1} = 1>\lambda_{2}$, $\{\mathbf{v}_{j}:j = 1,2,...\}$ is a complete orthonormal basis of $L^{2}(\mathbb{R}^{d})$, and $\mathbf{v}_{1}$ is a real-valued function. In particular, it holds that 
    \begin{align} \label{non_zero}
        \langle \mathbf{v}_{1}, \mathfrak{V}_{\beta} \rangle_{L^{2}(\mathbb{R}^{d})} \neq 0.
    \end{align}
    \item Recall that both of $V_{\beta}(x)$ and $\mathcal{G}^{0}(x)$are defined in (\ref{definition_V_and_Yukawa}), and $\mathbf{L}_{\beta}$ is given in (\ref{definition_hamiltonian_L}). Then $\mathbf{L}_{\beta} h_{\mathbf{v}_{1}} = 0$ in the sense of distribution and $h_{\mathbf{v}_{1}} \not \in L^{2}(\mathbb{R}^{d})$, where
    \begin{align} \label{definition_h_again}
    h_{\mathbf{v}_{1}}(x) := \int_{\mathbb{R}^{d}} \mathcal{G}^{0}(x-y) \mathfrak{V}_{\beta}(y) \mathbf{v}_{1}(y) dy.
    \end{align}
    In particular, $h_{\mathbf{v}_{1}}(\cdot) (1+|\cdot|)^{-s} \in L^{2}(\mathbb{R}^{d})$ for any $s > 1/2$ such that  
    \begin{align} \label{equation_zero_ro reson}
        h_{\mathbf{v}_{1}}(x) = \int_{\mathbb{R}^{d}} \mathcal{G}^{0}(x-y) V_{\beta}(y) h_{\mathbf{v}_{1}}(y) dy.
    \end{align}
\end{enumerate}
\end{lemma}
\begin{proof} 
We begin by proving all properties in the part (i) of Lemma \ref{lemma_Hilbert_Schmidt_operator} while assuming that 
\begin{align} \label{assumption_v_1}
    \langle \mathbf{v}, \mathfrak{V}_{\beta} \rangle_{L^{2}(\mathbb{R}^{d})} \neq 0 \quad \text{for all $\mathbf{v} \in L^{2}(\mathbb{R}^{d})$ such that } \mathbf{v} = \mathbf{T}^{0}_{\beta}\mathbf{v} \text{ and } ||\mathbf{v}||_{L^{2}(\mathbb{R}^{d})} = 1 \text{ when } \beta = \beta_{L^{2}}.
\end{align}
First, the non-negativity of the Birman-Schwinger operator $\mathbf{T}^{0}_{\beta}$ has been proved in Lemma \ref{lemma_properties_Birman-Schwinger_operators}. Second, due to (\ref{definition_Birman-Schwinger_operator}) and the fact that $V_{\beta} \in C_{c}^{\infty}(\mathbb{R}^{d};\mathbb{R}_{+})$, the following estimate shows that  $\mathbf{T}^{0}_{\beta}$ is a Hilbert-Schmidt operator:
\begin{align} \label{Hilbert_Schmidt_operator_check}
    \int_{\mathbb{R}^{2 \times d}} dzdz' \cdot V_{\beta}(z) \frac{1}{|z-z'|^{2}} V_{\beta}(z') < \infty,
\end{align}
where we have used the fact that $\mathcal{G}^{0}(z) = C/|z|$ for some positive constant $C$. Third, (\ref{identity_beta2+_=_betaL2}) and Lemma \ref{lemma_representation_beta_2_+} implies $\lambda_{1} = ||\mathbf{T}^{0}_{\beta}||_{L^{2}(\mathbb{R}^{d})} = 1$. Furthermore, if $\lambda_{k} = 1$ for some $k\geq 2$, then the following non-zero function clearly contradicts the assumption (\ref{assumption_v_1}):
\begin{align*}
    \widetilde{\mathbf{v}}(x) =  \biggbkt{\mathbf{v}_{k}(x) - \frac{\langle \mathbf{v}_{k}, \mathfrak{V}_{\beta} \rangle_{L^{2}(\mathbb{R}^{d})}}{\langle \mathbf{v}_{1}, \mathfrak{V}_{\beta} \rangle_{L^{2}(\mathbb{R}^{d})}} \mathbf{v}_{1}(x) }\biggl /\Biggl|\Biggl|\mathbf{v}_{k}(x) - \frac{\langle \mathbf{v}_{k}, \mathfrak{V}_{\beta} \rangle_{L^{2}(\mathbb{R}^{d})}}{\langle \mathbf{v}_{1}, \mathfrak{V}_{\beta} \rangle_{L^{2}(\mathbb{R}^{d})}} \mathbf{v}_{1}(x) \Biggr|\Biggr|_{L^{2}(\mathbb{R}^{d})},
\end{align*}
Then, one has $\lambda_{2} < \lambda_{1} = 1$. Finally, since $\lambda_{2} < \lambda_{1}$ and $\overline{\mathbf{v}}_{1} = \mathbf{T}^{0}_{\beta}\overline{\mathbf{v}}_{1}$, it follows that $\overline{\mathbf{v}}_{1} = \mathbf{v}_{1}$.

To complete the proof of the part (i) of Lemma \ref{lemma_Hilbert_Schmidt_operator}, it remains to prove (\ref{assumption_v_1}). We prove it by contradiction. Hence, we assume that $\langle \mathfrak{V}_{\beta}, \mathbf{v} \rangle_{L^{2}(\mathbb{R}^{d})} = 0$, where $\mathbf{v}$ satisfies the rest of the conditions in (\ref{assumption_v_1}). Now, we consider the auxiliary function $h_{\mathbf{v}}$, where it is defined as (\ref{definition_h_again}). Notice that $h_{\mathbf{v}} \not \equiv 0$ follows from $\mathbf{v}(x) = \mathbf{T}^{0}_{\beta}\mathbf{v}(x) = \mathfrak{V}_{\beta}(x) h_{\mathbf{v}}(x)$. The main idea to reach a contradiction is to show that $h_{\mathbf{v}} = c \cdot \widetilde{h}_{\mathbf{v}_{1}}$, where $c \in \mathbb{C}$ and $\widetilde{h}$ is a strictly positive function. Indeed, let $\widetilde{\mathbf{v}}(x) :=  \mathfrak{V}_{\beta}(x) \widetilde{h}_{\mathbf{v}}(x)$, we know that $\widetilde{\mathbf{v}}$ is strictly positive in the support of $\mathfrak{V}_{\beta}$, and thus, $\langle \mathfrak{V}_{\beta}, \widetilde{\mathbf{v}} \rangle_{L^{2}(\mathbb{R}^{d})} > 0$. This implies that $\langle \mathfrak{V}_{\beta}, \mathbf{v} \rangle_{L^{2}(\mathbb{R}^{d})} = c \cdot \langle \mathfrak{V}_{\beta}, \widetilde{\mathbf{v}} \rangle_{L^{2}(\mathbb{R}^{d})}\neq 0$.

With the above clarification, let us now show that, under the assumption $\langle \mathfrak{V}_{\beta}, \mathbf{v} \rangle_{L^{2}(\mathbb{R}^{d})} = 0$, where $\mathbf{v}$ has the rest of the conditions in (\ref{assumption_v_1}), $h_{\mathbf{v}}  = c \cdot \widetilde{h}_{\mathbf{v}}$ for some strictly positive function $\widetilde{h}_{\mathbf{v}}$. The key is to apply \cite[Theorem 11.8]{lieb2001analysis}. Namely, our goal is to show that when $\beta = \beta_{L^{2}}$, the following holds: ($\underline{1}$) $h_{\mathbf{v}} \in H^{1}(\mathbb{R}^{d})$; ($\underline{2}$) $V_{\beta} \cdot |h_{\mathbf{v}}|^{2} \in L^{1}(\mathbb{R}^{d})$; ($\underline{3}$) $\langle h_{\mathbf{v}},\mathbf{L}^{\beta} h_{\mathbf{v}} \rangle_{L^{2}(\mathbb{R}^{d})} = \sup \mathbf{L}^{\beta}$. To do this, applying \cite[Lemma 1.2.3]{albeverio2012solvable} with potential function $-2V_{\beta}(x)$, where the $\phi$ in \cite[(1.2.24)]{albeverio2012solvable} corresponds to $\mathbf{v}$, we know that $\sqrt{2}h_{\mathbf{v}}$ has the properties mentioned in the conclusion of \cite[Lemma 1.2.3]{albeverio2012solvable}. This shows the following: (a) $\nabla h_{\mathbf{v}} \in L^{2}(\mathbb{R}^{d})$; (b) $\mathbf{L}^{\beta} h_{\mathbf{v}} = 0$ in the sense of distribution; 
\begin{align*} 
    \text{(c)}\quad h_{\mathbf{v}} \in L^{2}(\mathbb{R}^{d}) \quad \text{if and only if} \quad \langle \mathfrak{V}_{\beta}, \mathbf{v} \rangle_{L^{2}(\mathbb{R}^{d})} = 0.
\end{align*}
Then, thanks to our assumption, $h \in L^{2}(\mathbb{R}^{d})$ follows from (c), and thus, combining this fact with (a) implies ($\underline{1}$). Moreover, ($\underline{2}$) follows from the fact $V_{\beta} \in C_{c}^{\infty}(\mathbb{R}^{d})$. To prove ($\underline{3}$), we show that $\langle h_{\mathbf{v}},\mathbf{L}^{\beta} h_{\mathbf{v}} \rangle_{L^{2}(\mathbb{R}^{d})} = 0 = \sup \mathbf{L}^{\beta}$. On the one hand, since $h_{\mathbf{v}} \in H^{1}(\mathbb{R}^{d})$ and $\mathbf{L}^{\beta} h_{\mathbf{v}} = 0$ in the sense of distribution, a standard density argument then implies that $\langle h_{\mathbf{v}},\mathbf{L}^{\beta} h_{\mathbf{v}} \rangle_{L^{2}(\mathbb{R}^{d})} = 0$. On the other hand, as the proof of Lemma \ref{lemma_continuity}, one can show that $\beta \mapsto \sup \mathbf{L}^{\beta}$ is continuous. Hence, thanks to (\ref{identity_beta2+_=_betaL2}), we see that $\sup \mathbf{L}^{\beta} = 0$ when $\beta = \beta_{L^{2}}$. As a result, ($\underline{3}$) is proved. Therefore, the proof of the part (i) is complete.

Finally, we explain the reason why the part (ii) follows from the part (i). First, we notice that the first property in the part (ii) is (b). Second, in the previous proof, $\langle \mathfrak{V}_{\beta}, \mathbf{v}_{1} \rangle_{L^{2}(\mathbb{R}^{d})} \neq 0$ is proved. Then, the second property in the part (ii) follows from (c). Third, (\ref{equation_zero_ro reson}) follows from (\ref{definition_h_again}) and the fact that $\mathbf{v}_{1}$ is a eigenvector of $\mathbf{T}_{\beta}^{0}$ with eigenvalue $1$. Finally, we show that $h_{\mathbf{v}_{1}}(\cdot) (1+|\cdot|)^{-s} \in L^{2}(\mathbb{R}^{d})$ for any $s > 1/2$. To this aim, we notice that $|h_{\mathbf{v}_{1}}(x)| \lesssim ||\mathbf{v}_{1}||_{L^{2}(\mathbb{R}^{d})} = 1$ by using Cauchy–Schwarz inequality and the fact $\mathfrak{V}_{\beta}\in C_{c}(\mathbb{R}^{d};\mathbb{R}_{+})$. Moreover, we observe that $|h_{\mathbf{v}_{1}}(x)| \lesssim \frac{1}{|x|-2r_{\phi}} ||\mathfrak{V}_{\beta}||_{L^{2}(\mathbb{R}^{d})} ||\mathbf{v}_{1}||_{L^{2}(\mathbb{R}^{d})}$ for every $|x| > 2r_{\phi}$, where (\ref{support_R}) has been used. With the above observations, $h_{\mathbf{v}_{1}}(\cdot) (1+|\cdot|)^{-s} \in L^{2}(\mathbb{R}^{d})$ can be proved if $s > 1/2$.

\end{proof}
\section{Proof of Lemma \ref{lemma_critical_second_moment_point_to_point}} \label{section_proof_lemma_critical_second_moment_point_to_point}
This section aims to prove Lemma \ref{lemma_critical_second_moment_point_to_point} by using \cite[Lemma 1.2.4]{albeverio2012solvable} in our framework. Throughout this section, we assume that $d =3$ and $\beta = \beta_{L^{2}}$. 

Before we proceed, we introduce the heuristics for Lemma \ref{lemma_critical_second_moment_point_to_point}. Fix $\Lambda > 0$ and $z\neq z'$.  To begin with, we will show that the Laplace transform of the rescaled Brownian exponential functional can be re-expressed in terms of the rescaled Birman-Schwinger operators $\mathbf{T}_{\beta}^{\varepsilon^{2} \cdot\Lambda}$ given in (\ref{definition_Birman-Schwinger_operator}):
\begin{align} \label{lemma_critical_second_moment_point_to_point_expression}
    \text{(L.H.S) of (\ref{part_i_lemma_critical_second_moment_point_to_point})}= \beta_{\varepsilon}^{2} \mathcal{G}^{\varepsilon^{2} \cdot \Lambda}(z-z') \beta^{2} + \biggl\langle \mathbf{a}_{\varepsilon;z}, \varepsilon \cdot \biggbkt{I-\mathbf{T}_{\beta}^{\varepsilon^{2} \cdot \Lambda}}^{-1} \mathbf{b}_{\varepsilon;z'} \biggr\rangle_{L^{2}(\mathbb{R}^{d})} \quad\forall \varepsilon>0.
\end{align}
Here, $\mathbf{a}_{\varepsilon;z}$ and $\mathbf{b}_{\varepsilon;z'}$ are functions in $L^{2}(\mathbb{R}^{d})$ defined as follows:
\begin{align*}
    \mathbf{a}_{\varepsilon;z}(x) := \beta^{2}\mathcal{G}^{\varepsilon^{2} \cdot \Lambda}(z-x) \cdot \mathfrak{V}_{\beta}(x) \quad \text{and} \quad \mathbf{b}_{\varepsilon;z'}(x) := \mathfrak{V}_{\beta}(x) \mathcal{G}^{\varepsilon^{2} \cdot \Lambda}(x-z') \beta^{2} \quad \forall \varepsilon \geq 0,
\end{align*}
where $\mathcal{G}^{\varepsilon^{2}\cdot \Lambda}(x)$ and $\mathfrak{V}_{\beta}(x)$ are given in (\ref{definition_V_and_Yukawa}). With the above notations, the key to prove Lemma \ref{lemma_critical_second_moment_point_to_point} is to show the following approximation when $\beta = \beta_{L^{2}}$:
\begin{align} \label{lemma_critical_second_moment_point_to_point_major_step}
    \lim_{\varepsilon \to 0^{+}} \varepsilon \cdot \biggbkt{I-\mathbf{T}_{\beta}^{\varepsilon^{2} \cdot \Lambda}}^{-1} = \frac{2\pi}{\sqrt{2\Lambda}} \cdot \mathbf{A}_{\beta}, \quad \text{where} \quad\mathbf{A}_{\beta}(x,x') := |\langle \mathbf{v}_{1},\mathfrak{V}_{\beta} \rangle_{L^{2}(\mathbb{R}^{d})}|^{-2} \cdot \mathbf{P}_{\beta;1},
\end{align}
where the above convergence is in the sense of the convergence of bounded linear operators on $L^{2}(\mathbb{R}^{d})$, $\mathbf{P}_{\beta;1}$ is the projection operator along the direction $\mathbf{v}_{1}$ defined in Lemma \ref{lemma_Hilbert_Schmidt_operator}, and the well-definedness of $\mathbf{A}_{\beta}$ follows from (\ref{non_zero}).

Let us now explain the reason why combining (\ref{lemma_critical_second_moment_point_to_point_expression}) and (\ref{lemma_critical_second_moment_point_to_point_major_step}) implies Lemma \ref{lemma_critical_second_moment_point_to_point}.
\begin{proof} [Proof of Lemma \ref{lemma_critical_second_moment_point_to_point}]
Due to $z \neq z'$, the first term on the right-hand side of (\ref{lemma_critical_second_moment_point_to_point_expression}) vanishes as $\varepsilon \to 0^{+}$. Since $\mathbf{a}_{\varepsilon;z} \to \mathbf{a}_{0;z}$ and $\mathbf{b}_{\varepsilon;z'} \to \mathbf{b}_{0;z'}$ in $L^{2}(\mathbb{R}^{d})$, applying
(\ref{lemma_critical_second_moment_point_to_point_major_step}) gives
\begin{align} 
    &\lim_{\varepsilon \to 0^{+}}\int_{0}^{\infty} dt \exp(-\Lambda \cdot t)
    \beta_{\varepsilon}^{2}  \mathbf{E}_{\varepsilon \cdot z}\bigbkt{\exp\bkt{\beta_{\varepsilon}^{2} \int_{0}^{t} R_{\varepsilon}(\sqrt{2}B(r)) dr} : B(t) = \varepsilon \cdot z'} \beta_{\varepsilon}^{2}\nonumber\\
    &= \frac{2\pi}{\sqrt{2\Lambda}} \cdot \langle \mathbf{a}_{0;z}, \mathbf{A}_{\beta} \mathbf{b}_{0;z'} \rangle_{L^{2}(\mathbb{R}^{d})} = \frac{2\pi}{\sqrt{2\Lambda}} \cdot \mathscr{C}(z,z'),
\end{align}
where $\mathscr{C}(z,z')$ is defined in (\ref{definition_C_z_z'}). Hence, we conclude the part (i) of Lemma \ref{lemma_critical_second_moment_point_to_point}. Moreover, due to (\ref{lemma_critical_second_moment_point_to_point_major_step}) and the facts that the kernel of $\mathbf{T}_{\beta}^{\varepsilon^{2} \cdot \Lambda}$ is non-negative and that $\mathbf{v}_{1}$ is real-valued, where both of them are proved in Lemma \ref{lemma_Hilbert_Schmidt_operator}, the following estimate holds:
\begin{align*}
    0&\leq \lim_{\varepsilon \to 0^{+}}\biggl\langle u, \varepsilon\cdot\sum_{k=0}^{\infty} (\mathbf{T}_{\beta}^{\varepsilon^{2} \cdot \Lambda})^{k} \widetilde{u} \biggr\rangle_{L^{2}(\mathbb{R}^{d})} = \lim_{\varepsilon \to 0^{+}}\biggl\langle u, \varepsilon \cdot \biggbkt{I-\mathbf{T}_{\beta}^{\varepsilon^{2} \cdot \Lambda}}^{-1} \widetilde{u} \biggr\rangle_{L^{2}(\mathbb{R}^{d})} \\
    &= \biggbkt{|\langle \mathbf{v}_{1},\mathfrak{V}_{\beta} \rangle_{L^{2}(\mathbb{R}^{d})}|^{-2} \cdot \frac{2\pi}{\sqrt{2\Lambda}} } \cdot \int_{\mathbb{R}^{2\times d}} dxdx' u(x)\widetilde{u}(x') \cdot \mathbf{v}_{1}(x)\mathbf{v}_{1}(x') \quad \forall u,\widetilde{u} \in C_{c}^{\infty}(\mathbb{R}^{d};\mathbb{R}_{+}).
\end{align*}
This estimate shows that $\mathbf{v}_{1}$ can be chosen to be a non-negative function in $L^{2}(\mathbb{R}^{d})$. Hence, the proof of the part (ii) of Lemma \ref{lemma_critical_second_moment_point_to_point} is complete.
    
\end{proof}

To complete the proof of Lemma \ref{lemma_critical_second_moment_point_to_point}, let us now present the proof of (\ref{lemma_critical_second_moment_point_to_point_expression}).
\begin{proof} [Proof of (\ref{lemma_critical_second_moment_point_to_point_expression})]
To begin with, observe that the following identity follows from the scaling property of Brownian motion:
\begin{align} \label{identity_laplace_0}
    &\text{(L.H.S) of (\ref{lemma_critical_second_moment_point_to_point_expression})} = \int_{0}^{\infty} dt \exp(-\varepsilon^{2} \Lambda  \cdot t)
    \beta_{\varepsilon}^{2}  \mathbf{E}_{z}\bigbkt{\exp\bkt{\beta^{2} \int_{0}^{t} R(\sqrt{2}B(r)) dr} : B(t) = z'} \beta^{2}.
\end{align}
Moreover, applying the fundamental theorem of calculus (\ref{identity_fundamental}) to the Brownian exponential functional on the right-hand side of the above identity, it holds that:
\begin{align} \label{identity_laplace_1}
    &\text{(R.H.S) of (\ref{identity_laplace_0})}\nonumber\\
    &= \int_{0}^{\infty} dt \exp(-\varepsilon^{2} \Lambda  \cdot t)
    \beta_{\varepsilon}^{2} \biggbkt{G_{t}(z-z')+\sum_{k=1}^{\infty} \int_{s_{0} = 0<s_{1}<...<s_{k}<t} \prod_{j=1}^{k} ds_{j} \cdot\int_{\mathbb{R}^{k\times d},\; x_{0} = z} \prod_{j=1}^{k} dx_{j} \cdot \nonumber\\
    &\biggbkt{\prod_{j=1}^{k} G_{s_{j}-s_{j-1}}(x_{j}-x_{j-1}) V_{\beta}(x_{j}) } \cdot G_{t-s_{k}}(z'-x_{k})}\beta^{2}\nonumber\\
    &= \beta_{\varepsilon}^{2} \mathcal{G}^{\varepsilon^{2}\cdot\Lambda}(z-z')\beta^{2}\nonumber\\
    &+\beta_{\varepsilon}^{2} \cdot \biggbkt{ \sum_{k=1}^{\infty} \int_{\mathbb{R}^{k\times d},\; x_{0} = z} \prod_{j=1}^{k} dx_{j} \cdot\biggbkt{\prod_{j=1}^{k} \mathcal{G}^{\varepsilon^{2}\cdot \Lambda}(x_{j}-x_{j-1}) V_{\beta}(x_{j}) } \cdot \mathcal{G}^{\varepsilon^{2}\cdot \Lambda}(z'-x_{k})} \cdot \beta^{2},
\end{align}
where $V_{\beta}(x)$ and $\mathcal{G}^{\lambda}(x)$ are defined in (\ref{definition_V_and_Yukawa}), and the following fact has been used in the last equality:
\begin{align*}
    \int_{0}^{\infty} dt \exp(-\lambda \cdot t) \int_{0}^{t} ds h(t-s) \cdot r(s) = \biggbkt{\int_{0}^{\infty} dt \exp(-\lambda \cdot t) h(t)} \cdot \biggbkt{\int_{0}^{\infty} dt \exp(-\lambda \cdot t) r(t) }.
\end{align*}
Therefore, re-expressing the right-hand side of (\ref{identity_laplace_1}) by using (\ref{definition_Birman-Schwinger_operator}) and (\ref{definition_V_and_Yukawa}) gives the right-hand side of (\ref{lemma_critical_second_moment_point_to_point_expression}).
\end{proof}

Finally, it remains to prove (\ref{lemma_critical_second_moment_point_to_point_major_step}). 
\begin{proof} [Proof of (\ref{lemma_critical_second_moment_point_to_point_major_step})]
To apply \cite[Lemma 1.2.4]{albeverio2012solvable}, we observe that 
\begin{align*}
    -\mathbf{T}^{\epsilon^{2} \cdot \Lambda}_{\beta}(x,x') = \lambda(\epsilon) \cdot (-1)\biggbkt{2 V_{\beta}(x)}^{1/2} \cdot \frac{1}{4\pi |x-x'|} \exp\biggbkt{ik\varepsilon \cdot |x-x'|} \cdot \biggbkt{2 V_{\beta}(x')}^{1/2},
\end{align*}
where the right-hand side is of the form \cite[(1.2.13)]{albeverio2012solvable} with potential function $-2V_{\beta}(x)$.  Here, we have set $\lambda(\epsilon) := 1$ and $k := i\sqrt{2\Lambda}$, and the fact $\mathcal{G}^{\lambda}(x) = \exp(-\sqrt{2\lambda}|x|)/(2\pi |x|)$ has been used. Moreover, thanks to Lemma \ref{lemma_Hilbert_Schmidt_operator}, we know that the $(\phi,\psi) := (\mathbf{v}_{1},\sqrt{2}h_{\mathbf{v}_{1}})$ satisfies the conditions in the case II in \cite[pp. 20]{albeverio2012solvable}. Consequently, applying \cite[Lemma 1.2.4]{albeverio2012solvable} shows that 
\begin{align*}
    \lim_{\varepsilon \to 0^{+}}\varepsilon \cdot \biggbkt{I+ (-\mathbf{T}^{\epsilon^{2} \cdot \Lambda}_{\beta})}^{-1} = \frac{4\pi}{ik} \cdot|\langle \sqrt{2} \mathfrak{V}_{\beta},\mathbf{v}_{1} \rangle_{L^{2}(\mathbb{R}^{d})}|^{-2} \cdot (-\mathbf{P}_{\beta;1}).
\end{align*}
Then, the proof is complete.
\end{proof}

\small
\bibliographystyle{plainurl}
\bibliography{mybibliography.bib}

\begin{thebibliography}{10}

\bibitem{Alberts_2013}
T.~Alberts, K.~Khanin, and J.~Quastel.
\newblock The continuum directed random polymer.
\newblock {\em Journal of Statistical Physics}, 154(1-2):305--326, 2013.
\newblock \href {https://doi.org/10.1007/s10955-013-0872-z} {\path{doi:10.1007/s10955-013-0872-z}}.

\bibitem{albeverio2012solvable}
S.~Albeverio, F.~Gesztesy, R.~Hoegh-Krohn, and H.~Holden.
\newblock {\em Solvable Models in Quantum Mechanics}.
\newblock Theoretical and Mathematical Physics. Springer Berlin Heidelberg, 2012.

\bibitem{teta2}
G.~Basti, C.~Cacciapuoti, D.~Finco, and A.~Teta.
\newblock Three-body hamiltonian with regularized zero-range interactions in dimension three.
\newblock {\em Annales Henri Poincar{\'e}}, 24(1):223--276, 2023.
\newblock \href {https://doi.org/10.1007/s00023-022-01214-9} {\path{doi:10.1007/s00023-022-01214-9}}.

\bibitem{betaclessthenbetaL2_3}
M.~Birkner and R.~Sun.
\newblock {Annealed vs quenched critical points for a random walk pinning model}.
\newblock {\em Annales de l'Institut Henri Poincaré, Probabilités et Statistiques}, 46(2):414--441, 2010.
\newblock \href {https://doi.org/10.1214/09-AIHP319} {\path{doi:10.1214/09-AIHP319}}.

\bibitem{brezis2011functional}
H.~Brezis.
\newblock {\em Functional analysis, Sobolev spaces and partial differential equations}.
\newblock Universitext. Springer New York, 2010.
\newblock \href {https://doi.org/10.1007/978-0-387-70914-7} {\path{doi:10.1007/978-0-387-70914-7}}.

\bibitem{UU}
F.~Caravenna, R.~Sun, and N.~Zygouras.
\newblock Universality in marginally relevant disordered systems.
\newblock {\em The Annals of Applied Probability}, 27(5):3050--3112, 2017.
\newblock \href {https://doi.org/10.1214/17-AAP1276} {\path{doi:10.1214/17-AAP1276}}.

\bibitem{caravenna2019moments}
F.~Caravenna, R.~Sun, and N.~Zygouras.
\newblock On the moments of the (2+ 1)-dimensional directed polymer and stochastic heat equation in the critical window.
\newblock {\em Communications in Mathematical Physics}, 372(2):385--440, 2019.
\newblock \href {https://doi.org/10.1007/s00220-019-03527-z} {\path{doi:10.1007/s00220-019-03527-z}}.

\bibitem{NK_KPZ}
F.~Caravenna, R.~Sun, and N.~Zygouras.
\newblock {The two-dimensional KPZ equation in the entire subcritical regime}.
\newblock {\em The Annals of Probability}, 48(3):1086--1127, 2020.
\newblock \href {https://doi.org/10.1214/19-AOP1383} {\path{doi:10.1214/19-AOP1383}}.

\bibitem{Heat_flow}
F.~Caravenna, R.~Sun, and N.~Zygouras.
\newblock {The critical $2D$ Stochastic Heat Flow}.
\newblock {\em Inventiones mathematicae}, 233:325--460, 2023.
\newblock \href {https://doi.org/10.1007/s00222-023-01184-7} {\path{doi:10.1007/s00222-023-01184-7}}.

\bibitem{chen}
Y.-T. Chen.
\newblock {Delta-Bose gas from the viewpoint of the two-dimensional stochastic heat equation}.
\newblock {\em The Annals of Probability}, 52(1):127 -- 187, 2024.
\newblock \href {https://doi.org/10.1214/23-AOP1649} {\path{doi:10.1214/23-AOP1649}}.

\bibitem{comets2017directed}
F.~Comets.
\newblock {\em Directed polymers in random environments}, volume 2175 of {\em Lecture Notes in Mathematics}.
\newblock Springer International Publishing, 2017.
\newblock \href {https://doi.org/10.1007/978-3-319-50487-2} {\path{doi:10.1007/978-3-319-50487-2}}.

\bibitem{GFDDP}
C.~Cosco and S.~Nakajima.
\newblock {Gaussian fluctuations for the directed polymer partition function in dimension $d\ge 3$ and in the whole ${L^{2}}$-region}.
\newblock {\em Annales de l'Institut Henri Poincaré, Probabilités et Statistiques}, 57(2):872--889, 2021.
\newblock \href {https://doi.org/10.1214/20-AIHP1100} {\path{doi:10.1214/20-AIHP1100}}.

\bibitem{COSCO2022127}
C.~Cosco, S.~Nakajima, and M.~Nakashima.
\newblock {Law of large numbers and fluctuations in the sub-critical and $L^{2}$ regions for SHE and KPZ equation in dimension $d\geq 3$}.
\newblock {\em Stochastic Processes and their Applications}, 151:127--173, 2022.
\newblock \href {https://doi.org/10.1016/j.spa.2022.05.010} {\path{doi:10.1016/j.spa.2022.05.010}}.

\bibitem{cosco2023moments}
C.~Cosco and O.~Zeitouni.
\newblock {Moments of partition functions of $2D$ gaussian polymers in the weak disorder regime - I}.
\newblock {\em Communications in Mathematical Physics}, 403(1):417--450, 2023.
\newblock \href {https://doi.org/10.1007/s00220-023-04799-2} {\path{doi:10.1007/s00220-023-04799-2}}.

\bibitem{CCO}
C.~Cosco and O.~Zeitouni.
\newblock {Moments of partition functions of $2D$ Gaussian polymers in the weak disorder regime - II}.
\newblock {\em Electronic Journal of Probability}, 29(none), 2024.
\newblock \href {https://doi.org/10.1214/24-EJP1148} {\path{doi:10.1214/24-EJP1148}}.

\bibitem{ferretti2023contact}
F.~Daniele.
\newblock Contact interactions for many-particle quantum systems in dimension three.
\newblock 2023.
\newblock URL: \url{https://iris.uniroma1.it/handle/11573/1669960}.

\bibitem{Teta}
G.~F. Dell'Antonio, R.~Figari, and A.~Teta.
\newblock Hamiltonians for systems of {N} particles interacting through point interactions.
\newblock {\em Annales de l'I.H.P. Physique th\'eorique}, 60(3):253--290, 1994.
\newblock URL: \url{http://www.numdam.org/item/AIHPA_1994__60_3_253_0/}.

\bibitem{efimov2}
V~N Efimov.
\newblock Weakly bound states of three resonantly interacting particles.
\newblock {\em Yadern. Fiz.}, 12:1080--1091, 1970.
\newblock URL: \url{https://www.osti.gov/biblio/4068792}.

\bibitem{evans2010partial}
L.C. Evans.
\newblock {\em Partial Differential Equations}, volume~19 of {\em Graduate studies in mathematics}.
\newblock American Mathematical Society, second edition, 2010.
\newblock URL: \url{https://bookstore.ams.org/GSM-19-R/}.

\bibitem{DBG_summary}
R.~Figari and A.~Teta.
\newblock On the hamiltonian for three bosons with point interactions.
\newblock {\em Quantum and Stochastic Mathematical Physics}, pages 127--145, 2023.
\newblock \href {https://doi.org/10.1007/978-3-031-14031-0_6} {\path{doi:10.1007/978-3-031-14031-0_6}}.

\bibitem{frank2024hardy}
R.L Frank, T.~Hoffmann-Ostenhof, A.~Laptev, and J.P. Solovej.
\newblock Hardy inequalities for large fermionic systems.
\newblock {\em Journal of Spectral Theory}, 14(2), 2024.
\newblock \href {https://doi.org/10.4171/JST/511} {\path{doi:10.4171/JST/511}}.

\bibitem{griesemer2023weakness}
M.~Griesemer and M.~Hofacker.
\newblock On the weakness of short-range interactions in fermi gases.
\newblock {\em Letters in Mathematical Physics}, 113(1):1, 2023.
\newblock \href {https://doi.org/10.1007/s11005-022-01624-0} {\path{doi:10.1007/s11005-022-01624-0}}.

\bibitem{Gu_2d}
Y.~Gu, J.~Quastel, and L.-C. Tasi.
\newblock {Moments of the $2D$ SHE at criticality}.
\newblock {\em Probability and Mathematical Physics}, 2(1):179 -- 219, 2021.
\newblock URL: \url{https://doi.org/10.2140/pmp.2021.2.179}.

\bibitem{Guzu}
D.~Guzu.
\newblock {Many-Particle Hardy Type Inequalities}.
\newblock {\em Journal of Mathematical Sciences}, 239(3), 2019.
\newblock \href {https://doi.org/10.1007/s10958-019-04305-x} {\path{doi:10.1007/s10958-019-04305-x}}.

\bibitem{hall2013quantum}
B.C. Hall.
\newblock {\em Quantum theory for mathematicians}, volume 267.
\newblock Springer Science \& Business Media, 2013.
\newblock \href {https://doi.org/10.1007/978-1-4614-7116-5} {\path{doi:10.1007/978-1-4614-7116-5}}.

\bibitem{griesemer2023weakness_thesis}
M.~Hofacker.
\newblock {\em From Short-Range to Contact Interactions in Many-Body Quantum Systems}.
\newblock PhD thesis. University of Stuttgart, 2022.

\bibitem{hoffmann2008many}
M.~Hoffmann-Ostenhof, T.~Hoffmann-Ostenhof, A.~Laptev, and J.~Tidblom.
\newblock Many-particle hardy inequalities.
\newblock {\em Journal of the London Mathematical Society}, 77(1):99--115, 2008.
\newblock \href {https://doi.org/10.1112/jlms/jdm091} {\path{doi:10.1112/jlms/jdm091}}.

\bibitem{GHS}
S.~Janson.
\newblock {\em {Gaussian Hilbert spaces}}.
\newblock Cambridge Tracts in Mathematics. Cambridge University Press, 1997.
\newblock \href {https://doi.org/10.1017/CBO9780511526169} {\path{doi:10.1017/CBO9780511526169}}.

\bibitem{junk_new_new}
S.~Junk.
\newblock {Fluctuations of partition functions of directed polymers in weak disorder beyond the $L^2$-phase}, 2023.
\newblock \href {https://doi.org/10.48550/arXiv.2202.02907} {\path{doi:10.48550/arXiv.2202.02907}}.

\bibitem{junk2024taildistributionfunctionpartition}
Stefan Junk and Hubert Lacoin.
\newblock The tail distribution function of the partition function for directed polymer in the weak disorder phase.
\newblock 2024.
\newblock URL: \url{https://arxiv.org/abs/2405.04335}, \href {https://arxiv.org/abs/2405.04335} {\path{arXiv:2405.04335}}.

\bibitem{KPZ}
M.~Kardar, G.~Parisi, and Y.-C. Zhang.
\newblock Dynamic scaling of growing interfaces.
\newblock {\em Phys. Rev. Lett.}, 56:889--892, 1986.
\newblock \href {https://doi.org/10.1103/PhysRevLett.56.889} {\path{doi:10.1103/PhysRevLett.56.889}}.

\bibitem{lieb2001analysis}
E.~H. Lieb and M.~Loss.
\newblock {\em Analysis}.
\newblock American Mathematical Society, second edition, 2001.
\newblock \href {https://doi.org/10.1090/gsm/014} {\path{doi:10.1090/gsm/014}}.

\bibitem{Feynman-Kac}
B.~Lorenzo and C.~Nicoletta.
\newblock {The stochastic heat equation: Feynman-Kac formula and intermittence}.
\newblock {\em Journal of Statistical Physics}, 78, 1995.
\newblock \href {https://doi.org/10.1007/BF02180136} {\path{doi:10.1007/BF02180136}}.

\bibitem{Nikos3}
D.~Lygkonis and N.~Zygouras.
\newblock {Edwards–Wilkinson fluctuations for the directed polymer in the full ${L^{2}}$-regime for dimensions $d\ge 3$}.
\newblock {\em Annales de l'Institut Henri Poincaré, Probabilités et Statistiques}, 58(1):65--104, 2022.
\newblock \href {https://doi.org/10.1214/21-AIHP1173} {\path{doi:10.1214/21-AIHP1173}}.

\bibitem{DDP_regime}
D.~Lygkonis and N.~Zygouras.
\newblock {Moments of the $2D$ Directed Polymer in the Subcritical Regime and a Generalisation of the Erdös–Taylor Theorem}.
\newblock {\em Communications in Mathematical Physics}, 401(1):2483 -- 2520, 2023.
\newblock \href {https://doi.org/10.1007/s00220-023-04694-w} {\path{doi:10.1007/s00220-023-04694-w}}.

\bibitem{weakandstrong}
C.~Mukherjee, A.~Shamov, and O.~Zeitouni.
\newblock {Weak and strong disorder for the stochastic heat equation and continuous directed polymers in $d\geq 3$}.
\newblock {\em Electronic Communications in Probability}, 21, 2016.
\newblock \href {https://doi.org/10.1214/16-ECP18} {\path{doi:10.1214/16-ECP18}}.

\bibitem{reed1975ii}
M.~Reed and B.~Simon.
\newblock {\em II: Fourier Analysis, Self-Adjointness}.
\newblock Number v. 2 in Methods of Modern Mathematical Physics. Elsevier Science, 1975.
\newblock URL: \url{https://shop.elsevier.com/books/ii-fourier-analysis-self-adjointness/reed/978-0-08-092537-0}.

\bibitem{reed1981functional}
M.~Reed and B.~Simon.
\newblock {\em I: Functional Analysis}.
\newblock Methods of Modern Mathematical Physics. Elsevier Science, 1981.
\newblock URL: \url{https://shop.elsevier.com/books/i-functional-analysis/reed/978-0-08-057048-8}.

\bibitem{TAMURA1991433}
H.~Tamura.
\newblock \textup{The Efimov effect of three-body Schrödinger operators}.
\newblock {\em Journal of Functional Analysis}, 95(2):433--459, 1991.
\newblock \href {https://doi.org/10.1016/0022-1236(91)90038-7} {\path{doi:10.1016/0022-1236(91)90038-7}}.

\bibitem{me}
T.-C. Wang.
\newblock {On the space-time fluctuations of the SHE and KPZ equation in the entire $L^{2}$-regime for spatial dimensions $d \geq 3$}.
\newblock 2024.
\newblock URL: \url{https://arxiv.org/abs/2402.06874}, \href {https://arxiv.org/abs/2402.06874} {\path{arXiv:2402.06874}}.

\bibitem{efimov0}
D.R. Yafaev.
\newblock On the theory of the discrete spectrum of the three-particle schrödinger operator.
\newblock {\em Mat. Sb. (N.S.)}, 94(136):567--593, 1974.
\newblock \href {https://doi.org/10.1070/sm1974v023n04abeh001730} {\path{doi:10.1070/sm1974v023n04abeh001730}}.

\end{thebibliography}

\end{document}